\documentclass[11pt]{gthesis2}
\usepackage{thesis2}
\usepackage{latexsym}
\usepackage{epsfig}
\usepackage{amssymb}
\usepackage{amsmath}
\usepackage{longtable,lscape}
\usepackage{array}
\usepackage{verbatim}
\usepackage{amssymb}
\usepackage{graphicx}
\usepackage{color}
\usepackage{url}
\usepackage{eufrak}
\input xy
\xyoption{all}
\DeclareMathOperator{\epse}{\epsilon}

\let\Book=1
\newcommand{\vecx}{\mathbf{x}}
\newcommand{\matX}{\mathbf{X}}

\newcommand{\jordanalgebra}{\mathbf{M}_3(\mathbb{O})}
\newcommand{\tr}[1]{\mathbf{tr}\hspace{.000em}\left(#1\right)}
\newcommand{\ad}[1]{\mathbf{ad}\hspace{.000em}\left(#1\right)}
\newcommand{\mydet}[1]{\mathbf{det}\hspace{.000em}\left(#1\right)}
\newcommand{\p}{\mathcal{P}}
\newcommand{\n}{\mathcal{N}}
\newcommand{\g}{\varrho}
\renewcommand{\l}{\ell}
\newcommand{\F}{\mathbb{F}}

\newcommand{\SO}{so}

\newcommand{\dotBtz}{\dot B^{(2)}_{tz} - \dot B^{(3)}_{tz}}
\newcommand{\jordanproduct}{\circ}
\newcommand{\myendofproof}{\hspace{5.40in}$\square$}
\newcommand{\und}{\underline{\hspace{.66cm}}, \hspace{.08cm} }
\newcommand{\unde}{\underline{\hspace{.66cm}} }

\newcommand{\dotproduct}{\textrm{ {\tiny $\bullet$} }}
\renewcommand{\subsubsection}[1]{{\Large #1}}
\newcommand{\mynarrowlabel}{\mbox{\hspace{-.05in} $D_2 = so(4)$ \hspace{.20in} $A_2 = su(3)$ \hspace{.24in} $B_2 = so(5)$ \hspace{.20in} $C_2 = sp(2\cdot 2)$ \hspace{.10in} $G_2 = Aut(\mathbb{O})$} }

\setmyname{Aaron D. Wangberg} 
\setmyday{August 8, 2007}            
\setmycomyear{2008}                   
\setmyadvisor{Tevian Dray}       
\setmytitle{The Structure of $E_6$}

\setmymajor{Mathematics}                 
\setmydept{Mathematics}                  
\setmyschool{Oregon State University}    
\setmydegree{Doctor of Philosophy}       
\setmythesis{Thesis}                     

\begin{document}




%
\clearpage 
%
%
%
\begin{center}
AN ABSTRACT OF THE THESIS OF
\end{center}
\noindent {\raggedright \underline{\myname} for the degree of \underline{\mydegree} in \underline{\mymajor} presented on
\underline{\myday}. \\


Title: \underline{The Structure of $E_6$}}


\hfil\strut\\
Abstract approved: \hspace{0.5cm} \hrulefill\\
\phantom{Abstract approved:\ }\hfil\myadvisor\hfil\break

\vspace{-0.5em}

\noindent  A fundamental question related to any Lie algebra is to know its subalgebras.  This is particularly true in the case of $E_6$, an algebra which seems just large enough to contain the algebras which describe the fundamental forces in the Standard Model of particle physics.
  In this situation, the question is not just to know which subalgebras exist in $E_6$ but to know how the subalgebras fit inside the larger algebra and how they are related to each other.

\noindent In this thesis, we present the subalgebra structure of $sl(3,\mathbb{O})$, a particular real form of $E_6$ chosen for its relevance to particle physics through the connection between its associated Lie group $SL(3,\mathbb{O})$ and generalized Lorentz groups.  Given the complications related to the non-associativity of the octonions $\mathbb{O}$ and the restriction to working with a real form of $E_6$, we find that the traditional methods used to study Lie algebras must be modified for our purposes.  We use an explicit representation of the Lie group $SL(3,\mathbb{O})$ to produce the multiplication table of the corresponding algebra $sl(3,\mathbb{O})$.  Both the multiplication table and the group are then utilized to find subalgebras of $sl(3,\mathbb{O})$.  In particular, we identify various subalgebras of the form $sl(n,\mathbb{F})$ and $su(n,\mathbb{F})$ within $sl(3,\mathbb{O})$ and we also find algebras corresponding to generalized Lorentz groups.  Methods based upon automorphisms of complex Lie algebras are developed to find less obvious subalgebras of $sl(3,\mathbb{O})$.  While we focus on the subalgebra structure of our real form of $E_6$, these methods may also be used to study the subalgebra structure of any other real form of $E_6$.  A maximal set of simultaneously measurable observables in physics corresponds to a maximal set of Casimir operators in the Lie algebra.  We not only identify six Casimir operators in $E_6$, but produce a nested sequence of subalgebras and Casimir operators in $E_6$ containing both $su(3)\oplus su(2) \oplus u(1)$ corresponding to the Standard Model and the Lorentz group of special relativity.

\clearpage

\DisplayCopyright 
\DisplayTitlePage 
%
%
%
\thispagestyle{empty} {\baselineskip=14.5pt
\def\ruleline{\hbox to \hsize{\hrulefill}\\[-2ex]}
\noindent \underline{\mydegree} thesis of \underline{\myname} presented on \underline{\myday}
\strut\\
\strut\\
APPROVED:\\
\strut\\
\strut\\
\ruleline
Major Professor, representing \mymajor\\
\strut\\
\strut\\
\ruleline
Chair of the Department of \mydept\\
\strut\\
\strut\\
\ruleline
Dean of the Graduate School\\
\strut\\
\strut\\
\strut\\
\strut\\
I understand that my thesis will become part of the permanent collection of Oregon State University libraries.  My signature below authorizes
release of my thesis to any reader upon request.
\strut\\
\strut\\
\ruleline \hbox to \textwidth{\hfil \myname, Author \hfil} }
%
%
%
\clearpage
%
%
\begin{center}
ACKNOWLEDGEMENTS
\end{center}

\vspace{3 cm}

\noindent I wish to thank my advisor, Tevian Dray, and Corinne Manogue for their participation and advice as I worked to finish this thesis.  I also want thank them for allowing me to pursue interesting questions and pursuits, whether they be related to mathematics, teaching, or exploring the Pacific Northwest.

\noindent I also greatly appreciate the support and encouragement of my wife, Robyn.  This would not have been completed without her.

\noindent Thank you.

%
%
%
%

%
\clearpage
\tableofcontents 
\clearpage
\listoffigures

\clearpage

\listoftables 
\clearpage


\thispagestyle{lotappendixone}
\setcounter{reducedpagenumber}{\value{page}}

\clearpage

\newpage

\DisplayTitle \vspace*{3em} \pagestyle{myheadings}
\setcounter{page}{1} \thispagestyle{empty}

\newpage{}

\part{Introduction}
\label{part:intro}

\section{Motivation}

A fundamental idea in classifying Lie algebras is not only to list the possible algebras, but to find each algebra's subalgebras.  The traditional approach to this problem follows the work of Eugene Dynkin \cite{dyn2}, who considered complex Lie algebras.  Many properties of an algebra are more easily seen when it is considered as a complex algebra, and it is possible to list the complex algebras which exist within, and which contain, a given complex Lie algebra.
  It is also possible to determine the real subalgebras of a complex Lie algebra.  However, for applications, such as to physics, it is important not only to know the subaglebras of a real algebra but also to describe them in terms of an explicit representation.  Dynkin's methods for classifying complex algebras are not helpful in this situation.

There are four infinite families of complex Lie algebras and exactly five algebras, called the {\it exceptional Lie algebras}, which do not fit into those families \cite{killing, cartan}.  The infinite families of Lie algebras correspond to certain well-known groups over the three division algebras $\mathbb{R}$, $\mathbb{C}$, and $\mathbb{H}$ \cite{baez}.  It is known that the exceptional Lie algebras, named $G_2$, $F_4$, $E_6$, $E_7$, and $E_8$, are related to the octonions $\mathbb{O}$, the largest division algebra \cite{baez}.  The smallest algebra, $G_2$, is the automorphism group of the octonions.  The other four exceptional Lie algebras appear as the bottom row and right column of the Magic Square of Lie algebras, constructed first by Freudenthal and Tits in 1964 and 1966, respectively \cite{freudenthal_magic_square, tits_magic_square}.  Each division algebra $\mathbb{F} = \mathbb{R}, \mathbb{C}, \mathbb{H}, \mathbb{O}$ labels both a column and row in the $4 \times 4$ Magic Square, and each entry is a Lie algebra built over $\mathbb{F}_1 \times \mathbb{F}_2$.  The Magic Square idea naturally yields several nested chains of subalgebras associated with the division algebras.  For example, in \cite{gursey_magic_square}, G\"ursey showed that
$$ SU(2) \times U(1) \subset SU(2) \times U(1) \times SU(2) \times U(1) $$
$$ \hspace{3cm} \subset SU(4) \times SU(2) \times U(1) \subset SO(10) \times SO(2) \subset E_6$$
We will construct additional subalgebra chains contained within $E_6$ in this thesis.

As discussed in more detail below, the exceptional Lie group $E_6$ is related to the exceptional Jordan algebra of $3 \times 3$ hermitian octonionic matrices.  The corresponding Dynkin diagram for $E_6$ shows that it contains the group $SO(8)$, indicating that $E_6$ may share some of the triality characteristics exhibited by $SO(8)$.  While the Dynkin diagram can be used to reconstruct the commutation relations for $E_6$, and while these commutation relations can in turn be used to find subalgebras of $E_6$, we note that this is done entirely in the complex setting.  The subalgebras identified by these constructions are complex Lie algebras, and a real form of the algebra might not contain a real form of the complex subalgebra.

The division algebras also seem to play an important role in the physics of string theory and supersymmetry.  In 1983, Kugo and Townsend \cite{kugo_townsend} related spinors in spacetimes of dimension $3, 4,$ and $6$ with the sequence of Lorentz groups $SL(2,\mathbb{R})$, $SL(2,\mathbb{C})$, and $SL(2,\mathbb{H})$.  In \cite{manogue_fairlie_composite_string}, Fairlie and Manogue extended this to $SL(2,\mathbb{O})$ and constructed a parametrization of the bosonic string coordinates in space-times of dimension $3$, $4$, $6$, and $10$.  This parametrization exists in those dimensions due to the identity $|z||z^\prime| = |zz^\prime|$ which holds for the four division algebras.  In \cite{manogue_fairlie_covariant_superstring}, Fairlie and Manogue gave a  parametrization of the covariant superstring that both encodes the bosonic and fermionic components of the superstring and satisfies three of the four equations of motion.  Division algebras again played an integral part in this construction, leading Fairlie and Manogue to assert that the special dimensionality associated with sypersymmetric gauge theories is just a reflection of the properties of division algebras.  

A number of results specifically relate the octonions to the superstring and supersymmetry.  In \cite{manogue_sudbery_gen_sol}, Manogue and Sudbery use an octonionic formalism for ten-dimensional vectors and spinors to solve the equations of motion for the Green-Schwarz Lagrangian for the superstring. By representing a vector in  $(9+1)$-dimensional Minkowski space with a $2 \times 2$ hermitian octonionic matrix, they construct a lightlike ten-vector by {\it squaring} a column of two octonions, which is a spinor.  In \cite{schray}, Schray combines the spinor and vector into a $3 \times 3$ Grassmann, octonionic Jordan matrix and expresses the Lorentz and supersymmetry transformations of fermionic and bosonic variables in terms of Jordan products.  This connection between the Jordan algebra and superstrings has been pursued by others \cite{corrigan1989}.    In fact, the octonionic Jordan algebra has another nice connection to physics.  In \cite{gunaydin_1993}, G\"unaydin notes that the reduced structure group of the Jordan algebra of Hermitian octonionic $3 \times 3$ matrices is a real form of the complex Lie algebra $E_6$ with signature $(52,26)$.  Other approaches also relate $E_6$ to the octonionic Jordan algebra, see \cite{toppan, lin}.  The exceptional Lie group $E_6$ seems to be the right size to contain the Lie group $SU(3,\mathbb{C}) \times SU(2,\mathbb{C}) \times U(1,\mathbb{C})$, which describes the fundamental forces in the Standard Model of particle physics, and the Lorentz group $SO(3,1,\mathbb{R})$ of special relativity \cite{e6_and_decays, e6_and_physics}.  In the Standard Model, the strong force is associated with the Lie group $SU(3,\mathbb{C})$, the weak force is tied to $SU(2,\mathbb{C})$, and the electromagnetic force is represented by $U(1,\mathbb{C})$.

\section{Summary}

Groups are used in quantum mechanics and particle physics in a fundamental way.  In quantum mechanics, operators describe the energy and time evolution of a particle.  The transformations which leave the Hamiltonian operator invariant in the Schr\"odinger equation form a group.  The particle is only measured with certain values, and these operator eigenvalues may be labeled with a representation of the group.  In traditional quantum mechanics, the representations are affiliated with complex Lie algebras.  This complexification process involves a square root of -1, usually labeled $i$, and causes some difficulty.  Myrheim notes in \cite{myrheim} that the standard complex notation identifies the role played by several different square roots of -1 in the field equation and Hilbert space of quantum mechanics.  When working with the octonions, it is even more difficult to identify the object or objects which should play the various roles of $\sqrt{-1}$.  If $\sqrt{-1}$ lies within the octonions, which one is 'special' and why was it chosen?  If $\sqrt{-1}$ lies outside the octonions, where does it come from and how does it interact with the octonions?  By avoiding issues related to the complexification of Lie algebras, we can leave these questions unanswered while still learning a great deal about the real structure of $sl(3,\mathbb{O})$, a real form of $E_6$.

Given the connection between the Lie group $E_6$, Jordan matrices, and superstrings, several authors have discussed the mathematical properties of the exceptional Lie algebra $E_6$ and its subalgebra $F_4$ \cite{schafer,lin,allenAndFerrar}.  However, these discussions typically involve complex Lie algebras instead of real Lie algebras.
In \cite{manogue_dray}, Manogue and Dray construct the Lie group $SL(3,\mathbb{O})$, and the Lie algebra $sl(3,\mathbb{O})$ corresponding to this group is a real form of the Lie algebra $E_6$.  The present work builds upon their results to provide the algebraic subalgebra structure of $sl(3,\mathbb{O})$.  We provide the multiplication table for the real Lie algebra $sl(3,\mathbb{O})$ and give towers of subalgebras in $sl(3,\mathbb{O})$.  We use our knowledge of the multiplication table of $sl(3,\mathbb{O})$ to label the towers with the appropriate Casimir operators and give explicit bases for the various subalgebras.  We use the Casimir operators and explicit bases to distinguish between the physically interesting subalgebras of $sl(3,\mathbb{O})$ and other physicially uninteresting, but isomorphic, subalgebras.

Even though the material presented in Chapters $2$ and $3$ summarizes known results, the geometric treatment of root and weight diagrams in Chapter $2$ is not readily available in the literature.  While the interpretation of the group $E_6$ as $SL(3,\mathbb{O})$ is due to Manogue and Dray \cite{manogue_dray}, the detailed analysis of the Lie algebra $sl(3,\mathbb{O})$ presented here is new.  In particular, the notion of continuous type transformations in Section \ref{ch:E6_basic_structure.Type_Transformation}, the use of involutive automorphisms in Chapter \ref{ch:E6_further_structure}, and the detailed towers of subalgebras with nested Casimir operators in Section \ref{ch:E6_further_structure.chains_of_subalgebras} are due to the author.

\subsection{Thesis Organization}

This thesis is arranged into seven chapters.  The second chapter provides a general introduction to Lie groups and Lie algebras.  Section \ref{ch:Lie_groups_Lie_algebras.Basic_Definitions} gives the basic definitions we use for the study of Lie groups and Lie algebras.  We introduce an example in Section \ref{ch:Lie_groups_Lie_algebras.Basic_Definitions.An_Example} and use this example to illustrate the definitions in the further sections.  The definitions needed for Lie groups are given in Section \ref{ch:Lie_groups_Lie_algebras.Basic_Definitions.Lie_Groups}, while definitions for Lie algebras are contained within Section \ref{ch:Lie_groups_Lie_algebras.Basic_Definitions.Lie_Algebras}.
  Properties of a Lie group can often be studied using its associated Lie algebra (and vice versa) using the tools of differentiation, exponentiation, and commutation; these tools are defined in Section \ref{ch:Lie_groups_Lie_algebras.Basic_Definitions.Diff_Exp_Comm}.  The regular representation of a Lie algebra and the algebra's Killing form can be used to identify the signature of a simple Lie algebra; these tools are described in Section \ref{ch:Lie_groups_Lie_algebras.Basic_Definitions.regular_rep_Killing_form}.  We finish Section \ref{ch:Lie_groups_Lie_algebras.Basic_Definitions} with a discussion in Section \ref{ch:Lie_groups_Lie_algebras.Basic_Definitions.conventions} of the different conventions employed by mathematicians and physicists when defining the Lie algebras of the classical matrix groups; These groups and algebras are listed in Section \ref{ch:Lie_groups_Lie_algebras.Basic_Definitions.classical_matrix_groups}.

Section \ref{ch:Joma_Paper} includes verbatim a non-technical paper \cite{joma_paper} which covers the classification of the simple complex Lie algebras using their root and weight diagrams, and concludes with the complex Lie subalgebra structure of $E_6$.  This paper is intended for a non-expert audience and the definitions are kept deliberately informal; we therefore precede the paper here with precise definitions in Section \ref{ch:Joma_Paper.definitions}.
The paper itself begins in Section \ref{ch:Joma_Paper.Introduction} with a brief introduction.
  In Section \ref{ch:Joma_Paper.Root_Weight_Construction}, we show how to construct the root and weight diagrams of any simple complex algebra from its Dynkin diagram.  This section concludes with the root and weight diagrams of the rank $3$ complex simple Lie algebras.  In Section \ref{ch:Joma_Paper.subalgebras} we show how root and weight diagrams of algebras $g$ and $h$ can be used to show $g \subset h$.  In particular, we develop two methods, called {\it slicings} and {\it projections}, which we use to identify the complex subalgebras of algebras whose root and weight diagrams are of any dimension.  In Section \ref{ch:Joma_Paper.subalgebras_greater_than_3}, we apply these two methods to the four-dimensional root diagram of $F_4$.  We discuss how these methods may be applied to the six-dimensional diagrams of $E_6$, and produce a diagram showing the nesting of complex subalgebras contained within the complex Lie algebra $E_6$.  This diagram will be utilized and expanded upon in Chapters \ref{ch:E6_basic_structure} and \ref{ch:E6_further_structure}, where we produce the nesting of the real subalgebras contained within the real Lie algebra $sl(3,\mathbb{O})$, which is a real form of the complex Lie algebra $E_6$.  The conclusion of this paper, in Section \ref{ch:Joma_Paper.conclusion}, is also incorporated in the conclusion (Chapter \ref{ch:Conclusion}) of this thesis.

Chapter \ref{ch:Division_Algebras_Applications} discusses the division algebras and their application to structures pertinent to geometry and physics.  The first two sections of this chapter summarize the work of Manogue and Schray, given in \cite{manogue_schray}. Section \ref{ch:Division_Algebras_Applications.Quaternions_Octonions} discusses the algebraic and geometric properties of the four normed division algebras.  Their algebraic properties are covered in Section \ref{ch:Division_Algebras_Applications.Quaternions_Octonions.Division_Algebras}, and in Section \ref{ch:Division_Algebras_Applications.Quaternions_Octonions.Conjugation_Reflecitons_Rotations}, we construct rotations and reflections in $\mathbb{R}^4$ and $\mathbb{R}^8$ using conjugation maps.  These maps will play a vital role in the construction of the generalized Lorentz transformations of $SO(9,1,\mathbb{R})$, which is discussed in Section \ref{ch:Division_Algebras_Applications.Lorentz_Transformations.Lorentz_Transformations_for_Spacetimes_involving_Division_Algebras}, and in our construction of $SL(3,\mathbb{O})$ in Chapter \ref{ch:E6_basic_structure}.  We discuss triality as it pertains to the octonions in Section \ref{ch:Division_Algebras_Applications.Quaternions_Octonions.Triality}.  Section \ref{ch:Division_Algebras_Applications.Lorentz_Transformations} discusses Lorentz transformations and spacetime.  In particular, Section \ref{ch:Division_Algebras_Applications.Lorentz_Transformations.General_Lorentz_Transformations} discusses~$(k+1)$-dimensional spacetime and Lorentz transformations for any finite~$k$, Section \ref{ch:Division_Algebras_Applications.Lorentz_Transformations.Division_Algebras_and_Spacetimes} shows how spacetimes where~$k = 2, 3, 5, 9$ can be identified with vectors and spinors by using the division algebras, and Lorentz transformations utilizing the division alagebras are constructed for these particular spacetimes in Section \ref{ch:Division_Algebras_Applications.Lorentz_Transformations.Lorentz_Transformations_for_Spacetimes_involving_Division_Algebras}.  We conclude Chapter \ref{ch:Division_Algebras_Applications} by constructing the exceptional Jordan Algebra $M_3(\mathbb{O})$, which is also known as the {\it Albert Algebra}.  We describe the Lie group $SL(3,\mathbb{O})$ as the group which preserves the determinant of $M_3(\mathbb{O})$.  By letting $SL(3,\mathbb{O})$ act on the Albert Algebra in Chapter \ref{ch:E6_basic_structure}, we are able to construct the commutation table of the Lie algebra $sl(3,\mathbb{O})$.

Chapter \ref{ch:E6_basic_structure} incorporates the work of Manogue and Dray \cite{manogue_dray} to construct the Lie group $SL(3,\mathbb{O})$ and explicitly extends their results to the Lie algebra $sl(3,\mathbb{O})$, which is a real form of $E_6$.  We extend the construction of $SL(2,\mathbb{O}) = SO(9,1,\mathbb{R})$ which was given by Manogue and Schray in \cite{manogue_schray} and reviewed in Section \ref{ch:Division_Algebras_Applications.Lorentz_Transformations} from the $2 \times 2$ case to the $3 \times 3$ case.  This work uses the exceptional Jordan Algebra $M_3(\mathbb{O})$ and a discrete {\it type} map, as described in \cite{manogue_dray},  and builds upon the work of Manogue and Dray who showed that $SL(3,\mathbb{O})$ contains three $SO(9,1,\mathbb{R})$ subgroups.  We expand the results of Manogue and Dray by using the explicit representation of the Lie group $SL(3,\mathbb{O})$ to construct its Lie algebra $sl(3,\mathbb{O})$ in Section \ref{ch:E6_basic_structure.ConstructingE6algebra}.  In particular, we find an explicit basis for the underlying vector space of the algebra and use the group transformations to define the product in the algebra.  These computations were done using Maple and {\sc Perl}.  We also extend the work of Manogue and Dray by giving explicit bases for the three subalgebra chains of $so(9,1,\mathbb{R})$ within $sl(3,\mathbb{O})$ and identifying six Casimir operators for the exceptional Lie algebra $E_6$.  We discuss triality as it relates to the Lie groups $G_2$ and $SO(8,\mathbb{R})$ in Section \ref{ch:E6_basic_structure.Triality}, and introduce the notion of {\it strong} triality.  We also relate triality to our notion of {\it type} in this section.  In Section \ref{ch:E6_basic_structure.Type_Transformation}, we generalize our discrete notion of type to a continuous transformation in the group $SL(3,\mathbb{O})$.  We identify this transformation in a particular subgroup $SO(3,\mathbb{R})$ of $SL(3,\mathbb{O})$, and, using its corresponding algebra $so(3,\mathbb{R})$, we find that any algebra containing this algebra $so(3,\mathbb{R})$ must be {\it type independent}.  We use this idea to separate the algebra $sl(3,\mathbb{O})$ into chains of subalgebra of the form $g_1 \oplus g_2$.  In Section \ref{ch:E6_basic_structure.Reduction_of_O} we give an explicit basis for each of the subaglebras $su(n,\mathbb{F})$ and $sl(n,\mathbb{F})$ for $\mathbb{F} = \mathbb{C}, \mathbb{H}, \mathbb{O}$ and $n = 2,3$.  We again create a chain of subalgebras $g_1 \oplus g_2$ of $sl(3,\mathbb{O})$ where $g_1$ is either $su(n,\mathbb{F})$ or $sl(n,\mathbb{F})$ and $g_2$ is the maximal subalgebra of $sl(3,\mathbb{O})$ which commutes with $g_1$.  In Section \ref{ch:E6_basic_structure.fix_l} we identify a subalgebra $F_\l$ of $sl(3,\mathbb{O})$ which is the stabilizer of the type $1$ octonionic unit $\l$ in $M_3(\mathbb{O})$.  This subalgebra is neither simple nor semi-simple, and consists of $so(8,1,\mathbb{R})$ together with an abelian ideal of dimension $16$.  Finally, in Section \ref{ch:E6_basic_structure.Gell-Mann}, we give an isomorphism between a basis for our preferred octonionic representation of the subalgebra $su(3,\mathbb{C})^C \subset G_2$ and the $3 \times 3$ Gell-Mann matrices.

In Chapter \ref{ch:E6_further_structure}, we find additional subalgebras of $sl(3,\mathbb{O})$ by adapting methods involving automorphisms of Lie algebras.  We review the properties of automorphisms and involutive automorphisms of complex Lie algebras in \ref{ch:E6_further_structure.Automorphisms_of_Lie_Algebras}.  The Killing form is used to show that an involutive automorphism of a complex Lie algebra $g^\mathbb{C}$ will map one real form of $g^\mathbb{C}$ to a possibly different real form of $g^\mathbb{C}$.  In Section \ref{ch:E6_further_structure.Automorphisms_of_Lie_Algebras.id_algebras} we provide standard theorems and lemmas from \cite{cornwell, gilmore} that help us establish the identity and ``$\mathbb{R}$-simpleness'' of a real subalgebra of $sl(3,\mathbb{O})$.  One such technique which we repeatedly employ identifies a real subalgebra of $sl(3,\mathbb{O})$ using only the algebra's dimension, rank, and signature.  We construct three involutive automorphisms of $E_6$ in Section \ref{ch:E6_further_structure.three_automorphisms} and identify new subalgebras of $sl(3,\mathbb{O})$ using these automorphisms.  In particular, if $\phi$ is such an involutive automorphism, then both $sl(3,\mathbb{O}) \cap \phi\left(sl\left(3,\mathbb{O}\right)\right)$ and the pre-image of the compact part of $\phi\left(sl\left(3,\mathbb{O}\right)\right)$ are subalgebras of $sl(3,\mathbb{O})$.  We use the methods from Section \ref{ch:E6_further_structure.Automorphisms_of_Lie_Algebras.id_algebras} to determine the $\mathbb{R}$-semi-simple decompositions of certain real subalgebras of $sl(3,\mathbb{O})$.  We expand this approach in Section \ref{ch:E6_further_structure.more_applications_of_three_automorphisms} to find additional subalgebras of $sl(3,\mathbb{O})$ by using the fact that the set of involutive automorphisms form a group.  We use the action of the individual automoprhisms to partition the basis of $sl(3,\mathbb{O})$ into sets and combine various sets to form subalgebras of $sl(3,\mathbb{O})$.  This approach yields a number of subalgebras of $sl(3,\mathbb{O})$ which could not be readily found using the methods of Chapter \ref{ch:E6_basic_structure}.  We find a real form of $F_4$ with signature $(36,16)$ and multiple real forms of $A_5$, $C_4$ and $C_3$.  We provide explicit bases for these real forms of these subalgebras and identify their Casimir operators.  Finally, in Section \ref{ch:E6_further_structure.chains_of_subalgebras}, we produce chains of subalgebras of $sl(3,\mathbb{O})$ in which each larger algebra is obtained by extending the basis of a subalgebra.  In particular, the subalgebras in these chains use the same choice of Casimir operators.  We find this structure to be much more intricate than the structure indicated for the complex Lie algebra $E_6$ indicated in Chapter \ref{ch:Joma_Paper.subalgebras_greater_than_3}.

We summarize this work in Chapter \ref{ch:Conclusion} after identifying some open questions in Chapter \ref{ch:Open_Questions}.

\newpage{}

\part{Lie Groups and Lie Algebras}

\label{ch:Lie_groups_Lie_algebras}

This chapter begins with a brief introduction to the principle ideas needed for the study of Lie groups and Lie algebras.  In Section 2.1.1, we begin with an example of a Lie group.  We will develop this example throughout Section 2.1 to help illustrate the definitions of a group, algebra, differentiation and exponentiation, as found in \cite{boothby}.  Although we give the general definitions of differentiation and exponentiation, we shall also give equivalent definitions that are helpful for the study of the classical matrix groups, which is of particular interest to this thesis. After discussing the adjoint representation of a Lie algebra and the Killing form in Section 2.1.7, we include a paper in Section 2.2 that discusses the classification of Lie algebras.  This paper is to appear in the Journal of Online Mathematics, and is included without modification.  It concludes with a tower of complex Lie algebras contained within $E_6$.  Additional information about Lie groups and Lie algebras can be found in \cite{gilmore, cornwell, jacobsen}.

\section{Basic Definitions}

\label{ch:Lie_groups_Lie_algebras.Basic_Definitions}

\subsection{An Example}

\label{ch:Lie_groups_Lie_algebras.Basic_Definitions.An_Example}

We begin with a simple example.  The solid ball of radius $r = \pi$, in three dimensions, is an example of a Lie group of rotations.  Each point $p$ in the ball can be used to define an axis of rotation extending from the origin through the point $p$, while the distance $|p|$ from the origin to $p$ can be used to specify the amount of rotation about that axis.  Hence, the point at the origin represents a rotation about any axis through $0$ radians, or the identity transformation, while the point $p = (x = \pi, y = 0, z = 0)$ represents a rotation of $\pi$ radians about the $x$-axis. Every rotation $p$ has an inverse rotation $-p$, and we note that antipodal points on the boundary of the ball may be identified since they are the same rotation.  If we compose one rotation $p$ about one axis with another rotation $q$ about a second axis, the result is a rotation $s = q \circ p$ about a third axis.  It is important to note that the resulting axis and amount of rotation depends continuously on both $p$ and $q$.  Thus, we have given the ball of radius $r = \pi$ in three dimensions a group structure.  

\subsection{Lie Groups}
\label{ch:Lie_groups_Lie_algebras.Basic_Definitions.Lie_Groups}

Informally, a Lie group is both a group and a manifold.  Before we give a precise definition, we must first define the concepts of a topological and differentiable manifold as well as a group. We follow the treatment of Lie groups given in \cite{boothby}, where additional information may be found.

A {\it topological manifold} $M$ of dimension $n$ is a topological space which is Hausdorff, locally Euclidean of dimension $n$, and has a countable basis of open sets.  A {\it differentiable manifold} is a topological manifold $M$ with a $C^\infty$ (smooth) structure, which is a family $\mathcal{U} = \lbrace \left( u_\alpha, \phi_\alpha \right) \rbrace$ of coordinate neighborhoods such that
\begin{enumerate}
  \item $\cup\ u_\alpha$ covers $M$
  \item $u_\alpha \cap u_\beta \ne \emptyset \implies \phi_\alpha \circ \phi_\beta^{-1} $ and $\phi_\beta \circ \phi_\alpha^{-1}$ are diffeomorphisms 
  \item Any coordinate neighborhood $\left( u_\gamma, \phi_\gamma \right)$ compatible with every $ \left( u_\alpha, \phi_\alpha \right) \in \mathcal{U}$ is in $\mathcal{U}$
\end{enumerate}

In general, a differentiable manifold may require more than one coordinate neighborhood. The first condition above guarantees that there is a coordinate neighborhood for each part of the manifold, while the second condition states that we can smoothly change from one coordinate neighborhood to another whenever their intersection is non-empty.  The third condition guarantees that our atlas $\mathcal{U}$ of coordinate neighborhoods is complete.

A {\it group} is a set $G$ together with a map $\cdot : G \times G \to G$, referred to as multiplication, such that
\begin{enumerate}
\item $\exists \, e \in G$ such that $x \cdot e = x = e \cdot x \; \forall \, x \in G$
\item $\forall \, x \in G, \; \exists \, x^{-1} \in G$ such that $x \cdot x^{-1} = e = x^{-1} \cdot x$
\item $(x \cdot y) \cdot z = x \cdot (y \cdot z) \; \forall \, x,y,z \in G$
\end{enumerate}
That is, there is an identity element in the group, there is an inverse for every element in the group, and the multiplication is associative.

A {\it Lie group} is a differentiable manifold $G$ which is also a group whose composition map $G \times G \to G$ defined by $(x,y) = x\cdot y$ and inverse map $G \to G$ defined by $x \to x^{-1}$ are both $C^\infty$ mappings.

In the example from Section \ref{ch:Lie_groups_Lie_algebras.Basic_Definitions.An_Example}, the ball of radius $r = \pi$ is a $3$-dimensional topological manifold.  To treat this solid ball as a differentiable manifold, we must be able to put coordinates on the ball by mapping it to $\mathbb{R}^3$.  We already saw that we could interpret the solid ball as a group in Section \ref{ch:Lie_groups_Lie_algebras.Basic_Definitions.An_Example}.  The solid ball became a Lie group once it was assigned a smooth map giving the inverse of each rotation and a smooth map composing two rotations.  We note that the identity $e$ of this Lie group is the origin of the ball.  We stress that while we are primarily interested in thinking of Lie groups as rotation groups for this thesis, there are many other examples of Lie groups.

\subsection{Lie Algebras}
\label{ch:Lie_groups_Lie_algebras.Basic_Definitions.Lie_Algebras}

We review here the definition of a Lie algebra in its general formalism.  This treatment again is based upon \cite{boothby}.  Additional information about Lie algebras may also be found in \cite{cornwell, gilmore, jacobsen, onishchik, bourbaki}.

A {\it real Lie algebra}  $g$ is a vector space $V$ over $\mathbb{R}$ along with a product $\left[ \hspace{.15cm} , \hspace{.15cm} \right] : g \times g \to g$, called a commutator or bracket, which is bilinear,
$$ \left[ax + by, z \right] = a\left[x, z\right] + b\left[ y, z\right]$$
anti-commutative,
$$ \left[ x, y\right] = - \left[ y, x \right]$$
and satisfies the Jacobi identity
$$\left[ \left[ x,y \right], z\right] + \left[ \left[ y,z \right], x \right] + \left[ \left[ z, x \right], y \right] = 0$$
for any $x,y,z \in g$ and $a,b \in \mathbb{R}$.  

A {\it complex Lie algebra} $g$ is a vector space $V$ over $\mathbb{C}$ with a bilinear product which satisfies the three properties above, but now with $a,b \in \mathbb{C}$.

{\bf Example 1:}  Given a Lie group $G$, let $L_p : G \to G$ be left translation by $p$ (more formally a diffeomorphism from $G$ to itself given by $L_p(r) = pr$ for $p,r \in G$).  A vector field $X$ on G is {\it left-invariant} if, for any $p, r \in G$, it satisfies $(dL_{rp^{-1}})(X_p) = X_r$ where $X_p$ is the value of $X$ at $p$.  
Let $\mathfrak{X}(G)$ denote the set of all $C^\infty$ left-invariant vector fields defined on $G$.
This is a vector space over $\mathbb{R}$, and it may be equipped with a multiplication
$\left[ \hspace{.15cm}, \hspace{.15cm} \right] : \mathfrak{X}(G) \times \mathfrak{X}(G) \to \mathfrak{X}(G)$ given by
$$ \left[ X, Y \right]_p f = (XY - YX)_p f = X_p(Yf) - Y_p(Xf), \hspace{.75cm} f \in C^\infty(p)$$
called the commutator of $X$ and $Y$.  As shown in \cite{boothby}, the commutator of left-invariant vector fields is again left-invariant.  Hence, the vector space $\mathfrak{X}(G)$ along with the product $\left[ X, Y \right]$ is a real Lie algebra whose dimension is the dimension of the Lie group $G$. 

{\bf Example 2:} The Lie algebra of a classical matrix group $G$ is the tangent space $T_eG$ at the identity $e$ of the group with the product $\left[ \hspace{.15cm}, \hspace{.15cm} \right] : T_eG \times T_eG \to T_eG$ defined by 
$$ \left[ X, Y \right] = X Y - Y X$$
for any $X, Y \in T_eG$ where $X Y$ denotes ordinary matrix multiplication, as shown in \cite{cornwell}.  This concept of a Lie algebra of a matrix Lie group will be useful in Chapter \ref{ch:E6_basic_structure}.

There is a $1-1$ correspondence between left-invariant vector fields on a Lie group $G$ and the vectors of $T_eG$, as any vector from $T_eG$ may be used to generate a left-invariant vector field on $G$ and each left-invariant vector field on $G$ is completely determined by its value at the identity element $e$ of $G$.

In our example using the solid ball of rotations as the Lie group, the corresponding Lie algebra consists of the vector space $\mathbb{R}^3$.  The tangent vector to each rotation at the identity (origin) in the group is parallel to the axis of rotation.  The commutator of two tangent vectors $x$ and $y$ can be given by their cross product, $\left[ x, y \right] = x \times y$ in $\mathbb{R}^3$.

Finally, we introduce some standard terminology for Lie algebras:
\begin{itemize}
\item A {\it basis} for the Lie algebra $g$ is a basis of the underlying vector space.
\item Given subsets $g_1, g_2$ of a Lie algebra $g$, the notation $\left[ g_1, g_2 \right]$ denotes the set 
$$ \{ \left[ z_1, z_2 \right] | z_1 \in g_1, z_2 \in g_2 \} $$
\item A {\it Lie subalgebra} $g_1$ of a Lie algebra $g$ is a subset of $g$ that forms a Lie algebra with the same commutator and field as that of $g$.
\item Given a basis $B = \lbrace b_1, \cdots, b_n \rbrace$ of a real (complex) Lie algebra $g$, the {\it structure constants of $g$ with respect to the basis $B$ } are the $n^3$ real (complex) numbers $c^r_{pq}$ satisfying
$$ [ b_p, b_q] = \sum_{r=1}^{n} c^r_{pq} b_r \hspace{.75cm} p,q = 1,2,\cdots, n$$
\item The {\it complexification $V^\mathbb{C}$ of a real vector space $V$} is the complex vector space 
$$V \otimes_\mathbb{R} \mathbb{C} = V \oplus i V$$
The {\it complexification $g^\mathbb{C}$ of a real Lie algebra $g$} is a complexification of the vector space $g$ with the commutator
$$ \left[ x_1 + i y_1, x_2 + i y_2 \right] = \left[ x_1, x_2 \right] - \left[ y_1, y_2 \right] + i\left( \left[ x_1, y_2 \right] + \left[ y_1, x_2 \right] \right), \hspace{.75cm} x_1, x_2, y_1, y_2 \in g$$
extending the commutator of $g$.\footnote{Care must be taken when considering real Lie algebras whose elements are complex matrices, see \cite{big_cornwell}.}  Further information may be found in \cite{onishchik}.
\item A {\it real form} of a complex Lie algebra $g$ is a real Lie algebra $h$ whose complexification $h^\mathbb{C}$ is isomorphic to $g$ as a complex Lie algebra.
\item A Lie algebra $g$ is {\it abelian} if its commutator is identically zero, i.e. $\left[ g, g \right] = 0$.
\item A subalgebra $g_1$ of $g$ is called  {\it invariant} if if $[a, b] \in g_1$ for all $a \in g_1$ and $b \in g$.  Equivalently, we have $\left[ g_1, g \right] \subset g_1$.
\item A real or complex Lie algebra $g$ is called {\it simple} if its complexification $g^\mathbb{C}$ is not abelian and does not possess a proper invariant Lie subalgebra.
\item A real or complex Lie algebra $g$ is called {\it semi-simple} if its complexification $g^\mathbb{C}$ does not possess an abelian invariant subalgebra.  Equivalently, a semi-simple Lie algebra is a direct sum of simple Lie algebras.
\end{itemize}

Throughout this thesis, we often write $C = R$ where $C$ is a complex Lie algebra and $R$ is a real Lie algebra whose complexification is isomorphic to $C$.  Given our primary interest in the structure of the real subalgebras of a real form of $E_6$, we define two non-standard definitions for real Lie algebras corresponding to simple and semi-simple Lie algebras:
\begin{enumerate}
\item If $g$ is a real Lie algebra that is not abelian and does not possess a proper invariant Lie subalgebra, then we say $g$ is {\it $\mathbb{R}$-simple}.
\item If the real Lie algebra $g$ does not possess an abelian invariant subalgebra,  then it is called {\it $\mathbb{R}$-semi-simple}.  An $\mathbb{R}$-semi-simple Lie algebra $g$ is the direct sum of $\mathbb{R}$-simple real Lie algebras.
\end{enumerate}
These terms are intended to comment only on the real structure of a real Lie algebra.  With this terminology, the real Lie algebra $so(3,1,\mathbb{R})$ is $\mathbb{R}$-simple and not decomposable into the form $g_1 \oplus g_2$, even though $so(3,1,\mathbb{R})$ and its complexification, $D_2 = A_1 \oplus A_1$, are not simple.  
We will encouter the difference between simple, semi-simple, and abelian Lie algebras again in Section \ref{ch:Lie_groups_Lie_algebras.Basic_Definitions.regular_rep_Killing_form}. 

Finally, we note that derivations may be used to construct a Lie algebra from any algebra $A$ over a field $\mathbb{F} = \mathbb{R}$ or $\mathbb{C}$.    If $A$ is an algebra over the field $\mathbb{F}$, then an endomorphism\footnote{An {\it endomorphism} is an $\mathbb{F}$-linear map from a vector space $V$ to itself, where $V$ is a vector space over the field $\mathbb{F}$.  We use $End_\mathbb{F} V$ to denote the associative algebra of all $\mathbb{F}$-linear maps from $V$ to $V$ where $V$ is a vector space over $\mathbb{F}$. }
 $D$ of the vector space $A$ is called a {\it derivation} if $A(xy) = A(x)y + xA(y)$ for all $x,y \in A$.  If $D_1$ and $D_2$ are derivations of $A$, then 
$$[D_1, D_2] = D_1 D_2 - D_2 D_1$$
 is also a derivation of $A$.  The set of derivations of $A$ is a subalgebra of $gl(A)$.  This construction will be useful in Section \ref{ch:Lie_groups_Lie_algebras.Basic_Definitions.regular_rep_Killing_form} in the particular case that $A$ is a Lie algebra.

\subsection{Correspondence between Lie groups and Lie algebras}
\label{ch:Lie_groups_Lie_algebras.Basic_Definitions.Diff_Exp_Comm}

There is a very nice correspondence between elements of a Lie group $G$ and elements of its corresponding Lie algebra $g$.  
Roughly speaking, algebra elements are produced from group elements by differentiation, while exponentiation produces a group element from an element in the Lie algebra.  As seen in the two examples above, elements of the Lie algebra of a Lie group $G$ may be described either as the set of left-invariant vector fields of $G$ or as tangent vectors at the identity of $G$.
  Using these two processes, the commutator of Lie algebra elements can be described using the Lie group.  The following treatment again follows \cite{boothby}.

We begin with {\it differentiation}, the process used to pass from the Lie group to the associated Lie algebra.  Given a Lie group $G$ and the Lie group $\mathbb{R}$, the additive group of real numbers, if $\alpha : \mathbb{R} \to G$ is a homomorphism of $\mathbb{R}$ into $G$, then the image $\alpha(\mathbb{R})$ is called a {\it one-parameter subgroup of $G$}.  Note that the homomorphism property requires both that $\alpha(0) = e$, the identity element in $G$, and that
$$ \alpha(s_1) \alpha(s_2) = \alpha(s_1 + s_2) \hspace{.75cm} s_1,s_2 \in \mathbb{R}$$
The corresponding left-invariant vector field $X$ is determined by its value $X_e$ at the identity of $G$, which is given by
$$ X_e = \dot \alpha = \frac{d \alpha(s)}{d s}|_{s = 0}$$
It is an element of the Lie algebra $g$ associated to the Lie group $G$.
We will often use the notation $\dot \alpha$ to denote the tangent vector at the identity of the one-parameter subgroup $\alpha : \mathbb{R} \to G$.

Having created an element of the Lie algebra from a one-parameter subgroup of a Lie group $G$, we now describe the inverse process, called {\it exponentiation}.  Let $G$ be a Lie group and $X$ a $C^\infty$ vector field on $G$.  
  For each $p \in G$, there is a neighborhood $V \subset G$, a real number $\delta > 0$, and a $C^\infty$ mapping
$$ \gamma^v : I_\delta \times V \to G$$
which satisfies 
$$ \frac{\partial \gamma^v(s,q)}{\partial s} = X_{\gamma^v(s,q)}$$
and
$$ \gamma^v(0,q) = q $$
for all $q \in V$, where $I_\delta = \left( -\delta, \delta \right)$ is a real open interval.
This theorem is the manifold version of the existence theorem for ordinary differential equations.
We use this theorem in the special case where $X$ is a left-invariant vector field of $G$, in which case the interval $I_\delta$ may be extended to all of $\mathbb{R}$ and we may take $V = G$.  Let $p = e \in G$.
  As noted in Section \ref{ch:Lie_groups_Lie_algebras.Basic_Definitions.Lie_Algebras}, each left-invariant vector field $X$ on $G$ may be identified by its value $X_e$ at the point $e$ in $G$, that is, with a tangent vector $v = X_e$ in $T_eG$.  Then the mapping $\gamma^v: \mathbb{R} \times G \to G$ given above is a global action of $\mathbb{R}$ on $G$, and the map $\gamma^v_e : \mathbb{R} \to G$ given by 
$$ \gamma^v_e(s) = \gamma^v(s,e)  \hspace{.75cm} s \in \mathbb{R}$$
is a one-parameter subgroup of $G$.  We note that $\gamma^v_e$ is defined for each particular $v \in T_eG$.  Indeed, the map $\gamma_e : \mathbb{R} \times T_eG \to G$ given by 
$$\gamma_e(s,v) = \gamma^v_e(s)$$
is a $C^\infty$ map.  This map satisfies $\gamma_e(s,s_1 v) = \gamma_e(s_1 s, v)$ for all $s, s_1 \in \mathbb{R}$.  The {\it exponential map}, $\exp{}: T_eG \to G$, is defined by the formula 
$$\exp{(v)} = \gamma_e(1, v)$$.

The exponential map derives its name from the case of its use with the classical matrix groups.  Consider the Lie group $G = GL(n,\mathbb{F})$, with $\mathbb{F} = \mathbb{R} $ or $\mathbb{C}$, which consists of all non-singular $n\times n$ matrices over $\mathbb{F}$.  The Lie algebra of this Lie group is $g = gl_n(\mathbb{F})$, which is isomorphic to $M_n(\mathbb{F})$, the set of all $n \times n$ matrices over $\mathbb{F}$.  Let $A \in gl_n(\mathbb{F})$.  The flow of the left-invariant vector field defined by $A$ is obtained by integrating $\frac{d X}{d s} = X \cdot A$ with $s \in \mathbb{R}$ and $X(0) = Id$.  If $\mathbb{F}^n$ is equipped with a norm and $M_n(\mathbb{F})$ with the associated endomorphism norm, then 
$$ || A^n || \le ||A||^n$$
Hence, the series 
$$ \sum_{r = 0}^\infty \left( \frac{s^r}{r!} \right) A^r$$
converges, as do all of its derivatives.  Thus, 
$$ \exp{(sA)} = \sum_{r = 0}^\infty \left( \frac{s^r}{r!} \right) A^r$$
giving the exponential map its name.

We can illustrate differentiation and exponentiation using our example of the solid ball of rotations from Section \ref{ch:Lie_groups_Lie_algebras.Basic_Definitions.An_Example}.  We begin with differentiation.  Let $G$ be the Lie group under consideration.  Choose a point $p$ in the solid ball, and consider the straight line path $\alpha: \mathbb{R} \to G$ extending through $-p$ and $p$, which is given by
$$ \alpha(s) = s\vec{p} \hspace{.75cm} s \in \mathbb{R}$$
where $\vec{p}$ is the vector pointing from the origin to $p$.  
This map is a homomorphism, and the line traced out by $\alpha$ is a one-parameter subgroup of $G$.
We note that $\alpha(1) = \vec{p}$ and $\alpha(-1) = -\vec{p}$ are rotations about the same axis but in opposite directions.  Then
$$ \dot \alpha \equiv \frac{d \alpha(s)}{d s}|_{s=0} = \vec{p}$$
is the tangent vector to the path $\alpha$ at the origin.  It is an element of the Lie algebra $g$ associated to the Lie group $G$.
Our example also illustrates the exponential of a Lie algebra element.  The Lie group $G$ is of dimension three, implying that its Lie algebra $g = T_eG$  is also of dimension three.  Let $\vec{v} \in T_eG$ be any vector in $\mathbb{R}^3$ and define the path $\gamma^{\vec{v}}_e : \mathbb{R} \to G$ in $G$ by
$$ \gamma^{\vec{v}}_e(s) = s\vec{v} \hspace{.75cm} s \in \mathbb{R}$$
The image of this map is a line extending through the origin in $G$.  Hence, it is a one-parameter subgroup of $G$.  The line also depends continuously on our choice of $\vec{v} \in \mathbb{R}^3$, so we may define $ \gamma_e : \mathbb{R} \times T_eG \to G$ by 
$$\gamma_e(s,\vec{v}) = \gamma^{\vec{v}}_e(s) = s\vec{v}$$
  The exponential map $\exp{} : T_eG \to G$ is then given by the formula 
$$\exp{(\vec{v})} = \gamma_e(1,\vec{v}) = \vec{v} = v$$
 where we identify $\vec{v}$ with the point $v \in G$.  We emphasize that the image of 
$$\exp{(s\vec{v})} = \gamma_e(1,s\vec{v}) = s\vec{v} = \gamma^{\vec{v}}_e(s) \hspace{.75cm} s \in \mathbb{R}$$
 is not only a line through the origin containing $v$, but it is a one-parameter subgroup of $G$ whose tangent vector at the identity is $\vec{v} \in T_eG$.

It will be useful in Section \ref{ch:E6_basic_structure.ConstructingE6algebra} to have an expression for the commutator of two Lie algebra elements in terms of the corresponding group elements.  Given a Lie algebra $g$ of a Lie group $G$, we identify $g$ with the tangent space $T_eG$ of $G$.  Let $u,v \in T_eG$, and let $\exp{(su)}$ and $\exp{(sv)}$ denote their corresponding one-parameter subgroups in $G$.  We define the curve $\sigma : \mathbb{R} \to G$ in $G$ by 
$$ \sigma{(s)} = \exp{(-\frac{s}{2}u)} \exp{(-\frac{s}{2}v)} \exp{(\frac{s}{2}u)} \exp{(\frac{s}{2}v)} \hspace{.75cm} s \in \mathbb{R}$$
We note that $\sigma{(0)} = e$ is the identity in $G$, but that the image of $\sigma{(\mathbb{R})}$ is not in general a one-parameter subgroup of $G$.  Nevertheless, its tangent vector is the commutator $\left[ u, v \right]$ in the algebra.\footnote{ The first-order derivative of $\sigma{(s)}$ vanishes at $s = 0$.  However, the second-order derivative $\frac{d^2 \sigma{(s)}}{ds^2}|_{s=0}$ will be tangent to the curve $\sigma{(s)}$ at $s = 0$.  Further information may be found in \cite{varadarajan}.}

We finish this section with a note on the direct products of Lie groups and the direct sum of Lie algebras.  In particular, we show that if $g_1$ and $g_2$ are the Lie algebras of the Lie groups $G_1$ and $G_2$, then the Lie algebra of the direct product $G_1 \times G_2$ is $g_1 \oplus g_2$.  Further information regarding this fact may be found in \cite{cornwell}.

Let $G_1$, $G_2$ be two Lie groups with identities $e_1$ and $e_2$, respectively.  Consider the direct product of the manifolds $G_1$ and $G_2$, denoted $G_1\times G_2$, with the product map
$$  (x_1, y_1) (x_2, y_2) \to (x_1 x_2, y_1 y_2) \textrm{ for } x_1, x_2 \in G_1, y_1, y_2 \in G_2$$
and inverse map
$$  (x_1, y_1) \to ((x_1)^{-1}, (y_1)^{-1}) \textrm{ for } x_1, \in G_1, y_1 \in G_2 $$
where $(x_1)^{-1}$ and $(y_1)^{-1}$ is the inverse of $x_1$ and $y_1$ in $G_1$ and $G_2$, respectively.  Since $G_1$ and $G_2$ are Lie groups, the product defined in this way also makes $G_1\times G_2$ into a Lie group.  Its identity element is $(e_1,e_2)$.  As a manifold, its dimension is $|G_1|+|G_2|$.

The direct sum of Lie algebras $g_1$ and $g_2$ over $\mathbb{R}$ is denoted $g_1\oplus g_2$.  It is a direct sum of the underlying vector spaces with the product
$$           [(x_1, y_1), (x_2, y_2) ] = [ [x_1,x_2], [y_1,y_2]]  \textrm{ for } x_1, x_2 \in g_1, y_1, y_2 \in g_2 $$
which is anti-commutative, satisfies the Jacobi identity, and is bilinear.  The dimension of $g_1 \oplus g_2$ is $|g_1|+|g_2|$.

Cornwell \cite{cornwell} uses the fact that every finite-dimensional Lie algebra has a finite-dimensional faithful linear representation to prove that if $G_1$ and $G_2$ are two matrix Lie groups of dimensions $|G_1|$ and $|G_2|$, then $G_1 \times G_2$ is a matrix Lie group of dimension $|G_1|+|G_2|$, and that if $g_1$ and $g_2$ are the real Lie algebras of $G_1$ and $G_2$, respectively, then the real Lie algebra of $G_1 \times G_2$ is isomorphic to $g_1 \oplus g_2$.  
Noting that the Lie algebra of $G_1\times G_2$ will have dimension $|G_1|+|G_2|$, we shall prove this theorem for general finite-dimensional Lie groups $G_1$ and $G_2$ by constructing a basis for $g_1 \oplus g_2$ from $G_1 \times G_2$.  Let $\mathbb{B}$ be a basis for $G_1$.  For each $x_i \in \mathbb{B}$, the image of $(\exp{(sx_i)}, e_2)$ with $s \in \mathbb{R}$ is a one-parameter subgroup of $G_1 \times G_2$.  Hence, we obtain $|G_1|$ one-parameter subgroups of $G_1 \times G_2$.  Similarly, we construct $|G_2|$ one-parameter subgroups in $G_1 \times G_2$ of the form $(e_1, \exp{(sy_i)})$ for each $y_i$ in a basis of $G_2$.  Hence, we have $|G_1| + |G_2|$ subgroups of $G_1 \times G_2$.  Differentiating, we obtain a set of $|G_1| + |G_2|$ linearly independent elements of $g_1 \oplus g_2$, each of the form $(x_i, 0)$ or $(0, y_i)$.  Hence, we have constructed a basis of $g_1 \oplus g_2$ from $G_1 \times G_2$.  Since their dimension are the same, we have shown that the Lie algebra of $G_1 \times G_2$ is $g_1 \oplus g_2$.

\subsection{Regular Representation and Killing Form}
\label{ch:Lie_groups_Lie_algebras.Basic_Definitions.regular_rep_Killing_form}

In this section, we describe the adjoint representation of a Lie algebra and the Killing form.  We begin with the adjoint representation of the Lie group, and show how it produces the adjoint representation of the Lie algebra associated with the Lie group.  We then discuss the Killing form of the Lie algebra, a tool which is fundamental in classifying real and complex Lie algebras.  Finally, we include a note on general representations of Lie algebras and include definitions which will be useful in Section \ref{ch:Joma_Paper.Root_Weight_Construction}.  The treatment of the adjoint representation of a Lie group, representations of a Lie algebra, and the Killing form may be found in a number of references, including \cite{cornwell, boothby, knapp}.

We first define a representation of a Lie group $G$ on a finite-dimensional vector space $V$ before considering the case where $V = g$, the Lie algebra of the Lie group.  
A {\it finite-dimensional representation of a topological group $G$} is a continuous homomorphism $\Psi : G \to GL_{\mathbb{C}}(V)$ from $G$ into the group of invertible linear transformations on a finite-dimensional complex vector space $V$.
In the case that $G$ is a Lie group, it is often useful to let $V = g$, the Lie algebra of $G$.
  Identify $g$ with the tangent space $T_eG$ to $G$ at the identity.  Then, for any element $p \in G$, the map $\Psi_p : G \to G$ given by $\Psi_p(r) = prp^{-1}$ is a smooth isomorphism.  
If $\alpha : \mathbb{R} \to G$ is a smooth curve in $G$ such that $\alpha(0) = e \in G$ and $ \dot \alpha \in g$, then $\Psi_p(\alpha(t)) = p \alpha(t) p^{-1}$ is also a smooth curve in $G$ which passes through $e$ at $t = 0$.  Hence, the tangent vector to this curve, $p \dot \alpha p^{-1}$, is in $g$.  Thus, the corresponding isomorphism $Ad(p) := d\Psi_p : T_eG \to T_eG$ is a representation $Ad : G \to Aut(T_eG)$ of $G$ on its tangent space $T_eG$.  Identifying $T_eG$ with $g$, we call this isomorphism the {\it adjoint representation of the group on its algebra}.  Further information on representations may be found in \cite{knapp}.

Having defined a representation of a Lie group, we now discuss the related concept for a Lie algebra.
A {\it representation of a Lie algebra $g$ on a vector space $V$} is a homomorphism of Lie algebras $\rho : g \to (End_{\mathbb{F}}(V))^\mathbb{K}$ where $V$ is a vector space over $\mathbb{F} = \mathbb{R} \textrm{ or } \mathbb{C}$.
The map $\rho$ must be $\mathbb{K}$-linear and satisfy
$$ \rho([x,y]) = \rho(x)\rho(y) - \rho(y)\rho(x) \textrm{ for all } x,y \in g$$
We note that $\mathbb{K} = \mathbb{R}$ if $g$ is a real Lie algebra, and $\mathbb{K} = \mathbb{C}$ if $g$ is a complex Lie algebra.  
\footnote{The dimension of the space of endomorphisms $(End_{\mathbb{R}}(V))^\mathbb{R}$ is twice that of $(End_{\mathbb{C}}(V))^\mathbb{R}$.}
Ado has proven that every finite-dimensional Lie algebra has a finite-dimensional faithful linear representation \cite{bourbaki}.

As shown in \cite{knapp}, if $\Psi$ is a representation of a Lie group on a vector space $V$, then the differential at $e \in G$ gives a representation of the real Lie algebra $g$ of $G$ on the space $V$.  We illustrate this using the adjoint representation of a matrix Lie group $G$ on its Lie algebra $g$.
Let $\alpha: \mathbb{R} \to G$ be a smooth curve in $G$ such that $\alpha(0) = e \in G$ and $\dot \alpha = y \in g$.  
The map $t \to Ad_{\alpha(t)}(x)$, with $x \in g$, is a smooth map into $g$.  Computing its derivative, we find that
$$
\begin{array}{ccl}
 \frac{d}{dt} Ad_{\alpha(t)}(x) & = & \frac{d}{dt} [ \alpha(t) x \alpha^{-1}(t) ] \\
                                & = & \frac{d \alpha(t)}{dt} x \alpha^{-1}(t) + \alpha(t) x \frac{d \alpha^{-1}(t)}{dt} \\
                                & = & \frac{d \alpha(t)}{dt} x \alpha^{-1}(t) - \alpha(t) x \alpha^{-1}(t) \frac{d \alpha(t)}{dt} \alpha^{-1}(t) \\
\end{array}
$$
which, when evaluated at $t = 0$, gives
$$ \begin{array}{ccl}
  \frac{d}{dt} Ad_{\alpha(t)}(x)|_{t = 0} & = & \dot \alpha x - x \dot \alpha \\
                                          & = & yx - xy
\end{array}
$$
which is also in $g$.  The homomorphism $Ad$ induces a homomorphism of Lie algebras, denoted by $ad: g \to End(g)$, given by $ad_x(y) = [x,y]$.  This homomorphism is known as the {\it adjoint representation of a Lie algebra $g$}.

The adjoint representation of a Lie algebra may also be defined without any reference to its Lie group.  
It is important to note that, for any Lie algebra $g$, the linear map $ad : g \to End(g)$ given by 
$$(ad\; x)(y) = [x,y]$$
is a representation of the Lie algebra $g$.  This homomorphism satisfies
$$ ad\; [x,y] = (ad\; x) (ad\; y) - (ad\; y) (ad\; x)$$
which can be seen by rewriting the Jacobi identity for $g$ as
$$ [ [x,y], z] = [x, [y,z]] - [y, [x,z]]$$
However, the Jacobi identity of $g$ may also be rewritten as
$$ [ x, [y,z]] = [[x,y],z] + [y,[x,z]]$$
or 
$$ (ad\; x)([y,z]) = [(ad\; x)(y), z] + [y, (ad\; x)(z)]$$
showing that $ad$ is a derivation of $g$.
Hence, as noted in Section \ref{ch:Lie_groups_Lie_algebras.Basic_Definitions.Lie_Algebras}, the image $ad\; g$ of $g$ under this representation is a Lie algebra, and its commutator is given by
$$ [ad\; x, ad\; y] = (ad\; x) (ad\; y) - (ad\; y) (ad\; x)$$
where the product on the right hand side is the the product in $End(g)$.  Further information may be found in \cite{knapp}.

We find that the adjoint representation of a Lie algebra $g$ is particularly useful.  The representation is fully defined by its action on a basis of $g$.  Hence, given a basis $\lbrace b_1, \cdots, b_n \rbrace$ of an $n$-dimensional Lie algebra $g$, we may express the linear transformation $ad_{b_i}$ as an $n \times n$ matrix $\ad{b_i}$ whose components are given by
$$ [ b_i, b_j ] = \sum_{k = 1}^n \ad{b_i}_{kj}b_k$$
We note that the components of $\ad{b_i}$ are the structure constants of $g$.  The coefficients in column $j$ of $\ad{a}$ are the coefficients in the expansion of $ad_a(b_j) = [a, b_j] $ using the basis previously chosen.  Hence, the components of the matrices $\ad{b_i}$, $i = 1, \cdots, n$, are the structure constants of $g$.  Even though the individual structure constants are basis dependent, the entire collection describes the entire algebra, regardless of the chosen basis, in the following way:  Given a real (complex) Lie algebra $g$, we define the {\it Killing form} to be a symmetric bilinear $B$ form over $\mathbb{R}$ ($\mathbb{C}$) on $g$ whose value is given by $B(x,y) = \tr{\ad{x} \ad{y}}$ for any $x,y \in g$.  The Killing form is invariant under all automorphisms of $g$.  Further information on the Killing form may be found in \cite{knapp}.

We now list some important results related to the Killing form which will be particularly useful in Chapters \ref{ch:E6_basic_structure} and \ref{ch:E6_further_structure}.  Consult \cite{cornwell, gilmore} for a further discussion of these facts.

The Killing form provides useful criteria for proving that a Lie algebra $g$ is semi-simple, known as {\it Cartan's criterion for semi-simplicity}.  Given a Lie algebra $g$ with basis $\lbrace b_1, \cdots, b_n \rbrace$, we may use the Killing form to define an $n \times n$ matrix $\mathbf{B}$ whose elements are defined by $B_{ij} = B(b_i, b_j)$.  The Lie algebra $g$ is semi-simple if and only if the Killing form is non-degenerate, that is, if $\det{ \mathbf{B}} \ne 0$.
According to a theorem in \cite{cornwell}, it is possible to choose a basis $\lbrace b_1, \cdots, b_n \rbrace$ of $g$ such that $B(b_i, b_i) \in \lbrace -1, 0 +1 \rbrace$ and $B(b_i, b_j) = 0 $ for $i \ne j$. That is, $\mathbf{B}$ is a diagonal matrix with diagonal entries equal to $+1$, $-1$, or $0$.  Let $d_{+}, d_{-},$ and $d_{0}$ denote the number of $+1$, $-1$, and $0$ diagonal entries in $\mathbf{B}$ with such a basis chosen. While there are any number of ways to choose the basis of $g$, a theorem \cite{cornwell} on bilinear forms states that if a basis of $g$ is chosen so that $\mathbb{B}$ is diagonal with entries $+1, -1$ and $0$, then the values $d_{+}$, $d_{-}$, and $d_{0}$ are invariant.  
If $d_{0} > 0$, then $\mathbf{B}$ is degenerate and $g$ is not semi-simple.
On the other hand, if $d_{0} = 0$, then $g$ is semi-simple and the Killing form provides a further description of the Lie algebra $g$.
The elements $b_i$ such that $B(b_i, b_i) = -1$ correspond to {\it compact generators} or {\it rotations} in the Lie group, while the elements for which $B(b_i,b_i) = +1$ correspond to {\it non-compact generators} or {\it boosts}.  In this case, we say that $(d_{-}, d_{+})$ is the {\it signature} of $g$.  We note that a Lie algebra is compact if and only if its Killing form is negative definite.  Indeed, as noted in \cite{knapp}, a Lie algebra $g$ of a Lie group $G$ is compact if and only if $G$ is compact.

There are two other aspects of the Killing form which we will find useful in Chapters \ref{ch:E6_basic_structure} and \ref{ch:E6_further_structure}.  First, we say $x \in g$ is a {\it null rotation} if $B(x,x) = 0$.  We note that this definition neither implies that $\mathbf{B}$ is a diagonal matrix nor that $g$ has an abelian invariant subalgebra.  We will encounter null rotations in Section \ref{ch:E6_basic_structure.fix_l}.  Second, if $\phi$ is any automorphism of $g$, then 
$$B(\phi(x), \phi(y)) = B(x,y) \textrm{ for } x,y \in g$$
We will use this important property of the Killing form extensively in Chapter 5, 
Further information about the properties of the Killing form can be found in \cite{cornwell}.

In my work below, I often represent the adjoint representation of a Lie algebra in a {\it commutation table $C$}.  Having chosen a basis $\lbrace b_1, \cdots, b_n \rbrace$, the $ij$ entry in the table $C$ is given by $C_{ij} = [ b_i, b_j]$.  This entry contains the expansion of the $j$-th column of $\ad{b_i}$ with the basis $\lbrace b_1, \cdots, b_n \rbrace$.  Thus, each column of $C$ represents a complete $n\times n$ matrix in the adjoint representation.  The Killing form of $b_i$ and $b_j$ is computed in part by expanding the $i$-th and $j$-th columns of $C$ into their equivalent $n\times n$ adjoint representations of $b_i$ and $b_j$.  If $g$ is abelian, then $C$ will be the zero matrix. If $g$ has an invariant subalgebra $g^\prime \subset g$, we may expand a basis $\lbrace b_1, \cdots, b_k \rbrace$ for $g^\prime$ to a basis of $g$.  Then each of the entries in the first $k$ rows and first $k$ columns of $C$ will be linear combinations of $\lbrace b_1, \cdots, b_k \rbrace$ only.   Finally, if $g$ is semi-simple but not simple, we may choose a basis of $g$ over $\mathbb{C}$ such that $C$ is block diagonal.  This block-diagonal structure may not be possible if $g$ is a real Lie algebra and the basis must be chosen over $\mathbb{R}$.  However, if $g$ is a real Lie algebra which is $\mathbb{R}$-semi-simple, then we may choose a basis of $g$ over $\mathbb{R}$ such that $C$ is block diagonal.

\subsection{Different Conventions used by Physicists and Mathematicians}
\label{ch:Lie_groups_Lie_algebras.Basic_Definitions.conventions}
Mathematicians and physicists have slightly different conventions for the differentiation and exponentiation processes related to the classical matrix groups.  The conventions above are typically used by mathematicians.  For the classical matrix groups, these conventions lead to Lie algebra elements which are anti-hermitian.  In physics, one of the fundamental ideas in quantum mechanics is that observables are associated with hermitian matrices.  Hence, physicists insert a factor of $-i$ in the definition of differentiation, using $\sigma = -i \dot \alpha = -i \frac{d \alpha(s)}{d s}|_{s=0}$.  This factor of $-i$ is removed in the exponentiation process, using the definition $\exp(s\sigma) = \sum_{m=0}^\infty \frac{{i s \sigma}^m}{m!} = \sum_{m=0}^\infty \frac{( s \dot \alpha)^m}{m!}$, with $s \in \mathbb{R}$, for the classical matrix groups.  

It is important to note that these differences do not affect the content of the commutator in the algebra.  If $\lbrace \dot \alpha_x, \dot \alpha_y, \dot \alpha_z \rbrace$ are the mathematicians' Lie algebra elements, the corresponding physicists' conventions are $\lbrace -i \dot \alpha_x, -i\dot \alpha_y, -i\dot \alpha_z \rbrace$.  The mathematicians' commutator $[ \dot \alpha_x, \dot \alpha_y ] = \dot \alpha_z$ and physicists' commutator $[ -i \dot \alpha_x, -i \dot \alpha_y] = -i(-i \dot \alpha_z)$ are consistent, even if the physicists' commutator initially appears to have an extra factor of~$i$.

\subsection{Classical Matrix Groups and their Algebras}
\label{ch:Lie_groups_Lie_algebras.Basic_Definitions.classical_matrix_groups}
We list here the classical matrix groups, using their standard definitions as found in \cite{curtis}. We use $GL(n,\mathbb{F})$ with $\mathbb{F} = \mathbb{R}$ or $\mathbb{C}$ to denote the general linear group composed of $n\times n$ invertible matrices over $\mathbb{F}$ under the operation of matrix multiplication.  The notation $A^T$ is used to denote the transpose of the matrix $A$.  The hermitian conjugate of the matrix $A$ is denoted by $A^\dagger = \left( \overline{A} \right)^T$.  That is, if $A = \left(a_{ij}\right)$, then the hermitian conjugate of $A$ is given by $A^\dagger = \left( \overline{a_{ji}}\right)$.  Finally, the $n \times n$ identity matrix whose off-diagonal components are zero and whose diagonal components are one is denoted by $\mathbb{I}_n$, and $J_n = \left( \begin{array}{cc} 0 & \mathbb{I}_n \\ -\mathbb{I}_n & 0 \end{array} \right)$.

\begin{enumerate}
 \item Real special linear group of degree $n$: $$ SL(n,\mathbb{R}) = \lbrace g \in GL(n,\mathbb{R}) | \mydet{g} = 1 \rbrace$$
 \item Complex special linear group of degree $n$: $$SL(n,\mathbb{C}) = \lbrace g\in GL(n,\mathbb{C}) |\mydet{g} = 1 \rbrace$$
 \item Real orthogonal group of degree $n$: $$O(n,\mathbb{R}) = \lbrace g \in GL(n,\mathbb{R}) | g^{T} g = \mathbb{I}_n \rbrace$$
 \item Complex orthogonal group of degree $n$: $$O(n,\mathbb{C}) = \lbrace g \in GL(n,\mathbb{C}) | g^\dagger g = \mathbb{I}_n \rbrace$$
 \item Real (Complex) special orthogonal group of degree $n$: $$SO(n,\mathbb{R}) = O(n,\mathbb{R}) \cap SL(n,\mathbb{R})$$
$$SO(n,\mathbb{C}) = O(n,\mathbb{C}) \cap SL(n,\mathbb{C})$$
 \item Real symplectic group $$Sp(n,\mathbb{R}) = \lbrace g \in GL(2n,\mathbb{R}) | g^T J_n g = J_n \rbrace$$
 \item Complex symplectic group $$Sp(n,\mathbb{C}) = \lbrace g \in GL(2n,\mathbb{C}) | g^\dagger J_n g = J_n \rbrace$$
 \item Unitary group of degree $n$: $$U(n, \mathbb{C}) = \lbrace g \in GL(n,\mathbb{C}) | g^\dagger g = \mathbb{I}_n \rbrace$$
 \item Special unitary group of degree $n$: $$SU(n,\mathbb{C}) = U(n) \cap SL(n,\mathbb{C})$$
 \item Unitary symplecitic group of degree $n$: $$Sp(n) = U(2n,\mathbb{C})\cap Sp(n,\mathbb{C})$$
\end{enumerate}

For each of the classical matrix groups just given, we list the corresponding Lie algebras:

\begin{enumerate}
  \item $$sl(n,\mathbb{R}) = \lbrace x \in M(n,\mathbb{R}) | \tr{X} = 0 \rbrace$$
  \item $$sl(n,\mathbb{C}) = \lbrace x \in M(n,\mathbb{C}) | \tr{x} = 0 \rbrace$$
  \item $$o(n,\mathbb{R}) = \lbrace x \in M(n,\mathbb{R}) | x^T = -x \rbrace $$
  \item $$o(n,\mathbb{C}) = \lbrace x \in M(n,\mathbb{C}) | x^T = -x \rbrace $$
  \item $$so(n,\mathbb{R}) = o(n,\mathbb{R}) \cap sl(n,\mathbb{R})$$
        $$so(n,\mathbb{C}) = o(n,\mathbb{C}) \cap sl(n,\mathbb{C})$$
  \item $$sp(n,\mathbb{R}) = \lbrace x \in M(2n,\mathbb{R}) | x^T J_n + J_n x = 0 \rbrace$$
  \item $$sp(n,\mathbb{C}) = \lbrace x \in M(2n,\mathbb{C}) | x^T J_n + J_n x = 0 \rbrace$$
  \item $$u(n,\mathbb{C}) = \lbrace x \in M(n,\mathbb{C}) | X^\dagger = -x \rbrace$$
  \item $$su(n,\mathbb{C}) = \lbrace x \in M(n,\mathbb{C}) | x^\dagger = -x, \tr{x} = 0 \rbrace = u(n) \cap sl(n,\mathbb{C})$$
  \item $$sp(n,\mathbb{C}) = \lbrace x \in M(2n,\mathbb{C}) | x^\dagger J_n + J_n x = 0, x^\dagger = -x \rbrace = u(2n,\mathbb{C}) \cap sp(n,\mathbb{C})$$
\end{enumerate}

\section{Classifying Complex Lie Algebras}
\label{ch:Joma_Paper}

Lie algebras can be represented by various diagrams.
The remaining part of Chapter \ref{ch:Lie_groups_Lie_algebras} is a paper, to appear in the {\it Journal of Online Mathematics and its Applicataions}, which discusses the classification of complex Lie algebras.  It is included without modification, and includes some minor duplication of material previously covered in Chapter \ref{ch:Lie_groups_Lie_algebras}.  It is intended for a non-expert audience, and therefore many definitions have been kept deliberately informal in this non-technical paper.  Some of the terms with informal definitions have been defined previously in Chapter \ref{ch:Lie_groups_Lie_algebras}, and we have included the precise definitions of the other terms in Section \ref{ch:Joma_Paper.definitions}.  The actual text of the paper begins in Section \ref{ch:Joma_Paper.Introduction}.  

The paper is organized as follows.  Section \ref{ch:Joma_Paper.Root_Weight_Construction} discusses how the root and weight diagrams of a complex Lie algebra my be constructed from the algebra's Dynkin diagram.  This section includes a brief introduction to Lie algebras and other concepts covered previously in this thesis.  An algebra's root and weight diagrams may be used to identify its subalgebras.  This is shown for the three-dimensional case in Section \ref{ch:Joma_Paper.subalgebras} using the algebra $B_3 = so(7)$.  In particular, two methods are described which find all the sublagebras of the rank $3$ algebras.  These methods are generalized to higher dimensional diagrams in Section \ref{ch:Joma_Paper.subalgebras_greater_than_3}, including a detailed treatment of $F_4$.  This paper concludes with a tower of the complex Lie algebras contained within $E_6$.  The goal of the rest of this thesis is to expand upon the structure of this tower by considering the real subalgebras of $sl(3,\mathbb{O})$, our real form of $E_6$.

\subsection{Precise Definitions}
\label{ch:Joma_Paper.definitions}

We list here precise definitions for the terms which we left inadequately defined in the following paper.

\begin{enumerate}

\item  In the following paper, we often write $C = R$ where 
$$C = A_n, B_n, C_n, D_n, E_6, E_7, E_8, F_4, G_2$$
 is the Dynkin classification of the complex Lie algebra and $R$ is a real algebra whose complexification is $C$.  The equal sign is a short notation for the statement {\it The complexification of $R$ is equal as a complex Lie algebra to $C$.}

The {\it Cartan subalgebra} $h$ of a Lie algebra $g$ is a maximal nilpotent subalgebra of $g$.  Equivalently, it can be shown that the Cartan subalgebra $h$ of a Lie algebra $g$ is the maximal abelian subalgebra $h$ of a simple Lie algebra $g$.  See \cite{varadarajan} for additional information.

\item  A {\it representation} $p$ of a real (resp. complex) Lie algebra $g$ on a real (resp. complex) vector space $V$ is a mapping which assigns to each $x \in g$ a linear transformation $p(x)$ on $V$, and which is such that
$$ p(a x + by) = ap(x) + bp(y) \textrm{ for all } x,y \in g \textrm{ and scalars } a,b; $$
$$ p\left( [x,y] \right) = p(x)p(y) - p(y)p(x)$$
For linear transformations on a vector space, the expression $[X,Y]$ will mean $XY-YX$. Further information about these definitions may be found in \cite{Hausner_and_schwartz}.

Equivalently, a representation of a real (resp. complex) Lie algebra $g$ is a homomorphism $p : g \to gl_d(\mathbb{F})$ for some $d$.  This $d$, called the {\it degree} of the representation, is unrelated to the dimension $n$ of the Lie algebra $g$.  Two representations $p, p^\prime$ of a Lie algebra $g$ of degree $d$ are called {\it equivalent} if there exists a non-singular matrix $T$ over $\mathbb{F}$ such that $p^\prime(x) = T^{-1}p(x)T$ for all $x \in g$.  The degree of the adjoint representation of a Lie algebra $g$ is equal to the dimension of $g$.  Further information may be found in \cite{carter_macdonald_segal_taylor}.

We tacitly require that $p$ is not the trivial representation, i.e. $p(x) = 0$ for all $x \in g$.

\item If $p$ is a representation of a complex Lie algebra $g$ on a vector space $V$, then the linear functional $\lambda$ on $g$ is called a {\it weight} of $p$ if there exists a non-zero vector $v \in V$ such that $p(x)v = \lambda(x)v$ for each $x \in g$.  
 The vector $v$ is called a {\it weight vector} belonging to the weight $\lambda$.
A weight is called a {\it zero weight} if the linear functional is identially zero, and a {\it non-zero weight} if the linear functional is not identically zero.  
 Further information about these definitions may be found in \cite{Hausner_and_schwartz}.

We once refer to the adjoint representation as the $n \times n$ representation, in Section \ref{ch:Joma_Paper.Root_Weight_Construction}.

\item  In the case that $p$ is the adjoint representation of a complex Lie algebra $g$ on a vector space $V$, the weights are called {\it roots}.  We reserve the term {\it root vector} for another use in this paper, but refer to the weight vectors in the adjoint representation as {\it operators}.  The {\it zero roots} and {\it non-zero roots} are defined analogously.  We use the symbol $\Delta$ to denote the set of non-zero roots.
Further information about these definitions may be found in \cite{Hausner_and_schwartz}.

\item A {\it simple set of roots} is a set $\lbrace r_1, \cdots, r_k \rbrace$ of linearly independent roots of $\Delta$ such that every root $r \in \Delta$ may be written as $$r = \sum_{i = 1}^k c_i r_i$$ where the coefficients $c_i$ are unique and are integers which either  satisfy $c_i \ge 0$ or $c_i \le 0$ for all $i$ for each root $r$.  The ordered basis $\lbrace r_1, \cdots, r_k \rbrace$ of simple roots can be used to provide an ordering for the roots of $\Delta$, by writing $r > r^\prime$ if $r = \sum_{i}^{k} c_i r_i$ and $r^\prime = \sum_{i}^{k} c_i^\prime r_i$ and $c_i > c_i^\prime$ at the first $i$ for which $c_i \ne c_i^\prime$.  The operator $v$ is called a {\it raising operator} if its associated root $\lambda$ satisfies $\lambda > 0$, where $0$ is the zero linear functional.  It is called a {\it lowering operator} if its associated root $\lambda$ satisfies $\lambda < 0$. Further information about these definitions may be found in \cite{Hausner_and_schwartz}.

\item  Each weight $\lambda$ of a representation of a Lie algebra $g$ may be written uniquely as a real linear combination of the simple roots $\lbrace r_1, \cdots, r_k \rbrace$ of the adjoint representation of $g$.  The ordered basis $\lbrace r_1, \cdots, r_k \rbrace$ of simple roots provides an ordering for the weights, by writing $\lambda > \lambda^\prime$ if $\lambda = \sum_{i}^k c_i r_i$ and $\lambda^\prime = \sum_{i}^k c^\prime_i r_i$ and $c_i > c_i^\prime$ at the first $i$ for which $c_i \ne c_i^\prime$.  The {\it highest weight} $\lambda_0$ of a representation of a Lie algebra $g$ satisfies $\lambda_0 > \lambda$ for every weight $\lambda$.  Further information about these definitions may be found in \cite{Hausner_and_schwartz}.

\item A {\it root diagram} is a diagram in $\mathbb{R}^d$, where $d$ is the rank of the Lie algebra $g$, in which every root of the adjoint representation is plotted as a point $\lambda$ in $\mathbb{R}^d$, vectors are drawn extending from the origin to each $\lambda$, and a vector is drawn between two roots $\lambda_1$ and $\lambda_2$ precisely when their difference is another root in $\Delta$.  Further information about this definition may be found in \cite{kass}.

\item  We reserve the term {\it root vector} for either the arrow in space which points from the origin to the root $\lambda$ in a root or weight diagram of $g$ or a translation of that arrow which does not reverse its orientation.  The endpoint of a root vector is referred to as a {\it state}, to help differentiate the root from the vector pointing towards it.

\item A {\it weight diagram} is a diagram in $\mathbb{R}^d$, where $d$ is the rank of the Lie algebra $g$, in which every weight of a non-trivial representation is plotted as a point $\lambda$ in $\mathbb{R}^d$ and a vector is drawn between two weights $\lambda_1$ and $\lambda_2$ precisely when their difference is a root in $\Delta$.  Further information about this definitions may be found in \cite{kass}.

\item  A {\it minimal weight diagram} of a Lie algebra $g$ is a diagram with the least number of weights which still contains a vector for every root in $\Delta$.  There may be more than one minimal weight diagram in the sense that each arises from a different representation of the Lie algebra $g$.  However, all minimal weight diagrams are congruent.

\end{enumerate}

In the following paper, we often refer to the image of the representation $\Gamma$ of a Lie algebra as a representation of the Lie algebra.  We also refer to the weight and root diagrams as a representation of the Lie algebra.  We use this abuse of terminology not only because it is prevalent in the literature but also because the image of $\Gamma$ may be used to construct the weight or root diagrams, which in turn may be used to recreate the complexified algebra's structure constants.

The root and weight diagrams in this paper are best viewed in color.  In particular, we often refer to roots which have been colored red, blue, etc.  We have strived to include other descriptions which may be useful if this paper is printed in black and white only.

We describe in the paper how we may identify $g_1$ as a subalgebra of $g_2$ using root and weight diagrams, and simply state that $g_1$ is a subalgebra of $g_2$ if a root or weight diagram of $g_1$ may be identified as a subdiagram of $g_2$.  We also comment that both diagrams must use the same highest weight.  We list here the precise criteria that allows us to state that the complex Lie algebra $g_1$ is a subalgebra of the complex Lie algebra $g_2$:

  We say $D_1$ may be {\it embedded}  in $D_2$ if there is an isometry $r: \mathbb{R}^n \to \mathbb{R}^n$ such that the image of each vertex in $D_1$ is a vertex in $D_2$.    Each root or weight diagram $D_i$ of an algebra $g_i$ has a {\it highest shell of weights } $\overline{W_i}$, which consists of the weights which are the furthest from the center of the diagram $D_i$.  The diagram $D_1$ is a {\it subdiagram} of $D_2$ if $D_1$ may be embedded in $D_2$ and the image of $\overline{W_1}$ is contained within $\overline{W_2}$.   In this case, the algebra represented by the image of $D_1$ in $D_2$ will close, allowing us to say that the algebra $g_1$ is a subalgebra of $g_2$.

Now, let $p: \mathbb{R}^n \to \mathbb{R}^{n-k}$ denote a projection with $k < n$ which satisfies $|p(x) - p(y)| \le |x - y|$ for all $x,y \in \mathbb{R}^n$, and write $pD_2$ for the result of projecting the diagram $D_2$ from $\mathbb{R}^n $ to $\mathbb{R}^{n-k}$.  Considered as a diagram itself, $pD_2$ will contain a weight $pW^0_2$ such that $|pW^0_2| \ge |pW^i_2|$ for every other weight $pW^i_2$ of $pD_2$, where $|pW^i_2|$ denotes the distance from the origin to the weight $pW^i_2$. Let $\overline{pW_2}$ denote the set of weights which are a distance $|pW^0_2|$ away from the center of $pD_2$.  Then the algebra $g_1$ is a subalgebra of $g_2$ if $D_1$ may be embedded in $pD_2$.

\subsection{Introduction}
\label{ch:Joma_Paper.Introduction}

Lie algebras are classified using Dynkin diagrams, which encode the geometric structure of root and weight diagrams associated with an algebra.  This paper begins with an introduction to Lie algebras, roots, and Dynkin diagrams.  We then show how Dynkin diagrams define an algebra's root and weight diagrams, and provide examples showing this construction.  In Section~\ref{ch:Joma_Paper.subalgebras}, we develop two methods to analyze subdiagrams.  We then apply these methods to the exceptional Lie algebra ~$F_4$, and describe the slight modifications needed in order to apply them to ~$E_6$.  We conclude by listing all Lie subalgebras of ~$E_6$.

\subsection{Root and Weight Diagrams of Lie Algebras}
\label{ch:Joma_Paper.Root_Weight_Construction}
We summarize here some basic properties of root and weight diagrams.  Further information can be found in ~\cite{jacobsen}, ~\cite{cornwell}, and ~\cite{wiki_lie_algebra}.  A description of how root and weight diagrams are applied to particle physics is also given in ~\cite{cornwell}.

\subsubsection{Lie Algebras}

A Lie algebra $g$ of dimension $n$ is an $n$-dimensional vector space along with a product $[\ ,\ ]$:$ g \times g \to g$, called a commutator, which is anti-commutative ($[x,y] = -[y,x]$) and satisfies the Jacobi Identity
$$[x,[y,z]] + [y,[z,x]] + [z,[x,y]] = 0$$
for all $x,y,z \in g$.  A Lie algebra is called simple if it is non-abelian and contains no non-trivial ideals.  All complex semi-simple Lie algebras are the direct sum of simple Lie algebras.  Thus, we follow the standard practice of studying the simple algebras, which are the building blocks of the semi-simple algebras.

There are four infinite families of Lie algebras as well as five exceptional Lie algebras.  The algebras in the four infinite families correspond to special unitary matrices (or their generalizations) over different division algebras.  The algebras ~$B_n$ and ~$D_n$ correspond to real special orthogonal groups in odd and even dimensions, ~$SO(2n+1)$ and ~$SO(2n)$, respectively.  The algebras ~$A_n$ are the complex special unitary groups ~$SU(n+1)$, and the algebras ~$C_n$ correspond to unitary groups ~$SU(n,\mathbb{H})$ over the quaternions.

While Lie algebras are usually classified using their complex representations, there are particular real representations, based upon the division algebras, which are of interest in particle physics. Manogue and Schray ~\cite{manogue_schray} describe the use of quaternions ~$\mathbb{H}$ to construct ~$su(2,\mathbb{H})$ and ~$sl(2, \mathbb{H})$, which are real representations of ~$B_2=so(5)$ and ~$D_3 = so(5,1)$, respectively.  As they discuss, their construction naturally generalizes to the octonions ~$\mathbb{O}$, yielding the real representations ~$su(2,\mathbb{O})$ and ~$sl(2,\mathbb{O})$ of ~$B_4 = so(9)$ and ~$D_5 = so(9,1)$, respectively.  This can be further generalized to the ~$3 \times 3$ case, resulting in ~$su(3, \mathbb{O})$ and ~$sl(3, \mathbb{O})$, which preserve the trace and determinant, respectively, of a ~$3 \times 3$ octonionic hermitian matrix, and which are real representations of two of the exceptional Lie algebras, namely ~$F_4$ and ~$E_6$, respectively ~\cite{corrigan}.   The remaining three exceptional Lie algebras are also related to the octonions ~\cite{baez, okubo}.  The smallest, ~$G_2$, preserves the octonionic multiplication table and is ~$14$-dimensional, while ~$E_7$ and ~$E_8$ have dimensions ~$133$ and ~$248$, respectively.  A major step in describing all possible representations of ~$E_8$ was recently completed by the Atlas Project \cite{atlas_E8}.

In Sections \ref{ch:Joma_Paper.subalgebras} and \ref{ch:Joma_Paper.subalgebras_greater_than_3}, we label the Lie algebras using the standard complex Lie algebra label (e.g. ~$A_n, B_n$) and also with the name of a particular real form (e.g. ~$su(3)$, ~$so(7)$) of the algebra.  In section ~$4$, when discussing the subalgebras of ~$E_6$, we also give particular choices of real forms.

\subsubsection{Roots and Root Diagrams}

Every simple Lie algebra ~$g$ contains a maximal abelian subalgebra ~$h \subset g$, called a Cartan subalgebra, whose dimension is called the rank of ~$g$.  There exists a suitably normalized basis ~$\{h_{1}, \cdots, h_{l}\}$ of ~$h$ which can be extended to a basis 
$$\{ h_{1}, \cdots, h_{l}, g_1, g_{-1}, g_2, g_{-2}, \cdots, g_{\frac{n-l}{2}}, g_{-\frac{n-l}{2}} \} $$
 of ~$g$ satisfying:

\begin{enumerate}
  \item $[h_{i}, g_{j}] = \lambda_i^j g_{j}$ (no sum), $\lambda_i^j \in \mathbb{R}$
  \item $[h_{i}, h_{j}] = 0 $
  \item $[g_{j}, g_{-j}] \in h$
\end{enumerate}

The basis elements ~$g_j$ and ~$g_{-j}$ are referred to as raising and lowering operators.  Property ~$1$ associates every ~$g_{j}$ with an ~$l$-tuple of real numbers ~$r^{j} = \langle\lambda_1^j, \cdots, \lambda_l^j\rangle$, called roots of the algebra, and this association is one-to-one.  Further, if ~$r^{j}$ is a root, then so is ~$-r^{j} = r^{-j}$, and these are the only two real multiples of ~$r^{j}$ which are roots.  According to Property ~$2$, each ~$h_i$ is associated with the ~$l$-tuple ~$\langle0, \cdots, 0\rangle$.  Because this association holds for every ~$h_i \in h$, these ~$l$-tuples are sometimes referred to as zero roots.  For raising and lowering operators ~$g_j$ and ~$g_{-j}$, Property ~$3$ states that ~$r^{j} + r^{-j} = \langle0, \cdots, 0\rangle$.

Let ~$\Delta$ denote the collection of non-zero roots.  For roots ~$r^{i}$ and ~$r^{j} \ne -r^i$, if there exists ~$r^{k} \in \Delta$ such that ~$r^i + r^j = r^k$, then the associated operators for ~$r^i$ and ~$r^j$ do not commute, that is, $[ g_i, g_j ] \ne 0$.  In this case, ~$[g_i, g_j] = C^{k}_{ij}g_{k}$ (no sum), with ~$C^k_{ij} \in \mathbb{C}, C^i_{ij} \ne 0$.  If ~$r^i + r^j \not\in \Delta$, then ~$[g_i, g_j] = 0$.

When plotted in ~$\mathbb{R}^l$, the set of roots provide a geometric description of the algebra.  Each root is associated with a vector in ~$\mathbb{R}^l$.  We draw ~$l$ zero vectors at the origin for the ~$l$ zero roots corresponding to the basis ~$h_1, \cdots, h_l$ of the Cartan subalgebra.  For the time being, we then plot each non-zero root ~$r^i = \langle\lambda_1^i, \cdots, \lambda_l^i\rangle$ as a vector extending from the origin to the point ~$\langle\lambda_1^i, \cdots, \lambda_l^i\rangle$.  The terminal point of each root vector is called a state.  As is commonly done, we use ~$r^i$ to refer to both the root vector and the state.  In addition, we allow translations of the root vectors to start at any state, and connect two states ~$r^i$ and ~$r^j$ by the root vector ~$r^k$ when ~$r^k + r^i = r^j$ in the root system.  The resulting diagram is called a root diagram.

As an example, consider the algebra ~$su(2)$, which is classified as ~$A_1$.  The algebra ~$su(2)$ is the set of ~$2 \times 2$ complex traceless Hermitian matrices.  Setting 
\[
 \sigma_1 = \left[ \begin{array}{cc}
0 & 1 \\
1 & 0 \\
\end{array}
\right]
\hspace{.5cm}
\sigma_2 = \left[ \begin{array}{cc}
0 & -i \\
i & 0 \\
\end{array}
\right]
\hspace{.5cm}
\sigma_3 = \left[ \begin{array}{cc}
1 & 0 \\
0 & -1 \\
\end{array}
\right]
\] we choose the basis ~$h_1 = \frac{1}{2}\sigma_3$ for the Cartan subalgebra ~$h$, and use~$g_1 = \frac{1}{2}(\sigma_1+i\sigma_2)$ and~$g_{-1} = \frac{1}{2}(\sigma_1-i\sigma_2)$ to extend this basis for all of ~$su(2)$.  Then 
\begin{enumerate}
\item $[h_1, h_1] = 0$
\item $[h_1, g_1] = 1 g_1$
\item $[h_1, g_{-1}] = -1 g_{-1}$
\item $[g_{1}, g_{-1}] = h_1$
\end{enumerate}

\begin{figure}[htbp]
\begin{minipage}[t]{6in}
  \begin{center}
    \leavevmode
    \resizebox{3in}{!}{\includegraphics*[50,210][531,262]{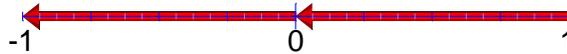}}
  \end{center}
\end{minipage}
\caption{Root diagram of $A_1 = su(2)$}
\label{fig:a1_root_diagram}
\end{figure}

By Properties ~$2$ and ~$3$, we associate the root vector ~$r^1 = \langle1\rangle$ with the raising operator ~$g_1$ and the root vector ~$r^{-1} = \langle-1\rangle$ with the lowering operator ~$g_{-1}$.  Using the zero root ~$\langle0\rangle$ associated with ~$h_1$, we plot the corresponding three points $(1)$, $(-1)$, and $(0)$ for the states ~$r^1$, ~$r^{-1}$, and ~$h_1$.  We then connect the states using the root vectors.  Instead of displaying both root vectors ~$r^1$ and ~$r^{-1}$ extending from the origin, we have chosen to use only the root vector ~$r^{-1}$, as ~$r^{-1} = - r^{1}$, to connect the states ~$(1)$ and ~$(0)$ to the states ~$(0)$ and ~$(-1)$, respectively.  The resulting root diagram is illustrated in Figure ~\ref{fig:a1_root_diagram}.

\subsubsection{Weights and Weight Diagrams}

An algebra ~$g$ can also be represented using a collection of ~$d \times d$ matrices, with ~$d$ unrelated to the dimension of ~$g$.  The matrices corresponding to the basis ~$h_1, \cdots, h_l$ of the Cartan subalgebra can be simultaneously diagonalized, providing ~$d$ eigenvectors.  Then, a list ~$w^m$ of ~$l$ eigenvalues, called a {\it weight}, is associated with each eigenvector.  Thus, the diagonalization process provides ~$d$ weights for the algebra ~$g$.  The roots of an ~$n$-dimensional algebra can be viewed as the non-zero weights of its ~$n \times n$ representation.

Weight diagrams are created in a manner comparable to root diagrams.  First, each weight ~$w^i$ is plotted as a point in ~$\mathbb{R}^l$.  Recalling the correspondence between a root ~$r^i$ and the operator ~$g^i$, we draw the root ~$r^k$ from the weight ~$w^i$ to the weight ~$w^j$ precisely when ~$r^k + w^i = w^j$, which at the algebra level occurs when the operator ~$g^k$ raises (or lowers) the state ~$w^i$ to the state ~$w^j$.

\begin{figure}[htbp]
\begin{center}
\begin{minipage}[t]{5cm}
  \begin{center}
    \leavevmode
    \resizebox{3cm}{!}{\includegraphics*[120,80][484,400]{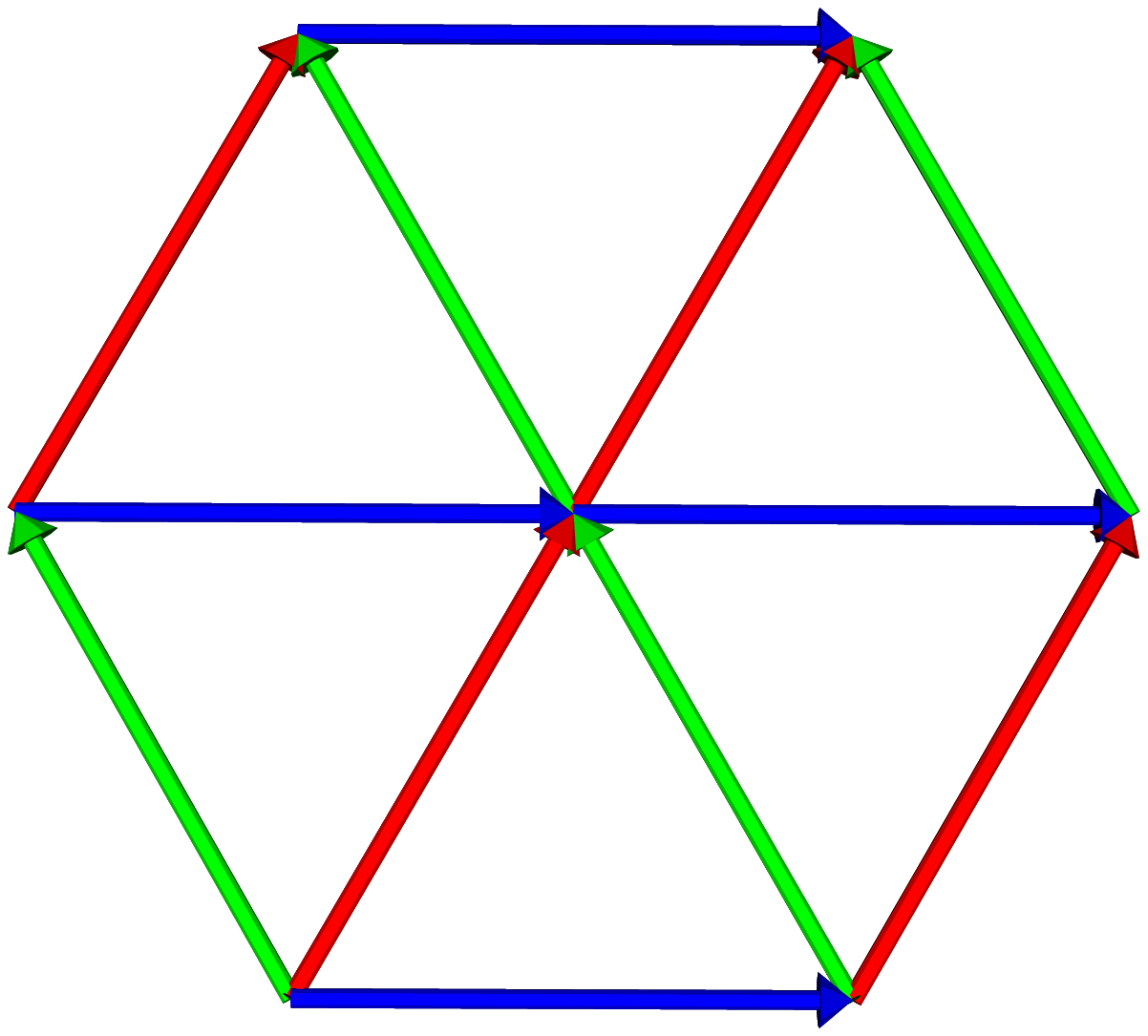}}
  \end{center}
\end{minipage}
\hspace{1cm}
\begin{minipage}[t]{5cm}
  \begin{center}
    \leavevmode
    \resizebox{3cm}{!}{\includegraphics*[200,100][380,260]{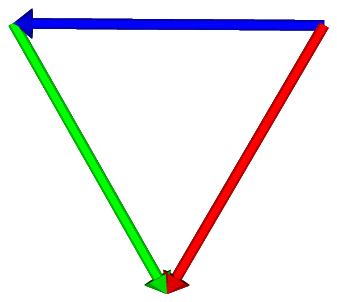}}
  \end{center}
\end{minipage}
\end{center}
\caption{Root and minimal weight diagrams of ~$A_2=su(3)$}
\label{fig:root_and_weight_diagrams_A2}
\end{figure}

The root and minimal non-trivial weight diagrams of the algebra ~$A_2 = su(3)$ are shown in Figure ~\ref{fig:root_and_weight_diagrams_A2}.  The algebra has three pairs of root vectors, which are oriented east-west (colored blue), roughly northeast-southwest (colored red), and roughly northwest-southeast (colored green).  
  The algebra's rank is the dimension of the underlying Euclidean space, which in this case is~$l=2$, and the number of non-zero root vectors is the number of raising and lowering operators.  The minimal representation contains three independent arrows, which provide six different roots (although we've only indicated one of each raising and lowering root pair).  Using the root diagram, the dimension of the algebra can now be determined either from the number of non-zero roots or from the number of root vectors extending in different directions from the origin.  Both diagrams indicate that the dimension of ~$A_2 = su(3)$ is ~$8 = 6+2$.  Root and (non-trivial) weight diagrams indicate the underlying algebra's dimension and rank, which is (almost) enough information to identify the algebra.
\footnote{For algebras of rank~$6$ and lower, the exceptions are that ~$B_n = so(2n+1)$ and ~$C_n=sp(2\cdot n)$ have the same dimension for each rank ~$n$, and that ~$B_6$, ~$C_6$, and ~$E_6$ all have dimension~$78$.}

\subsubsection{Constructing Root Diagrams from Dynkin Diagrams}

In 1947, Eugene Dynkin simplified the process of classifying complex semi-simple Lie algebras by using what became known as Dynkin diagrams ~\cite{dyn2}.  He realized that every root in a rank ~$l$ algebra can be expressed as an integer sum or difference of ~$l$ simple roots.  Further, the relative lengths and interior angle between pairs of simple roots fits one of four cases.  A Dynkin diagram records the configuration of an algebra's simple roots.

Each node in a Dynkin diagram represents one of the algebra's simple roots.  Two nodes are connected by zero, one, two, or three lines when the interior angle between the roots is ~$\frac{\pi}{2}$, ~$\frac{2\pi}{3}$, ~$\frac{3\pi}{4}$, or ~$\frac{5\pi}{6}$, respectively.  If two nodes are connected by ~$n$ lines, then the magnitudes of the corresponding roots satisfy the ratio ~$1:\sqrt{n}$.  An arrow is used in the Dynkin diagram to point toward the node for the smaller root.  If two roots are orthogonal, no direct information is known about their relative lengths.

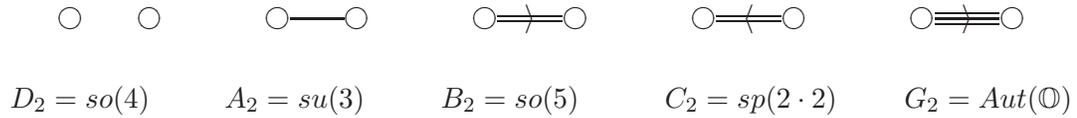
\begin{figure}[htbp]
\begin{minipage}[t]{6in}
  \begin{center}
    \leavevmode
    \xymatrixcolsep{10pt}
    \xymatrix@M=0pt@W=0pt{
*+=[o][F-]{\hspace{.35cm} } & & *+=[o][F-]{\hspace{.35cm} } & \hspace{.65cm} & *+=[o][F-]{\hspace{.35cm} } \ar@{-}[rr]&  & *+=[o][F-]{\hspace{.35cm} } &  \hspace{.65cm} & *+=[o][F-]{\hspace{.35cm} } \ar@2{-}[rr] & \rangle & *+=[o][F-]{\hspace{.35cm} } & \hspace{.65cm} & *+=[o][F-]{\hspace{.35cm} }  & \langle & *+=[o][F-]{\hspace{.35cm} } \ar@2{-}[ll] & \hspace{.65cm} & *+=[o][F-]{\hspace{.35cm} } \ar@3{-}[rr] & \rangle & *+=[o][F-]{\hspace{.35cm} } \\
    }\\ 
\vspace{1em}
\hspace{-.00in} $D_2 = so(4)$ \hspace{.29in} $A_2 = su(3)$ \hspace{.30in} $B_2 = so(5)$ \hspace{.35in} $C_2 = sp(2\cdot 2)$ \hspace{.25in} $G_2 = Aut(\mathbb{O})$
  \caption{Rank 2 Dynkin diagrams}
  \label{fig:rank_2_dynkin_diagrams}
  \end{center}
\end{minipage}
\end{figure}

\begin{figure}[htbp]
\begin{minipage}[t]{6in}
  \begin{center}
    \leavevmode
    \resizebox{5.0in}{!}{\includegraphics*[20,50][553,120]{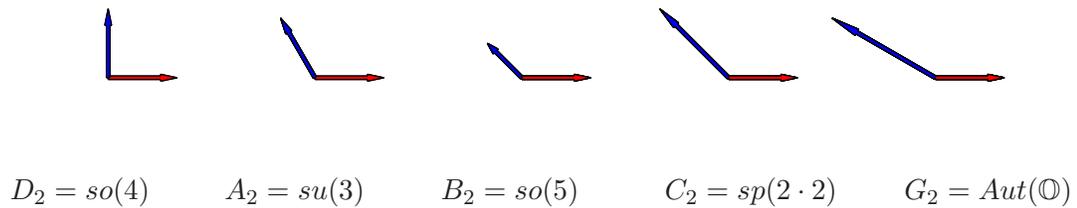}}
\\
\vspace{1em}
\hspace{-.00in} $D_2 = so(4)$ \hspace{.29in} $A_2 = su(3)$ \hspace{.30in} $B_2 = so(5)$ \hspace{.35in} $C_2 = sp(2\cdot 2)$ \hspace{.25in} $G_2=Aut(\mathbb{O})$
  \caption{Rank 2 simple roots}
  \label{fig:rank_2_simple_roots}
  \end{center}
\end{minipage}
\end{figure}

We give the Dynkin diagrams for the rank ~$2$ algebras in Figure ~\ref{fig:rank_2_dynkin_diagrams} and the corresponding simple root configurations in Figure ~\ref{fig:rank_2_simple_roots}.  For each algebra, the left node in the Dynkin diagram corresponds to the root ~$r^1$ of length ~$1$, colored red and lying along the horizontal axis, and the right node corresponds to the other root ~$r^2$, colored blue.

In ~$\mathbb{R}^l$, each root ~$r^i$ defines an ~$(l-1)$-dimensional hyperplane which contains the origin and is orthogonal to ~$r^i$.  A {\it Weyl reflection} for ~$r^i$ reflects each of the other roots ~$r^j$ across this hyperplane, producing the root ~$r^k$ defined by 
$$r^k = r^j - 2(\frac{r^j \dotproduct r^i}{|r^i|})\frac{r^i}{|r^i|} $$
According to Jacobsen ~\cite{jacobsen}, the full set of roots can be generated from the set of simple roots and their associated Weyl reflections.

\begin{figure}[htbp]
\begin{center}
  \begin{minipage}[t]{5in}
    \begin{center}
      \leavevmode
      Step ~$1$: Weyl reflection using root $r^1$\\
      \resizebox{5in}{!}{\includegraphics*[8,22][567,123]{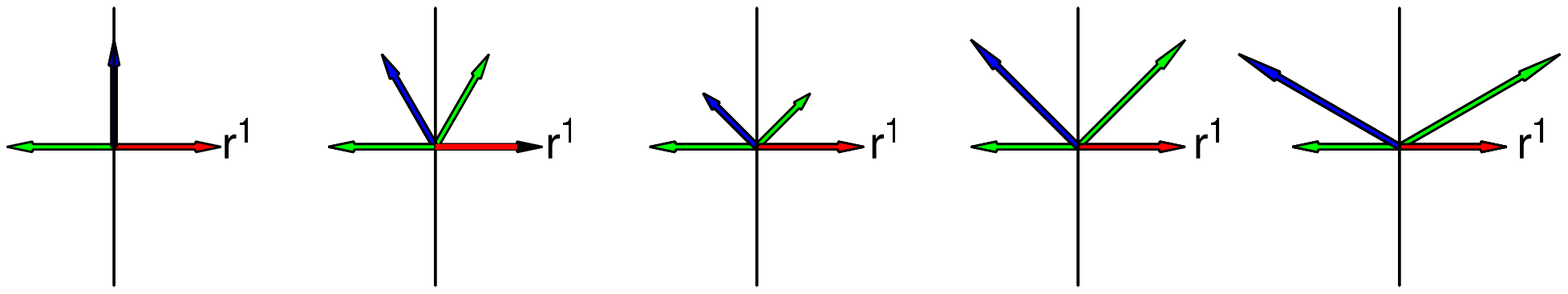}}
      \begin{minipage}[t]{5in}
        \mynarrowlabel
      \end{minipage}
    \end{center}
    \vspace{.6cm}
  \end{minipage}
  \begin{minipage}[t]{5in}
    \begin{center}
      \leavevmode
      Step ~$2$: Weyl reflection using root $r^2$\\
      \resizebox{5in}{!}{\includegraphics*[8,22][567,141]{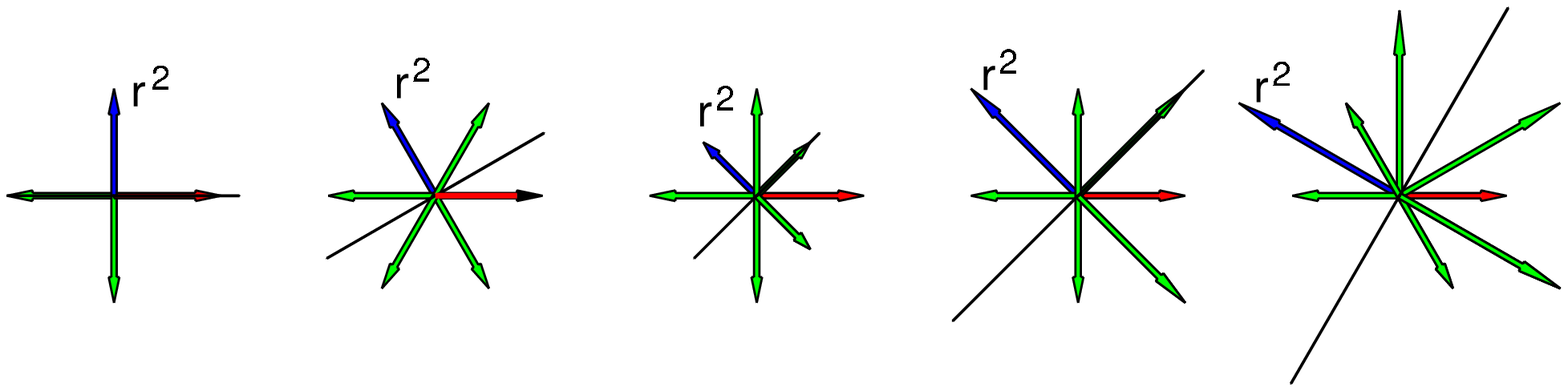}}
      \begin{minipage}[t]{5in}
        \mynarrowlabel
      \end{minipage}
    \end{center}
    \vspace{.6cm}
  \end{minipage}
  \begin{minipage}[t]{5in}
    \begin{center}
      \leavevmode
      Step ~$3$: Weyl reflection using root $r^1$\\
      \resizebox{5in}{!}{\includegraphics*[8,22][567,141]{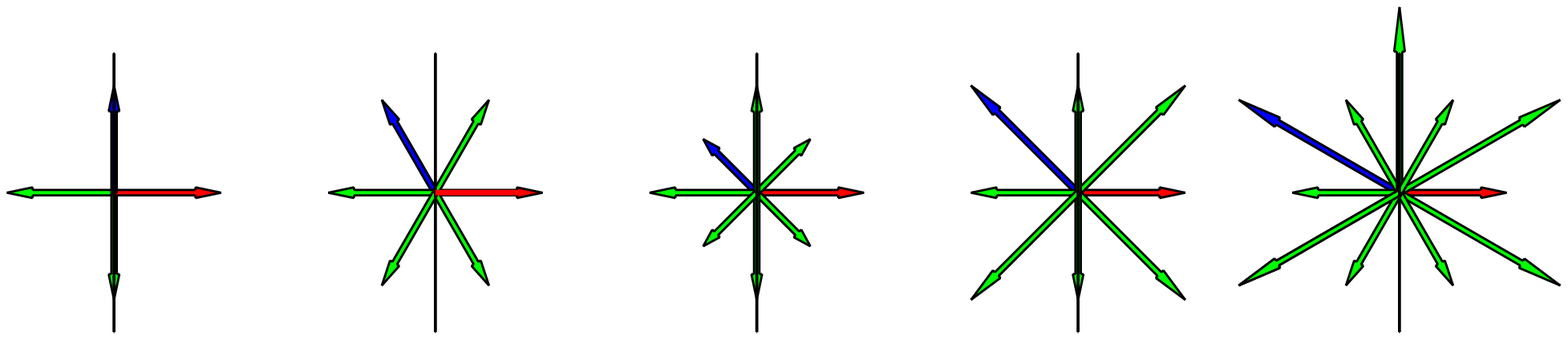}}
      \begin{minipage}[t]{5in}
        \mynarrowlabel
      \end{minipage}
    \end{center}
    \vspace{.6cm}
  \end{minipage}
  \begin{minipage}[t]{5in}
    \begin{center}
      \leavevmode
      Step ~$4$: Weyl reflection using root $r^2$\\
      \resizebox{5in}{!}{\includegraphics*[8,4][567,141]{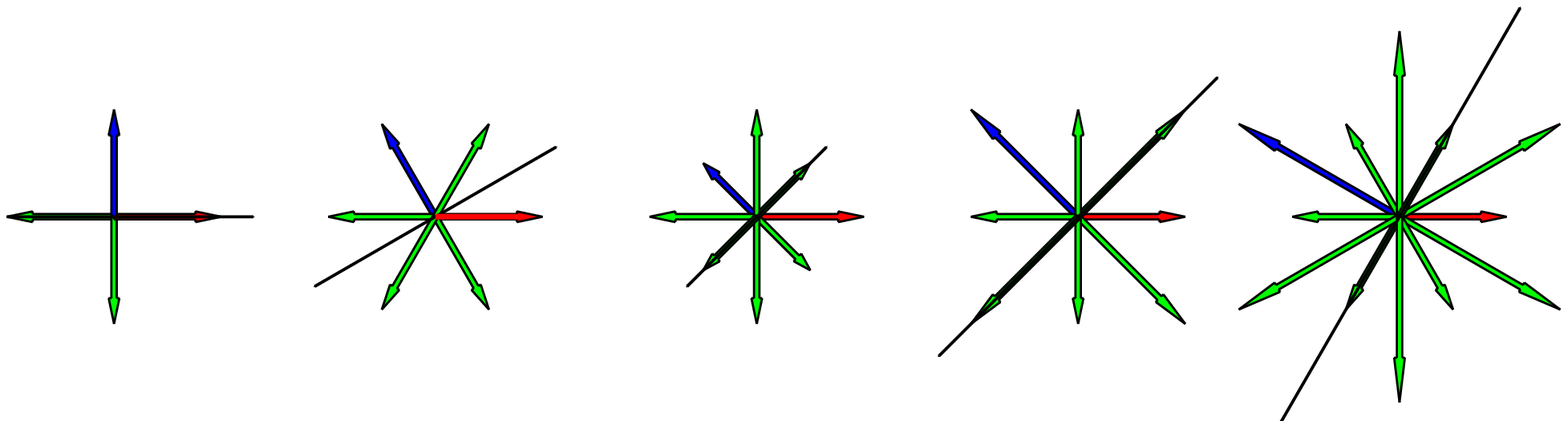}}
      \begin{minipage}[t]{5in}
        \mynarrowlabel
      \end{minipage}
    \end{center}
  \caption{Generating an algebra's full root system using Weyl reflections}
  \label{fig:rank_2_weyl_reflections_animation}
  \end{minipage}
\end{center}
\end{figure}

We illustrate how the full set of roots can be obtained from the simple roots using Weyl reflections in Figure ~\ref{fig:rank_2_weyl_reflections_animation}.  We start with the two simple roots for each algebra, as given in Figure ~\ref{fig:rank_2_simple_roots}.  For each algebra, we refer to the horizontal simple root, colored red, as ~$r^1$, and the other simple root, colored blue, as ~$r^2$.  Step ~$1$ shows the result of reflecting the simple roots using the Weyl reflection associated with ~$r^1$.  In this diagram, the black thin line represents the hyperplane orthogonal to ~$r^1$, and the new resulting roots are colored green.  Step ~$2$ shows the result of reflecting this new set of roots using the Weyl reflection associated with ~$r^2$.  At this stage, both ~$D_2 = so(4)$ and~$A_2 = su(3)$ have their full set of roots.  We repeat this process again in steps 3 and 4, using the Weyl reflections associated first with ~$r^1$ and then with ~$r^2$.  The full root systems for ~$B_2 = so(5)$ and ~$C_2 = sp(2\cdot 2)$ are obtained after the first three Weyl reflections.  Only ~$G_2$ requires all four Weyl reflections.

\begin{figure}[htbp]
\begin{center}
  \begin{minipage}[t]{2.4cm}
    \begin{center}
      \leavevmode
      \resizebox{2.1cm}{!}{\includegraphics*[202,144][389,333]{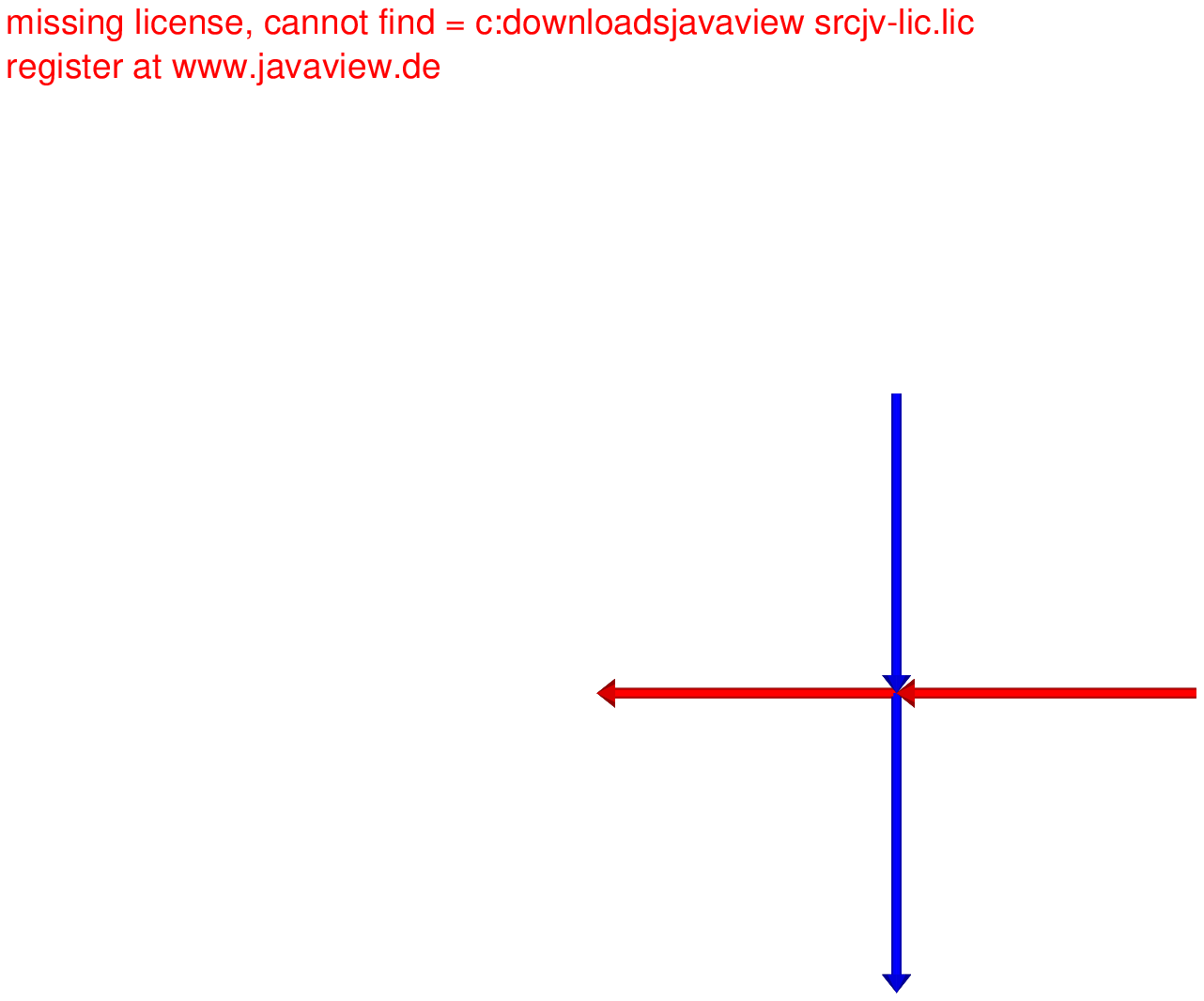}}\\
      $D_2 = so(4)$
    \end{center}
  \end{minipage}
  \begin{minipage}[t]{2.4cm}
    \begin{center}
      \leavevmode
      \resizebox{2.1cm}{!}{\includegraphics*[120,70][500,400]{adjoint_rep_A2}}\\
      $A_2 = su(3)$
    \end{center}
  \end{minipage}
  \begin{minipage}[t]{2.4cm}
    \begin{center}
      \leavevmode
      \resizebox{2.1cm}{!}{\includegraphics*[120,60][480,420]{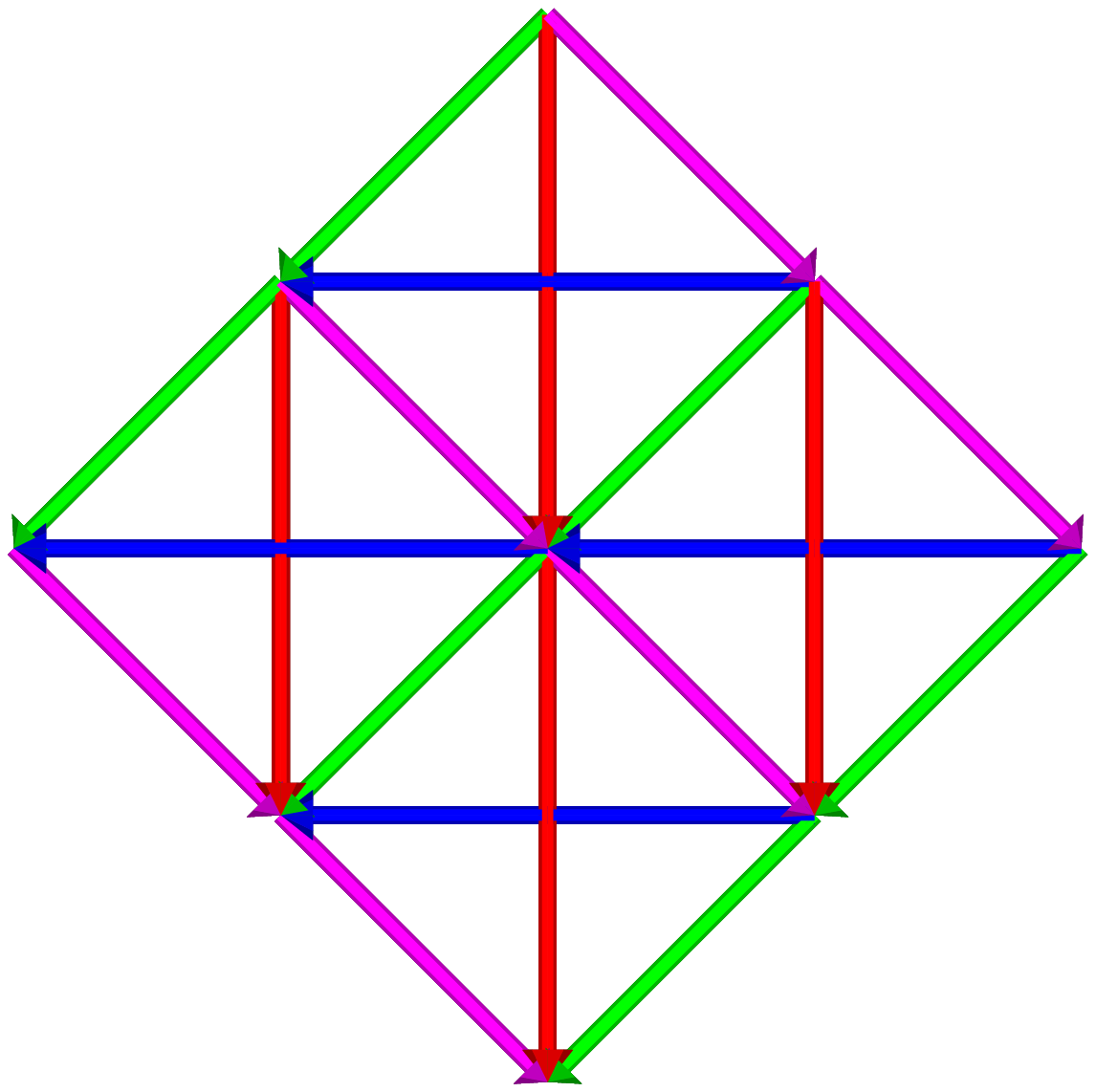}}\\
      $B_2=so(5)$
    \end{center}
  \end{minipage}
  \begin{minipage}[t]{2.4cm}
    \begin{center}
      \leavevmode
      \resizebox{2.1cm}{!}{\includegraphics*[230,170][375,315]{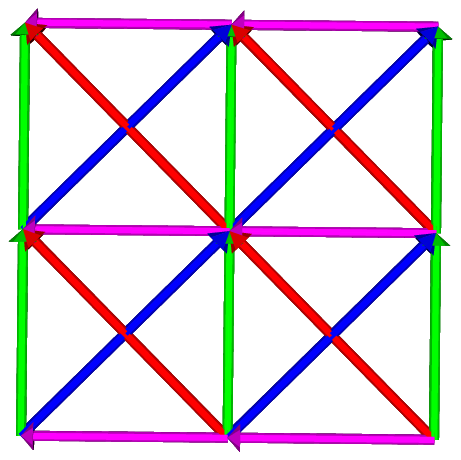}}\\
      $C_2=sp(2\cdot 2)$
    \end{center}
  \end{minipage}
  \begin{minipage}[t]{2.4cm}
    \begin{center}
      \leavevmode
      \resizebox{2.1cm}{!}{\includegraphics*[152,79][440,401]{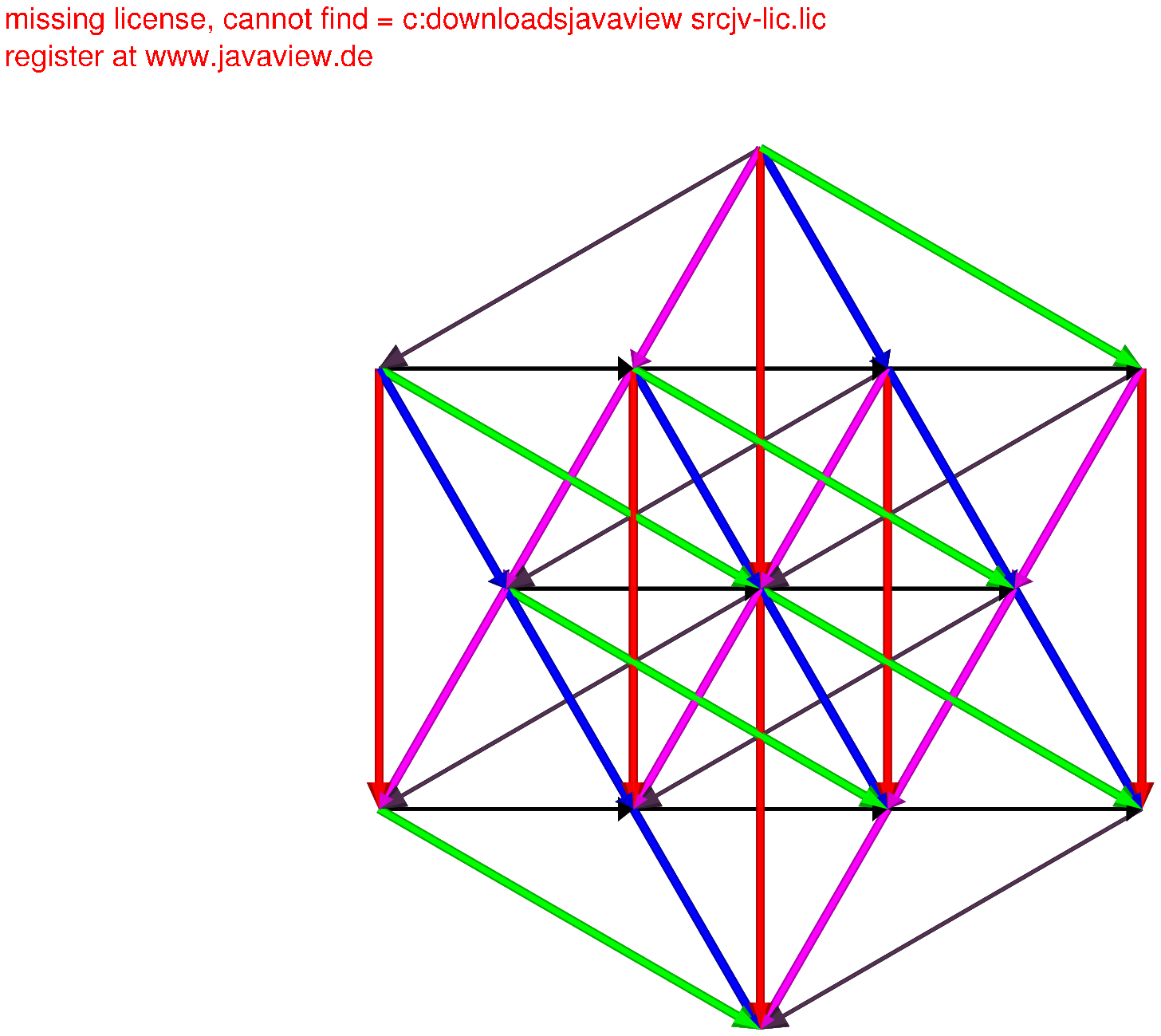}}\\
      $G_2 = Aut(\mathbb{O})$
    \end{center}
  \end{minipage}
  \caption{Root diagrams of simple rank ~$2$ algebras}
  \label{fig:rank_2_root_diagrams}
\end{center}
\end{figure}

The full set of roots have been produced once the Weyl reflections fail to produce additional roots.  The root diagrams are completed by connecting the tips of the roots ~$r^i$ and ~$r^j$ via the root ~$r^k$ precisely when ~$r^k + r^i = r^j$.  From the root diagrams in Figure ~\ref{fig:rank_2_root_diagrams}, it is clear that the dimension of~$A_2=su(3) \textrm{ is }8$ and the dimension of ~$G_2 \textrm{ is }14$, while both ~$B_2=so(5)$ and ~$C_2=sp(2\cdot 2)$ have dimension ~$10$.  Further, since the diagram of ~$B_2$ can be obtained via a rotation and rescaling of the root diagram for ~$C_2$, it is clear that ~$B_2$ and ~$C_2$ are isomorphic.

\subsubsection{Constructing Weight Diagrams from Dynkin Diagrams}

Root diagrams are a specific type of weight diagram.  While the states and roots can be identified with each other in a root diagram, this does not happen for general weight diagrams.  Weight diagrams are a collection of states, called weights, and the roots are used to move from one weight to another.  Although it has only one root diagram, an algebra has an infinite number of weight diagrams.

Like a root diagram, any one of an algebra's weight diagrams can be constructed entirely from its Dynkin diagram.  The Dynkin diagram records the relationship between the ~$l$ simple roots of a rank ~$l$ algebra.  Each simple root ~$r^1, \cdots, r^l$ is associated with a weight ~$w^1, \cdots, w^l$.  When all integer linear combinations of these weights are plotted in ~$\mathbb{R}^l$, they form an infinite lattice of possible weights.  This lattice contains every weight from every ~$d \times d$ ~$(d \ge 0)$ representation of the algebra.  A particular weight diagram is chosen by specifying a highest weight ~$W$ in this lattice.  Weyl reflections related to the weights ~$w^1, \cdots, w^l$ then define the boundary of weights ~$\overline{W}$ for that particular weight diagram, and root vectors are used to connect pairs of valid weights in ~$\overline{W}$ and within ~$\overline{W}$.\footnote{The distance between adjacent weights in the infinite lattice is less than or equal to the length of each root vector.  Root vectors do not always connect adjacent weights, but often skip over them.}  Choosing a new highest weight ~$W$ outside the original weight diagram's boundary ~$\overline{W}$ selects a new, larger, weight diagram, while choosing a new highest weight inside of ~$\overline{W}$ selects a new, smaller, weight diagram.

An equivalent method ~\cite{kass} constructs weight diagrams while avoiding the explicit use of Weyl reflections.  First, the  ~$l$ simple roots ~$r^1, \cdots, r^l$ are used to define the Cartan matrix ~$A$, whose components are ~$a_{ij} = 2 \frac{r^i \dotproduct r^j}{r^i \dotproduct r^i}$.  Equivalently, the components are:

\[
a_{ij} =\begin{cases} 2, &\text{if $i = j$;}\\
0, &\text{if $r^i$ and $r^j$ are orthogonal;}\\
-3, &\text{if the interior angle between $r^i$ and $r^j$ is $\frac{5\pi}{6}$, and $\sqrt{3}|r^i| = |r^j|$;}\\
-2, &\text{if the interior angle between roots $r^i$ and $r^j$ is $\frac{3\pi}{4}$, and $\sqrt{2}|r^i| = |r^j|$;}\\
-1, &\text{in all other cases;}\\
\end{cases}
\]

Cartan matrices are invertible.  Thus, the {\it fundamental weights} ~$w^1, \cdots, w^l$ defined by \[w^i = \sum_{k=1}^{l}(A^{-1})_{ki}r^k\] are linearly independent.  We create an infinite lattice ~$\mathbb{W} = \{ m_i w^i | m_i \in \mathbb{Z}\}$ of possible weights in ~$\mathbb{R}^l$, and then label each weight ~$W^i = m^i_j w^j \in \mathbb{W}$ using the ~$l$-tuple ~$M^i = \left[m^i_1, \cdots, m^i_l\right]$, called a {\it mark}.    We choose one of the infinite number of weight diagrams by specifying ~$l$ non-negative integers ~$m^0_1, \cdots, m^0_l$ for the mark of the highest weight ~$W^0$.

While the components of the Cartan matrix record the geometry of the simple roots, the components of each mark record the geometric configuration of the weights within the lattice.  Each weight ~$W^i$ is part of a shell of weights ~$\overline{W}$ which are equidistant from the origin.  
As discussed earlier, Weyl reflections associated with the fundamental weights can be used to find the weights in ~$\overline{W}$.  The diagram's entire set of valid weights can also be determined using the mark ~$M^i$ of each weight.   The positive integers ~$m^i_j \in M^i$ list the maximum number of times that the simple root ~$r^j$ can be subtracted from ~$W^i$ while keeping the new weights on or within the boundary ~$\overline{W}$.
Thus, weights ~$W^i - 1 r^j, \cdots, W^i - (m^i_j) r^j \textrm{ (no sum)}$ are valid weights occurring on or inside ~$\overline{W}$.  These new weights have marks ~$M^i - A_i^{\mathsf{T}}, M^i - 2 A_i^{\mathsf{T}}, \cdots, M^i - m^i_j A_i^{\mathsf{T}}$, where ~$A_i^{\mathsf{T}}$ is the transpose of the ~$i$th column of the Cartan matrix ~$A$.  Thus, given a weight diagram's highest weight ~$W^0$, this procedure selects all of the diagram's weights from the infinite lattice.  The diagram is completed by connecting any two weights ~$W^i$ and ~$W^k$ by the root ~$r^j$ whenever ~$W^i + r^j = W^k$.

We carry out this procedure for the algebra $B_3 = so(7)$, whose simple roots are

\[
\begin{array}{ccc}
 r^1 = \langle\sqrt{2},0,0\rangle, & r^2 = \langle-\sqrt{\frac{1}{2}}, -\sqrt{\frac{3}{2}}, 0\rangle, & r^3 = \langle0, \sqrt{\frac{2}{3}}, \sqrt{\frac{1}{3}}\rangle 
\end{array}
\]

\noindent
We produce the Cartan matrix ~$A$, find ~$A^{-1}$, and list the fundamental weights.

\hspace{-.5cm}
\begin{minipage}[t]{4.05cm}
\begin{center}
Cartan Matrix
\end{center}
$$ A = \left( \begin{array}{ccc}
2 & -1 & 0 \\
-1 & 2 & -1 \\
0 & -2 & 2 
\end{array}
\right)
$$
\end{minipage}
\hfill
\begin{minipage}[t]{4.05cm}
\begin{center}
Inverse Cartan Matrix
\end{center}
$$ A^{-1} = \left( \begin{array}{ccc}
1 & 1 & \frac{1}{2} \\
1 & 2 & 1 \\
1 & 2 & \frac{3}{2}
\end{array}
\right)
$$
\end{minipage}
\hfill
\begin{minipage}[t]{4.05cm}
\begin{center}
Fundamental Weights
\end{center}
$$
\begin{array}{c}
w^1 = 1r^1 + 1r^2 + 1r^3 \\
w^2 = 1r^1 + 2r^2 + 2r^3 \\
w^3 = \frac{1}{2}r^1 + 1r^2 + \frac{3}{2}r^3
\end{array}
$$
\end{minipage}

\noindent
Starting with the highest weight ~$W^0 = w^1$, whose mark is $\left[1,0,0\right]$, the above procedure generates the following weights:

\noindent
\begin{minipage}[htbp]{5in}
\begin{center}
\xymatrixcolsep{6pt}
\xymatrixrowsep{1pt}
\xymatrix@W=0pt{
 & \textrm{Marks} & \hfill & \textrm{Weights} \\
                    & \left[1, 0, 0\right]\ar[ddl]_{-r^1} &                    & W^0 = \left[1, 0, 0\right] \dotproduct \langle w^1, w^2, w^3 \rangle = 1r^1 + 1r^2 + 1r^3\\
  &  &  & \\
\left[-1, 1, 0\right]\ar[ddd]_{-r^2} &                    &                    & W^1 = \left[-1, 1, 0\right] \dotproduct \langle w^1, w^2, w^3 \rangle = 0r^1 + 1r^2 + 1r^3\\
  &  &  & \\
  &  &  & \\
\left[0, -1, 2\right]\ar[ddr]^{-r^3}&                    &  & W^2 = \left[0, -1, 2\right] \dotproduct \langle w^1, w^2, w^3 \rangle = 0r^1 + 0r^2 + 1r^3\\
  &  &  & \\
                    & \left[0, 0, 0\right]\ar[ddr]^{-r^3} &  & W^3 = \left[0, 0, 0\right] \dotproduct  \langle w^1, w^2, w^3 \rangle = 0r^1 + 0r^2 + 0r^3\\
  &  &  & \\
                    &                    & \left[0, 1, -2\right]\ar[ddd]_{-r^2} & W^4 = \left[0, 1, -2\right] \dotproduct  \langle w^1, w^2, w^3 \rangle = 0r^1 + 0r^2 -1r^3\\
  &  &  & \\
  &  &  & \\
                    &                    & \left[1, -1, 0\right]\ar[ddl]_{-r^1} &  W^5 = \left[1, -1, 0\right] \dotproduct  \langle w^1, w^2, w^3\rangle = 0r^1 - 1r^2 -1r^3\\
  &  &  & \\
                    & \left[-1, 0, 0\right]         &                   &  W^6 = \left[-1, 1, 0\right] \dotproduct  \langle w^1, w^2, w^3\rangle = -1r^1 - 1r^2 -1r^3
}
\end{center}
\end{minipage}

\begin{figure}[hbtp]
\begin{minipage}[hb]{7cm}
  \begin{center}
    \leavevmode
  \resizebox{3.5cm}{!}{\includegraphics*[100,90][338,320]{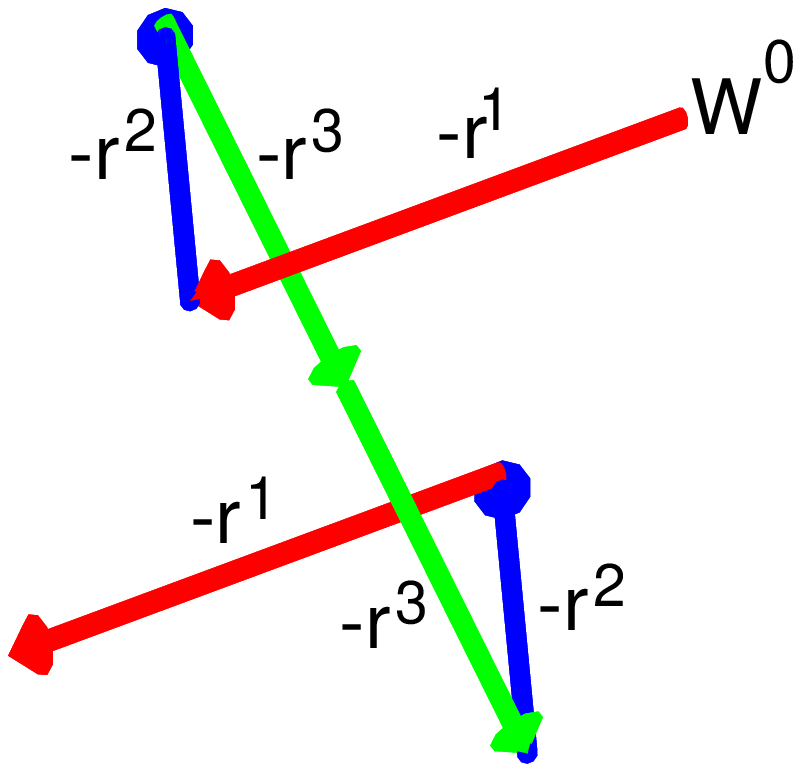}}
  \caption{\label{fig:weight_skeleton} \hbox{$B_3 = so(7)$ weight skeleton}}
  \end{center}
\end{minipage}
\hfill
\begin{minipage}[hb]{7cm}
  \begin{center}
    \leavevmode
  \resizebox{3.5cm}{!}{\includegraphics*[100,90][338,320]{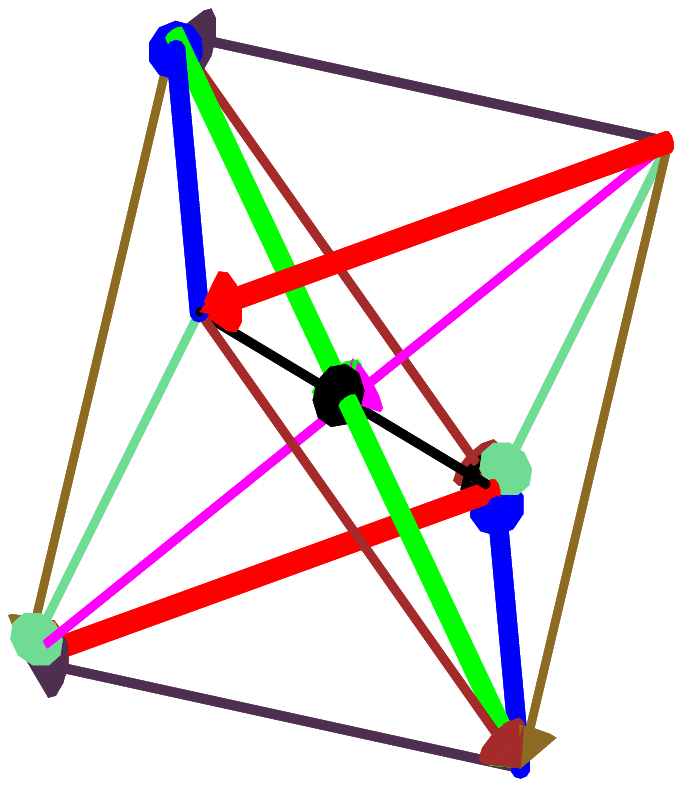}}
  \caption{\label{fig:weight_diagram}  \hbox{$B_3=so(7)$ weight diagram}}
  \end{center}
\end{minipage}
\end{figure}

We plot the weight ~$W^0, \cdots, W^6$ and the appropriate lowering roots ~$-r^1, \cdots, -r^3$ at each weight for this weight diagram of ~$B_3 = so(7)$ in Figure ~\ref{fig:weight_skeleton}.  The full set of roots are used to connect pairs of weights, giving the complete weight diagram of ~$B_3 = so(7)$ in Figure ~\ref{fig:weight_diagram}.

\subsubsection{Rank 3 Root and Weight Diagrams}

Using the procedures outlined above, we catalog the root and minimal weight diagrams for the rank ~$3$ algebras.  We did this by writing a Perl program to generate the weight and root diagram structures and return executable Maple code.  We then used Maple and Javaview to display the pictures.

\begin{figure}[htbp]
\begin{minipage}[t]{4.1cm}
  \begin{center}
    \leavevmode
  \resizebox{4.0cm}{!}{\includegraphics*[185,155][352,304]{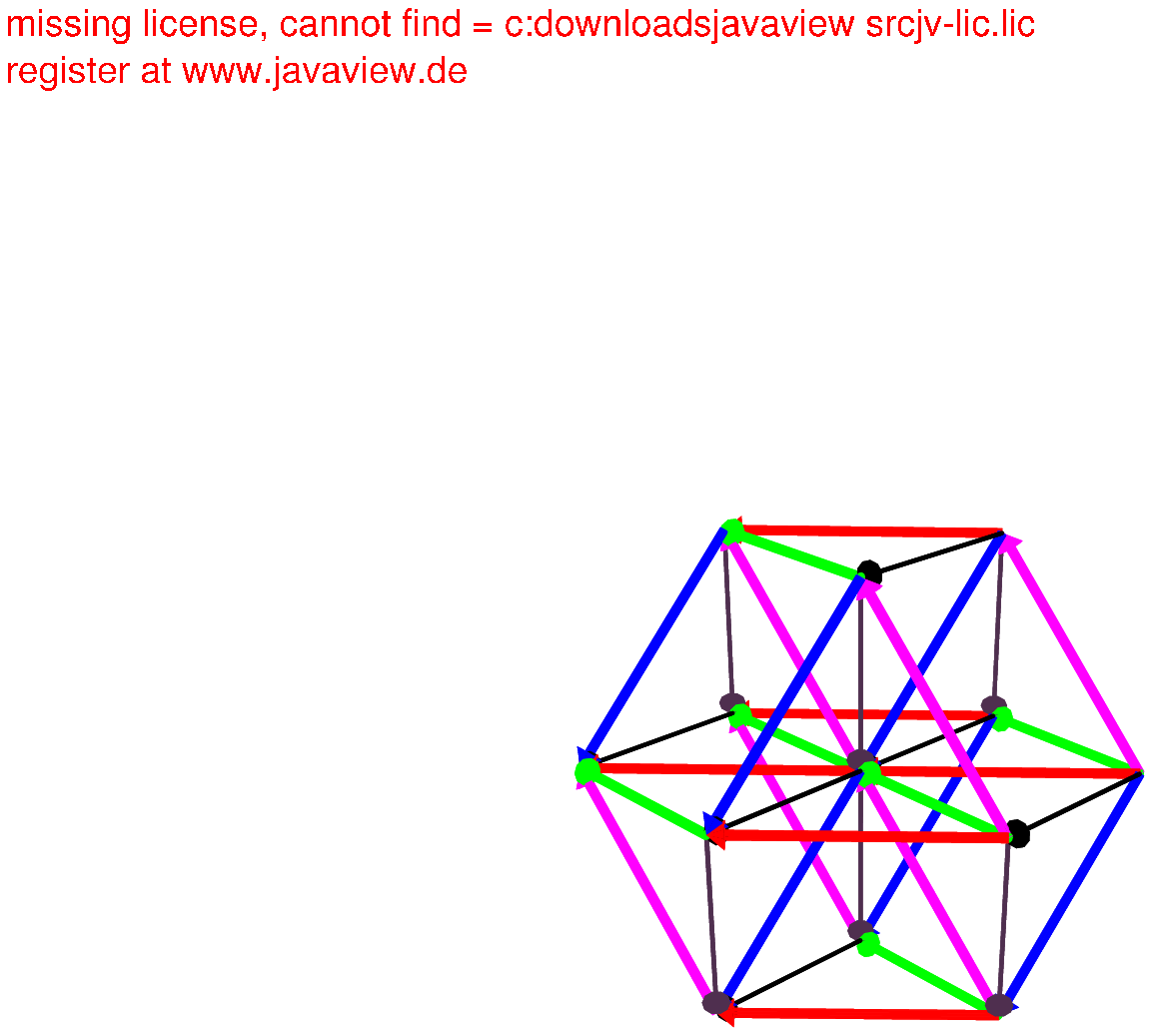}}\\
  $su(4) = A_3 = D_3 = so(6)$
  \end{center}
\end{minipage}
\hfill
\begin{minipage}[t]{4.1cm}
  \begin{center}
    \leavevmode
  \resizebox{4.0cm}{!}{\includegraphics*[204,160][370,308]{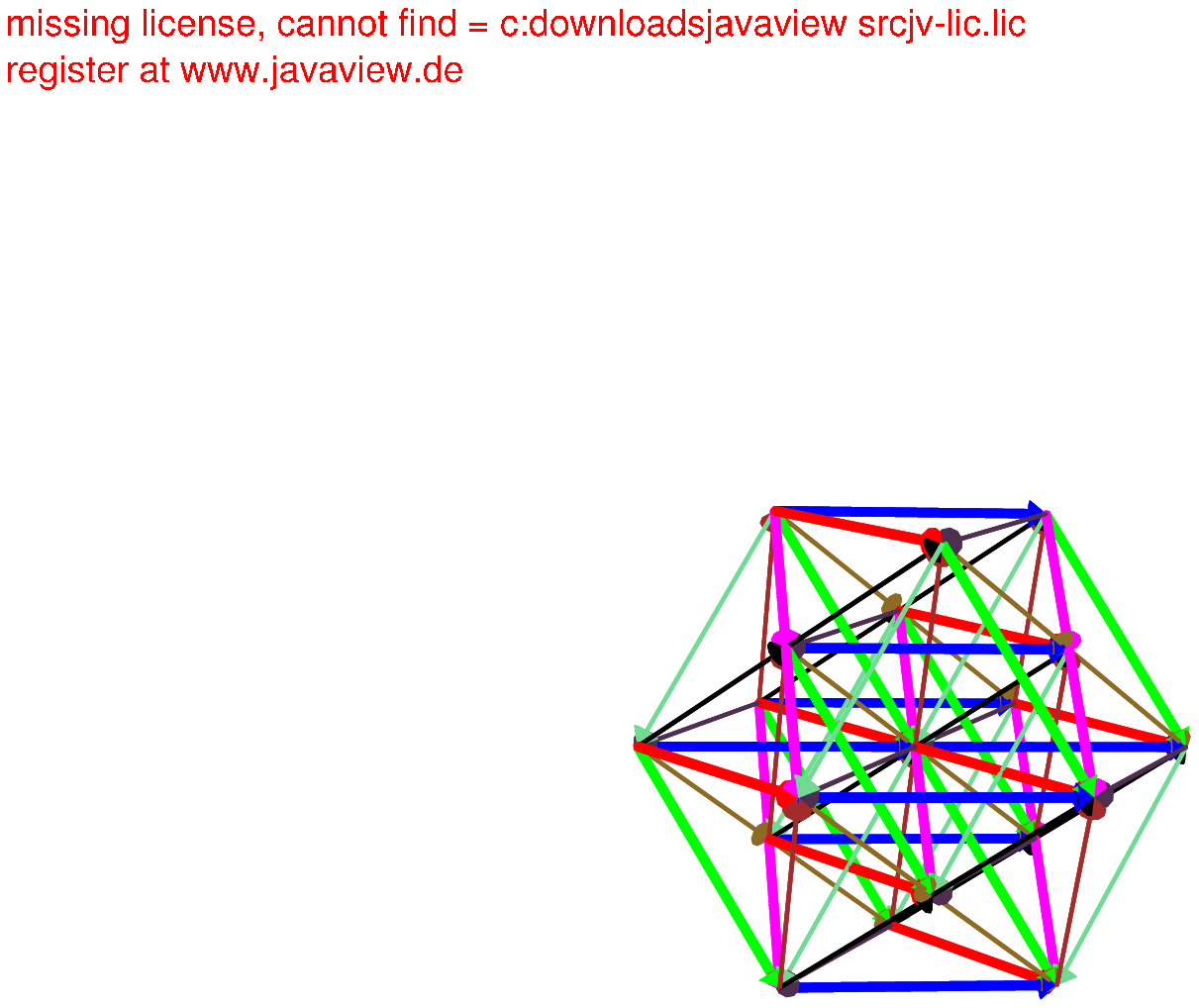}}\\
  $B_3 = so(7)$
  \end{center}
\end{minipage}
\hfill
\begin{minipage}[t]{4.1cm}
  \begin{center}
    \leavevmode
  \resizebox{4.0cm}{!}{\includegraphics*[159,109][440,362]{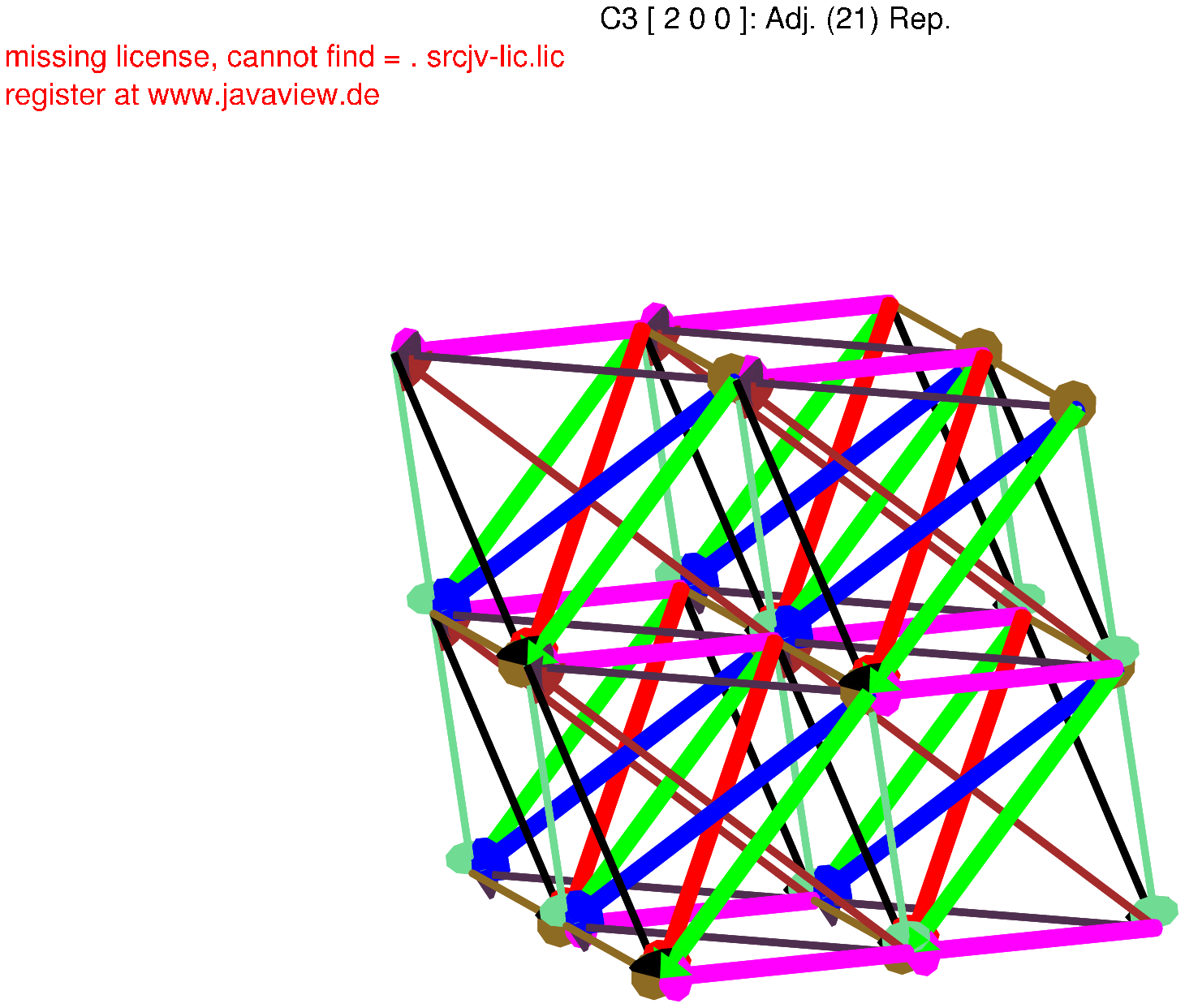}}\\
  $C_3 = sp(2\cdot 3)$
  \end{center}
\end{minipage}
\caption{Rank 3 root diagrams}
\label{fig:Rank_3_root_diagrams}
\end{figure}

The ~$3$-dimensional root diagrams of the rank ~$3$ algebras are given in Figure ~\ref{fig:Rank_3_root_diagrams}.  The root diagrams of ~$A_3=su(4)$ and $D_3=so(6)$ are identical.  Hence, $A_3 = D_3$.  These algebras have dimension ~$15$, as their root diagram contains ~$12$ non-zero states. The ~$B_3=so(7)$ and ~$C_3=sp(2\cdot 3)$ algebras both have dimension ~$21$, and their root diagrams each contain ~$18$ non-zero states.  Each of these algebras contain the ~$A_2=su(3)$ root diagram, which is a hexagon.  This can be seen at the center of each rank ~$3$ algebra by turning its root diagram in various orientations.   The modeling kit {\it ZOME} ~\cite{zome} can be used to construct the rank ~$3$ root diagrams.  The kit contains connectors of the right length and nodes with the correct configuration of connection angles to construct most of the rank ~$3$ root diagrams.\footnote{The root diagram for the algebra ~$C_3 = sp(2\cdot 3)$ can not be built to scale because it requires an angle not in the existing {\it ZOME} tools.}

\begin{figure}[htbp]
\begin{minipage}{4.1cm}
  \begin{center}
    \leavevmode
  \resizebox{3.4cm}{!}{\includegraphics*[242,130][370,240]{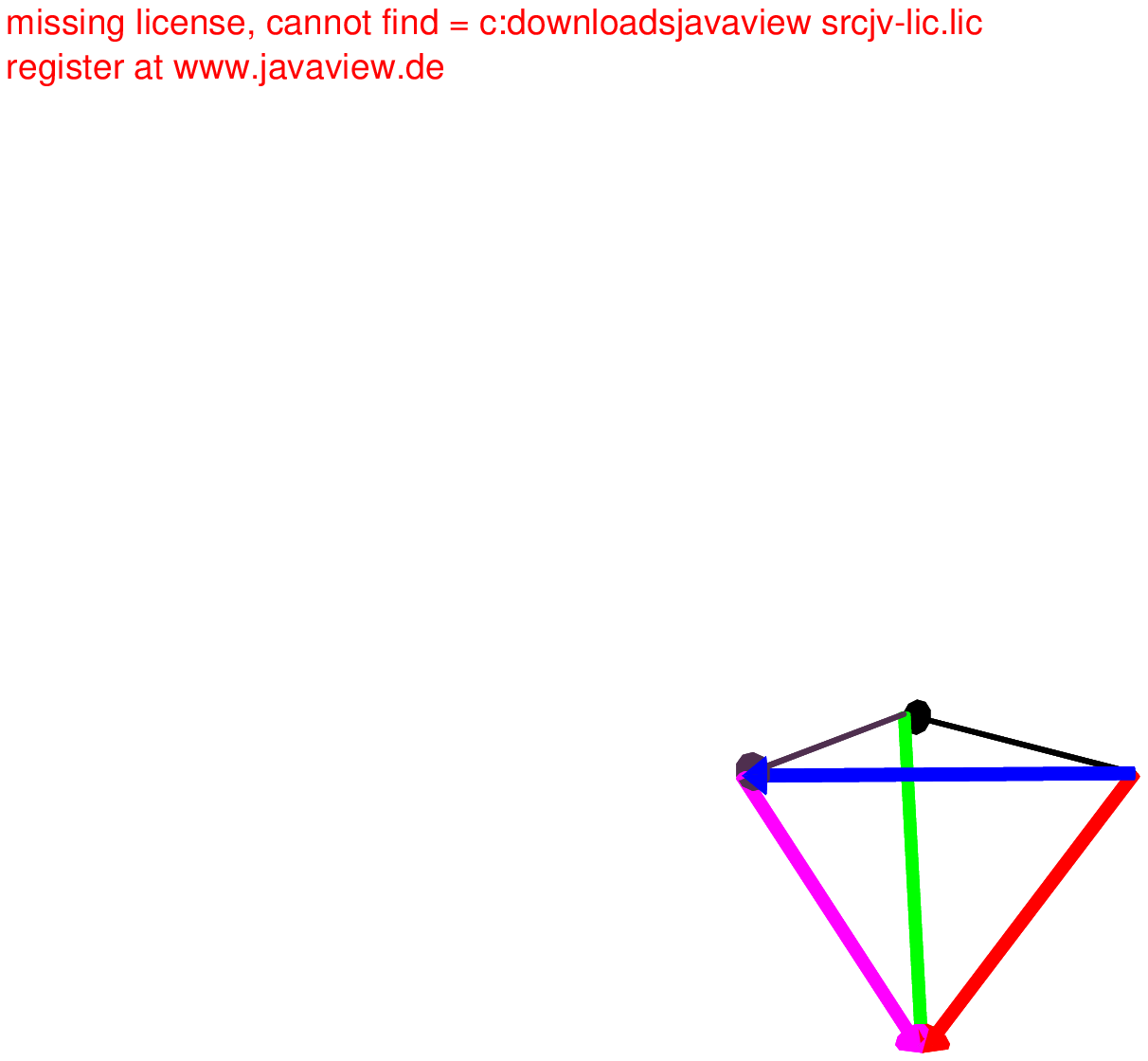}}\\
  $su(4) = A_3 = D_3 = so(6)$
  \end{center}
\end{minipage}
\hfill
\begin{minipage}{4.1cm}
  \begin{center}
    \leavevmode
  \resizebox{3.4cm}{!}{\includegraphics*[225,145][350,268]{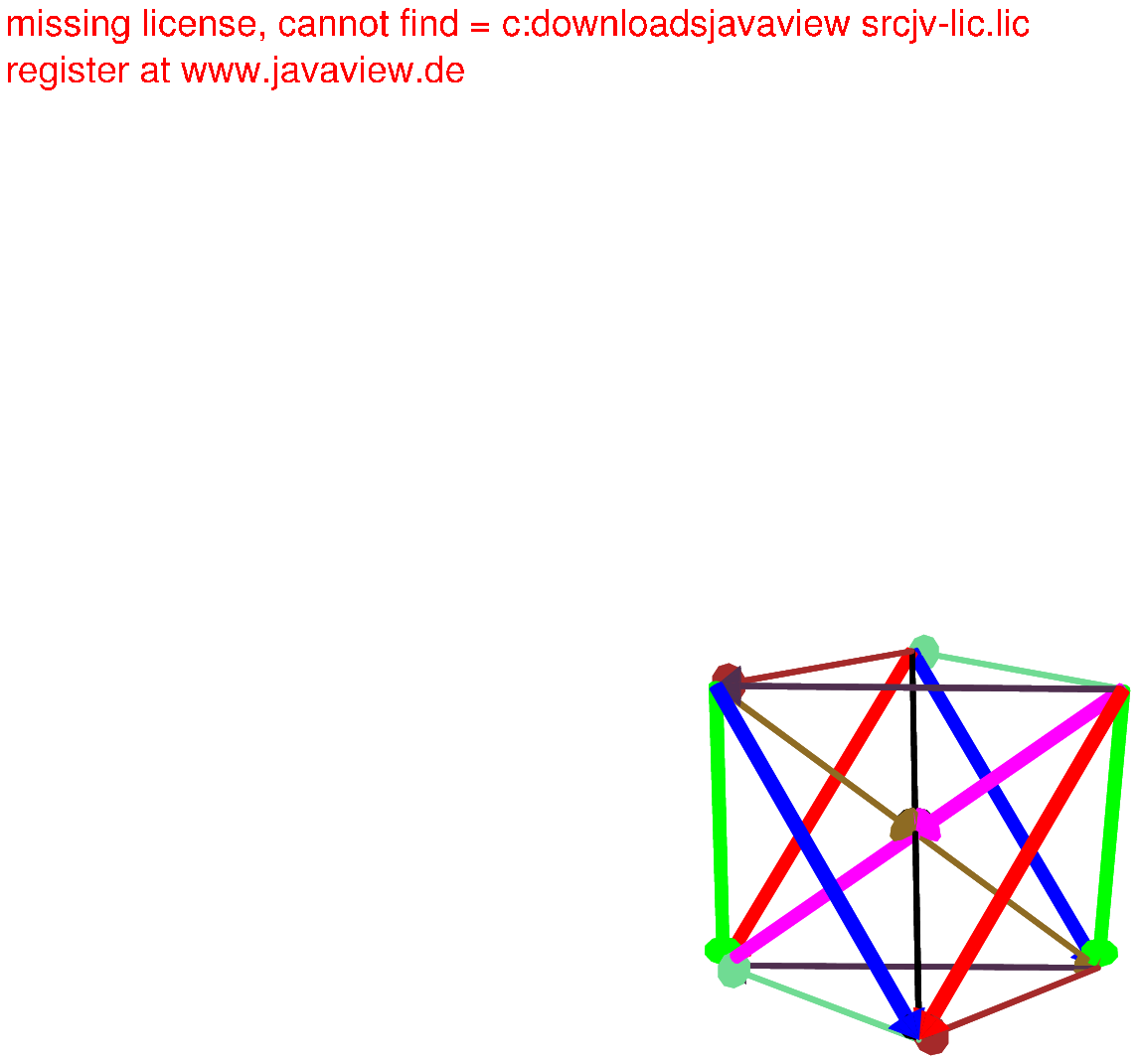}}\\
  $B_3 = so(7)$
  \end{center}
\end{minipage}
\hfill
\begin{minipage}{4.1cm}
  \begin{center}
    \leavevmode
  \resizebox{3.4cm}{!}{\includegraphics*[158,115][431,365]{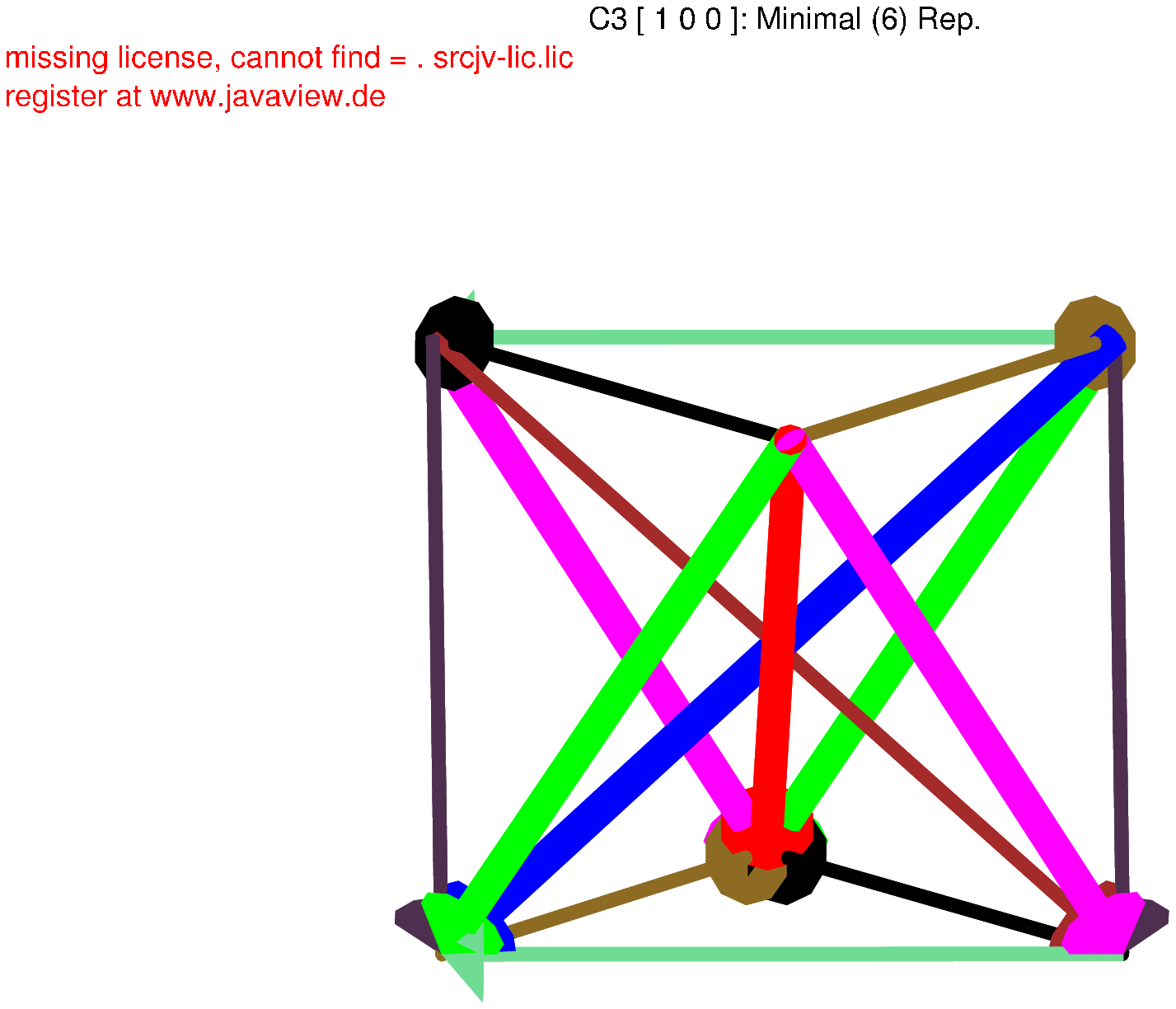}}\\
  $C_3=sp(2\cdot 3)$
  \end{center}
\end{minipage}
\caption{Rank 3 minimal weight diagrams}
\label{fig:rank_3_minimal_weight_diagrams}
\end{figure}

An algebra's minimal weight diagram has the fewest number of weights while still containing all of the roots.  Figure ~\ref{fig:rank_3_minimal_weight_diagrams} shows the  minimal weight diagrams for ~$A_3=D_3$, ~$B_3$, and ~$C_3$.  Each root occurs once in the diagram for ~$A_3=D_3$, while in the diagram for ~$B_3$ every root is used twice.  The minimal weight diagram for ~$C_3$ is centered about the origin, and the roots passing through the origin (colored red, blue, and brown) occur once, while the other roots occur twice.

\subsection{Subalgebras of Complex Lie Algebras}
\label{ch:Joma_Paper.subalgebras}

We have already noted that we can recognize $A_2 = su(3)$ as a subalgebra of $A_3 = su(4)$, $B_3 = so(7)$, and $C_3 = sp(2\cdot 3)$ by identifying its hexagonal root diagram in each of the larger root diagrams.  This section develops two further methods which help us identify subalgebras.

\subsubsection{Subalgebras of ~$B_3=so(7)$ using Root Diagrams}

Root and weight diagrams can be used to identify subalgebras.  For algebras of rank ~$l \le 3$, we can recognize subdiagrams corresponding to subalgebras by building the algebra's root and weight diagrams and rotating them in ~$\mathbb{R}^3$.

\begin{figure}[htbp]
\begin{center}
\begin{minipage}[t]{4.1cm}
  \begin{center}
    \leavevmode
\resizebox{4.0cm}{!}{\includegraphics*[150,175][320,310]{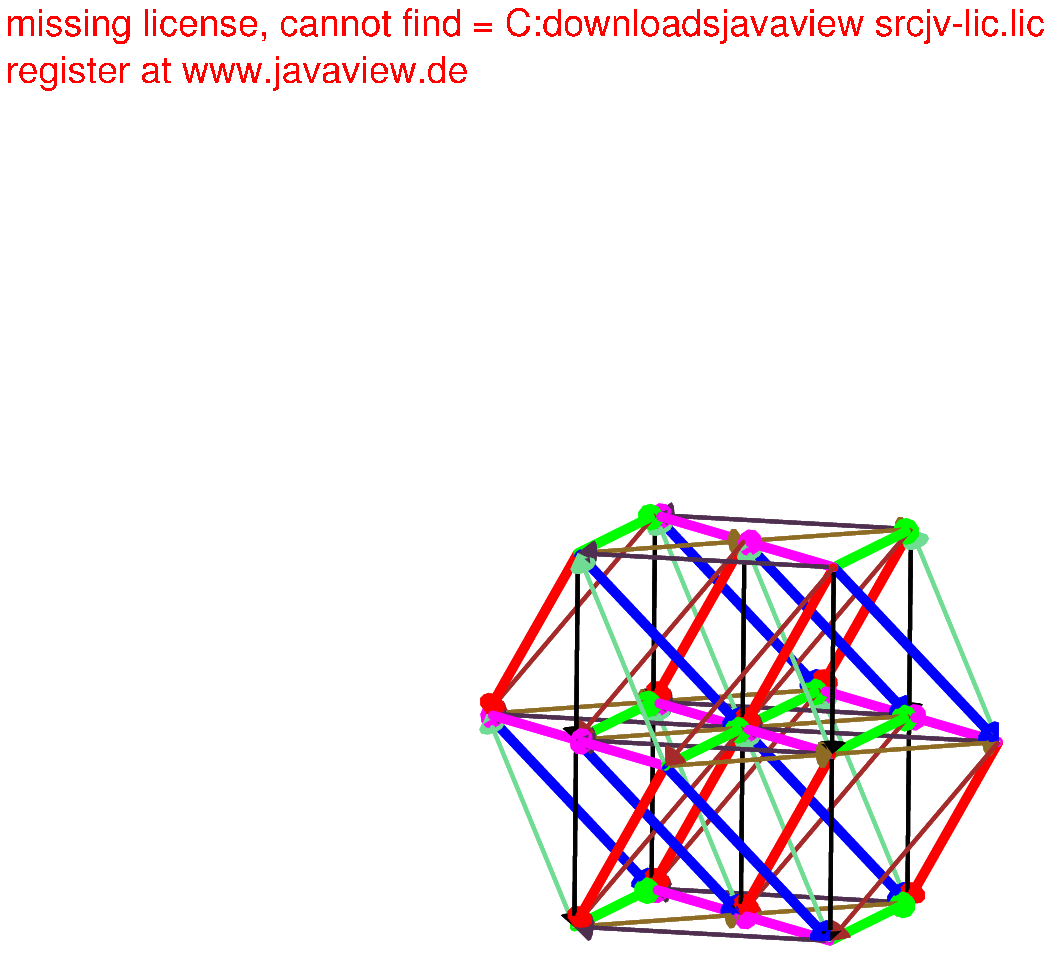}}
  \caption{ \hbox{$B_2 \subset B_3$}}
  \label{fig:B3_slice_to_B2}
  \end{center}
\end{minipage}
\hfill
\begin{minipage}[t]{4.1cm}
  \begin{center}
    \leavevmode
\resizebox{3.2cm}{!}{\includegraphics*[175,190][325,360]{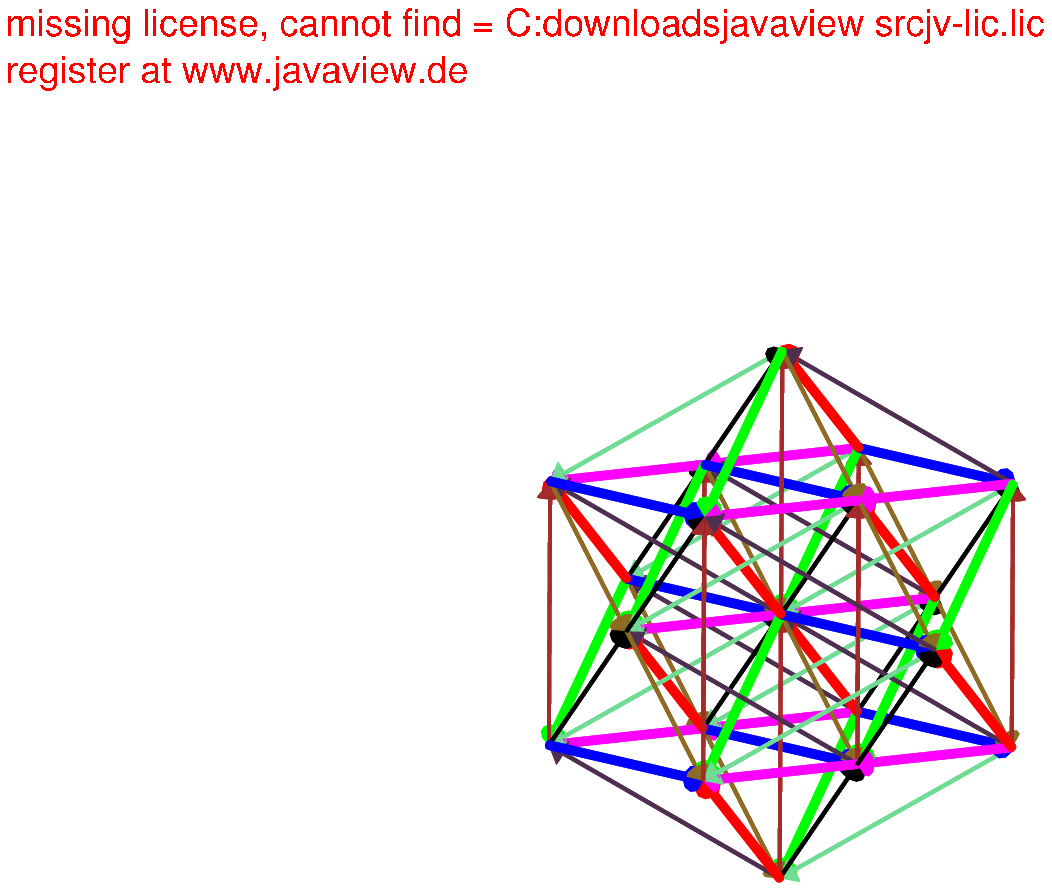}}
  \caption{ \hbox{$D_2 \subset B_3$}}
  \label{fig:B3_slice_to_D2}
  \end{center}
\end{minipage}
\hfill
\begin{minipage}[t]{4.9cm}
  \begin{center}
    \leavevmode
\resizebox{4.0cm}{!}{\includegraphics*[152,167][322,323]{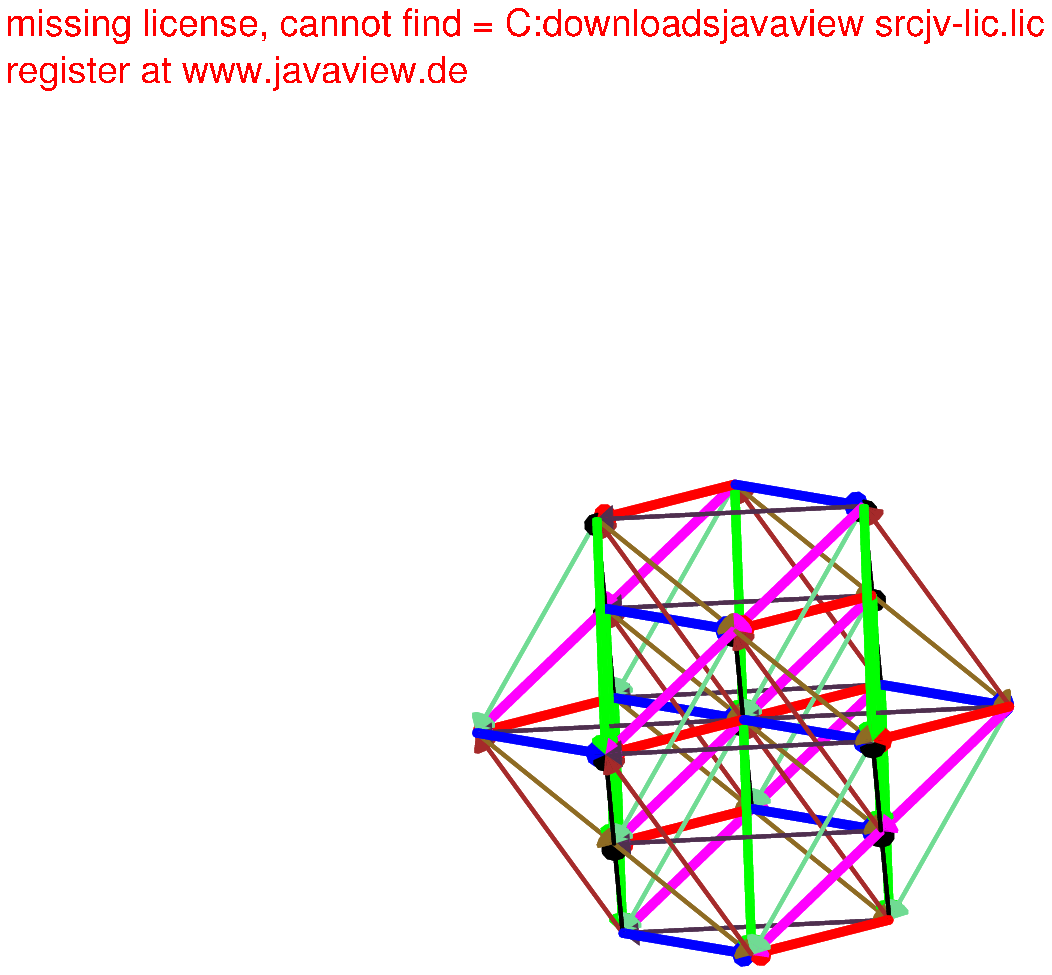}}
  \caption{ \hbox{$A_2 \subset D_3 \subset B_3$}}
  \label{fig:B3_slice_to_D3_and_A2}
  \end{center}
\end{minipage}
\\
\vspace{.5cm}
\begin{center}Identified Subalgebras of $B_3=so(7)$ using the Root Diagram\end{center}
\end{center}
\end{figure}

We illustrate the process of finding all subalgebras of ~$B_3=so(7)$ by using its root diagram.  When rotated to the position in Figure ~\ref{fig:B3_slice_to_B2}, the large horizontal ~$2$-dimensional rectangle through the origin (containing roots colored green and magenta) contains eight non-zero nodes.  This subdiagram is the root diagram of the ~$10$-dimensional algebra ~$B_2 = C_2$.  The smaller rectangular diagrams lying in parallel planes above and below this large rectangle contain all of the roots vectors of ~$B_2=C_2$.  These diagrams are minimal weight diagrams of ~$B_2 = C_2$.

Two additional rotations of the root diagram of ~$B_3=so(7)$ produce subalgebras.  In Figure ~\ref{fig:B3_slice_to_D2}, the horizontal plane through the origin contains only two orthogonal roots, which are colored blue and fuchsia.  These roots comprise the root diagram of ~$D_2=su(2)+su(2)$.  We have already identified the hexagonal ~$A_2=su(3)$ root diagram in the horizontal plane containing the origin in Figure ~\ref{fig:B3_slice_to_D3_and_A2}.  Each horizontal triangle above and below this plane is a minimal weight representation of ~$A_2$.  In addition, the subdiagram containing the bottom triangle, middle hexagon, and top triangle is the ~$15$-dimensional, rank ~$3$ algebra ~$A_3 = D_3$.  Thus, it is possible to identify both rank ~$l-1$ and rank ~$l$ subalgebras of a rank ~$l$ algebra.

Not all subalgebras of a rank ~$l$ algebra can be identified as subdiagrams of its root or weight diagram.  When extending a rank ~$l$ algebra to a rank ~$l+1$ algebra, each root ~$r^i = \langle \lambda_1^i, \cdots, \lambda_l^i \rangle$ is extended in ~$\mathbb{R}^{l+1}$ to the roots ~$r^{i_1}, \cdots, r^{i_m}$, where ~$r^{i_j} = \langle \lambda_1^i, \cdots, \lambda_l^i, \lambda_{l+1}^{i_j} \rangle$.  Here, ~$\lambda_{l+1}^{i_j}$ is one of ~$m \ge 1$ different eigenvalues values defined by the extension of the algebra.  Hence, although roots ~$r^{i_1}, \cdots, r^{i_m}$ are all distinct in ~$\mathbb{R}^{l+1}$, they are the same root when restricted to their first ~$l$ coordinates.  Thus, we can identify subalgebras of a rank ~$l+1$ algebra by projecting its root and weight diagrams along any direction.

\begin{figure}[hbtp]
\begin{center}
\begin{minipage}[t]{7cm}
  \begin{center}
    \leavevmode
  \resizebox{4cm}{!}{\includegraphics*[190,125][360,300]{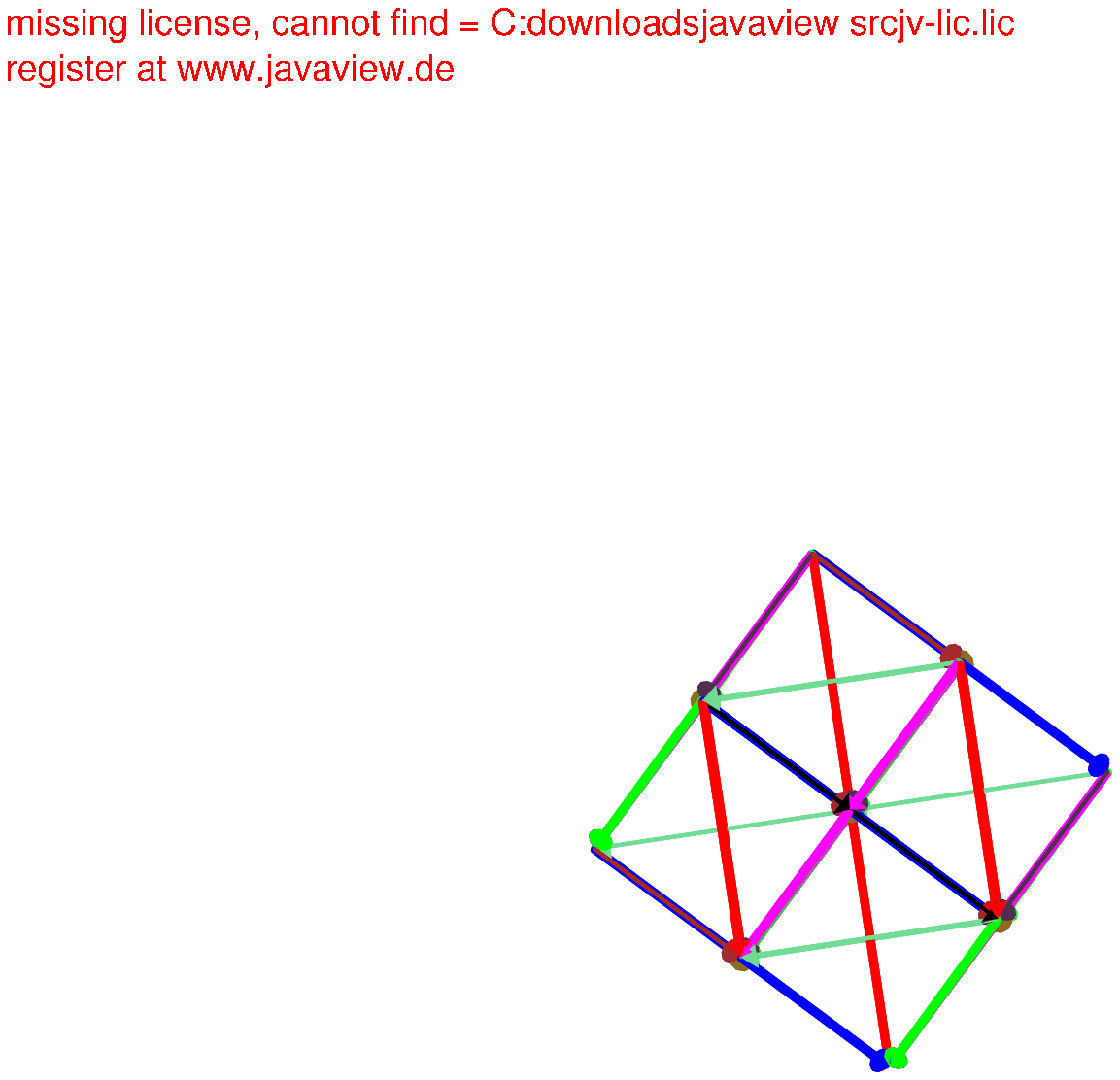}}
  \caption{\hbox{$B_2=so(5) \subset B_3=so(7)$}}
  \label{fig:B3_project_to_B2}
  \end{center}
\end{minipage}
\hfill
\hspace{-.5cm}
\begin{minipage}[t]{7cm}
  \begin{center}
    \leavevmode
\resizebox{3.5cm}{!}{\includegraphics*[250,180][400,350]{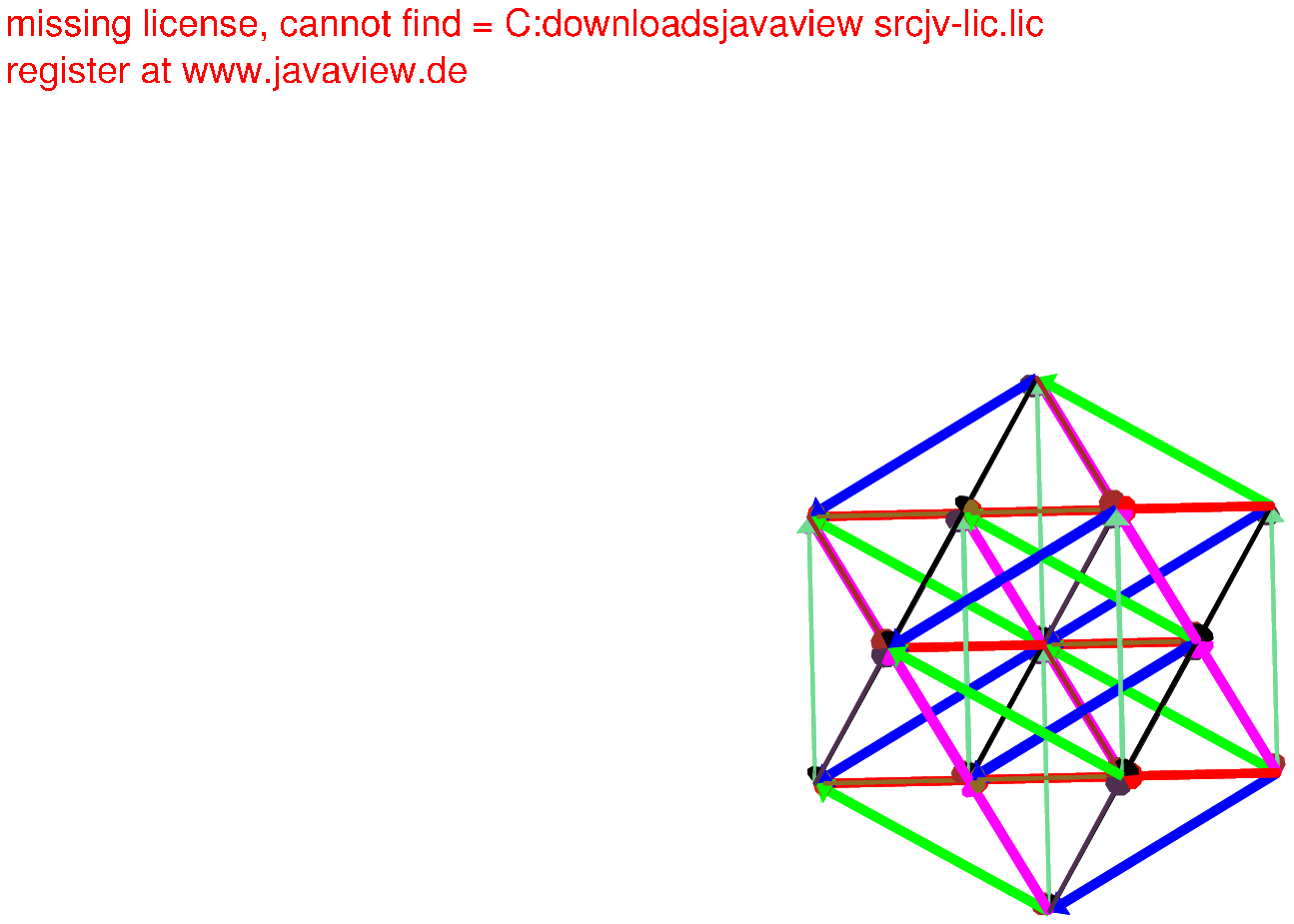}}
  \caption{$G_2 \subset B_3=so(7)$}
  \label{fig:B3_project_to_G2}
  \end{center}
\end{minipage}
\\
\vspace{.5cm}
\begin{center}Identify Subalgebras of ~$B_3=so(7)$ using Projections of its Root Diagram\end{center}
\end{center}
\end{figure}

Just as our eyes ``see'' objects in ~$\mathbb{R}^3$ by projecting them into ~$\mathbb{R}^2$, we can identify subalgebras of ~$B_3=so(7)$ by projecting its root and weight diagram into ~$\mathbb{R}^2$. In Figures ~\ref{fig:B3_project_to_B2} and ~\ref{fig:B3_project_to_G2}, we have rotated the root diagram of ~$B_3$ so that our eyes project one of the root vectors onto another root vector.  This method allows us to identify ~$B_2=so(5)$ and ~$G_2$ as subalgebras of ~$B_3$ by recognizing their root diagrams in Figures ~\ref{fig:B3_project_to_B2} and ~\ref{fig:B3_project_to_G2}, respectively.  Although we previously identified ~$B_2 \subset B_3$ by using a subdiagram confined to a plane, we could not identify ~$G_2$ as a subalgebra of ~$B_3$ using that method.

\subsubsection{Geometry of Slices and Projections}

The procedures in subsections ~$2.4$ and ~$2.5$ allow us to create root and weight diagrams for algebras of rank greater than ~$3$.  We now give precise definitions for the two processes we used in subsection ~$3.1$ to identify subalgebras of ~$B_3 = so(7)$.

We first generalize the method that allowed us to find subalgebras of ~$B_3=so(7)$ by identifying subdiagrams in specific cross-sections of its root diagram.  We refer to this method as ``slicing'', based on an analogy to slicing a loaf of bread.  A knife, making parallel cuts through the bread, creates several independent slices of the bread.  We use the same idea to slice an  algebra's root or weight diagram.  A set of ~$l-1$ linearly independent vectors ~$V = \{v^1, \cdots, v^{l-1}\}$ defines an ~$(l-1)$-dimensional hyperplane in ~$\mathbb{R}^l$.  Given ~$V$, a {\it slice}~$D_{\alpha}$ of an algebra's ~$l$-dimensional diagram ~$D$ is the subdiagram consisting of the vertex ~$\alpha$ and all root vectors and vertices in the hyperplane spanned by ~$V$ containing ~$\alpha$.  A {\it slicing} of a diagram~$D$ using ~$V$ separates ~$D$ into a finite set of disjoint slices.  The root vectors which connect vertices from two different slices are called {\it struts}.  We are interested in slices of ~$D$ which contain diagrams corresponding to the algebra's subalgebras.  Hence, the slices must contain root vectors of ~$D$, and in practice we choose ~$V$ to consist of integer linear combinations of simple roots.

A diagram's slices can tell us about its original structure.  When dealing with bread, we can obviously stack the slices on top of each other, in order, to recreate an image of the pre-cut loaf of bread - our mind removes the cuts made by the knife.     When dealing with root and weight diagrams, we do not have the benefit of using the shape of an ``outer crust'' to guide the stacking of the slices of the diagram.
Instead of severing the struts, we color them grey to make them less prominent.  This allows us to use the slices and struts to recreate the structure of the original root or weight diagram.

\begin{figure}[htbp]
\begin{center}
\begin{minipage}[t]{5in}
  \begin{center}
    \leavevmode
    \resizebox{8.5cm}{!}{\includegraphics*[40,160][555,310]{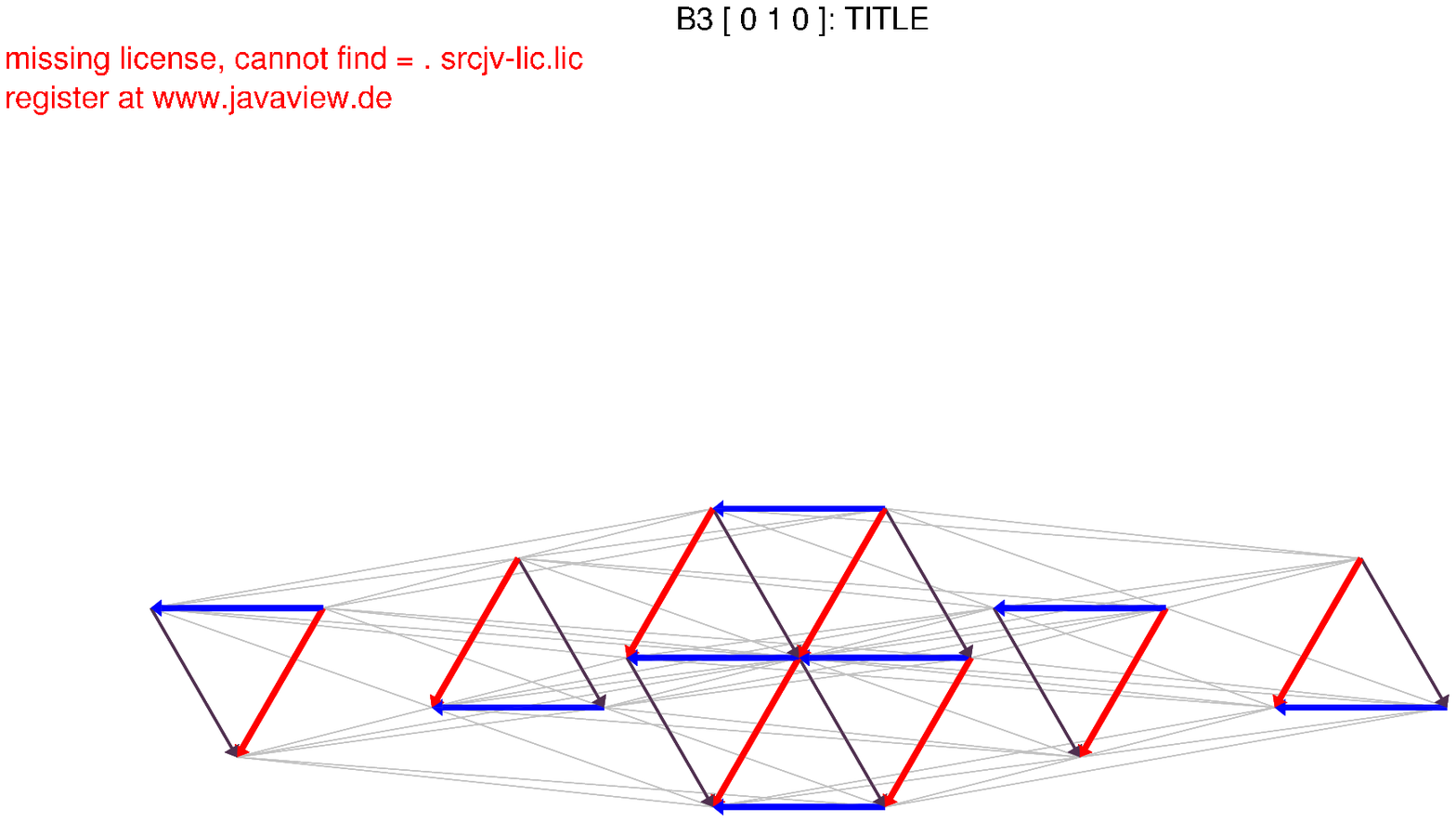}}
    \caption{  \label{fig:B3_slice_using_roots_1_and_2} Slicing of $B_3=so(7)$ using root $r^1$, colored red, and root $r^2$, colored blue.}
  \end{center}
\end{minipage}
\hfill
\begin{minipage}[t]{5in}
  \begin{center}
    \leavevmode
    \resizebox{11cm}{!}{\includegraphics*[86,188][505,290]{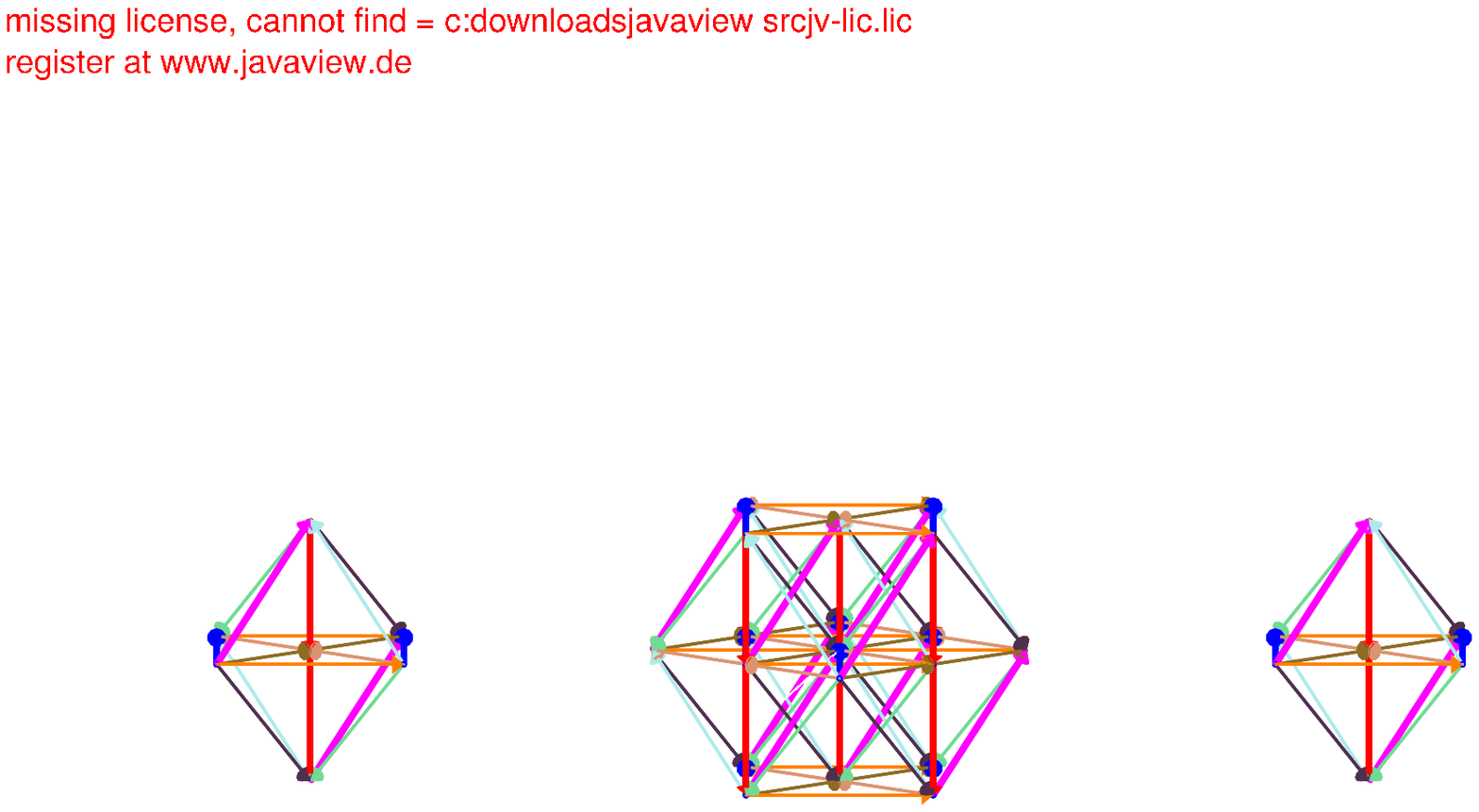}}
    \caption{\label{fig:B4_slice_showing_B3}  Slicing of ~$B_4=so(9)$ using roots ~$r^1$ (red), ~$r^2$ (green), and ~$r^3$ (blue).}
  \end{center}
\end{minipage}
\end{center}
\end{figure}

When the dimension of the diagram is greater than ~$3$, stacking slices on top of each other is not an effective means of recreating the root or weight diagram.  Instead, we lay the slices out along one direction, much as slices of bread are laid in order along a countertop to make sandwiches.  This allows us to display a~$3$-dimensional diagram in two dimensions, as we have shown for ~$B_3=so(7)$ in Figure ~\ref{fig:B3_slice_using_roots_1_and_2}, or a~$4$-dimensional diagram in three dimensions, as we have shown for the root diagram of ~$B_4 =so(9)$ in Figure ~\ref{fig:B4_slice_showing_B3}.  Of course, a~$5$-dimensional diagram can be displayed in three dimensions by first laying~$4$-dimensional slices along the ~$x$ axis, and then slicing each of these diagrams and spreading them along directions parallel to the ~$y$ axis.  This procedure generalizes easily to diagrams for rank six algebras, and can be modified to allow any compact ~$n$-dimensional diagram to be displayed in three dimensions.  

For the Lie algebras of rank ~$l = 6$ or less, we implement a slicing as follows.  We first build the algebra's root or weight diagram as described in subsection ~$2.5$.  Given three vectors ~$v^1$, ~$v^2$, and ~$v^3$ to define the slicing, we apply an orthonormal transformation so that the slices are contained within the first ~$3$~coordinates of ~$\mathbb{R}^l$.  We then use the projection (given here for ~$l=6$) ~$\mathbb{R}^6 \to \mathbb{R}^3: (x,y,z,u_1,u_2, u_3) \to (x + s_1\dotproduct u_1, y + s_2\dotproduct u_2, z + s_3 \dotproduct u_3)$, where ~$s_1$, ~$s_2$, and ~$s_3$ are separation factors used to separate the slices from each other when placed on our~$3$-dimensional countertop.\footnote{An equivalent projection ~$\mathbb{R}^6 \to \mathbb{R}^3: (x,y,z,u_1,u_2,u_3) \to (x+s_1 u_1 + s_2 u_2 + s_3 u_3, y, z)$ with~$s_1 > \eta_2 s_2 > \eta_3 s_3$ for sufficiently large ~$\eta_2$ and ~$\eta_3$ will string a ~$6$-dimensional diagram along one axis.  The ~$u_1$ coordinate will separate the different ~$5$-dimensional slices, with ~$s_1 u_1$ moving these different slices very far apart. The ~$u_2$ and ~$u_3$ coordinates will locally separate the ~$4$-dimensional and ~$3$-dimensional slices along the ~$x$ axis, but sufficiently small ~$s_2$ and ~$s_3$ will keep the subslices of one ~$5$-dimensional slice from interfering with another ~$5$-dimensional slice.  This method generalizes to~$n$-dimensional diagrams, but creates a very long string of the resulting diagrams.}
When helpful, we keep the grey colored struts in the sliced diagrams.\footnote{The large number of grey lines in a diagram can hide the important roots within each slice in addition to causing computational overload.}  While slicing preserves the length and direction of any root vector within a slice, laying the slices out along one direction obviously changes these characteristics for the grey struts.

The second method used to find subalgebras of ~$B_3=so(7)$ in subsection~$3.1$ involved projecting the ~$3$-dimensional diagram into a ~$2$-dimensional diagram.  A projection of an ~$l$-dimensional diagram to an ~$(l-1)$-dimensional diagram is accomplished by projecting along a direction specified by ~$p$.  To be useful, this projection must preserve the lengths of the roots, and the angles between them, when those roots are orthogonal to ~$p$.

Given a direction specified by a vector ~$p$, we create a linear transformation to change the basis from the standard basis ~$e^1, \cdots, e^l$ to a new orthonormal basis whose first basis vector is ~$\frac{p}{|p|}$.  This is accomplished by applying Gram-Schmidt orthonormalization to the ordered set of vectors ~$P = \{p, e^1, e^2, \cdots, e^l\}$, which is linearly dependent, and keeping the first ~$l$ non-zero vectors.  We then use a linear transformation to convert the standard basis to this new basis and apply it to the simple roots.  Finally, we throw away the first coordinate in the expression for each simple root.  This allows us to build the root or weight diagram following the procedure in subsection ~$2.5$, using the original simple roots to define the weights ~$W^0\textrm{,} \cdots \textrm{,} W^n$, which are then constructed using our projected simple roots.  It is faster computationally to apply the projection to the simple roots before building the diagram than to apply the projection to the entire diagram after it has been built.  As we can only display diagrams in three dimensions, when ~$l \ge 4$, we repeat this procedure ~$l-3$ times using ~$l-3$ projection directions ~$p^1, \cdots, p^{l-3}$.

The Gram-Schmidt process smoothly transforms a set of linearly independent vectors into a set of orthonormal vectors, and we choose our vectors in ~$P$ in a smooth way.  However, as ~$P$ is linearly dependent, our resulting change of basis transformation will not smoothly depend on ~$p$ if ~$p \in span(e^1, \cdots, e^{l-1})$.  Thus, we place the restriction that ~$|\frac{p}{|p|} \cdot e^l| > \epsilon$, for some small ~$\epsilon$.  In practice, we are usually interested in directions ~$p$ which are integer or half-integer linear combinations of the simple roots, and change ~$p$ to ~$p + 0.015e^1 + 0.015e^2 + \cdots + 0.015e^l$.  This assures that our projection smoothly depends upon ~$p$, or in the case ~$l \ge 4$, on ~$p^1, \cdots, p^{l-3}$.

By setting the separation factors ~$s_i = 0$, the slicing method can be used to produce another projection of the diagram ~$D$.  This choice for ~$s_i$ collapses the separate slices ~$D_{\alpha}$ onto one another, and centers them about the origin.  Given a hyperplane ~$V$, this {\it slice and collapse} method projects the slices along a direction perpendicular to ~$V$.  This provides less flexibility that the true projection method, which allows a projection along any direction ~$p$ when ~$\{p\} \cap V = 0$.  Nevertheless, by setting some ~$s_i = 0$, the slice and collapse method can produce useful projections of ~$5$-dimensional and ~$6$-dimensional diagram.

\begin{figure}[htbp]
\begin{center}
\begin{minipage}[t]{6cm}
  \begin{center}
    \leavevmode
    \resizebox{!}{4cm}{\includegraphics*[124,87][375,369]{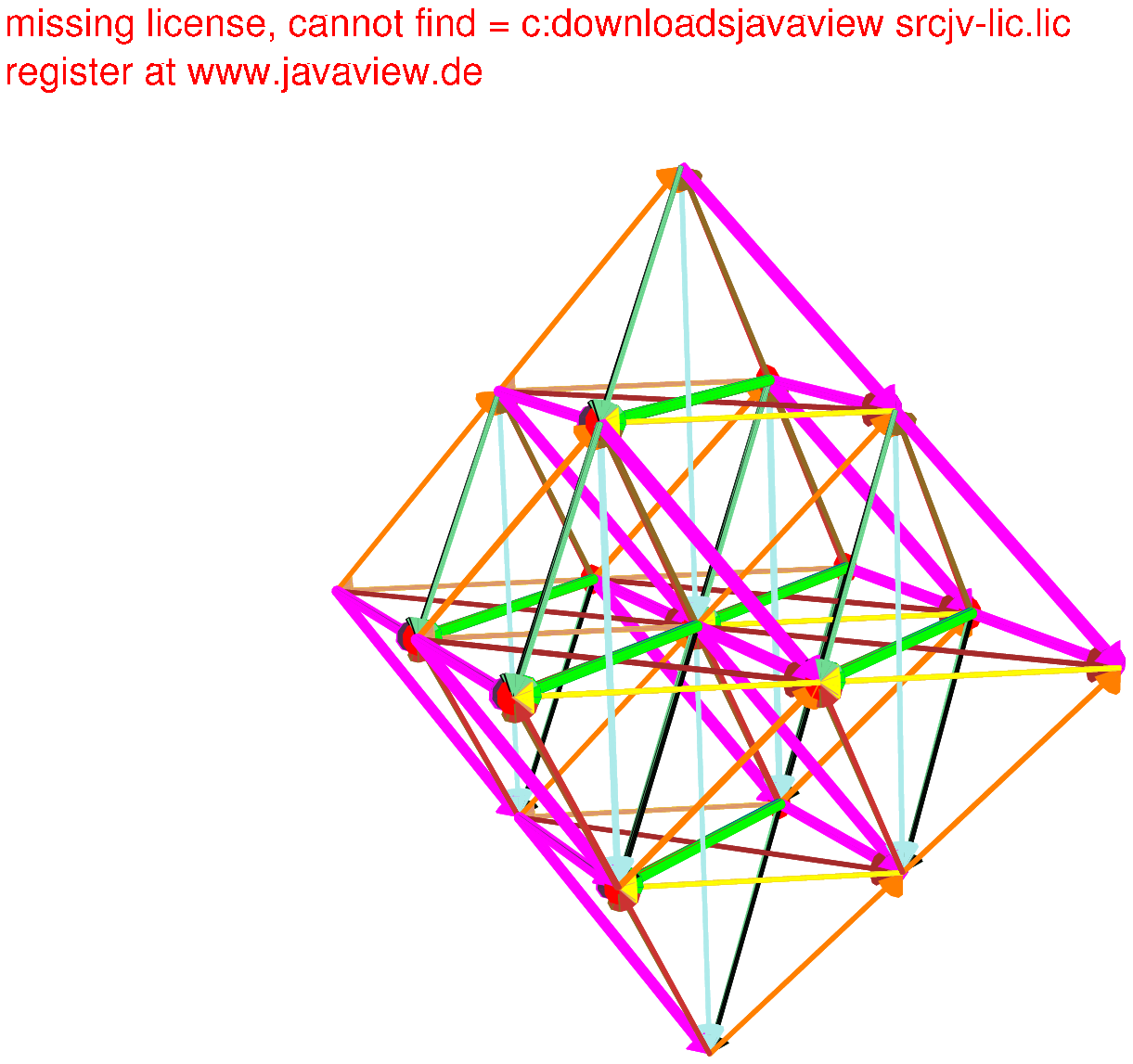}}
    \caption{\label{fig:C4_project_along_root_1} Projecting the ~$C_4=sp(2\cdot 4)$ root diagram along simple root ~$r^1$}
  \end{center}
\end{minipage}\hfill
\begin{minipage}[t]{6cm}
  \begin{center}
    \leavevmode
    \resizebox{!}{4cm}{\includegraphics*[164,84][394,327]{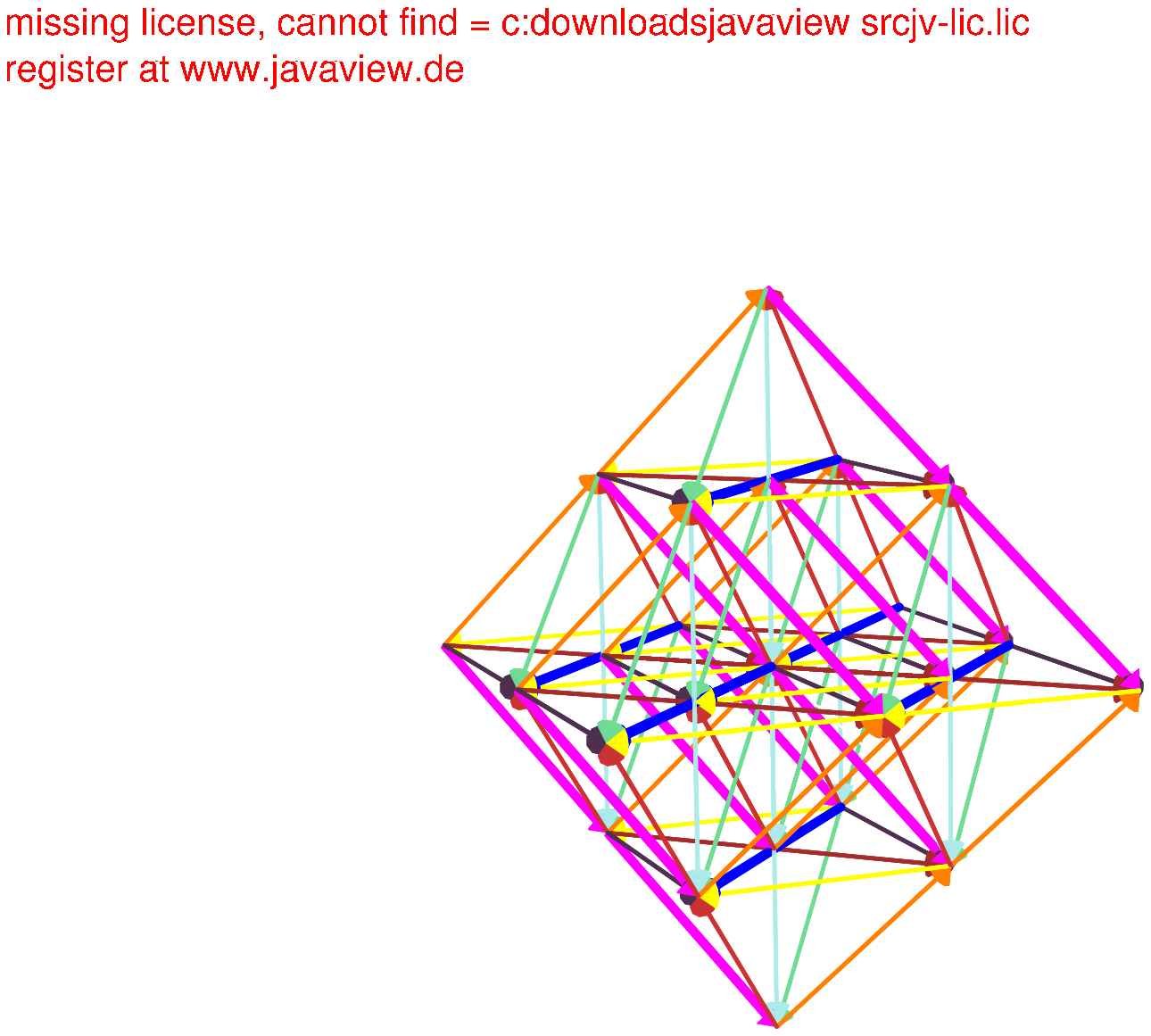}}
  \caption{\label{fig:C4_slice_using_roots_2_3_4}  Collapsing the slices of ~$C_4=sp(2\cdot 4)$ defined by simple roots ~$r^2$, ~$r^3$, and ~$r^4$ onto the origin.}
  \end{center}
\end{minipage}
\end{center}
\end{figure}

The slight difference between the projection method and slice and collapse method is illustrated in Figure ~\ref{fig:C4_project_along_root_1} and Figure~\ref{fig:C4_slice_using_roots_2_3_4}.  Figure ~\ref{fig:C4_project_along_root_1} projects the root diagram of ~$C_4=sp(2\cdot 4)$ along the simple root ~$r^1$.  The result is the root diagram of ~$C_3 = sp(2\cdot 3)$.  Figure ~\ref{fig:C4_slice_using_roots_2_3_4} collapses the slices of the root diagram of ~$C_4 = sp(2\cdot 4)$, defined using the simple roots ~$r^2, r^3$, and ~$r^4$, onto the origin.  While this diagram contains the ~$C_3$ root diagram, consisting of two large triangles on either side of a large hexagon, it also contains smaller triangles on either side of the hexagon.  These small triangles are part of an octahedron, which is one of the original slices of ~$C_4$.  In Figure ~\ref{fig:C4_slice_using_roots_2_3_4}, the octahedron is actually disjoint from the ~$C_3$ root diagram.\footnote{While the root vectors in the octahedron and the ~$C_3$ root diagram appear to intersect, they do not terminate or start from any common vertex.}  However, in the true projection, in Figure ~\ref{fig:C4_project_along_root_1},  the octahedron is placed by the projection either above or below the origin.  Whenever the slice and project method produces overlapping disjoint diagrams, a better projection can be obtained by translating the collapsed slice away from the origin.

\subsubsection{Identifying Subalgebras using Slicings and Projections}

Given an algebra ~$g$, the slicing and projection techniques described above produce subdiagrams of the algebra's root and weight diagrams.  We then compare these subdiagrams to a list of known Lie algebra diagrams.  If the subdiagram's vertices and root configuration exactly matches the configuration of a diagram for the Lie algebra ~$g^\prime$, then ~$g^\prime$ is a subalgebra of ~$g$.

While a projection along one of the diagram's root vectors will only allow an identification of subalgebras of rank ~$l-1$, slicings of a rank ~$l$ algebra's root diagram allow identifications of subalgebras of rank ~$l$ and ~$l-1$.  When successfully slicing root diagrams, the middle slice will be the root diagram of a rank ~$l-1$ subalgebra, and other slices will be different weight diagrams for that subalgebra~\cite{danielsson}.  We use a subdiagram consisting of some slices to identify a rank ~$l$ subalgebra.  In addition, that subdiagram must contain the original diagram's highest weight.  For instance, the short roots in the diagram of ~$C_3 = sp(2\cdot 3)$ look like they form the ~$A_3 = D_3$ algebra in Figure ~\ref{fig:Rank_3_root_diagrams}.  However, the highest weight of the ~$C_3$ diagram is the furthest away from the origin, at the tip of the longest root extending from the origin.  This weight is outside any embedding of the ~$A_3 = D_3$ root diagram into the root diagram of ~$C_3$.  In fact, ~$A_3 = D_3$ is not a subalgebra of ~$C_3$, as there are two operators $g_1$ and $g_2$ corresponding to the short roots whose commutator is an operator corresponding to one of the longer roots.  These rules are sufficient to identify rank ~$l$ and rank ~$l-1$ subalgebras of a rank ~$l$ algebra.

These methods only allow us to find subalgebras of the complex simple Lie algebras.  For the case of real subalgebras of real Lie algebras, these techniques can only indicate which containments are not possible (i.e. a real form of ~$C_3$ can not contain a real form of ~$A_3$).  See ~\cite{aaron_thesis} for information regarding the real Lie subalgebras of ~$E_6$.

\subsection{Applications to Algebras of Dimension Greater than~$3$}
\label{ch:Joma_Paper.subalgebras_greater_than_3}

We now show how these two techniques can be applied to find subalgebras of rank ~$l$ algebras, for ~$l \ge 4$.  We begin by using slices and projections to find subalgebras of the exceptional Lie algebra ~$F_4$.  We then show how to apply these techniques to algebras of higher rank.

\subsubsection{Subalgebras of ~$F_4$ using Slices}

\begin{figure}[htbp]
  \begin{center}
    \leavevmode
    \xymatrixcolsep{10pt}
    \xymatrix@M=0pt@W=0pt{
*+=[o][F-]{\hspace{.25cm} } \ar@{-}[rr]&  & *+=[o][F-]{\hspace{.25cm} } \ar@2{-}[rr] & \rangle & *+=[o][F-]{\hspace{.25cm} } \ar@{-}[rr] & \hspace{.25cm} & *+=[o][F-]{\hspace{.25cm} } 
    }
    \caption{$F_4$ Dynkin diagram}
    \label{fig:F4_dynkin_diagram}
  \end{center}
\end{figure}
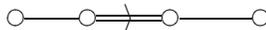

We apply the slice and projection techniques to the ~$52$-dimensional exceptional Lie algebra ~$F_4$, whose Dynkin diagram is shown in Figure ~\ref{fig:F4_dynkin_diagram}.  We number the nodes ~$1$ through ~$4$, from left to right, and use this numbering to label the simple roots ~$r^1, \dotsc, r^4$.  Thus, the magnitude of ~$r^1$ and ~$r^2$ is greater than the magnitude of ~$r^3$ and ~$r^4$.  We color these simple roots magenta ($r^1$), red ($r^2$), blue ($r^3$), and green ($r^4$).

\begin{figure}[htb]
\begin{center}
\begin{minipage}[t]{12cm}
  \begin{center}
    \leavevmode
    \resizebox{11cm}{!}{\includegraphics*[0,0][411,145]{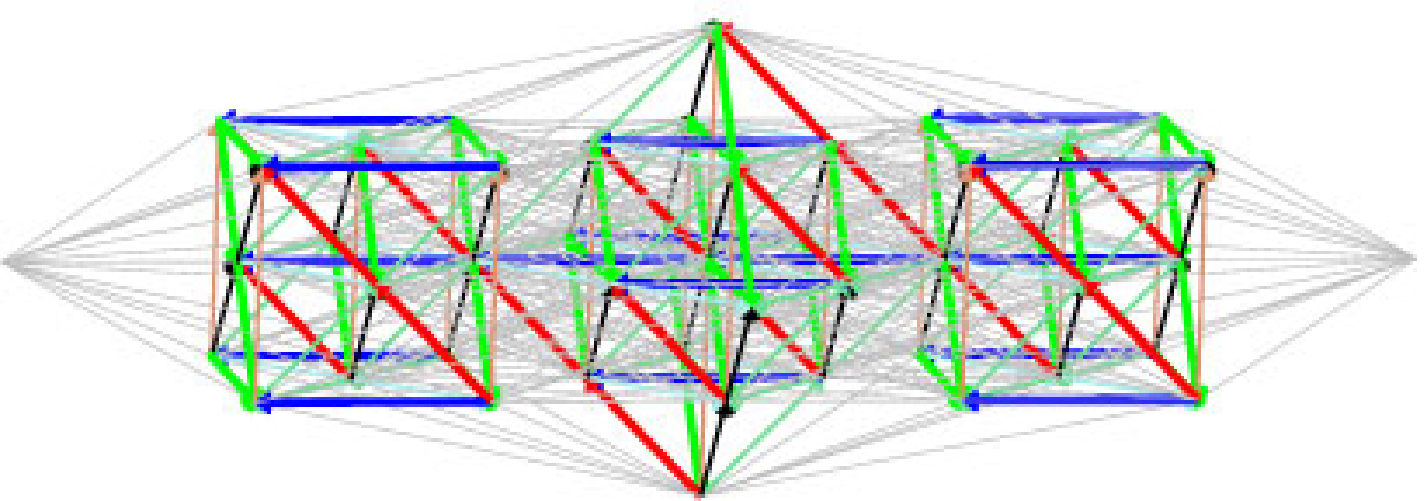}}
    \caption{\label{fig:F4_slice_showing_C3_greylines}  Slicing of ~$F_4$ using roots ~$r^2$ (red),  ~$r^3$ (blue), and ~$r^4$ (green). Grey colored struts connect vertices from different slices.}
  \end{center}
\end{minipage}\hfill
\begin{minipage}[t]{12cm}
  \begin{center}
    \leavevmode
    \resizebox{11cm}{!}{\includegraphics*[40,160][545,320]{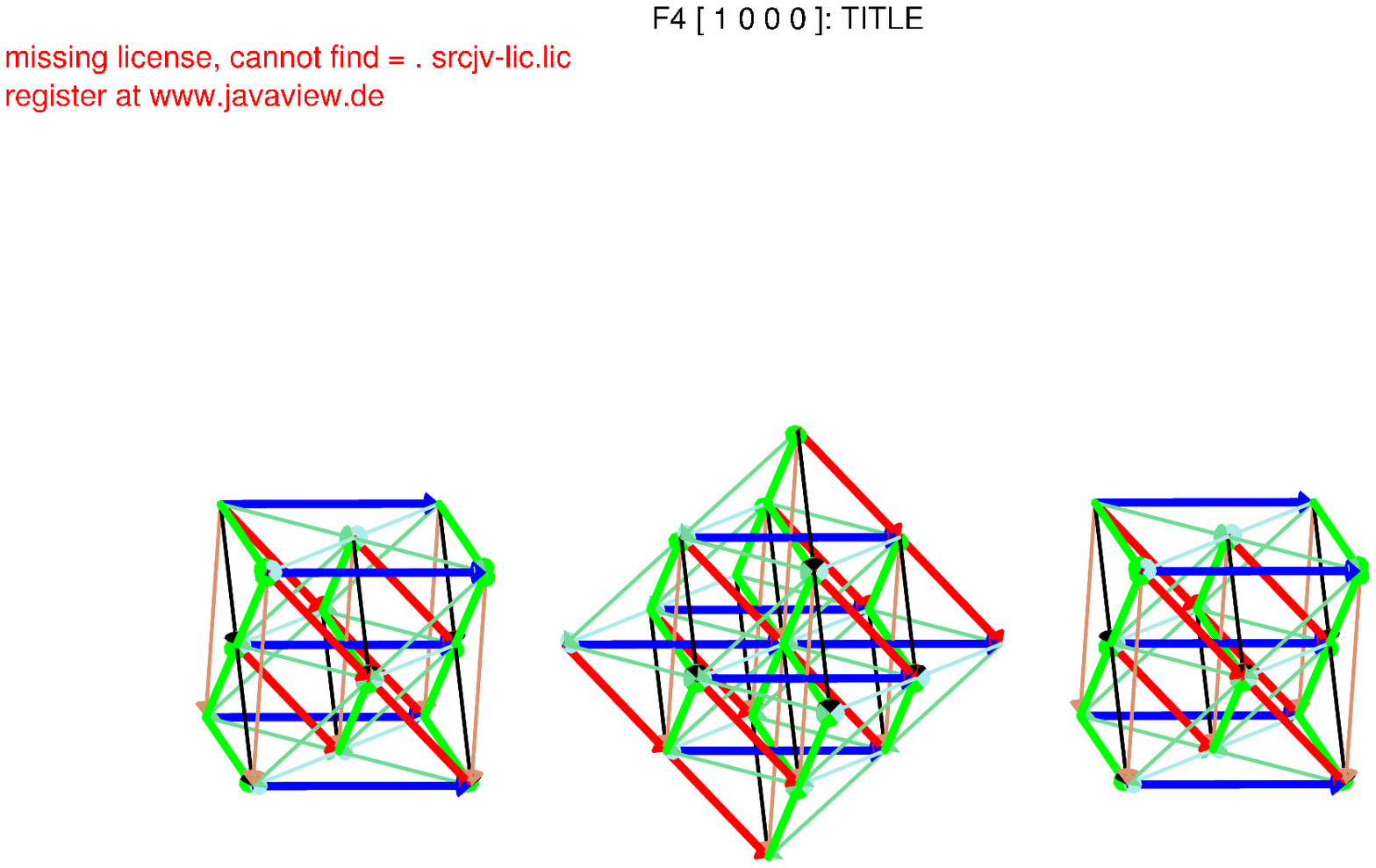}}
    \caption[Slicing of ~$F_4$ using roots ~$r^2$, ~$r^3$, and ~$r^4$. Eliminating the struts shows~$C_3=sp(2\cdot 3)~\subset~F_4$]{\label{fig:F4_slice_showing_C3}  \hbox{Slicing of ~$F_4$ using roots ~$r^2$, ~$r^3$, and ~$r^4$.} \hbox{Eliminating the struts shows~$C_3=sp(2\cdot 3)~\subset~F_4$}}
  \end{center}
\end{minipage}
\end{center}
\end{figure}

We consider the slicing of ~$F_4$ defined using roots ~$r^2$, ~$r^3$, and ~$r^4$.  Laying the slices along the ~$x$ axis, the large number of grey struts in the resulting diagram, Figure ~\ref{fig:F4_slice_showing_C3_greylines}, makes it difficult to observe the underlying structure of each slice, and so they are removed from the diagram in Figure ~\ref{fig:F4_slice_showing_C3}.  This diagram clearly contains three nontrivial rank ~$3$ root or weight diagrams.  Comparing this diagram to the root diagrams in Figure ~\ref{fig:Rank_3_root_diagrams}, we identify the middle diagram, containing ~$18$ non-zero vertices, as the root diagram of ~$C_3=sp(2\cdot 3)$.  The other two slices are identical non-minimal weight diagrams of ~$C_3$.  Because there are~$46$~non-zero vertices visible in Figure ~\ref{fig:F4_slice_showing_C3}, it is clear that two single vertices are missing from this representation of the root diagram of ~$F_4$, which has dimension ~$52$.

\begin{figure}[htbp]
\begin{center}
\begin{minipage}[t]{12cm}
  \begin{center}
    \leavevmode
    \resizebox{11cm}{!}{\includegraphics*[67,187][526,292]{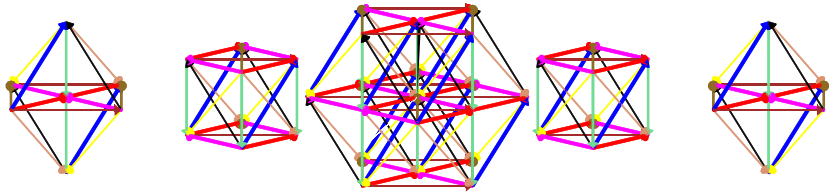}}
    \caption{\label{fig:F4_slice_showing_B3}  \hbox{Slice of ~$F_4$ showing ~$B_3=so(7) \subset B_4 =so(9) \subset F_4$}}
  \end{center}
\end{minipage}
\end{center}
\end{figure}

Figure ~\ref{fig:F4_slice_showing_B3} is the result of slicing the root diagram of ~$F_4$ using the simple roots ~$r^1$, ~$r^2$, and ~$r^3$.  The center diagram again contains ~$18$ non-zero weights, which we identify as ~$B_3=so(7)$ using Figure ~\ref{fig:Rank_3_root_diagrams}.  Hence, ~$B_3 \subset F_4$.  Furthermore, as all ~$48$ non-zero vertices are present and there are ~$5$ nontrivial slices in the root diagram, we compare this sliced root diagram of ~$F_4$ with that of ~$B_4=so(9)$, which is shown in Figure ~\ref{fig:B4_slice_showing_B3}, and see that ~$B_3 \subset B_4 \subset F_4$.  An additional slicing of ~$B_4$ shows ~$D_4 =so(8) \subset B_4 \subset F_4$.

\subsubsection{Subalgebras of ~$F_4$ using Projections}  Given the ~$4$-dimensional root diagram of ~$F_4$, we can observe its ~$3$-dimensional shadow when projected along any one direction.  However, as a single projection eliminates the information contained in one direction, it is not possible to understand the root diagram of ~$F_4$ using a single projection.  We work around this problem by creating an animation of projections, in which the direction of the projection changes slightly from one frame to the next.

\begin{figure}[htb]
\begin{center}
\begin{minipage}[t]{4.1cm}
  \begin{center}
    \leavevmode
    \resizebox{3.8cm}{!}{\includegraphics*[193,128][400,347]{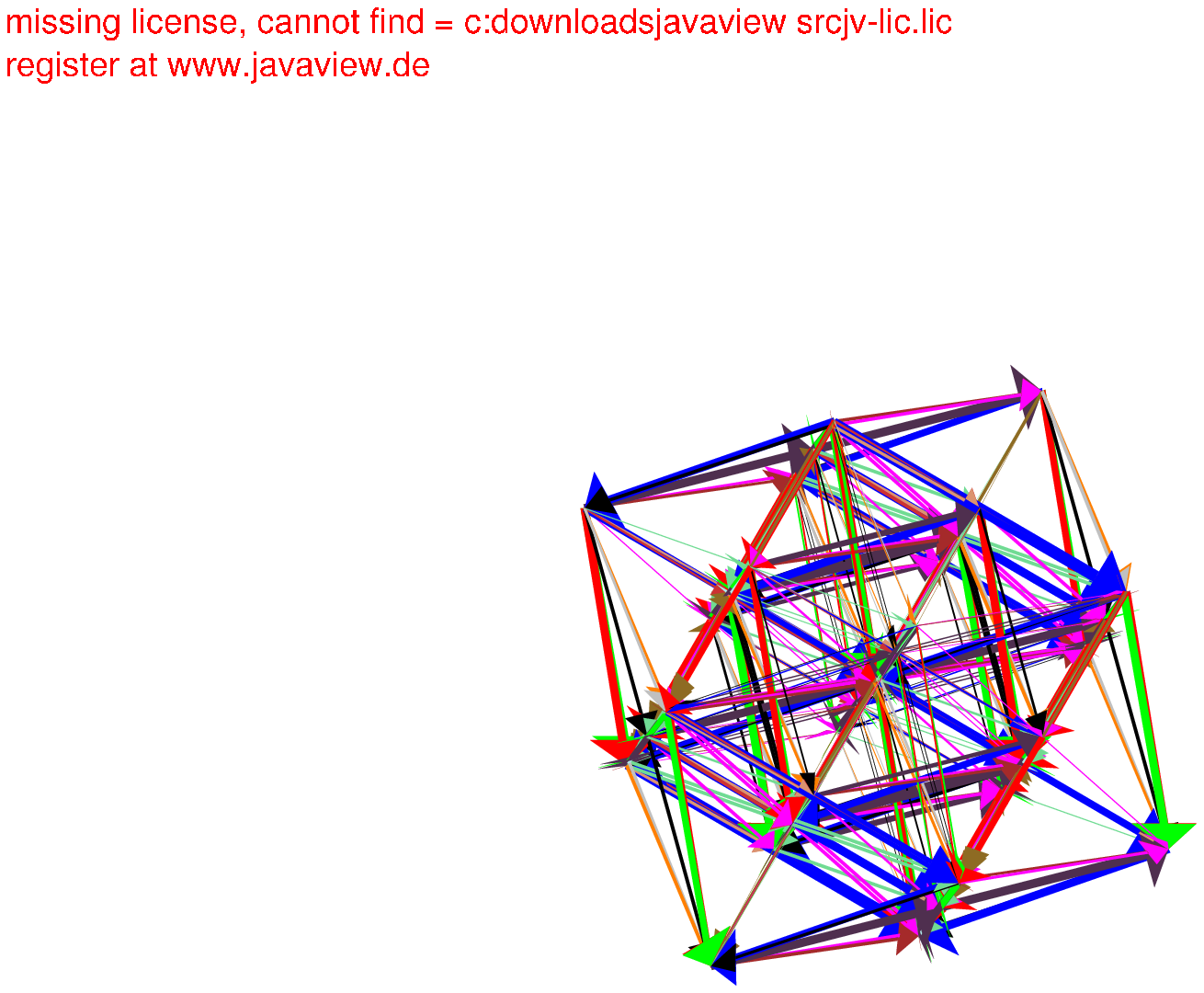}}\\
    $\theta = 0$
  \end{center}
\end{minipage}
\hfill
\begin{minipage}[t]{4.1cm}
  \begin{center}
    \leavevmode
    \resizebox{3.8cm}{!}{\includegraphics*[193,128][400,347]{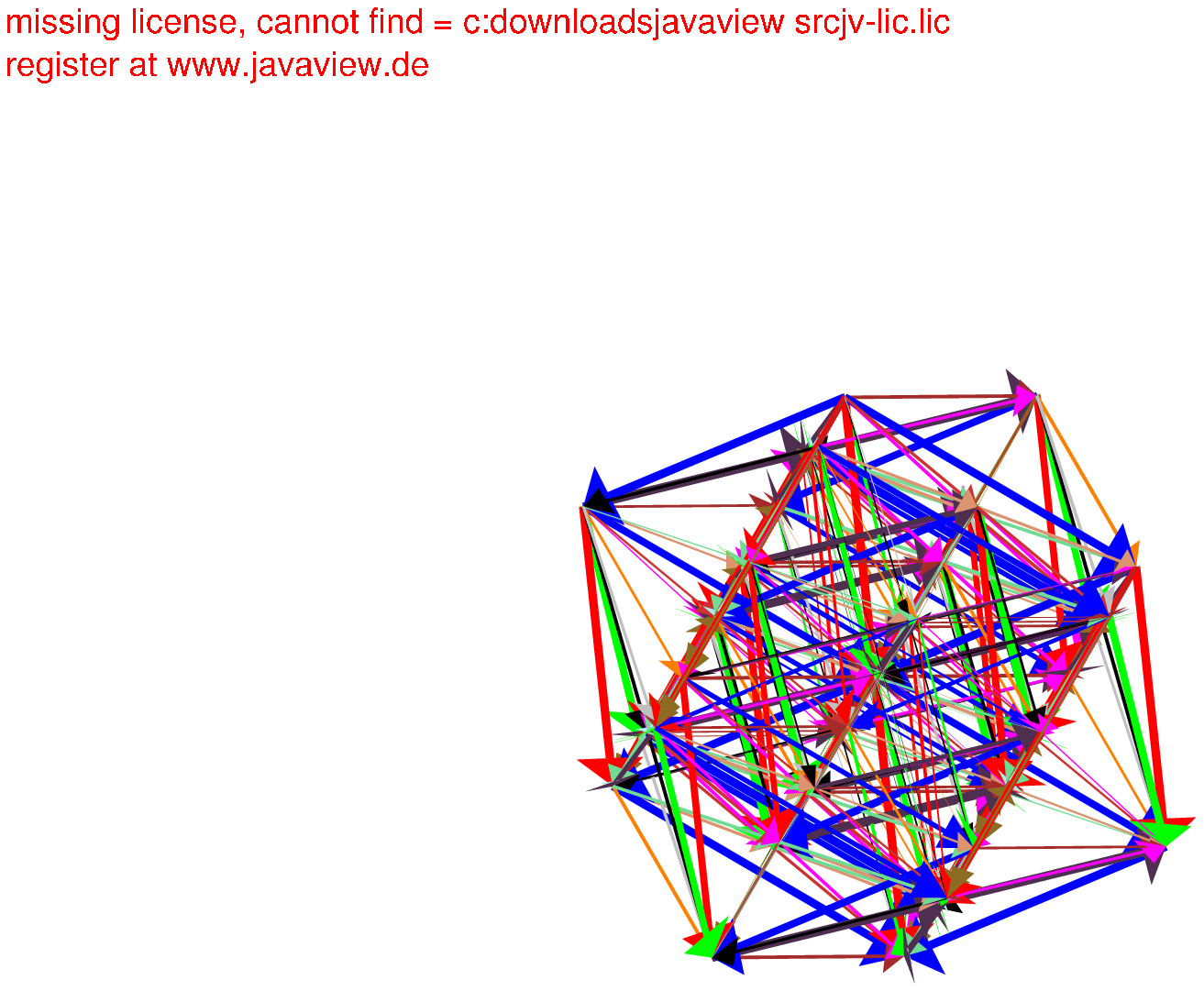}}\\
    $\theta = \frac{\pi}{8}$
  \end{center}
\end{minipage}
\hfill
\begin{minipage}[t]{4.1cm}
  \begin{center}
    \leavevmode
    \resizebox{3.8cm}{!}{\includegraphics*[193,128][400,347]{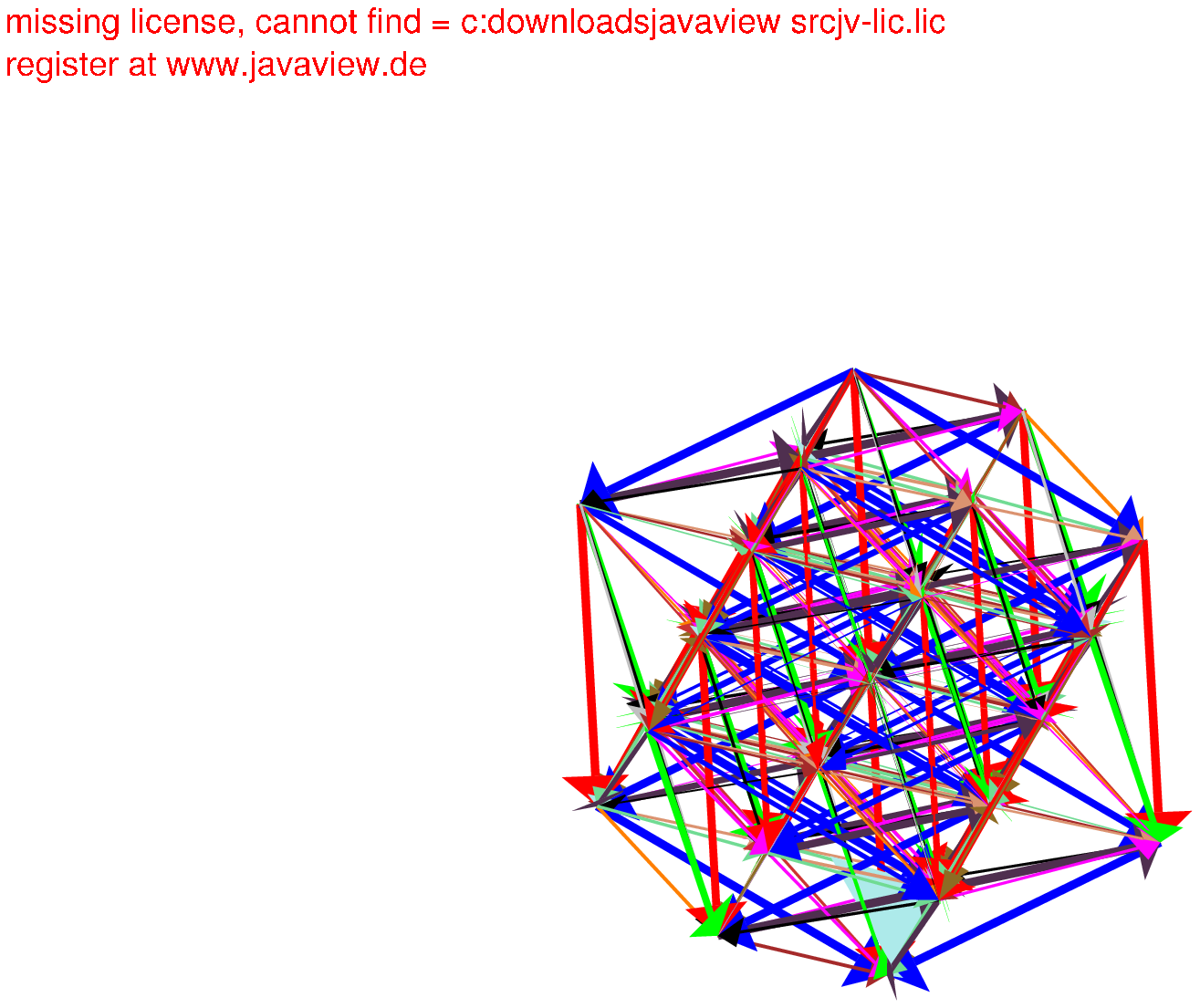}}\\
    $\theta = \frac{\pi}{4}$
  \end{center}
\end{minipage}
\hfill
\begin{minipage}[t]{4.1cm}
  \begin{center}
    \leavevmode
    \resizebox{3.8cm}{!}{\includegraphics*[193,128][400,347]{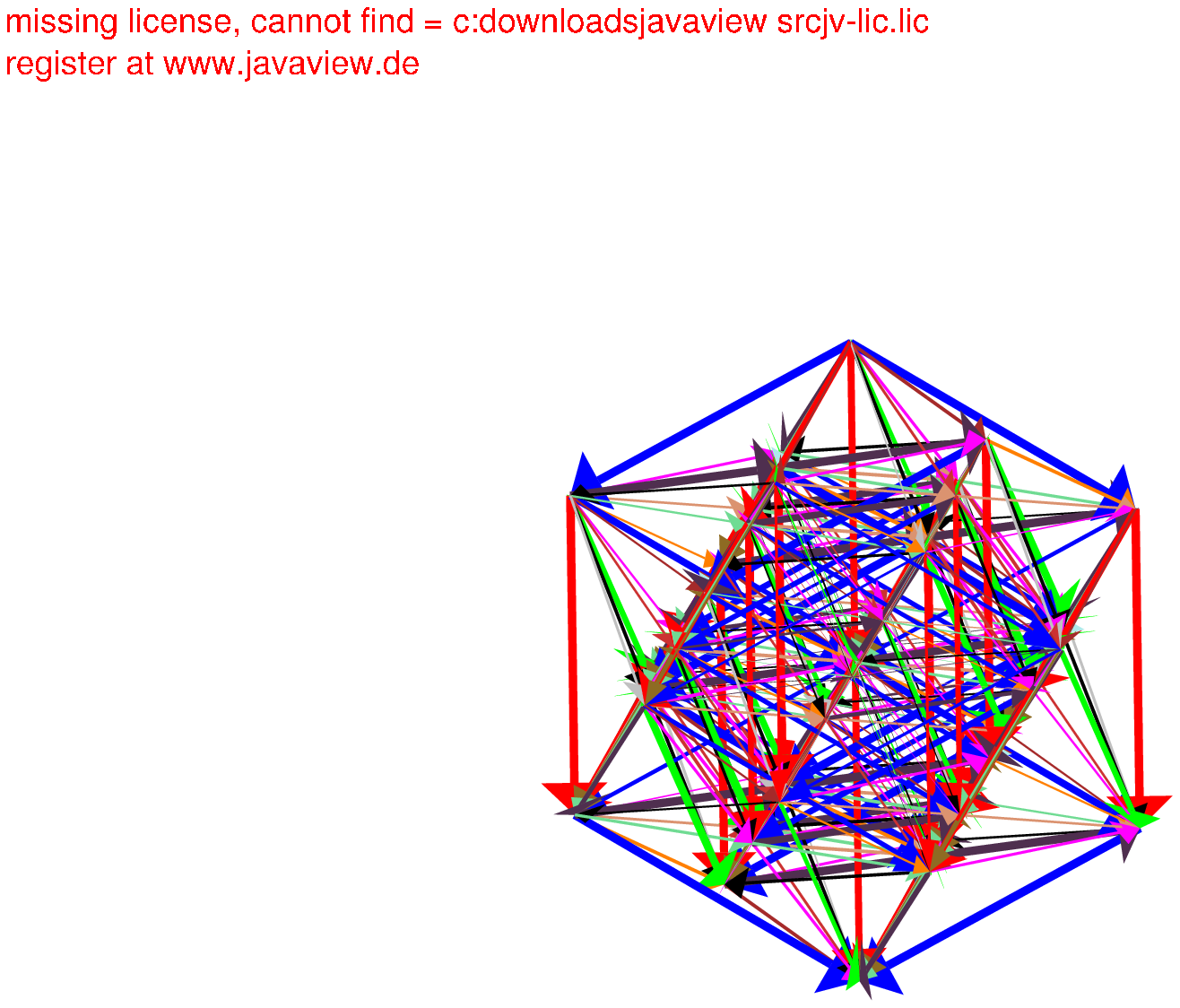}}\\
    $\theta = \frac{3\pi}{8}$
  \end{center}
\end{minipage}
\hfill
\begin{minipage}[t]{4.1cm}
  \begin{center}
    \leavevmode
    \resizebox{3.8cm}{!}{\includegraphics*[193,128][400,347]{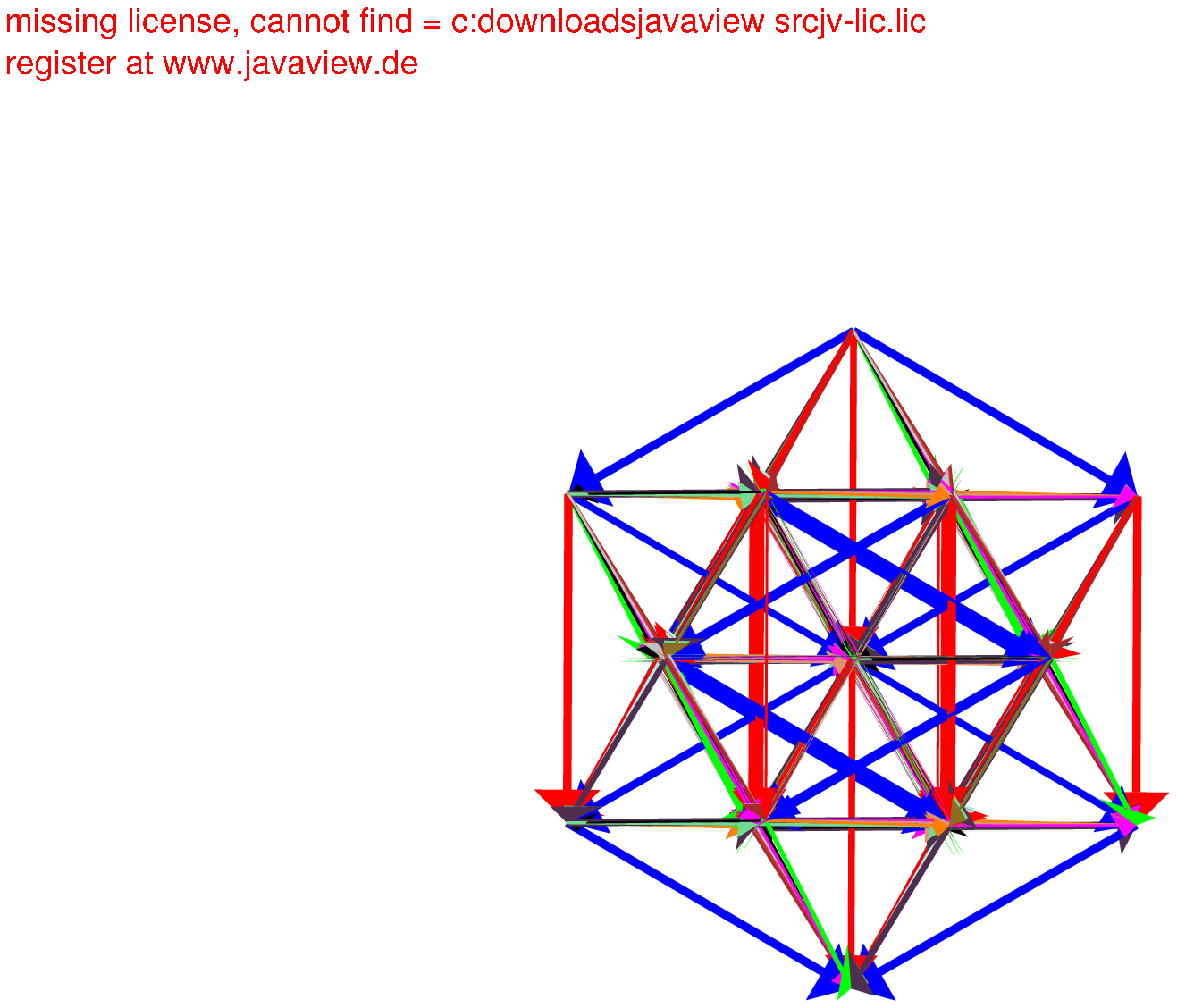}}\\
    $\theta = \frac{\pi}{2}$
  \end{center}
\end{minipage}
\hfill
\begin{minipage}[t]{4.1cm}
  \begin{center}
    \leavevmode
    \resizebox{3.8cm}{!}{\includegraphics*[199,145][394,338]{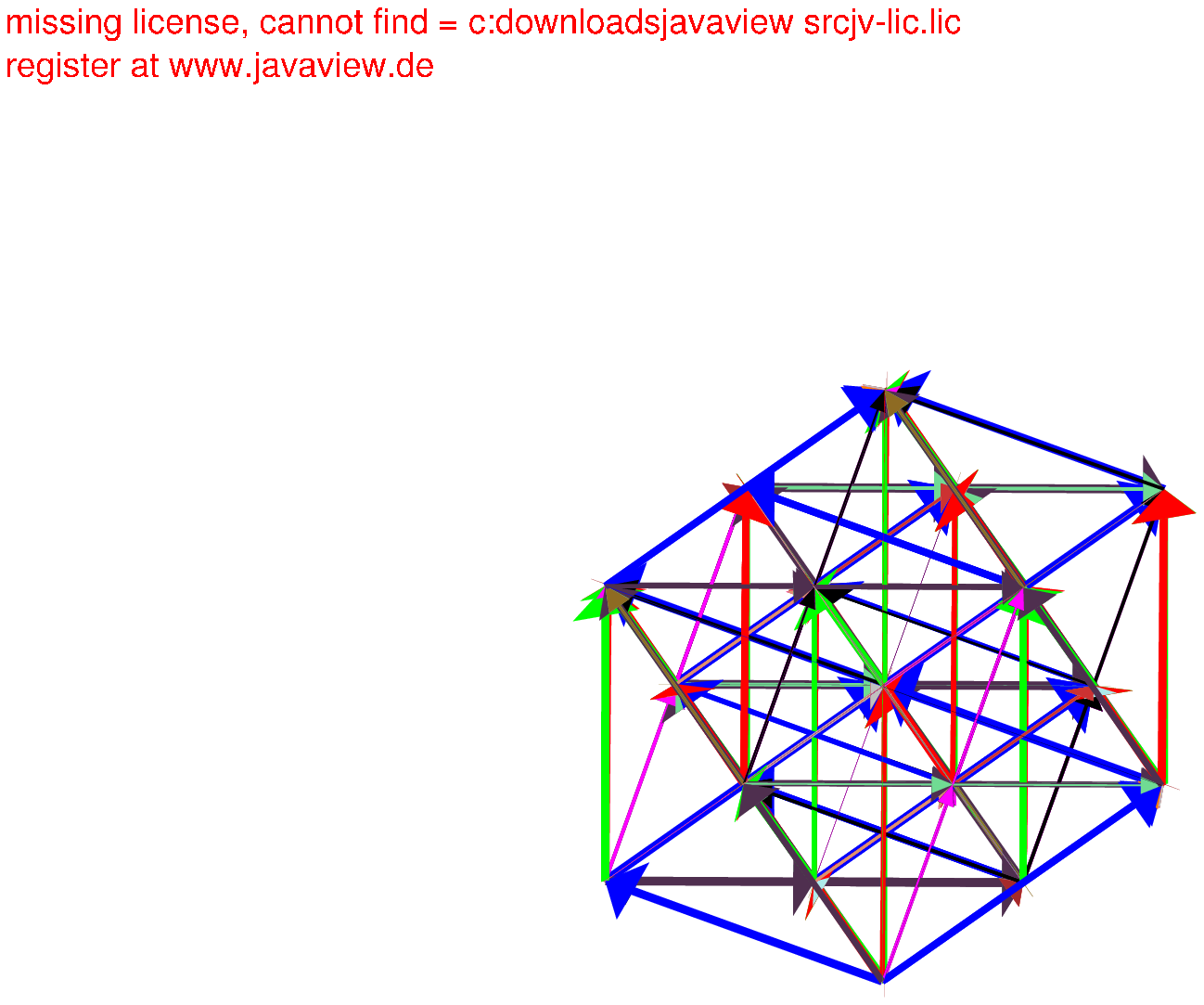}}\\
    Orthogonal View in $\mathbb{R}^4$
  \end{center}
\end{minipage}
\hfill
\caption{Animation of ~$F_4$, projecting along the \hbox{direction~$p_{\theta} = \cos(\theta)r^1 + \sin(\theta)r^4$} which is confined to a plane in~$\mathbb{R}^4$ containing ~$r^1$ and ~$r^4$.}
\label{animation_F4}
\end{center}
\end{figure}

The Dynkin diagram of ~$F_4$ reduces to the Dynkin diagram of ~$C_3=sp(2\cdot 3)$ or ~$B_3=so(7)$ by eliminating either the first or fourth node.  The simple roots ~$r^1$ and ~$r^4$ define a plane in ~$\mathbb{R}^4$, and we choose a projection vector ~$p_{\theta} = \cos \theta r^1 + \sin \theta r^4$ to vary discretely in steps of size ~$\frac{\pi}{18}$ from ~$\theta = 0$ to ~$\theta = \frac{\pi}{2}$ in this plane.  Each value of ~$\theta$ produces a frame of the animation sequence using the projection procedures of section~$3.2$.  The resulting animation is displayed in Figure ~\ref{animation_F4}.

The result of each projection of ~$F_4$ is a diagram in three dimensions.  We create the animation using Maple, and the software package Javaview is used to make a live, interactive applet of the animation.  The Javaview applet allows  the animation to be rotated in ~$\mathbb{R}^3$ as it plays.  In particular, when ~$\theta = 0$, we can rotate the diagram to show a weight diagram of ~$C_3=sp(2\cdot 3)$, which our eyes project down to the root diagram of ~$B_2 = C_2$.  Without rotating the diagram, the animation continuously changes the projected diagram as ~$\theta$ increases.  When ~$\theta = \frac{\pi}{2}$, our eyes project the root diagram of ~$B_3 = so(7)$ down to the root diagram of ~$G_2$.  However, it is also possible to rotate the animation to see the root diagram of ~$G_2$ at various other values of ~$\theta$.  The interactive animation makes it easier to explore the structure of ~$F_4$.

This interactive animation can also illustrate an obvious fact about planes in ~$\mathbb{R}^4$.  As ~$p_{\theta}$ is confined to a plane, there is a plane ~$P^{\perp}$ which is orthogonal to each of the projection directions.  Thus, the projection does not affect ~$P^{\perp}$, and it is possible to see this plane in ~$\mathbb{R}^3$ by rotating the animation to the view shown in the sixth diagram in Figure ~\ref{animation_F4}.  In this configuration, the roots and vertices in this diagram do not change as the animation varies from ~$\theta = 0$ to ~$\theta = \frac{\pi}{2}$.  While this is obvious from the standpoint of Euclidean space, it is still surprising this plane can be seen in ~$\mathbb{R}^3$ even as the projected diagram is continually changing.

\subsubsection{Modifications of methods for ~$E_6$}

Of particular interest is the exceptional Lie algebra~$E_6$, which preserves the determinant of elements of the Cayley plane.  As explained in Section~$2.1$, this allows us to write ~$E_6 = sl(3,\mathbb{O})$.  It therefore naturally contains the subalgebras ~$sl(2,\mathbb{O})$ and ~$su(2, \mathbb{O})$, which are identified as real forms of ~$D_5$ and ~$B_4$, respectively ~\cite{aaron_thesis}.

Because the root diagram of ~$E_6$ is a ~$6$-dimensional diagram, we apply  three consecutive slicings to determine subalgebras of ~$E_6$.  The first slicing creates a string of ~$5$-dimensional diagrams along the ~$x$ axis.  The second slicing, when applied to each of these diagrams, creates a string of ~$4\textrm{-dimensional}$ diagrams laid parallel to the ~$y$-axis.  The third slicing turns each ~$4$-dimensional diagram into a string of ~$3$-dimensional diagrams laid out parallel to the ~$z$ axis.  We apply the same process to rank ~$5$ algebras and use two slicings.  To identify a rank ~$5$ subalgebra ~$g$ of ~$E_6$, we need to identify its sliced root diagram as a subdiagram of the sliced ~$E_6$ diagram.  Of course, we must be careful to check that the highest weight of the ~$E_6$ root diagram is used in the subdiagram.  A similar technique is applied for finding rank ~$4$ subalgebras of both rank ~$5$ algebras and of ~$E_6$.

The projection technique can also be used to identify subalgebras of rank ~$E_6$.  In one version, we project the root diagram of ~$E_6$ along one direction, thereby creating a diagram that possibly corresponds to a rank ~$5$ algebra ~$g$.  We then apply the same pair of projections to our projected ~$E_6$ diagram and to the candidate root diagram of ~$g$.  If these two projections preserve the number of vertices in the ~$5$-dimensional diagrams, it is possible to compare the resulting diagrams in ~$\mathbb{R}^3$.  If we have identified the correct subalgebra of ~$E_6$, the resulting two diagrams should match for every pair of projections applied to the ~$5\textrm{-dimensional}$ diagrams.

Projections of rank ~$5$ and ~$6$ algebras can also be simulated using slicings of their root diagrams.  This is done using the slice and collapse technique, which collapses all the slices onto one another in a particular direction.  When using this technique, we draw the grey struts, as we are now interested in the root diagram's structure after the projection.  This technique provides clearer pictures compared to the pure projection method.

\subsubsection{Towers of Complex Lie Algebras}

We list in Figure ~\ref{fig:subalgebras_of_E6} the subalgebras of ~$E_6$ found using the slicing and projection techniques applied to an algebra's root diagram.  As mentioned in Section ~$2.1$, we list certain real representations of the subalgebras of the ~$sl(3,\mathbb{O})$ representation of ~$E_6$ in the diagram.  The particular real representations are listed below each algebra.

We use different notations to indicate the particular method that was used to identify subalgebras.  The notation ~$\xymatrix@1{A\ar@{.>}[r]&B}$ indicates the slicing method was used to identify ~$A$ as a subalgebra of ~$B$.  The notation ~$\xymatrix@1{A\ar@{->}[r]&B}$ indicates that ~$A$ was identified as a subalgebra of ~$B$ using the normal projection technique, while we indicate projections done by the slice and collapse method as ~$\xymatrix@1{A\ar@{^{}.>}[r]^{s+p}&B}$.  If both dotted and solid arrows are present, then ~$A$ can be found as a subalgebra of ~$B$ using both slicing and projection methods.  If~$A$ and ~$B$ have the same rank, only the slicing method allows us to identify the root diagram of ~$A$ as a subdiagram of ~$B$.  This case is indicated in the diagram using the notation ~$\xymatrix@1{*+++\txt{$A$}\ar@{^{(}.>}[r]&B}$.  Each of the subalgebra inclusions below can be verified online ~\cite{subalgebras_online}.

\begin{figure}[htbp]
\begin{center}
\begin{minipage}{6in}
\begin{center}
\xymatrix@M=3.5pt@H=1pt{
               &                           &   *++\txt{$E_6$\\ $sl(3,\mathbb{O})$}           &             \\
\\
               &   *++\txt{$F_4$\\ $su(3,\mathbb{O})$} \ar[uur]^(.4){s+p} &   *++\txt{$D_5 = so(10)$\\$sl(2,\mathbb{O})$} \ar@{.>}[uu]{}    & *++\txt{$A_5 = su(6, \mathbb{C})$\\$sl(3,\mathbb{H})$} \ar@{.>}[uul]{}  \\
              &                           &                                    &              \\
*++\txt{$C_4 = sp(2\cdot 4)$\\$su(3,1, \mathbb{H})$} \ar@/^4pc/[uuuurr]^{s+p} & *++\txt{$B_4 = so(9)$\\$su(2,\mathbb{O})$} \ar@{_{(}.>}[uu]^(.4){}  \ar[uur]^(.4){s+p}   &   D_4 = so(8) \ar@{.>}[uu]_(.3){} \ar@{_{(}.>}[l]_(.4){}  &   A_4 = su(5, \mathbb{C}) \ar@{.>}[uul]<1ex>{} \ar[uul]_(.6){s+p} \ar@{.>}[uu]<1ex>{} \ar[uu]_(.4){}  \\
              &                           &                                    &        \\
*++\txt{$C_3=sp(2\cdot 3)$\\$su(3,\mathbb{H})$} \ar@/_4pc/[uuuurrr]^(.6){} \ar@{.>}[uu]<1ex>{} \ar[uu]^(.4){} \ar@{.>}[uuuur]{} & B_3 = so(7) \ar@{.>}[uu]{}  \ar[uur]<.5ex>^(.6){} &                   *++\txt{$ A_3 = su(4,\mathbb{C})$\\$D_3=so(6)$\\$sl(2,\mathbb{H})$} \ar[uur]{} \ar@{_{(}.>}[l]^(.4){} \ar@{.>}[uur]<1ex>{} \ar@{.>}[uull]<-1ex>^(.4){} \ar@{.>}[uu]{} &  \\
               &                           &          &                \\
 *++\txt{$B_2 = so(5)$\\$C_2 = sp(2\cdot 2)$\\$su(2,\mathbb{H})$} \ar@/_2.3pc/[uurr]{} \ar@{.>}[uu]{}  \ar@{.>}[uur]<1ex> \ar[uur]^(.4){}  & G_2 = Aut(\mathbb{O}) \ar[uu]^(.4){}  & A_2 = su(3,\mathbb{C}) \ar@{.>}[uull]<0ex>{} \ar@{.>}[uu]<1ex>{} \ar[uu]_(.4){} \ar@{_{(}.>}[l]^(.4){} \\
               & \save[]+<0cm,.1cm>*++\txt{$D_2 = su(2,\mathbb{C}) \oplus su(2,\mathbb{C})$} \ar@{_{(}.>}[ul]^(.4){} \ar@{_{(}.>}[u]^(.4){} \restore &                  &                \\
               &                           &   A_1 = su(2,\mathbb{C}) \ar@{.>}[uu]<1ex>{} \ar[uu]_(.4){} \ar@{.>}[ul]<1ex>{} \ar[ul]^(.4){}     &\\
}
$\xymatrix@1{A\ar@{.>}[r]&B}$ indicates $A$ is realized using a slice of $B$.
\\
$\xymatrix@1{A\ar@{->}[r]&B}$ indicates A is a projection of B.
\\
$\xymatrix@1{A\ar@{->}[r]^{s+p}&B}$ indicates B projects to A with the slice and project method.
\\
$\xymatrix@1{*+++\txt{$A$} \ar@{^(.>}[r]&B}$ signals A and B have the same rank, but A is a subdiagram of B.

\caption{Subalgebras of ~$E_6$ together with some important real representations}
\label{fig:subalgebras_of_E6}
\end{center}
\end{minipage}
\end{center}
\end{figure}

Lists of subalgebra inclusions are found in ~\cite{danielsson}, which applies subalgebras to particle physics, and in ~\cite{gilmore}, which recreates the subalgebra lists of ~\cite{van_der_waerden}.  However, the list in ~\cite{gilmore} mistakenly has ~$C_4$ and ~$B_3$ as subalgebras of ~$F_4$, instead of ~$C_3$ and ~$B_4$.  Further, the list omits the inclusions ~$G_2 \subset B_3$, ~$C_4 \subset E_6$, ~$F_4 \subset E_6$, and ~$D_5 \subset E_6$.  The correct inclusions of ~$C_3 \subset F_4$ and ~$B_4 \subset F_4$ are listed in Section ~$8$ of ~\cite{van_der_waerden}, but the ~$B_n$ and ~$C_n$ chains are mislabeled in the final table of \cite{van_der_waerden} which was used by Gilmore in ~\cite{gilmore}.  Although \cite{van_der_waerden} uses root systems to determine subalgebra inclusions, it mistakenly claims that ~$D_n \subset C_n$ as a subalgebra in his Section ~$21$ of \cite{van_der_waerden}, which is not true since their root diagrams are based upon inequivalent highest weights.

In ~\cite{dyn2}, Dynkin classified subalgebras depending upon the root structure.  If the root system of a subalgebra can be a subset of the root system of the full algebra, the subalgebra is called a {\it regular} subalgebra.  Otherwise, the subalgebra is {\it special}.  A complete list of regular and special subalgebras is given in ~\cite{danielsson}.  All of the regular embeddings of an algebra in a subalgebra of ~$E_6$ can be found using the slicing method.  In many cases, the projection technique also identifies these regular embeddings of subalgebras, but there are regular embeddings which are not recognized as the result of projections.  The special embeddings of an algebra in a subalgebra of ~$E_6$ can only be found using the projection technique.  We conjecture that all the subalgebras of a complex Lie algebra may be found using our slicing and projection methods, as these methods find the same regular and special subalgebras of a complex Lie algebra as found by Dynkin.

\subsection{Conclusion}
\label{ch:Joma_Paper.conclusion}

We have presented here methods which illustrate how root and weight diagrams can be used to visually identify the subalgebras of a given Lie algebra.  While the standard methods of determining subalgebras rely upon adding, removing, or folding along nodes in a Dynkin diagram, we show here how to construct any of a Lie algebra's root or weight diagrams from its Dynkin diagram, and how to use geometric transformations to visually identify subalgebras using those weight and root diagrams.  In particular, we show how these methods can be applied to algebras whose root and weight diagrams have dimensions four or greater.    In addition to pointing out the erroneous inclusion of ~$C_4 \subset F_4$ in ~\cite{gilmore, van_der_waerden}, we provide visual proof that ~$C_4 \subset E_6$ and list all the subalgebras of ~$E_6$.  While we are primarily concerned with the subalgebras of ~$E_6$, these methods can be used to find subalgebras of any rank ~$l$ algebra.

\newpage{}

\part{Division Algebras and Applications}
\label{ch:Division_Algebras_Applications}


This chapter includes a review of the concepts we will use in the construction of the group~$E_6 = SL(3,\mathbb{O})$ in Chapter \ref{ch:E6_basic_structure}.  Readers familiar with the octonions and Lorentz transformations may want to skip this chapter, although they may want to quickly review of notion of triality in Section \ref{ch:Division_Algebras_Applications.Quaternions_Octonions.Triality} and Lorentz groups involving division algebras in Section \ref{ch:Division_Algebras_Applications.Lorentz_Transformations}.  In Section \ref{ch:Division_Algebras_Applications.Quaternions_Octonions}, we review the properties of the division algebras and point out that the complexes, quaternions, and octonions naturally describe Euclidean spaces of dimension $2$, $4$, and $8$.
  In Section \ref{ch:Division_Algebras_Applications.Quaternions_Octonions.Conjugation_Reflecitons_Rotations}, we review the work of Manogue and Schray \cite{manogue_schray} who show how to construct rotations and reflections in~$\mathbb{R}^4$ and~$\mathbb{R}^8$ using certain conjugation maps in the division algebras.  These maps will be used extensively in the construction of our real form of~$E_6$ in Chapter \ref{ch:E6_basic_structure}.  We finish our review of the division algebras with a summary of triality in Section \ref{ch:Division_Algebras_Applications.Quaternions_Octonions.Triality}.  This idea will be again visited in Section \ref{ch:E6_basic_structure.Triality}, where we find that~$sl(3,\mathbb{O})$, the Lie algebra corresponding to the Lie group~$SL(3,\mathbb{O})$, exhibits interesting characteristics due to triality.
We review the properties of Lorentz Transformations in Section \ref{ch:Division_Algebras_Applications.Lorentz_Transformations} and pay particular attention to the notion that the division algebras may be used to construct Lorentz Transformations in~$k+1$ dimensions, where \hbox{$k-1 = |\mathbb{K}|$} is the dimension of the division algebra~$\mathbb{K} = \mathbb{R}, \mathbb{C}, \mathbb{H}, \mathbb{O}$.  We conclude Chapter \ref{ch:Division_Algebras_Applications} in Section \ref{ch:Division_Algebras_Applications.Jordan_Algebras} by describing how the three smallest division algebras may be used to construct Jordan algebras while the octonions may be used to construct the exceptional Jordan algebra $M_3(\mathbb{O})$.  Octonionic Lorentz transformations and the exceptional Jordan algebra are utilized in Chapter \ref{ch:E6_basic_structure} to give a description of the Lie group $SL(3,\mathbb{O})$ and its associated Lie algebra $sl(3,\mathbb{O})$.

\section{Normed Division Algebras}
\label{ch:Division_Algebras_Applications.Quaternions_Octonions}
The complexes, quaternions and octonions are division algebras of dimension~$2$,~$4$, and~$8$ over the reals.  
Section \ref{ch:Division_Algebras_Applications.Quaternions_Octonions.Division_Algebras} includes a review of the basic algebraic and geometric properties of the normed division algebras.  These algebras provide a nice description of~$|\mathbb{K}|$-dimensional Euclidean space for~$\mathbb{K} = \mathbb{R}, \mathbb{C}, \mathbb{H},$ or~$\mathbb{O}$.  Certain multiplication maps in the division algebra can result in rotations or reflections in the corresponding Euclidean space, and this discussion, based on the work of \cite{manogue_schray}, is included in Section \ref{fig:octonionic_multiplication}.  Triality is a concept closely related to the multiplication properties of the division algebras as well as to representations of Lie algebras, and is discussed in Section \ref{ch:Division_Algebras_Applications.Quaternions_Octonions.Triality}.  Particular attention is paid to the triality related to the octonions and representations of~$so(8,\mathbb{R})$.

\subsection{Reals, Complexes, Quaternions and Octonions}
\label{ch:Division_Algebras_Applications.Quaternions_Octonions.Division_Algebras}
We review here the properties of the four division algebras, using the construction provided using the Cayley-Dickson process.  Additional information about this construction and the properties of the division algebras may be found in either \cite{dickson} or \cite{baez}.

\subsubsection{Real and Complex Numbers}
\label{ch:Division_Algebras_Applications.Quaternions_Octonions.Division_Algebras.Real_and_Complex_Numbers}

The Cayley-Dickson process can be used to create a~$2n$-dimensional algebra from an~$n$-dimensional associative algebra.
The real numbers,~$\mathbb{R}$, are associative, commutative, and have an identity, denoted~$1$.  Applied to~$\mathbb{R}$, the Cayley-Dickson process creates a two-dimensional algebra~$\mathbb{C} = \mathbb{R} \oplus \mathbb{R}i$ over~$\mathbb{R}$, using~$i$ to denote a square root of~$-1$.
  This algebra has multiplication
\[ (a, b) (c, d) = (ac-bd, ad+bc) \]
or equivalently,
\[ (a+bi) (c+di) = ac-bd + (ad+bc)i \]
In addition, this algebra has the property that ~$i^2 = -1$, or ~\hbox{$(0,1) (0,1) = (-1,0)$}.  For any ~$x,y,z \in \mathbb{C}$, this algebra is associative
\[ (xy)z = x(yz) \]
and commutative
\[ x y = y x \]
There is a conjugation map taking~$a+bi$ to~$\overline{a+bi} = a-bi$.  Using this map, we can define the norm 
$$|a+bi|^2 = (a+bi) \overline{(a+bi)} = a^2 + b^2$$
 of any complex number~$a+bi$.  Notice that the only complex number of norm~$0$ is the complex number~$0$ and that ~$z^{-1} = \frac{\overline{z}}{|z|^2}$ is the inverse of any non-zero complex number~$z$.  We call~$Re(z) = \frac{1}{2} \left( z + \overline{z} \right)$ and~$Im(z) = \frac{1}{2}\left( z - \overline{z}\right)$ the {\it real part} and {\it imaginary part} of~$z \in \mathbb{C}$.

One nice characteristic of the complex numbers is their ability to describe points in a~$2$-dimensional plane.  Each complex number~$a+bi \in \mathbb{C}$ may be identified with the point~$(a,b) \in \mathbb{R}^2$.  Euler's formula
$$ e^{i\alpha} = \cos \alpha + i \sin \alpha$$
may be used to write every complex number~$z = a+bi$ in the form~$z = |z|e^{i \alpha}$.  The distance from the origin to~$(a,b)$ is denoted by~\hbox{$|z| = \sqrt{a^2 + b^2}$}, while~$\alpha$ is the angle between the positive real axis and the ray extending from the origin to~$(a,b)$.  The quantity~$e^{i\alpha}$ is called a {\it phase}.  We note that complete set of complex phases describes~$S^1$.

\subsubsection{Quaternion Numbers}
\label{ch:Division_Algebras_Applications.Quaternions_Octonions.Division_Algebras.Quaternion_Numbers}

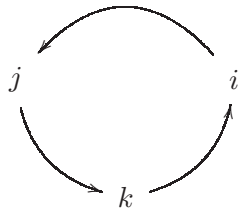
\begin{figure}[htbp]
\vspace{1cm}
\begin{center}
\begin{minipage}{3in}
\begin{center}
\mbox{
\xymatrix@M=3pt@H=10pt{
 *++\txt{$j$}  \ar@/_1pc/[dr] &              & *++\txt{$i$}  \ar@/_2.3pc/[ll] \\
               & *++\txt{$k$} \ar@/_1pc/[ur] &    \\
}
}
\end{center}
\caption{Quaternionic multiplication}
\label{fig:quaternionic_mult}
\end{minipage}
\end{center}
\end{figure}

As the complex numbers are associative, we can again use the Cayley-Dickson process and produce the quaternions~$\mathbb{H} = \mathbb{C} \oplus \mathbb{C}j$, a four-dimensional division algebra over~$\mathbb{R}$, by using~$j$, another square root of~$-1$.  The multiplication is again given by the rule 
\[ (a, b) (c, d) = (ac-bd, ad+bc) \]
where~$a, b, c, d \in \mathbb{C}$.  Setting~$a = q_1 + q_2 i$ and~$b = q_3 + q_4 i$, we will often write the quaternionic number ~$(a,b) = a + b j = q_1 + q_2 i + q_3 j + q_4 k$, where~$k$ is yet another square root of~$-1$ and is the product~$k = i j$.
  This multiplication rule leads to the products 
$$
\begin{array}{ccc}
i\; j = k  & j \; k = i  & k \; i = j \\
j \; i = -k & k \; j = -i & i \; k = -j 
\end{array}
$$
This multiplication is nicely summarized in Figure \ref{fig:quaternionic_mult}, where the arrow indicates whether the product of two imaginary unit quaternions will carry a plus (with the arrow) or minus sign.  The quaternions are not commutative, but they are still associative.

Although they are not commutative, the quaternions have similar norm properties to the complex numbers.
  Again, there is a conjugation map~$q \to \overline{q}$ where
$$\overline{q_1 + q_2 i + q_3 j + q_4 k} = q_1 - q_2 i - q_3 j - q_4 k$$
allowing the definition of a norm~$|q|^2 = q_1^2 + q_2^2 + q_3^2 + q_4^2$ for any quaternion~$q \in \mathbb{H}$.  The inverse of any non-zero~$q \in \mathbb{H}$ is given by~$q^{-1} = \frac{\overline{q}}{|q|^2}$.  Again,~$Re(q) = \frac{1}{2} \left( q + \overline{q} \right)$ and~$Im(q) = \frac{1}{2}\left( q - \overline{q}\right)$ are called the {\it real part} and {\it imaginary part}  of~$q \in \mathbb{H}$.

Just as the complex numbers can be used to describe a plane, the quaternions may be used to describe points in a~$4$-dimensional space.  The quaternion~$q = q_1 + q_2 i + q_3 j + q_4 k$ may be identified with the point~$(q_1, q_2, q_3, q_4) \in \mathbb{R}^4$.  We note that the imaginary quaternions may be identified with vectors in~$\mathbb{R}^3$, in which case the quaternionic multiplication corresponds to the ordinary cross-product.  The set of imaginary unit quaternions form the sphere~$S^{2}$ in~$\mathbb{R}^3$.  For any~$s \in S^{2}$, we may form the complex subalgebra with basis~$\lbrace 1, s \rbrace$.  Hence, Euler's formula may be used to write any quaternion~$q \in \mathbb{H}$ in the form~$q = |q|e^{s \alpha}$, where~$s$ is the unit imaginary quaternion pointing toward~$q$.

\subsubsection{Octonion Numbers}
\label{ch:Division_Algebras_Applications.Quaternions_Octonions.Division_Algebras.Octonion_Numbers}

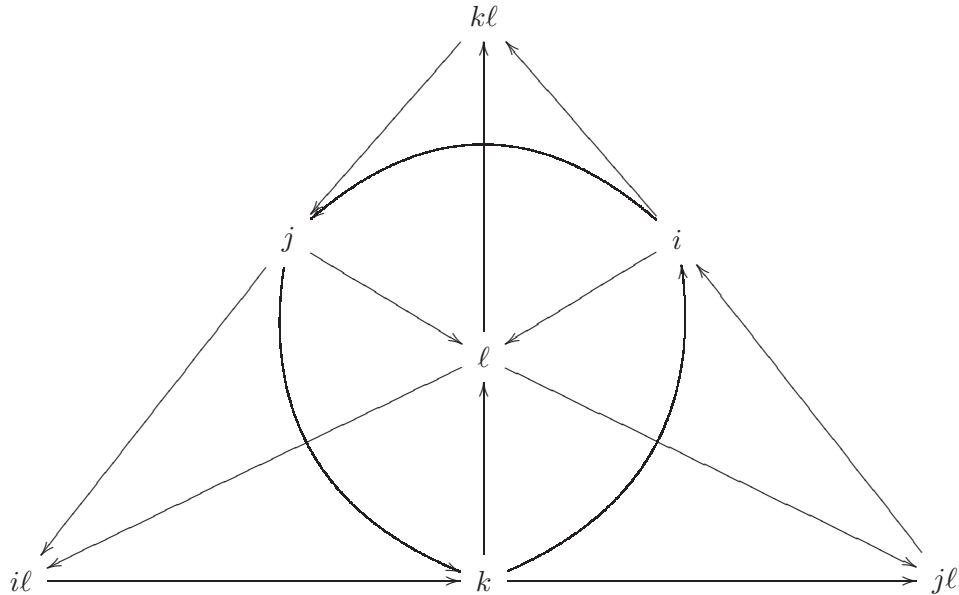
\begin{figure}[htb]
\begin{center}
\hspace{-.25in}
\begin{minipage}{5.0in}
\begin{center}
\mbox{
\xymatrix@M=3pt@H=10pt{
 &   &  &  &   *++\txt{$k\l$}   \ar[ddll]  & &  & &   \\
\\
 &   &  *++\txt{$j$}  \ar@/_2.7pc/[dddrr] \ar[dddll] \ar[drr] & &  & & *++\txt{$i$}  \ar@/_3pc/[llll]  \ar[uull] \ar[dll] \\
 &   &  &  &  *++\txt{$\l$} \ar[ddrrrr] \ar[ddllll] \ar[uuu]  & &  & &   \\
\\
*++\txt{$i\l$}  \ar[rrrr]   &  \hspace{1.0cm}      &                  &     & *++\txt{$k$} \ar@/_2.6pc/[uuurr]  \ar[uu] \ar[rrrr] & &        &     \hspace{1.0cm}    &  *++\txt{$j\l$} \ar[uuull] \\
}
}
\end{center}
\caption{Octonionic multiplication}
\label{fig:octonionic_multiplication}
\end{minipage}
\end{center}
\end{figure}

When applied to the quaternions, the Cayley-Dickson process produces the octonions,~$\mathbb{O} = \mathbb{H} + \mathbb{H}\l$, an eight-dimensional division algebra over~$\mathbb{R}$.  Here,~$\l$ is yet another square root of~$-1$, which is orthogonal to~$\mathbb{H}$.  The multiplication is again given by the rule
$$(a, b) (c, d) = (ac - bd, ad+bc)$$
with~$a,b,c,d \in \mathbb{H}$.  
We will often write the octonion~$(a,b) = a+b\l$ as
$$q_1 + q_2 i + q_3 j + q_4 k + q_5 k\l + q_6 j\l + q_7 i\l + q_8 \l$$
where~$a = q_1 + q_2 i + q_3 j + q_4 k$,~$b = q_8 + q_7 i + q_6 j + q_5 k$ and~$k\l, j\l,$ and~$i\l$ are the products of~$k, j$, and~$i$ with~$\l$, respectively.  We refer to~$\lbrace i, j, k, k\l, j\l, i\l, \l \rbrace$ as the standard basis of unit imaginary octonions.  It is convenient to encode the multiplication of the octonions as shown in Figure \ref{fig:octonionic_multiplication}, which contains seven directed loops, six of which are shown as directed lines instead of loops.  Each loop contains three octonions~$p$,~$q$ and~$r$.  The product of any two of these octonions~$p$ and~$q$ is~$\pm r$, where the positive sign is chosen if the order of multiplication follows the arrow and the negative sign is chosen if the order goes against the arrow.  Hence, we see that~$i \; \l = i\l$ and~$\l \; i\l = i$ but~$\l \; i = - i\l$.  We note that the octonions are not commutative, since they are constructed from the quaternions, and not associative, as~$i \; (j \; \l) = i \; (j\l) = -k\l$ but~$(i \; j) \; \l = (k) \; \l = k\l$.  However, the octonions are {\it alternative}, since 
$$ x(xy) = (xx)y \hspace{2cm}  (yx)x = y(xx)$$
for any $x,y \in \mathbb{O}$.
  The conjugation map~$q \to \overline{q}$ changes the sign of every imaginary basis unit in~$q$.  With this modification, the norm~$|q|^2 = q \overline{q}$, inverse~$q^{-1} = \frac{q}{|q|}$, real~$Re(q) = \frac{1}{2}(q + \overline{q})$ and imaginary~$Im(q) = \frac{1}{2}(q - \overline{q})$ parts of an octonion~$q \in \mathbb{O}$ are similar to those for the other division algebras.

Just as the complex numbers and quaternions could be used to describe~$\mathbb{R}^2$ and~$\mathbb{R}^4$, the octonions may be used to describe points in~$\mathbb{R}^8$ using the obvious identification.  The imaginary octonions describe~$\mathbb{R}^7$, and the set of imaginary unit octonions form the six-sphere~$S^6$ in~$\mathbb{R}^7$.  Of course, Euler's formula again allows us to write any octonion~$q \in \mathbb{O}$ in the form~$q = |q|e^{s \alpha}$ where now~$s \in S^6$ points in the same direction as~$Im(q)$.  
According to a theorem by Artin, the subalgebra generated by any two elements in an alternative algebra is associative \cite{schafer_book}.  
For the octonions, two imaginary orthogonal units~$s_1, s_2 \in S^6$ define a quaternionic subalgebra spanned by~$\langle 1, s_1, s_2, s_1 \;  s_2 \rangle$.  This result may be generalized to the case where~$s_1$ and~$s_2$ are not orthogonal (but not parallel) by finding the orthogonal projection of~$s_2$ onto~$s_1$ and re-normalizing the resulting octonions.  
We note that each triple of imaginary unit octonions in the multiplication diagram in Figure \ref{fig:octonionic_multiplication} define a quaternionic subalgebra.  When working with one or two octonions, it is advantageous to consider the complex or quaternionic subalgebra generated by those octonions.

\subsection{Conjugation, Reflections, and Rotations}
\label{ch:Division_Algebras_Applications.Quaternions_Octonions.Conjugation_Reflecitons_Rotations}

As mentioned in the previous section, each division algebra~$\mathbb{K} = \mathbb{R}, \mathbb{C}, \mathbb{H}, \mathbb{O}$ may be identified with~$\mathbb{R}^{|\mathbb{K}|}$.  This identification allows us to produce a geometric transformation in~$\mathbb{R}^{|\mathbb{K}|}$ using multiplication in ~$\mathbb{K}$.  In this section, we describe the transformation corresponding to the conjugation map~$f_x : \mathbb{K} \to \mathbb{K}$ given by~$f_x(y) = x y \overline{x}$ for each~$x \in \mathbb{K}$.  Noting that~$\mathbb{R}$ and ~$\mathbb{C}$ are commutative, implying the conjugation map is the identity map.  Hence, we only discuss the quaternionic or octonionic case.  The material in this section is a summary of the treatment given in \cite{manogue_schray}, and will be fundamental to the construction of~$SL(2,\mathbb{O})$ in Section \ref{ch:Division_Algebras_Applications.Lorentz_Transformations} and the construction of ~$SL(3,\mathbb{O})$ in Section \ref{ch:E6_basic_structure.Lorentz_Transformations}.

Consider the case~$\mathbb{K} = \mathbb{H}$, and let~$x = e^{i \alpha/2} \in \mathbb{H}$,~$\alpha \in \mathbb{R}$, be a phase.  We shall see that the conjugation map~$f_x : \mathbb{H} \to \mathbb{H}$ produces a rotation in one plane in~$\mathbb{R}^4$ through an angle~$\alpha$; this is most easily seen using an explicit example.  Let~$q = q_1 + q_2 i + q_3 j + q_4 k \in \mathbb{H}$.  Then~$f_x(q) = e^{i\alpha / 2} q e^{-i\alpha / 2}$.  The phase~$e^{i\alpha / 2}$ commutes with the basis elements~$1$ and~$i$.  However, for~$j$ or~$k$, we have
$$  j \; e^{-i\alpha / 2} = e^{i \alpha /2}\; j \hspace{1cm} k\; e^{-i\alpha / 2} = e^{i \alpha /2}\; k $$
Hence, the conjugation map~$f_x$ fixes the~$(1,i)$ plane and rotates the~$(j,k)$ plane through an angle of~$\alpha$ radians:
$$
\begin{array}{ccl}
 e^{i\alpha / 2}(q)e^{-i\alpha / 2}& =& e^{i\alpha / 2}(q_1 + q_2 i + q_3 j + q_4 k)e^{-i \alpha / 2}\\
 &= & e^{i\alpha / 2} e^{-i\alpha / 2}(q_1 + q_2 i) + e^{i\alpha / 2} e^{i \alpha / 2}(q_3 j + q_4 k) \\
 &= & q_1 + q_2 i + e^{i\alpha}(q_3 j + q_4 k)
\end{array}
$$
This calculation may be generalized.  For any imaginary unit $s \in \mathbb{H}$ and $\alpha \in \mathbb{R}$, the phase~$e^{s\alpha / 2} \in \mathbb{H}$ defines a unique plane~$P$ in~$\mathbb{R}^4$ perpendicular to~$1$ and~$s$ whenever $e^{s\alpha / 2} \ne \pm 1$.
  Then, analogous to the example given above, the conjugation map~$f_{e^{s\alpha / 2}} : \mathbb{H} \to \mathbb{H}$ will fix the~$(1,s)$ plane and produce a rotation in the plane~$P$ through~$\alpha$ radians.

In the case~$\mathbb{K} = \mathbb{O}$, the conjugation map~$f_x : \mathbb{O} \to \mathbb{O}$ no longer rotates one plane in~$\mathbb{R}^8$.  Indeed, if~$s$ is an imaginary unit octonion and~$\alpha \in \mathbb{R}$, then~$e^{s\alpha/2}$ is a phase in~$\mathbb{O}$.  Assume that~$Im(e^{s\alpha}) \ne 0$, or equivalently, that~$\alpha$ is not an integer multiple of~$2\pi$.  Then, orthogonal to the~$(1,s)$ plane in~$\mathbb{R}^8$, there are three planes~$P_1, P_2, P_3$ which are also pair-wise orthogonal.  We will now show that the conjugation map produces a rotation in three planes, not one!
Given an imaginary unit octonion~$s \in \mathbb{O}$, pick an orthonormal basis~$\lbrace 1, s, p_1, s_1, p_2, s_2, p_3, s_3 \rbrace$ for~$\mathbb{O}$ over~$\mathbb{R}$ with the property that~$s = p_a \; s_a$ for each~$a = 1, 2, 3$.  Let~$P_a$ denote the~$(p_a, s_a)$ plane.  By the construction of the basis,~$P_1, P_2, P_3$ are pairwise orthogonal and perpendicular to the~$(1,s)$ plane.
Expand~$q \in \mathbb{O}$ in terms of the orthonormal basis.  Conjugating~$q$ by the phase~$e^{s\alpha /2}$, we see that~$e^{s\alpha/2}$ will commute with~$1$ and~$e^{s\alpha/2}$, but
$$ p_a e^{- s \alpha/2} = e^{s \alpha/2} p_a \hspace{1cm} s_a e^{- s \alpha/2} = e^{s \alpha/2} s_a$$
for each~$a = 1,2,3$.
  Hence, just as with the quaternions, the plane spanned by~$p_a$ and~$s_a$ is rotated by~$\alpha$ radians.  We note, however, that in the case~$\mathbb{K} = \mathbb{O}$, the conjugation map rotates all three planes~$P_a$ through an angle of~$\alpha$ radians while fixing the~$(1,s)$ plane.

Manogue and Schray \cite{manogue_schray} note that a rotation of a single plane in~$\mathbb{R}^8$ may be constructed as the composition of two {\it flips}, which are reflections of~$\mathbb{R}^8$ across a two-dimensional plane containing the origin.  In terms of octonions, a flip is accomplished with the conjugation map~$f_x : \mathbb{O} \to \mathbb{O}$ where again $s \in \mathbb{O}$ is a unit imaginary octonion but now the phase angle~$\alpha$ in~$x = e^{s \alpha / 2}$ is chosen to be~$\alpha = \pi$.  This again fixes the~$(1,s)$ plane, but causes the three pairwise orthogonal planes~$P_a$ for~$a = 1, 2, 3$ perpendicular to the~$(1,s)$ plane in~$\mathbb{R}^8$ to rotate by~$\pi$ radians about the~$s$-axis.  Each rotation is indeed a reflection across the~$(1,s)$ plane.  Hence, in order to rotate the~$(r,s)$ plane through~$\alpha$ radians, where~$r,s \in \mathbb{O}$ are orthogonal imaginary unit octonions, we compose two reflections using the map~$f_{r,s,\alpha /2} : \mathbb{O} \to \mathbb{O}$
given by 
$$f_{r,s,\alpha/2}(q) = \left( \cos \left(\alpha/2 \right) r + \sin \left(\alpha/2 \right) s \right) \left( r  q \overline{r}\right) \left( \overline{ \cos \left(\alpha/2 \right) r + \sin \left(\alpha/2 \right) s } \right)$$
The first conjugation~$q \to r q \overline{r}$ in the map above reflects~$\mathbb{R}^8$ about the direction~$r$, while the  second conjugation map reflects all of~$\mathbb{R}^8$ back across the direction~$\cos \left(\alpha/2 \right)r + \sin \left(\alpha/2\right) s$.  Together, these two reflections cause the~$(r,s)$ plane to be rotated by~$ \alpha$ radians.  In the directions orthogonal to this plane, the composition of the two reflections is the identity transformation.  We shall use this map extensively in Chapter \ref{ch:E6_basic_structure} to construct a basis for~$SL(3,\mathbb{O})$.

Finally, notice that since the octonions~$r$ and~$\cos \left(\alpha /2 \right) r + \sin \left( \alpha /2 \right) s$ are imaginary, then ~$\overline{r} = -r$ and~$\overline{\cos \left(\alpha /2 \right) r + \sin \left( \alpha /2 \right) s} = -\cos \left(\alpha /2 \right) r - \sin \left( \alpha /2 \right) s$.  The fact that a rotation in a plane spanned by two imaginary orthogonal unit octonions involves two sign changes will be significant in future discussions of triality in Section \ref{ch:E6_basic_structure.Triality}.

\subsection{Triality}
\label{ch:Division_Algebras_Applications.Quaternions_Octonions.Triality}

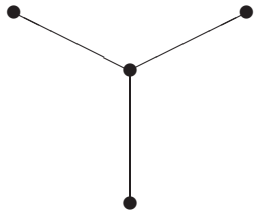
\begin{figure}[htbp]
\begin{center}
\begin{picture}(200,80)
\put(100,50){\circle*{5}}

\put(100,50){\line(2,1){44.72}}
\put(144,72){\circle*{5}}

\put(100,50){\line(0,-1){50}}
\put(100,0){\circle*{5}}

\put(100,50){\line(-2,1){44.72}}
\put(56,72){\circle*{5}}
\end{picture}
\caption{Dynkin Diagram of~$so(8,\mathbb{R})$}
\label{fig:dynkin_diagram_so8}
\end{center}
\end{figure}

In this section, we discuss triality as it relates to octonions and three different representations of~$SO(8,\mathbb{R})$.
  In 1925, Cartan \cite{cartanTriality} noted that the Dynkin diagram~$D$ of~$so(8,\mathbb{R})$, shown in Figure \ref{fig:dynkin_diagram_so8}, contained a three-fold symmetry~$\tau$ which satisfies~$\tau^2 = \tau^{-1}$ and~\hbox{$\tau^3 = Id$}.  He referred to these symmetries as ``triality'', and noted that~$\tau$ was an outer automorphism, not an inner automorphism, of~$so(8,\mathbb{R})$.
  In \cite{baez}, Baez mentions that each of the exterior nodes in the Dynkin diagram of~$so(8,\mathbb{R})$ may be identified with a spinor, dual-spinor, and vector representation of~$SO(8,\mathbb{R})$.  Hence, the triality~$\tau$ should provide an automorphism between these three representations of~$SO(8,\mathbb{R})$.
  Baez, however, discusses triality for the most part as a property of the division algebras.  After looking at that version of triality, we spend the remainder of this section discussing the more specific version of triality related to octonions given by Conway in \cite{conway}.  This idea of triality will be visited again in Chapter \ref{ch:E6_basic_structure.Triality} as we discuss our real form of~$E_6$.
  Additional information regarding triality and octonions may be found in \cite{baez, conway, schafer_book, jacobsen}, although we note that Jacobsen \cite{jacobsen} and Schafer \cite{schafer_book} give an infinitesimal version of triality which is less helpful for this work and hence we do not discuss their treatments of triality.

Given inner product spaces~$V_1, V_2,$ and~$V_3$, all of the same dimension, Baez \cite{baez} defines a {\it normed triality} as a trilinear map
$$ t: V_1 \times V_2 \times V_3 \to \mathbb{R}$$
where
$$ |t(v_1, v_2, v_3)| \le ||v_1|| ||v_2|| ||v_3|| $$
and such that for all~$v_1, v_2$ there exists~$v_3 \ne 0$ for which this bound is attained - and similarly for cyclic permutations of~$1, 2, 3$.
  This map may be dualized, resulting in a bilinear map
$$ m : V_1 \times V_2 \to V_3^*$$
which Baez conveniently calls a ``multiplication''. Using this dualized normed triality map, left multiplication by~$v_1 \in V_1, v_1 \ne 0$ gives an isomorphism $V_2 \to V_3^*$ while right multiplication by~$v_2 \in V_2, v_2 \ne 0$ gives an isomorphism $V_1 \to V_3^*$.  Similarly, appropriate results hold for cyclic permutations of~$1,2,3.$ Baez notes that with a normed triality, if unit vectors are picked from any two of the three spaces $V_1$, $V_2$, or $V_3$, then all three spaces may be identified giving a normed division algebra.  Evidently, every normed division algebra has a normed triality.

We are especially interested in the triality associated with the octonions, and find the particular treatment given in Conway \cite{conway} to be useful.  Conway discusses triality in relation to octonions and $SO(8,\mathbb{R})$, and so we must first indicate how octonions may be used to construct three different representations of~$SO(8,\mathbb{R})$.
  For any unit~$a \in \mathbb{O}$, let~$L_a$,~$R_a$, and~$S_a$ denote the left, right, and symmetric multiplication maps from~$\mathbb{O}$ to~$\mathbb{O}$ given by
$$ L_a(x) = ax \hspace{1cm} R_a(x) = xa \hspace{1cm} S_a(x) = axa \hspace{1cm} (x \in \mathbb{O})$$
As~$\mathbb{O}$ may be identified with~$\mathbb{R}^8$, and these multiplications are invertible and preserve the norm of~$x$, it can be shown \cite{conway} that~$\left( \mathbb{O}, L \right)$,~$\left( \mathbb{O}, R \right)$, and~$\left( \mathbb{O}, S \right)$ are representations of~$SO(8,\mathbb{R})$.  These representations are often called the ``spinor,'' ``dual-spinor,'' and ``vector'' representations due to the multiplication utilized in each space.  Denote a general transformation on~$\mathbb{O}$ induced by a sequence of left, right, and symmetric multiplications of unit octonions by~$L_\alpha$,~$R_\beta$, and~$S_\gamma$, respectively.  That is,
$$ L_\alpha = L_{a_1}  \cdots L_{a_{m_1}} \hspace{1cm} R_\beta = R_{b_1}  \cdots R_{b_{m_2}} \hspace{1cm} S_\gamma = S_{c_1} \cdots S_{c_{m_3}}~$$
for~$a_i, b_i, c_i \in \mathbb{O}$ with~$|a_i| = |b_i| = |c_i| = 1$.  
We use~$L_{-\alpha}$ to denote the transformation given by~$L_{-a_1} \cdots L_{-a_{m_1}}$. Similarly,~$R_{-\beta} = R_{-b_1} \cdots R_{-b_{m_2}}$ and~$S_{-\gamma} = S_{-c_1} \cdots S_{-c_{m_3}}$, implying~$S_{-\gamma} = S_{\gamma}$.

We are now prepared to define the specific version of triality which is helpful in the octonionic case, as presented in \cite{conway}.
Suppose~$x,y,z \in \mathbb{O}$ satisfy~$xy = z$.  If~$S_\gamma$ is any transformation in~$\left( \mathbb{O}, S \right)$, 
 then there exists 
\hbox{$L_\alpha \in \left( \mathbb{O}, L \right)$} and~$R_\beta \in \left( \mathbb{O}, R \right)$
 for which~$L_\alpha(x) R_\beta(y) = S_\gamma(z)$.  Moreover,~$\left( L_\alpha, R_\beta \right)$ and~$\left(L_{-\alpha}, R_{-\beta} \right)$ are the only choices from~$(\mathbb{O},L)\times (\mathbb{O},R)$ which satisfy this equation.  Hence, for each representation of~$SO(8,\mathbb{R})$, we instead consider~$PSO(8,\mathbb{R})$ by writing~$[L_\alpha] = \lbrace L_\alpha, L_{-\alpha} \rbrace$ to signify that~$L_\alpha$ is identified with ~$L_{-\alpha}$, and similarly for~$[R_\beta]$ and~$[S_\gamma]$. 
  In \cite{conway}, Conway states that $PSO(8, \mathbb{R})$ has an outer automorphism~$\tau$ of order 3 such that~$\tau \left( [L_\alpha], [R_\beta], [S_\gamma] \right) = \left( [R_\beta], [S_\gamma], [L_\alpha] \right)$.  This means that the general transformation~$[L_\alpha] \in PSO(8, \mathbb{R})$ may be written as either a composition of symmetric multiplications~$[S_\gamma] \in PSO(8,\mathbb{R})$ or, using~$\tau^2$, as a composition of right multiplications~$[R_\beta] \in PSO(8,\mathbb{R})$.  This triality map~$\tau$ on the three representations of~$PSO(8,\mathbb{R})$ induces a map
$$ SO(8,\mathbb{R}) \to SO(8,\mathbb{R}) \to SO(8,\mathbb{R})$$
between the spaces
$$ (\mathbb{O}, L) \to (\mathbb{O},R) \to (\mathbb{O},S) $$
generated by 
$$ L_\alpha \to R_\beta \to S_{\gamma}$$
In the case of~$(\mathbb{O}, L), (\mathbb{O}, R)$, and~$(\mathbb{O},S)$, this not only implies that each left multiplication may be written as the composition of right multiplications, but Conway shows that the general left multiplication requires an expression involving exactly seven right right multiplications \cite{conway}.
This form of triality will prove useful in Chapter \ref{ch:E6_basic_structure}, and in particular in Section \ref{ch:E6_basic_structure.Triality}, as we analyze the structure of~$SL(3,\mathbb{O})$.

\section{Lorentz Transformations}
\label{ch:Division_Algebras_Applications.Lorentz_Transformations}

Informally, a~$(k+1)$-dimensional spacetime consists of one time-like direction and~$k$ space-like directions.   Spacetimes for certain values of~$k$ naturally match physical phenomena, for instance,~$(3+1)$-dimensional spacetime, called {\it Minkowski spacetime}, is a natural place for working with Einstein's theory of special relativity.  As noted in \cite{kugo_townsend, manogue_fairlie_composite_string, manogue_fairlie_covariant_superstring, manogue_sudbery_gen_sol}, each of the division algebras~$\mathbb{K} = \mathbb{R}, \mathbb{C}, \mathbb{H}, \mathbb{O}$ give rise to particular \hbox{$(k+1)$-dimensional} spacetimes applicable to physics, where~$k-1 = |\mathbb{K}|$.  In this section, we first show how~$(k+1)$-dimensional spacetimes may be represented using~$2 \times 2$ Hermitian matrices over~$\mathbb{K}$ for~\hbox{$k-1 = |\mathbb{K}|$}.  We then show how the division algebras may be used to construct isometries in these spacetimes.  All of the concepts in Section \ref{ch:Division_Algebras_Applications.Lorentz_Transformations.General_Lorentz_Transformations} through Section \ref{ch:Division_Algebras_Applications.Lorentz_Transformations.Lorentz_Transformations_for_Spacetimes_involving_Division_Algebras}
 follow the treatment given in \cite{manogue_schray}, which is fundamental to the construction of~$SL(3,\mathbb{O})$ in Chapter \ref{ch:E6_basic_structure}.

\subsection{Spacetime and Lorentz Transformations}
\label{ch:Division_Algebras_Applications.Lorentz_Transformations.General_Lorentz_Transformations}

Let~$V$ be a~$(k+1)$-dimensional real vector space, and let
~$\lbrace e_1, \cdots, e_{k+1} \rbrace$ denote the standard ordered orthonormal basis on~$V$.
A~$(k+1)$-dimensional {\it spacetime}~$(V,g)$ is a \hbox{$(k+1)$-dimensional} real vector space~$V$ equipped with a non-degenerate, symmetric, bilinear form~$g$ with signature~$(+, -, -, \cdots, -)$.  
We use~$t$ to signify the light-light coordinate of~$e_1$, and write
$$\vecx = \left(\begin{array}{cccccc} t & x_2 & x_3 & \cdots & x_{k} & x_{k+1} \end{array} \right)^T \in V$$
Then the squared length of~$\vecx$ is given by 
$$ |\vecx |^2 = x^{T} g x = t^2 - x_2^2 - x^3_2 - \cdots - x_{k+1}^2$$
We note that $\vecx \in V$ is called a {\it null vector} if $|\vecx| = 0$.

  A {\it Lorentz Transformation} is a map~$M : V \to V$ which preserves the squared length of~$\vecx$, that is,
$$ |\vecx |^2 = |(M\vecx) |^2$$
Writing~$\vecx$ as a~$(k+1)$ column vector, it is convenient to use matrix multiplication to represent a Lorentz transformation.  Indeed, if~$M_{e_a,e_b} : V \to V$ is a rotation through an angle~$\alpha$ in the plane spanned by~$e_a,e_b$,
then the components of~$M_{e_a,e_b}$ are given by either (assuming~$a < b$)
$$ (M_{e_1,e_b})_{(1,1)} = (M_{e_1,e_b})_{(b,b)} = \cosh \alpha \hspace{1cm} 
(M_{e_1,e_b})_{(1,b)} = (M_{e_1,e_b})_{(b,1)} = \sinh \alpha~$$
$$(M_{e_1,e_b})_{(c,c)} = 1 \textrm{ if } c \ne 1 \textrm{ and } c \ne b$$
 or, when~$a \ne 1$, by 
$$ (M_{e_a,e_b})_{(a,a)} = (M_{e_a,e_b})_{(b,b)} = \cos \alpha \hspace{1cm} 
(M_{e_a,e_b})_{(a,b)} = -(M_{e_a,e_b})_{(b,a)} = -\sin \alpha~$$
$$(M_{e_a,e_b})_{(c,c)} = 1 \textrm{ if } c \ne a \textrm{ and } c \ne b$$
where all other entries are~$0$.
Hence, we see that~$M_{e_1,e_3}$ and~$M_{e_2,e_{k+1}}$ are given by

$$
M_{e_1,e_3} = \left( \begin{array}{cccccc}
\cosh \alpha &    0   & \sinh \alpha & 0 & \cdots  & 0 \\
      0      &    1   &       0      & 0 & \cdots  & 0 \\
\sinh \alpha &    0   & \cosh \alpha & 0 & \cdots  & 0 \\
      0      &    0   &       0      & 1 & \cdots  & 0 \\
   \vdots    & \vdots &    \vdots    & \vdots & \ddots & \vdots \\
      0      &    0   &       0      & 0 & \cdots     & 1 \\
\end{array} \right)$$

$$ M_{e_2,e_{k+1}} = \left( \begin{array}{cccccc} 
      1      &    0          &       0      & \cdots & 0 & 0 \\
      0      & \cos \alpha   &       0      & \cdots & 0 & -\sin \alpha \\
      0      &    0          &       1      & \cdots & 0 & 0 \\
    \vdots   &  \vdots       &    \vdots    & \ddots & \vdots & \vdots \\
      0      &    0          &       0      & \cdots & 1 & 0 \\
      0      & \sin \alpha   &       0      & \cdots & 0 & \cos \alpha \\
\end{array} \right)
$$
Any Lorentz transformation may be written as some composition of the matrices
$$M_{e_{1},e_{2}}, M_{e_{1},e_3}, \cdots, M_{e_1,e_{k+1}}, M_{e_2,e_3}, \cdots, M_{e_k,e_{k+1}}$$
using appropriate angles.  

\subsection{Spacetimes related to Division Algebras}
\label{ch:Division_Algebras_Applications.Lorentz_Transformations.Division_Algebras_and_Spacetimes}

For reasons expressed above, we are particularly interested in~$(k+1)$-dimensional spacetimes where~$k = 2, 3, 5, 9$ is one more than the dimension of~\hbox{$\mathbb{K} = \mathbb{R}, \mathbb{C}, \mathbb{H}, \mathbb{O}$}.  For such values of~$k$,
 it is convenient to write an element of~$(k+1)$-dimensional spacetime using a {\it vector}, which is a~$2 \times 2$ hermitian matrix
$$ \matX = \left( \begin{array}{cc} t+z & \overline{q} \\ q & t-z \end{array} \right)$$
with~$q \in \mathbb{R}, \mathbb{C}, \mathbb{H}, \mathbb{O}$, respectively.
In the complex case, we may identify the element
$$\vecx = \left( \begin{array}{cccc} t & x & y & z \end{array} \right)^T$$
of~$(3+1)$ spacetime dimensions with the vector 
$$ \matX = \left( \begin{array}{cc} t+z & x-iy \\ x+iy & t-z \end{array} \right)$$
by writing
$$\matX = t \sigma_t + x \sigma_x + y \sigma_y + z \sigma_z$$
where~$\sigma_t$ is the~$2 \times 2$ identity matrix and~$\sigma_x$,~$\sigma_y$, and~$\sigma_z$ are the Pauli spin matrices
$$ \sigma_x = \left( \begin{array}{cc} 0 & 1 \\ 1 & 0 \end{array} \right) \hspace{1cm}
\sigma_y = \left( \begin{array}{cc} 0 & -i \\ i & 0 \end{array} \right)\hspace{1cm}
\sigma_z = \left( \begin{array}{cc} 1 & 0 \\ 0 & -1 \end{array} \right)$$
Equipped with a metric~$\langle \; , \; \rangle$ whose value is given by~$\langle \sigma_a, \sigma_b \rangle = \frac{1}{2}\tr{\sigma_a \sigma_b}$, it can be shown that the Pauli matrices along with~$\sigma_t$ form an orthonormal basis of a four-dimensional vector space over~$\mathbb{R}$.
This identification between~$(3+1)$-dimensional spacetime and complex vectors may be extended to the quaternionic or octonionic case, involving spacetimes of dimension~$(5+1)$ and~$(9+1)$, respectively, by 
relabeling~$\sigma_y$ as~$\sigma_i$ and including ~$\sigma_q$, which is the matrix~$\sigma_i$ with~$i$ replaced with~$q = i,j,k$ or~$q = i, \cdots, i\l, \l$, respectively.  

We now show that the squared length of an element in~$(k+1)$-dimensional spacetime may be expressed as the determinant of its corresponding~$2 \times 2$ hermitian matrix, where~$k-1$ is the dimension of one of the division algebras.  We show this for the case involving octonions and~$(9+1)$-dimensional spacetime, and note that the other cases correspond to subalgebras and subspaces of this case.

Let~$q = q_x + q_i i + q_j j + q_k k + q_{k\l}k\l + q_{j\l}j\l + q_{i\l}i\l + q_\l \l \in \mathbb{O}$.  The general octonionic vector
$$ \matX = \left( \begin{array}{cc} t+z & \overline{q} \\ q & t-z \end{array} \right), q \in \mathbb{O}, t,z \in \mathbb{R}$$
satisfies its characteristic equation
$$ \matX^2 - (\tr \matX) \matX + (\det{\matX}) I = 0$$
where~$\tr \matX$ denotes the trace of~$\matX$.  We may solve this equation, giving an expression for~$\det{\matX}$ (See Section \ref{ch:Division_Algebras_Applications.Jordan_Algebras}).  As~$\matX$ is a $2 \times 2$ hermitian matrix and the components of~$\matX$ lie in a complex subalgebra of~$\mathbb{O}$, the expression for the determinant of~$\matX$ takes the familiar form
$$\det{\matX} = (t+z)(t-z) - |q|^2$$
Hence, using the identification between an element
$$\left( \begin{array}{ccccccc} t & q_x & q_i & \cdots & q_{i\l} & q_\l & z \end{array} \right)^T$$
of~$(9+1)$-dimensional spacetime and the octonionic vector~$\matX$ given above, we see that
$$\begin{array}{ccl} |( \begin{array}{ccccccc} t & q_x & q_i & \cdots & q_{i\l} & q_\l & z \end{array})^T|^2 &=& t^2 - q_x^2 - q_i^2 - \cdots - q_{i\l}^2 - q_\l^2 - q^2 \\
&=&  \det{\matX} \end{array}$$
We note that the time-like variable~$t$ is found using 
$$t = \frac{1}{2}\tr{\matX}$$
Similar results can be found for the real, complex, or quaternionic cases involving spacetimes of dimension~$(2+1)$,~$(3+1)$, and~$(5+1)$, respectively.

Finally, we can not discuss vectors for spacetimes associated to the division algebras without mentioning spinors and dual-spinors.  
  A {\it spinor} or {\it Weyl spinor} is an element of~$\mathbb{K} \oplus \mathbb{K}$, for~$\mathbb{K} = \mathbb{R}, \mathbb{C}, \mathbb{H}, \mathbb{O}$, which we conveniently write as a two-component column~$v = \left( \begin{array}{c} b \\ c \end{array} \right)$ with $b,c \in \mathbb{K}$.  A {\it dual spinor} is an element of the vector space dual to~$\mathbb{K} \oplus \mathbb{K}$.  It is convenient to write a dual spinor $v^\dagger = \left( \begin{array}{cc} \overline{b} & \overline{c} \end{array} \right)$ as the hermitian conjugate, or dagger, of a spinor.

Let~$\mathbb{K}$ be one of the four division algebras and let the vector~$v \in \mathbb{K} \oplus \mathbb{K}$ and spinor~$v^\dagger \in (\mathbb{K} \oplus \mathbb{K})^*$ be written as above.  There are two meaningful products which relate spinors to vectors.  The product~$v v^\dagger$ is the~$2 \times 2$ Hermitian matrix
$$ v v^\dagger = \left( \begin{array}{cc} |b|^2 & b\overline{c} \\ c\overline{b} & |c|^2 \end{array} \right)$$
called the {\it square of a spinor}.  Using the correspondence previously described in this chapter, we see that~$v v^\dagger$ is an element of a~$(k+1)$-dimensional spacetime where~$k-1 = |\mathbb{K}|$.  However, as 
$$ \begin{array}{ccl}
\mydet{ v v^\dagger } &=& \det{ \left( \begin{array}{cc} |b|^2 & b\overline{c} \\ c\overline{b} & |c|^2 \end{array} \right)} \\
& = & |b|^2 |c|^2 - (b\overline{c})\overline{(b\overline{c})} \\
& = & |b|^2 |c|^2 - |bc|^2  \\
& = & 0\\
\end{array}$$
we see that a squared spinor corresponds to a null vector in~$(k+1)$-dimensional spacetime.  The other product,~$v^\dagger v = |b|^2 + |c|^2$, gives twice the time-like coordinate of this null vector.

\subsection{Lorentz Transformations related to Division Algebras}
\label{ch:Division_Algebras_Applications.Lorentz_Transformations.Lorentz_Transformations_for_Spacetimes_involving_Division_Algebras}

For~$k=2,3,5,9$, the elements of~$(k+1)$-dimensional spacetime may be represented using real, complex, quaternionic, or octonionic vectors.  The final step is to construct Lorentz transformations.  As we shall see, such transformations will make use of the division algebras and the geometric transformations associated with its multiplication, as reviewed in Section \ref{ch:Division_Algebras_Applications.Quaternions_Octonions.Conjugation_Reflecitons_Rotations}.  This treatment again follows the work of \cite{manogue_schray}.  The constructions and notation in this section are fundamental to the construction of~$SL(3,\mathbb{O})$ in Chapter \ref{ch:E6_basic_structure}.

For the general case, let~$\matX$ be an octonionic vector.  Since vectors are hermitian matrices, a Lorentz transformation must preserve both the hermiticity of~$\matX$ as well as~$\det{\matX}$.  For general~$M$, neither one-sided multiplication~\hbox{$\matX \to M \matX$} nor~\hbox{$\matX \to \matX M$} preserves the hermiticity of~$\matX$.  However, the result of a transformation using conjugation
$$\matX \to M \matX M^\dagger$$
is a hermitian matrix for any matrix~$M$, provided the result is well-defined.
As shown in \cite{manogue_schray}, the product $M \matX M^\dagger := \left(M \matX \right) M^\dagger = M \left( \matX M^\dagger \right)$ is well-defined if either the components of $M$ lie in a complex subspace of $\mathbb{O}$ or the columns of the imaginary part of $M$ must be real multiples of each other.
In addition, the squared length of~$\vecx$ corresponding to~$\matX$ will be preserved if~$\det{\matX} = \mydet{M \matX M^\dagger}$.  For such~$M$, we also have
$$\mydet{M \matX M^\dagger} = \mydet{M M^\dagger} \mydet{\matX}$$
requiring that the matrices~$M$ also satisfy the condition\footnote{Over the complex numbers, we can multiply~$M$ by an overall phase~$e^{i\alpha}$ without changing~$|\det{M}|$ or~$\det{(M\matX M^{\dagger})}$, and hence we require that~$\det{M} = 1$.  Thus, the Lorentz group~$SO(3,1)$ is essentially~$SL(2,\mathbf{C})$.}
$$\mydet{M M^\dagger} = 1$$

In the case of complex vectors~$\matX \in M_{2}(\mathbb{C})$, we choose basis transformations 
$$R_{ab}: \mathbb{R} \times M_2(\mathbb{C}) \to M_2(\mathbb{C})$$
where
$$R_{ab}(\alpha)(\matX) = M_{a,b}(\frac{\alpha}{2}) \; \matX \; M_{a,b}^\dagger (\frac{\alpha}{2})$$
rotates the~$(a,b)$ plane through an angle~$\alpha$, for~$a,b \in \lbrace t,x,y,z \rbrace$.  For rotations in the~$(x,y)$,~$(y,z)$, and~$(z,x)$ planes, the matrices~$M_{a,b}(\alpha)$ are given by:

\[
M_{x,y}(\alpha) = \left( \begin{array}{cc} e^{i\alpha} & 0 \\ 0 & e^{-i\alpha}\end{array}\right)
\hspace{1cm}
M_{y,z}(\alpha) = \left( \begin{array}{cc} \cos \alpha & i\sin \alpha \\ i \sin \alpha & \cos \alpha\end{array}\right)
\]
\[ M_{z,x}(\alpha) = \left( \begin{array}{cc} \cos \alpha & \sin \alpha \\ -\sin \alpha & \cos \alpha\end{array}\right) \]

\noindent while the~$M_{a,b}(\alpha)$ matrices with~$a=t$ are given by

\[
M_{t,z}(\alpha) = \left( \begin{array}{cc} e^{\alpha} & 0 \\ 0 & e^{-\alpha}\end{array}\right)
\hspace{1cm}
M_{t,x}(\alpha) = \left( \begin{array}{cc} \cosh \alpha & \sinh \alpha \\ \sinh \alpha & \cosh \alpha\end{array}\right) 
\]
\[M_{t,y}(\alpha) = \left( \begin{array}{cc} \cosh \alpha & i\sinh \alpha \\ -i\sinh \alpha & \cosh \alpha\end{array}\right) \]

\noindent  
We find it convenient to call the Lorentz transformations which change the time-like coordinate {\it boosts}.  We label those transformations not with the label~$R_{t,b}(\alpha)$ but with the label~$B_{t,b}(\alpha)$.

These transformations also work in the quaternionic case~$\matX \in M_2(\mathbb{H})$ with a slight modification.  In this larger setting, we are now working with three imaginary units,~$i$,~$j$, and~$k$, instead of just one,~$i$.  Hence, it becomes helpful to relabel the Lorentz transformations~$R_{xy}$,~$R_{yz}$, and~$R_{ty}$ as~$R_{xi}$,~$R_{iz}$, and~$R_{ti}$.  We create six new Lorentz transformations by replacing the~$i$ in the subscripts for the three transformations~$R_{ab}(\alpha)$ and~$M_{a,b}(\alpha)$ with~$j$ and~$k$.  This gives 12 Lorentz transformations, each of which rotates a plane containing the~$t$,~$x$, or~$z$ direction.  Finally, it is only necessary to include the rotations which rotate a plane spanned by two imaginary quaternions in~$5+1$ dimensions.  However, in Section \ref{ch:Division_Algebras_Applications.Quaternions_Octonions.Conjugation_Reflecitons_Rotations}, we saw that the conjugation map
$$q \to e^{i\alpha/2} q e^{-i\alpha/2}$$
fixes the~$(1,i)$ plane and produces a rotation through an angle~$\alpha$ in the~$(j,k)$ plane.  Thus, for~$s = p \, q$ for~$p,q \in \lbrace i,j,k \rbrace, p\ne q$, the Lorentz transformation~\hbox{$R_{p,q}: \mathbb{R} \times M_2(\mathbb{C}) \to M_2(\mathbb{C})$} 
given by 
$$ R_{p,q}(\alpha)(\matX) = \left( \begin{array}{cc} \exp({s \; \alpha/2}) & 0 \\ 0 & \exp({s \; \alpha/2}) \end{array} \right)
\matX \left( \begin{array}{cc} \exp({s \; \alpha/2}) & 0 \\ 0 & \exp({r\; \alpha/2}) \end{array} \right)^\dagger$$
produces a rotation in the~$(p,q)$ plane through an angle of~$\alpha$.  The transformations involving two imaginary units are called {\it transverse rotations}.

In the octonionic case, Manogue and Schray \cite{manogue_schray} provided an explicit description of vector Lorentz transformations in~$(9+1)$ spacetime dimensions by generalizing the above results.  The boosts and non-transverse rotations generalize to this case in the obvious way, by extending the transformations~$B_{tq}(\alpha), R_{xq}(\alpha),$ and~$R_{zq}(\alpha)$ from~$q = i,j,k$ to now~$q = i,j,k,k\l,j\l,i\l,\l$.  More care must be taken to generalize the transverse rotations, as we saw in Section \ref{ch:Division_Algebras_Applications.Quaternions_Octonions.Conjugation_Reflecitons_Rotations} that the conjugation map~$q \to e^{i\alpha} q e^{-i\alpha}$ rotates three planes simultaneously instead of just one plane in the octonionic setting.  As discussed in Section \ref{ch:Division_Algebras_Applications.Quaternions_Octonions.Conjugation_Reflecitons_Rotations}, nested flips may be used to produce a rotation spanned by two orthogonal pure imaginary units~$p$ and~$q$.  That is, the transformation~$R_{p,q}: \mathbb{R} \times M_2(\mathbb{O}) \to M_2(\mathbb{O})$ given by
$$ R_{p,q}(\alpha) = M_{2(p,q)}\left(\alpha/2\right) \left( M_{1(p)} \matX M_{1(p)}^\dagger \right) M_{2(p,q)}\left(\alpha/2\right)^\dagger$$
where 
$$\begin{array}{ccl}  M_{1(p)} & = & \left( \begin{array}{cc} -p & 0 \\ 0 & -p \end{array} \right) \\
                      M_{2(p,q)}(\alpha/2) & = & \left( \begin{array}{cc} \cos(\alpha/2)\; p + \sin(\alpha/2) \;q & 0 \\
                                                                      0 & \cos(\alpha/2)\; p + \sin(\alpha/2)\; q 
                                                    \end{array} \right)
\end{array}
$$
is the transverse rotation in the plane defined by two orthogonal imaginary units~$p, q \in \mathbb{O}$.
  There are~$21$ different ways to pick pairs of these units from the seven orthogonal imaginary units in~$\mathbb{O}$, adding~$21$ transverse rotations to the~$15$ general rotations and~$9$ boosts.

Manogue and Schray \cite{manogue_schray} identified these~$45$ transformations as generators of the group $SO(9,1)$ and found several important subgroups using their finite Lorentz transformations.  The~$21$ transverse rotations produce all rotations among the various plains spanned by the any two of the seven imaginary orthogonal units of~$\mathbb{O}$, generating the group~$SO(7)$.  The group~$SO(8)$ is obtained by adding the seven~$R_{x,q}$ transformations, where~$q = i,j,k,k\l,j\l,i\l,\l$, while adding the remaining (non-boost) rotations gives~$SO(9)$.

In Table \ref{table:finite_octonionic_lorentz_transformations}, we list the finite Lorentz transformations as categorized by Manogue and Schray.  The Category 1 transformations are boosts, while the Category 2 transformations are rotations in a plane containing either the~$x$ or~$z$ axis.  The Category 3 transformations are the transverse rotations in the plane defined by two imaginary units~$p, q \in \mathbb{O}$ where~$p \perp q$.  We note that each of these basic transformations $R_{ab}$ only affect the~$(a,b)$ plane in~$9+1$-dimensional spacetime.  Hence,~$R_{ab}(\alpha)R_{ab}(\gamma) = R_{ab}(\alpha + \gamma)$ for any~$\alpha, \gamma \in \mathbb{R}$.  That is, each of the~$45$ basic transformations of~$SL(2,\mathbb{O})$ listed in Table \ref{table:finite_octonionic_lorentz_transformations} are generators of a one-parameter subgroup of~$SL(2,\mathbb{O})$.

\begin{table}[hbtp]
\begin{center}
\begin{minipage}{4.75in}
\[
\matX \to \left\{ \begin{array}{ccl}
M\matX M^\dagger, &\hspace{1cm} & \textrm{ for Categories 1 and 2} \\
M_2 \left( M_1 \matX M_1^\dagger \right) M_2^\dagger, & & \textrm{ for Category 3}
\end{array}
\right.
\]

$$
\begin{array}[t]{|l|c|}
\hline
\begin{minipage}[t]{0.85in}\textrm{Category 1:\\ Boosts}\end{minipage} &
\begin{array}{ccl}
B_{tz} & t \longleftrightarrow z & M = \left(
  \begin{array}{cc} \exp \left(\frac{\alpha}{2} \right) & 0 \\
                                 0         & \exp \left(-\frac{\alpha}{2} \right) \\
 \end{array}
 \right) \\
B_{tx} & t \longleftrightarrow x & M = \left( 
    \begin{array}{cc} \cosh \left( \frac{\alpha}{2} \right) & \sinh \left( \frac{\alpha}{2} \right) \\
                      \sinh \left( \frac{\alpha}{2} \right) & \cosh \left( \frac{\alpha}{2} \right)
    \end{array}                             
\right) \\
B_{tq} & t \longleftrightarrow q & M = \left( 
\begin{array}{cc} \cosh \left( \frac{\alpha}{2} \right)  & q \sinh \left( \frac{\alpha}{2} \right) \\
                  -q \sinh \left( \frac{\alpha}{2} \right) & \cosh \left( \frac{\alpha}{2} \right) 
\end{array}
\right) \\
\end{array} \\
\hline
\hline
\begin{minipage}[b]{0.85in} Category~2:\\ Rotations \end{minipage} & 
\begin{array}{ccl}
R_{xq} & x \longleftrightarrow q & M = \left(
\begin{array}{cc}
\exp \left( q \frac{\alpha}{2} \right) &          0 \\
                  0                      & \exp \left( -q \frac{\alpha}{2} \right)
\end{array} \right) \\
R_{xz} & x \longleftrightarrow z & M = \left(
\begin{array}{cc}
\cos \left( \frac{\alpha}{2} \right) & -\sin \left(\frac{\alpha}{2} \right) \\
\sin \left( \frac{\alpha}{2}\right) & \cos \left(\frac{\alpha}{2} \right) 
\end{array} \right) \\
R_{zq} & q \longleftrightarrow z & M = \left(
\begin{array}{cc}
\cos \left( \frac{\alpha}{2} \right) & q \sin \left(\frac{\alpha}{2} \right) \\
q \sin \left( \frac{\alpha}{2}\right) & \cos \left(\frac{\alpha}{2} \right) 
\end{array} \right) \\
\end{array} \\
\hline
\hline
\begin{minipage}[c]{0.85in}Category 3:\\ Transverse\\[-.66em] Rotations \\ \end{minipage} & 
\begin{array}{ccl}
R_{p,q} & p \longleftrightarrow q & M_1 = -p \; \mathbb{I}_2 \\
        &                & M_2 = \left( \; \cos \left(\frac{\alpha}{2}\right) p + \sin \left( \frac{\alpha}{2} \right) q \; \right) \mathbb{I}_2
\end{array} \\
\hline
\end{array}
$$
\end{minipage}
\caption{Finite octonionic Lorentz transformations}
\label{table:finite_octonionic_lorentz_transformations}
\end{center}
\end{table}

Finally, we note that Manogue and Schray carefully constructed their finite Lorentz transformations so that they were compatible, meaning that the parentheses may be re-arranged when it is applied to a squared spinor~$v v^\dagger$.  More precisely, a transformation~$\matX \to M \matX M^\dagger$ is {\it compatible } if
$$M (v v^\dagger) M^\dagger = (M v) (v^\dagger M^\dagger) = (M v) (M v)^\dagger$$
  This re-arrangement of parenthesis is possible for octonionic spinors~$v$ provided that the elements of~$M$ lie in a complex subalgebra \cite{manogue_schray}.

\noindent We use compatible vector finite Lorentz transformations to give spinor and dual spinor transformations simply by using the same matrices but applying only left-multiplication or right-multiplication, respectively, instead of conjugation.  That is, if~$R: \mathbb{R} \times M_2(\mathbb{O}) \to M_2(\mathbb{O})$ is a transformation given by
$$R(\alpha)(\matX) = M_n(\alpha)\left( \cdots \left(M_1(\alpha) \matX M_1^\dagger(\alpha) \right) \cdots \right) M_n^\dagger(\alpha)$$
then the spinor and dual spinor versions of~$R$ are then $R: \mathbb{R} \times (\mathbb{O}\oplus \mathbb{O}) \to (\mathbb{O}\oplus \mathbb{O})$ and $R: \mathbb{R} \times (\mathbb{O}\oplus \mathbb{O})^* \to (\mathbb{O}\oplus \mathbb{O})^*$,
respectively, where
$$R(\alpha)(\theta) = M_n(\alpha)\left( \cdots \left(M_1(\alpha) \; \theta \right) \cdots \right) $$
and
$$R(\alpha)(\theta^\dagger) = \left( \cdots \left(\theta^\dagger \; M_1^\dagger(\alpha) \right) \cdots \right) M_n(\alpha)^\dagger$$
We deliberately give the spinor and dual spinor versions of the transformations the same name as the vector Lorentz transformation, as the context will make it clear which version is appropriate.  We utilize the vector, spinor, and dual spinor versions of each of the finite Lorentz transformations in our construction of~$E_6 = SL(3,\mathbb{O})$ in Chapter \ref{ch:E6_basic_structure}.

\section{Jordan Algebras and Albert Algebras}
\label{ch:Division_Algebras_Applications.Jordan_Algebras}

We give here the definition of the Albert algebra, which is the exceptional Jordan algebra.  One particularly nice characteristic of the Albert algebra is that it naturally contains both spinors and vectors.  We develop the invariants for the Albert algebra corresponding to the trace and determinant of a matrix, and discuss how these invariants are tied to certain Lie groups.  Additional information on Jordan and Albert algebras, including their role in physics, may be found in \cite{okubo, manogue_dray_exceptional_eigenvalue_problem}.

A {\it Jordan algebra } is a commutative algebra~$A$ with product denoted by~$\circ$ that satisfies the Jordan identity
$$ (x^2 \circ y) \circ x = x^2\circ (y\circ x)$$
for any~$x,y \in A$, where~$x^2 \equiv x \circ x$.  Given an algebra~$A$ with a bi-linear product~$xy$, we can create a new algebra~$(A,\circ)$ by introducing the commutative product
$$ x \circ y = \frac{1}{2} \left(xy + yx \right)$$
The algebra~$(A,\circ)$ may or may not satisfy the Jordan identity.  If~$A$ is an associative algebra, then~$(A,\circ)$ will always satisfy the Jordan identity and is called a {\it special Jordan algebra}.  In 1934, Jordan, von Neumann, and Wigner \cite{jordan_von_neumann_wigner} proved that all finite-dimensional simple Jordan algebras are special Jordan algebras except for one case, known as the exceptional Jordan algebra.

The exceptional Jordan algebra~$\jordanalgebra$, also known as the Albert algebra, consists of~$3 \times 3$ octonionic hermitian matrices with the product
$$ X \circ Y = \frac{1}{2}\left(XY + YX \right)$$
where~$XY$ is matrix multiplication,~$X,Y \in \jordanalgebra$.
Although this product is commutative, it is not associative as
$$
\begin{array}{l}
(X \circ Y) \circ Z - X \circ (Y \circ Z) =  \\
\hspace{.6cm} \frac{1}{4} \left( (XY)Z + (YX)Z + Z(XY) + Z(YX) \right. \\ 
\hspace{1.2cm}  \left. -X(YZ) - X(ZY) - (YZ)X - (ZY)X 
\right)
\end{array}
$$
Nevertheless, the exceptional Jordan algebra does satisfy the Jordan identity.
Now, for $X \in \jordanalgebra$, define~$X^2 := (X \circ X)$ and~$X^3 := X^2 \circ X \equiv X \circ X^2$.

A general element of~$\chi = \jordanalgebra$ may be written in components as

\[
\chi = \left( \begin{array}{ccc}
p & \overline{a} & c \\
a & m & \overline{b} \\
\overline{c} & b & n \\
\end{array} \right)
\]

\noindent 
where~$p, m, n \in \mathbb{R}$ and~$a,b,c \in \mathbb{O}$.
  There are three natural ways to identify a~$2 \times 2$ Hermitian octonionic vector~$\matX$, 2-component octonionic spinor~$\theta$, and 2-component octonionic dual spinor~$\theta^\dagger$ in~$\chi$.  We will find it convenient in Chapter \ref{ch:E6_basic_structure} to use the following identification of~$\matX$,~$\theta$, and~$\theta^\dagger$ in~$\chi$:
\[
\chi = \left( \begin{array}{c|c}
\matX & \theta \\
\hline
\theta^\dagger & n
\end{array} \right) = 
\left( \begin{array}{cc|c}
p & \overline{a} & c \\
a & m & \overline{b} \\
\hline
\overline{c} & b & n
\end{array} \right)
\]

Standard expressions for the determinant of a matrix rely upon the associativity of the underlying division algebra, but a general element of the exceptional Jordan algebra contains three independent octonionic directions~$a$,~$b$, and~$c$.  Hence, except for very particular choices of~$a, b, c \in \mathbb{O}$, standard expressions of the determinant of a matrix fail to hold for \hbox{$\chi \in M_3(\mathbb{O})$}. 
Nevertheless, octonionic~$3 \times 3$ matrices satisfy their characteristic equation.  We now produce an equivalent expression of the determinant for~$3 \times 3$ octonionic hermitian matrices.

  The characteristic equation of a~$1 \times 1$ matrix~$A$ is 
$$ A - \sigma_1(A) = 0$$
Solving for~$\sigma_1(A)$ and taking the trace, we see that~$\tr{A} = \sigma_1(A)$ and note that for~$1 \times 1$ matrices,~$\det{A} = \sigma_1(A)$.

For the case when~$A$ is a~$2 \times 2$ matrix, the characteristic equation is
$$ A^2 - A \; \tr{A} + \sigma_2(A) \; \mathbb{I}_2 = 0$$
Again, we take the trace giving 
$$ \tr{A^2} - (\tr{A})^2 + 2\sigma_2(A) = 0$$
Solving for~$\sigma_2(A)$, we obtain an expression for the~$2 \times 2$ invariant~$\sigma$:
$$ \sigma(A) := \frac{1}{2}\left( \left( \tr{A}\right)^2 - \tr{A^2}\right)$$
We note that for~$2 \times 2$ matrices, the determinant of~$A$ is~$\sigma_2(A)$.

If~$A$ is a~$3 \times 3$ matrix, then the characteristic equation is given by 
$$ A^3 - A^2 \; \tr{A} + A \; \sigma(A) - \det{A} \;\mathbb{I}_3 = 0$$
Again, taking the trace and solving for~$\det{A}$ finds
$$\det(A) = \frac{1}{3}\left( \tr{A^3} - \frac{3}{2}\tr{A^2}\tr{A} + \frac{1}{2}\left(\tr{A}\right)^3 \right)$$ 
where we have substituted in the expression for~$\sigma(A)$ previously found.
This expression for~$\det(A)$ gives the determinant of a~$3 \times 3$ matrix in the special case that~$A \in \jordanalgebra$. \footnote{
The general pattern for the characteristic equation of an~$n \times n$ matrix is given by 
$$ A^n - A^{n-1}\epse_1(A) + A^{n-2}\epse_2(A) - A^{n-3}\epse_3(A) + \cdots \pm A^{1}\epsilon_{n-1}(A) \mp \det(A)\mathbb{I}_n = 0$$
where~$\epse_n(A)$ is the expression found when solving the~$n\times n$ case for~$\det{A}$.}

Equivalently, we find that the determinant of 
\[
\chi = \left( \begin{array}{ccc}
p & \overline{a} & c \\
a & m & \overline{b} \\
\overline{c} & b & n \\
\end{array} \right)
\in \jordanalgebra
\]
is given by 
$$ \det{\chi} = pmn - (p|b|^2 + m|c|^2 + n|a|^2) + 2Re(\overline{a} \overline{b} \overline{c})$$
This may be rewritten in the more convenient form
$$ \det{ \left( \chi \right) } = \det{ \left( \begin{array}{cc} p & \overline{a} \\ a & m \end{array} \right)} n + 2\left( \begin{array}{cc} p & \overline{a} \\ a & m \end{array} \right) \cdot \left[ \left( \begin{array}{c} c \\ \overline{b} \end{array} \right) \left( \begin{array}{cc} \overline{c} & b \end{array} \right) \right] $$
where $A \cdot B = \tr{A \jordanproduct B} - \tr{A} \tr{B}$ is the Lorentzian inner product of vectors $A$ and $B$.  
Of course, the trace of~$\chi$ is~$\tr{\chi} = p + m + n$.
Just as in the case of~$n \times n$ real, complex, or quaternionic matrices, the expressions for the trace and determinant of a Jordan matrix will be helpful in defining certain Lie groups.  In particular, we will see in the next chapter that the exceptional Jordan algebra is closely tied to the exceptional Lie groups~$G_2$,~$F_4$, and~$E_6$, as well as a property of octonions called triality.

\newpage{}

\part{The Basic Structure of~$E_6$ }
\label{ch:E6_basic_structure}

In this chapter, we explore the Lie group~$SL(3,\mathbb{O})$ and its associated Lie algebra~$sl(3,\mathbb{O})$.  We provide our methods for producing and studying our representation of the~$78$-dimensional exceptional Lie group~$E_6$ and its associated Lie algebra.  This group preserves the determinant of the exceptional Jordan Algebra~$\jordanalgebra$, justifying the interpretation
$$E_6 \equiv SL(3,\mathbb{O})$$
We show how to realize this group by extending the description of
$$SL(2,\mathbb{O}) \equiv SO(9,1,\mathbb{R})$$
given in~\cite{manogue_schray} from the~$2 \times 2$ case to the~$3 \times 3$ case.
  This process produces~$135$ group generators, of which only~$78$ are independent.  We follow the common practice of using linear dependencies among the associated Lie algebra elements to reduce our list to the~$78$ independent group generators, giving us a preferred basis for~$E_6$.

After constructing the group~$SL(3,\mathbb{O})$ in Section~\ref{ch:E6_basic_structure.LieGroup} and discussing how to construct its algebra in Section~\ref{ch:E6_basic_structure.E6_Lie_Algebra}, we study the characteristics of the group and algebra, including the relevant subgroups and subalgebras, in the remainig sections of this chapter.  In particular, we discuss the group~$SL(3,\mathbb{O})$ in Section~\ref{ch:E6_basic_structure.LieGroup}, presenting our construction of the group elements in Section~\ref{ch:E6_basic_structure.Lorentz_Transformations} and giving an improved basis of the~$3\times 3$ transformations generalizing the~$2 \times 2$ transverse rotations in Section~\ref{ch:E6_basic_structure.Lorentz_Transformations.BetterBasis}.  We study the Lie algebra corresponding to this group in Section~\ref{ch:E6_basic_structure.E6_Lie_Algebra}.  In particular, we discuss our methods which allow us to construct the algebra commutator for~$sl(3,\mathbb{O})$ in Section~\ref{ch:E6_basic_structure.ConstructingE6algebra}.  In Section~\ref{ch:E6_basic_structure.ConstructingE6algebra.LinearDependencies}, we use linear dependencies in the algebra to choose a preferred basis, and provide a basis for each of the three natural~$SO(9,1,\mathbb{O})$ subgroups of~$SL(3,\mathbb{O})$ in Section~\ref{ch:E6_basic_structure.Lorentz_Transformations.BetterBasisSubgroups}.  We discuss triality as it relates to~$E_6$ in Section~\ref{ch:E6_basic_structure.Triality}.  In Section~\ref{ch:E6_basic_structure.Type_Transformation}, we give some continuous extensions to the discrete type map used in our construction of~$SL(3,\mathbb{O})$ and find subgroups which contain or are independent of this type map in Section~\ref{ch:E6_basic_structure.Type_Transformation.Subgroups}.  Although our construction relies upon the octonions, we obtain other subgroups of~$SL(3,\mathbb{O})$ by restricting~$\mathbb{O}$ to the other division algebras.  These findings are summarized in Section~\ref{ch:E6_basic_structure.Reduction_of_O}.  Finally, we give an isomorphism between our preferred basis of~$su(3,\mathbb{C}) \subset G_2$ and the Gell-Mann matrices in Section~\ref{ch:E6_basic_structure.fix_l}.

\section{The Lie Group~$E_6 = SL(3,\mathbb{O})$ }
\label{ch:E6_basic_structure.LieGroup}

As mentioned in Section~\ref{ch:Joma_Paper.subalgebras_greater_than_3}, it is known that~$E_6$ is the group which preserves the determinant of the Jordan matrix~$\chi \in M_3(\mathbb{O})$.  It is in this sense that~$E_6$ is~$SL(3,\mathbb{O})$.  Dray and Manogue~\cite{manogue_dray} show~$E_6$ is the appropriate union of three copies of~\hbox{$SO(9,1,\mathbb{R}) = SL(2,\mathbb{O})$}.  In Section~\ref{ch:E6_basic_structure.Lorentz_Transformations}, we explicitly construct three different~$3 \times 3$ transformations from each of the~$2 \times 2$ finite Lorentz transformations which, as Manogue and Schray showed \cite{manogue_schray}, form a basis for~$SL(2,\mathbb{O})$.  These~$2 \times 2$ transformations were also discussed in Section~\ref{ch:Division_Algebras_Applications.Lorentz_Transformations.Lorentz_Transformations_for_Spacetimes_involving_Division_Algebras}.  This construction provides~$135$ transformations in~$E_6$, which has dimension~$78$.  In Section~\ref{ch:E6_basic_structure.Lorentz_Transformations.BetterBasis}, we choose an explicit basis for a preferred subgroup~$SO(7,\mathbb{R})$ of each~$SL(2,\mathbb{O})$.  We find this new basis will be helpful when discussing triality and the subgroup~$G_2$ in relation to~$SL(3,\mathbb{O})$, and will use this basis for the remainder of this thesis after Section~\ref{ch:E6_basic_structure.E6_Lie_Algebra}.

\subsection{Lorentz Transformations in~$3 \times 3$ case}
\label{ch:E6_basic_structure.Lorentz_Transformations}

We first describe the three natural ways to identify a vector~$\matX \in M_2(\mathbb{O})$ in an element~$\chi$ of the exceptional Jordan algebra~$M_3(\mathbb{O})$.  We write
$$\chi = \left( \begin{array}{ccc} t+z & \overline{a} & c \\ a & t-z & \overline{b} \\ \overline{c} & b & n \end{array} \right) \in M_3(\mathbb{O})$$
$$ t,z,n \in \mathbb{R} \hspace{1cm} a, b, c \in \mathbb{O}$$
where we use 
$$  a = a_x + a_i i + a_j j + a_k k + a_{k\l} k\l + a_{j\l} j\l + a_{i\l} i\l + a_\l \l \in \mathbb{O}$$
and similar expressions for~$b$ and~$c$ to denote a general element of~$\mathbb{O}$.  Recall from Section~\ref{ch:Division_Algebras_Applications.Lorentz_Transformations.Division_Algebras_and_Spacetimes} that an octonionic vector~$\matX$ is a~$2 \times 2$ hermitian matrix, an octonionic spinor~\hbox{$\theta \in (\mathbb{O} \oplus \mathbb{O})$} may be written as a two-component column matrix,  and a dual spinor~\hbox{$\theta^\dagger \in (\mathbb{O} \oplus \mathbb{O})^*$} may be written as the hermitian conjugate, or dagger, of a spinor.  There are three natural ways to identify a~$2 \times 2$ hermitian matrix in~$\chi$, as shown in Table~\ref{table:three_types_of_vector_locations}.  In each case, this identification breaks~$\chi$ into a block structure leaving a~$1\times 2$ octonionic submatrix and a~$2\times 1$ octonionic submatrix which are hermitian conjugates of each other. We identify them with the octonionic spinor and dual spinor, respectively.

\begin{table}[tbp]
\begin{center}
\begin{tabular}{ccccc}
\textrm{Location 1} & & \textrm{Location 2} & & \textrm{Location 3}\\
$\left( \begin{array}{c|c} \matX & \theta \\ \hline \theta^\dagger & \cdot \end{array}\right)$ & &
$\left( \begin{array}{c|c} \cdot & \theta^\dagger \\ \hline \theta & \matX \end{array}\right)$ & &
$\left( \begin{array}{c|c|c} \matX_{2,2} & \theta_2 & \matX_{2,1} \\ \hline \overline{\theta_2} & \cdot & \overline{\theta_1} \\ \hline \matX_{1,2} & \theta_{1} & \matX_{1,1} \end{array}\right)$ \\
$\parallel$ & & $\parallel$ & & $\parallel$ \\
$\left( \begin{array}{cc|c} t+z & a & \overline{c} \\ \overline{a} & t-z & b \\ \hline c & \overline{b} & n \end{array} \right)$ & &
$\left( \begin{array}{c|cc} t+z & a & \overline{c} \\ \hline \overline{a} & t-z & b \\ c & \overline{b} & n \end{array} \right)$ & &
$\left( \begin{array}{c|c|c} t+z & a & \overline{c} \\ \hline \overline{a} & t-z & b \\ \hline c & \overline{b} & n \end{array} \right)$ \\
\end{tabular}
\caption{\noindent Three natural locations of~$\matX$,~$\theta$, and~$\theta^\dagger$ in~$\chi$}
\label{table:three_types_of_vector_locations}
\end{center}
\end{table}

Each of the~$2 \times 2$ Lorentz transformations~$R \in SL(2,\mathbb{O})$ given in~\cite{manogue_schray} and used in Section~\ref{ch:Division_Algebras_Applications.Lorentz_Transformations.Lorentz_Transformations_for_Spacetimes_involving_Division_Algebras} are expressed using a nesting of~$2 \times 2$ matrices~$M_1, M_2, \cdots, M_n$.  That is, if~\hbox{$\matX \in M_2(\mathbb{O})$} is a vector, then~$R$ is given by
$$ R(\matX) = M_n \left( \cdots \left( M_2 \left( M_1 \; \matX M_1^\dagger \right) M_2^\dagger \right) \cdots \right) M_n^\dagger$$
As noted in Section~\ref{ch:Division_Algebras_Applications.Lorentz_Transformations.Lorentz_Transformations_for_Spacetimes_involving_Division_Algebras}, the spinor (dual-spinor) version of the Lorentz transformation is constructed from nesting matrices using left (right) multiplication by matrices:
$$ R(\theta) = M_n \left( \cdots \left( M_1 \; \theta \right) \cdots \right) $$
$$ R(\theta^\dagger) = \left( \cdots \left( \theta^\dagger M_1^\dagger \right) \cdots \right) M_n^\dagger$$
If we treat the upper left~$2 \times 2$ submatrix of~$\chi$ as the vector, as in Location 1 in Table~\ref{table:three_types_of_vector_locations}, and~$M$ is a~$2 \times 2$ matrix which we embed into the upper left~$2\times 2$ submatrix of a~$3 \times 3$ matrix~$\mathbb{M}$ which is block diagonal and~$\mathbb{M}_{3,3} = 1$, then we note that conjugating~$\chi$ by~$\mathbb{M}$ gives the following block structure

$$
\begin{array}{ccc}
\chi & \to & \mathbb{M} \chi \mathbb{M}^\dagger \\
\left( \begin{array}{c|c} \matX & \theta \\ \hline \theta^\dagger & n \end{array} \right) & \to & 
\left( \begin{array}{c|c} M  & 0 \\ \hline 0 & 1 \end{array} \right)
\left( \begin{array}{c|c} \matX & \theta \\ \hline \theta^\dagger & n \end{array} \right) 
\left( \begin{array}{c|c} M  & 0 \\ \hline 0 & 1 \end{array} \right)^\dagger \\
   &   & =  \left( \begin{array}{c|c} M \matX M^\dagger & M\theta \\ \hline (M \theta)^\dagger & n \end{array} \right) 
\end{array}
$$
where we must worry about~$M \matX M^\dagger$ being well-defined.
But each of the matrices~$M_a = M_1, \cdots, M_n$ involved in the  expression for the Lorentz transformation~$R$ satisfy certain conditions, outlined in~\cite{manogue_schray}, which ensures~$M \matX M^\dagger$ is well-defined.  
  In particular, constructing the~$3 \times 3$ matrix~$\mathbb{M}_a$ for each~$2 \times 2$ matrix~$M_a$ in this way and conjugating~$\chi$ successively with matrices~$\mathbb{M}_1, \cdots, \mathbb{M}_n$ gives 
$$\chi \to \mathbb{M}_n \left( \cdots \left( \mathbb{M}_1 \chi \mathbb{M}_1^\dagger \right) \cdots \right) \mathbb{M}_n^\dagger$$
The resulting matrix has the block structure
$$
\begin{array}{c}
  \left( \begin{array}{c|c} M_n \left( \cdots \left( M_1 \matX M_1^\dagger\right) \cdots \right) M_n^\dagger & M_n \left( \cdots \left( M_1 \theta \right) \cdots \right) \\ \hline \left(M_n \left( \cdots \left( M_1 \theta \right) \cdots \right)\right)^\dagger & n \end{array} \right) \\
= \left( \begin{array}{c|c} R(\matX) & R(\theta) \\ \hline R(\theta^\dagger) & n \end{array} \right)
\end{array}
$$
in which the vector, spinor, and dual spinor have all transformed according to their respective versions of the Lorentz transformation~$R$.

We now generalize this previous construction for all three natural locations of the vector~$\matX$ in~$\chi$.  Let~$M_1, \cdots, M_n$ be the matrices involved in the expression of the~$2 \times 2$ Lorentz transformation.  Let~$T^{(a)}: M_2(\mathbb{O}) \to M_3(\mathbb{O})$ be an embedding of a~$2 \times 2$ matrix~$M$ into a~$3 \times 3$ matrix according to
$$T^{(1)}(M) = \left( \begin{array}{c|c} M & 0 \\ \hline 0 & 1 \end{array} \right) \hspace{2cm} T^{(2)}(M) = \left( \begin{array}{c|c} 1 & 0 \\ \hline 0 & M  \end{array} \right) $$
$$T^{(3)}(M) = \left( \begin{array}{c|c|c} M_{2,2} & 0 & M_{2,1} \\ \hline 0 & 1 & 0 \\ \hline M_{1,2} & 0 & M_{1,1} \end{array} \right) $$
for~$a = 1,2,3$.  We call this the {\it discrete}\footnote{This type map is expanded to a continuous map in Section~\ref{ch:E6_basic_structure.Type_Transformation}} {\it type map}.  Given a~$2 \times 2$ matrix~$M$, we abbreviate~$T^{(a)}(M)$ as~$M^{(a)}$ and refer to~$M^{(1)}, M^{(2)},$ and~$M^{(3)}$ as type~$1, 2$, and~$3$ versions of~$M$.  We note that
$$M^{(2)} = \mathcal{T} M^{(1)} \mathcal{T}^\dagger \hspace{1.5cm} M^{(3)} = \mathcal{T} M^{(2)} \mathcal{T}^\dagger \hspace{1.5cm} M^{(1)} = \mathcal{T} M^{(3)} \mathcal{T}^\dagger$$
is a cyclic permutation of the three types of~$M$, where 
$$\mathcal{T} = \left( \begin{array}{ccc} 0 & 0 & 1 \\ 1 & 0 & 0 \\ 0 & 1 & 0\end{array} \right)$$
If~$R \in SL(2,\mathbb{O})$ is a vector Lorentz transformation which is constructed from the nesting of conjugations of~$2 \times 2$ matrices~$M_1, \cdots, M_n$, then it can be shown that
$$ \chi \to M_n^{(a)} \left( \cdots \left( M_1^{(a)} \chi M_1^{(a)^\dagger} \right) \cdots \right)M_n^{(a)\dagger}$$
respects the block structure of~$\chi$, in the sense that identifying a vector~$\matX$, spinor~$\theta$, and dual spinor~$\theta^\dagger$ in~$\chi$ according to the Location~$a$ given in Table~\ref{table:three_types_of_vector_locations} results in the transformations
$$ \matX \to R(\matX) \hspace{1cm} \theta \to R(\theta) \hspace{1cm} \theta^\dagger \to R(\theta^\dagger)$$
for the appropriate vector, spinor, and dual-spinor versions~$R$ of the Lorentz transformation.  We call~$R^{(a)}$ the transformation~$R$ of {\it Type~$a$}.

{\bf Lemma: }  For any Lorentz transformation~$R \in SL(2,\mathbb{O})$, the corresponding~$3 \times 3$ transformations~$R^{(a)}: M_3(\mathbb{O}) \to M_3(\mathbb{O})$ preserves the determinant of elements of~$M_3(\mathbb{O})$ for~$a = 1,2,3$.

{\bf Proof: } We first prove this for the case~$a = 1$:
By definition, the finite Lorentz transformation~$R \in SL(2,\mathbb{O})$ preserves the determinant of the~$2 \times 2$ vector~$\matX \in M_2(\mathbb{O})$.  That is,~$\det{ \left( \matX \right) } = \det{\left( R( \matX) \right)}$.  
Since~$R$ is a compatible transformation, then~\hbox{$R(\theta \theta^\dagger) = R(\theta) \; R(\theta^\dagger)$} for the appropriate vector, spinor, and dual-spinor versions of the Lorentz transformation, giving
$$\det{ \left(\theta \; \theta^\dagger \right)} = \det{ \left( R(\theta) \; R(\theta^\dagger) \right)}$$    
As shown in Section~\ref{ch:Division_Algebras_Applications.Jordan_Algebras}, the determinant of~$\chi \in M_3(\mathbb{O})$ is given by
$$ \det{ \left( \chi \right) } = \det{\left( \begin{array}{cc}\matX & \theta \\ \theta^\dagger & n \end{array} \right)} = \left( \det{ \matX } \right) n + 2\matX \cdot \theta \theta^\dagger$$
where~$A \cdot B = \tr{A \jordanproduct B} - \tr{A} \tr{B}$ is the Lorentzian inner product of~$A$ and~$B$.  
But
Lorentz transformations preserve the Lorentzian inner product of~$\matX$ and~$\theta \; \theta^\dagger$, so that

$$\begin{array}{ccc}
R(\matX) \cdot R(\theta \; \theta^\dagger) & = & \matX \cdot (\theta \; \theta^\dagger) 
\end{array}
$$
Hence, we see that 
$$ \begin{array}{ccl}
\det{ \left(R^{(1)}(\chi)\right) } & = & \det{\left( \begin{array}{c|c} R(\matX) & R(\theta) \\ \hline R(\theta^\dagger) & n \end{array} \right)}\\ 
& = &\det{\left(R(\matX) \right)} n + 2 R(\matX) \cdot \left(R(\theta) R(\theta^\dagger) \right)\\
& = &\det{\left(R(\matX) \right)} n + 2 R(\matX) \cdot R(\theta \theta^\dagger)\\
& = &\det{ \left( \matX \right) } n + 2\matX \cdot (\theta \theta^\dagger) \\
& = & \det{(\chi)}
 \end{array} $$
where the first two equalities are true by definition, the third equality is a result of compatibility, and the fourth equality results from the use of Lorentz transformations.  Hence, the transformation~$R^{(1)}$ preserves the determinant of~$\chi \in M_3(\mathbb{O})$ for each~$R \in SL(2,\mathbb{O})$.

For the case~$a = 2$ and~$a = 3$, we note that the determinant of~$\chi \in M_3(\mathbb{O})$ may be written
$$ \det{\chi} = (t+z)(t-z)n - ((t+z)|b|^2 + (t-z)|z|^2 + n|a|^2) + 2Re(\overline{a} \overline{b} \overline{c})$$
as was shown in Section~\ref{ch:Division_Algebras_Applications.Jordan_Algebras} with~$t+z = p$ and~$t-z = m$.  We see that it is cyclic in~$a,b,c$ and~$t+z,t-z,n$.  That is,~$\det{( \mathcal{T} \chi \mathcal{T}^\dagger) } = \det{ \chi}$.  Hence,
$$\det{ \chi} = \det{(R^{(1)}(\chi))} = \det{(R^{(2)}(\chi))} = \det{(R^{(3)}(\chi))}$$
showing that each of the three types of~$3 \times 3$ transformation constructed from \hbox{$R \in SL(2,\mathbb{O})$} are in~$SL(3,\mathbb{O})$.

\myendofproof

We conclude with a result showing that the~$3 \times 3$ transformations we just constructed will form a one-parameter subgroup of~$SL(3,\mathbb{O})$ provided that the original~$2 \times 2$ transformation is both compatible and forms a one-parameter subgroup of~$SL(2,\mathbb{O})$.

{\bf Lemma: } If the finite Lorentz transformation~$R: \mathbb{R} \times M_2(\mathbb{O}) \to M_2(\mathbb{O})$ is both compatible and a one-parameter subgroup of~$SL(2,\mathbb{O})$, then~$R^{(a)}: \mathbb{R}\times M_3(\mathbb{O}) \to M_3(\mathbb{O})$ is a one-parameter subgroup of~$SL(3,\mathbb{O})$ for~$a = 1,2,3$.

{\bf Proof: }  Our construction of~$R^{(a)}$ fixes a block structure in~$\chi$ in which the vector, spinor, and dual-spinor version of the transformation~$R$ act on the appropriate blocks of~$\chi$.  By hypothesis, the vector transformation~$R$ is a one-parameter subgroup.  It remains to show that the spinor and dual-spinor transformations are one-parameter subgroups.  
But if the vector version of~$R$ is compatible and a one-parameter subgroup, then for any~\hbox{$\alpha, \beta \in \mathbb{R}$}, we see that
$$\begin{array}{ccl}
\left( R(\alpha)\left(R(\beta)(\theta)\right) \right) \; \left (R(\alpha)\left(R(\beta)(\theta^\dagger)\right) \right) & = & R(\alpha)\left( R(\beta)(\theta) \; R(\beta)(\theta^\dagger) \right) \\
 & = & R(\alpha)\left( R(\beta)(\theta \; \theta^\dagger) \right) \\
 & = & R(\alpha + \beta)(\theta \; \theta^\dagger) \\
 & = & \left( R(\alpha + \beta)(\theta) \right) \; \left( R(\alpha+\beta)(\theta^\dagger)\right) \\
\end{array}
$$
where~$\theta$ is a spinor and~$\theta^\dagger$ is a dual spinor, and we have used the appropriate vector, spinor, or dual spinor version of the vector Lorentz transformation~$R$.  Hence, the products of spinor and dual-spinor transformations are one-parameter subgroups.  But this requires the spinor and dual-spinor transformations themselves to be one-parameter subgroups.
Since each of the vector, spinor, and dual-spinor versions of the transformation~$R$ are one-parameter subgroups, then our construction produces a new transformation~$R^{(a)}$ which is also a one-parameter subgroup.

\myendofproof

This lemma shows that since each of the basis transformations~$R$ constructed by Manogue and Schray~\cite{manogue_schray} which were discussed in~\ref{ch:Division_Algebras_Applications.Lorentz_Transformations.Lorentz_Transformations_for_Spacetimes_involving_Division_Algebras} are compatible one-parameter subgroups of~$SL(2,\mathbb{O})$, then~$R^{(a)}$ is a one-parameter subgroup of~$SL(3,\mathbb{O})$ for~$a = 1,2,3$.  This result will be useful in Section~\ref{ch:E6_basic_structure.ConstructingE6algebra}.

\subsection{A Basis for the Category 3 Transformations}
\label{ch:E6_basic_structure.Lorentz_Transformations.BetterBasis}

The type maps provide three different copies of~$SO(9,1)$ transformations which act on~$\chi$.  Hence, there are at least three~$SO(9,1)$ groups, one of each type, as well as subgroups~$SO(n,\mathbb{R})$ and~$SO(n,1,\mathbb{R})$ for~$n \le 9$, in the group~$E_6 = SL(3,\mathbb{O})$.  This gives~$3 \times 45 = 135$ generators in~$E_6$, which only has dimension~$78$.  In the next section, we will use linear dependencies among the corresponding Lie algebra elements for these transformations to reduce the~$135$ generators down to a basis of~$78$ independent generators.  In the remainder of this section, we introduce an alternate basis for the Category~$3$ transformations which were given by Manogue and Schray \cite{manogue_schray} and reviewed in Section \ref{ch:Division_Algebras_Applications.Lorentz_Transformations.Lorentz_Transformations_for_Spacetimes_involving_Division_Algebras}.  This new basis will prove useful in determining the dependencies among the~$135$ algebra generators.
  
In the~$2 \times 2$ setting, each Category~$3$ transformation~$R_{p,q}$ rotates a single plane spanned by the imaginary octonions~$p$ and~$q \perp p$.  
  Rather than simply rotating one plane at a time, we choose a new basis of transformations, where each transformation simultaneously rotates two or three planes.  For each basis octonion, say~$q = i$, there are three pairs of basis octonions~$\lbrace j, k \rbrace, \lbrace k\l, j\l \rbrace, \lbrace \l, i\l \rbrace$ which generate a quaternionic subalgebra containing~$q = i$.  Note that we have chosen the ordering of the pairs in such a way that the product of the pair is~$q = i$, e.g. ~$(j) (k) = i$.  
   Hence, we have used the three triples (each pair along with~$q = i$) to create a right-handed, three-dimensional coordinate frame.  The planes spanned by~$\lbrace j, k \rbrace$,~$\lbrace k\l, j\l \rbrace$, and~$\lbrace \l, i\l \rbrace$ are orthogonal to each other, and the Category~$3$ transformations~$R_{j,k}(\alpha)$,~$R_{k\l,j\l}(\alpha)$, and~$R_{\l,i\l}(\alpha)$ rotate the appropriate plane counter-clockwise about the~$q=i$ axis through the angle~$\alpha$.  From these rotations, we create the new class of transformations~$A_q$,~$G_q$, and~$S_q$ where, for~$q = i$, we have \footnote{Because the rotations are in orthogonal planes, our definition does not depend on the order of the rotations, i.e.~$A_i(\alpha) = R_{j,k}(\alpha) \circ R_{k\l, j\l}(-\alpha) = R_{k\l, j\l}(-\alpha) \circ R_{j,k}(\alpha)$}
$$ A_i(\alpha) = R_{j,k}(\alpha) \circ R_{k\l, j\l}(-\alpha) \hspace{1.5cm} G_i(\alpha) = R_{j,k}(\alpha) \circ R_{k\l, j\l}(\alpha)  \circ R_{\l,i\l}(-2\alpha)$$
$$ S_i(\alpha) = R_{j,k}(\alpha) \circ R_{k\l, j\l}(\alpha)  \circ R_{\l,i\l}(\alpha)$$
where~$\circ$ is used to denote composition of the group transformations.
The~$A_i$ transformation rotates the~$(j,k)$ plane counter-clockwise and the~$(k\l, j\l)$ plane clockwise about the~$q = i$ axis.  The~$G_i$ transformation rotates two planes counter-clockwise about the~$q=i$ axis through an angle of~$\alpha$ and one plane clockwise about the~$q=i$ axis through an angle of~$2\alpha$.
The~$S_i$ transformation rotates all three planes counter-clockwise about the~$q = i$ axis through the angle~$\alpha$.\footnote{  We choose the labels~$A_q$ since the seven~$A_q$ tranformations along with~$G_\l$ are isomorphic to an~$SU(3,\mathbb{C})$ normally represented by Gell-Mann matrices, which are typically labeled using the greek symbol~$\Lambda$, and the Roman symbol~$A$ closely resembles~$\Lambda$.  The~$G_q$ label is chosen for the relationship between these transformations and the group~$G_2$.  The~$S_q$ label is chosen to remind us that these are the symmetric transformations.}  
Combined, the seven transformations of each class~$A_q$,~$G_q$, and~$S_q$ form a new basis for~$SO(7,\mathbb{R})$.  Finally, we expand these~$2 \times 2$ transformations to the~$3 \times 3$ case by using the type map~$T^{(a)}, a = 1,2,3$ which is applied to each matrix used in the transformations~$A_q, G_q$, and~$S_q$ to produce the transformations~$A^{(a)}_q$,~$G^{(a)}_q$, and~$S^{(a)}_q$ for~$a = 1,2,3$.  

{\bf Lemma: }  The~$3 \times 3$ transformations~$A^{(a)}_q$,~$G^{(a)}_q$, and~$S^{(a)}_q$ built from the~$2 \times 2$ one-parameter transformations~$R_{p,q}$ are also one parameter transformations.

{\bf Proof: } The~$2 \times 2$ transformations~$R_{p,q}$ are one-parameter transformations, and rotate one plane spanned by~$p$ and~$q$.  Our new transformations rotate two or three planes at once, but because these planes are pair-wise orthogonal, the effect of one transformation in the plane spanned by~$p$ and~$q$ does not affect the transformation in the plane spanned by~$u$ and~$v$.  Hence, the rotations which comprise~$A_q, G_q,$ and~$S_q$ commute with each other.  But then we can rearrange all of the rotations in~$A_q(\alpha) \circ A_q(\beta)$ to form~$A_q(\alpha+\beta)$, and similarly for the transformations~$G_q$ and~$S_q$.  Hence, the~$2 \times 2$ transformations~$A_q$,~$G_q$, and~$S_q$ are one-parameter transformations.  That their~$3 \times 3$ generalizations~$A^{(a)}_q, G^{(a)}_q$, and~$S^{(a)}_q$ are also one-parameter transformations follows from the last lemma of Section~\ref{ch:E6_basic_structure.Lorentz_Transformations}.

\myendofproof

We list here our conventions for our new basis for the~$SO(7)$ transformations.  Each of the~$A_q$,~$G_q$, and~$S_q$ transformations are labeled by a unit octonion~$q \in \lbrace i,j,k,k\l,j\l, i\l, l \rbrace$ about which the planes rotate.  For each~$q$, we have three pairs of unit octonions which multiply to~$q$.  
  Our conventions for the pairs are listed in Table~\ref{table:AGS_conventions}.  We note that~$\l$ only appears in the third pair, if at all.

\begin{table}[htbp]
\begin{center}
\begin{tabular}{|c|c|c|c|}
\hline
$q$ & First pair & Second pair & Third pair \\
\hline
\hline
$i$  & $(j , k )$ & $(kl , jl)$ & $(l , il)$ \\
$j$  & $(k , i )$ & $(il , kl)$ & $(l , jl)$ \\
$k$  & $(i , j )$ & $(jl , il)$ & $(l , kl)$ \\
$kl$ & $(jl , i)$ & $(j , il) $ & $(k , l) $ \\
$jl$ & $(i , kl)$ & $(il , k) $ & $(j , l) $ \\
$il$ & $(kl , j)$ & $(k , jl) $ & $(i , l) $ \\
$l$  & $(il , i)$ & $(jl , j) $ & $(kl , k)$ \\
\hline
\end{tabular}
\caption{\noindent Quaternionic subalgebras chosen for~$A_q$,~$G_q$, and~$S_q$}
\label{table:AGS_conventions}
\end{center}
\end{table}

In the~$2 \times 2$ case, we may reproduce the original Category~$3$ transformations~$R_{q,r}(\alpha)$ simply by combining our new transformations with appropriate angles.    For instance,

$$R_{j,k}(6\alpha) = A_i(3\alpha) \circ G_i(\alpha) \circ S_i(2\alpha)$$

We use the anlges~$3 \alpha$,~$\alpha$, and~$2\alpha$ so that the clockwise and counter-clockwise rotations in the various planes cancel each other out, leaving only the rotation by~$6\alpha$ in the~$j,k$ plane.  These dependencies also extend to the~$3 \times 3$ transformations.

\section{The Lie Algebra~$sl(3,\mathbb{O})$ of~$E_6 = SL(3,\mathbb{O})$}
\label{ch:E6_basic_structure.E6_Lie_Algebra}

We have expanded the~$45$ finite Lorentz transformations from the~$2 \times 2$ case to the~$3 \times 3$ case, giving us~$3 \times 45 = 135$ transformations.  Each of these transformations preserves the determinant of the Jordan matrix~$\chi \in M_3(\mathbb{O})$.  Hence, they are a subset of~$SL(3,\mathbb{O}) = E_6$, which is a~$78$-dimensional exceptional Lie group.  The Lie algebra of this group is~$sl(3,\mathbb{O})$, which is the topic of this section.  

In Section~\ref{ch:E6_basic_structure.ConstructingE6algebra}, we present the techniques which allow us to conclude that our transformations generate the entire group.  Normally, we would do this by associating each transformation with its tangent vector in the Lie algebra~$sl(3,\mathbb{O})$.  Then, using linear dependencies in the Lie algebra, we could reduce our~$135$ generators to~$78$ independent generators to show that we have the entire group~$SL(3,\mathbb{O})$.  However, the nesting of matrices involved in our representation of the transformations makes this difficult.  Instead, we use the local action of our representation of~$SL(3,\mathbb{O})$ on~$M_3(\mathbb{O})$ to give a homomorphism of~$sl(3,\mathbb{O})$ into the Lie algebra of all vector fields on~$M_3(\mathbb{O})$.  The tangent vectors to the group orbits at a point~$\chi$ in~$M_3(\mathbb{O})$ constitute a homomorphic copy of the Lie algebra~$sl(3,\mathbb{O})$.  We also use the tangent vectors to the group orbits to construct the complete commutator multiplication in the Lie algebra~$sl(3,\mathbb{O})$.  In Section~\ref{ch:E6_basic_structure.ConstructingE6algebra.LinearDependencies}, we use the linear dependencies among these tangent vectors to find linear dependencies among the~$135$ group generators constructed in the previous section.  These results are combined with our construction of~$SL(3,\mathbb{O})$ in Section~\ref{ch:E6_basic_structure.Lorentz_Transformations.BetterBasis} to give an explicit basis for subgroups~$SO(9,1,\mathbb{R})$ and~$SO(n,1,\mathbb{R})$ for~$n \le 9$, as well as for the exceptional groups~$G_2$ and~$F_4 = SU(3,\mathbb{O})$.

In what follows, we find it cumbersome to always refer to {\it tangent vectors to the group orbits} and {\it the homomorphic copy of the Lie algebra~$sl(3,\mathbb{O})$}.  Hence, we often refer to the group orbit of a one-parameter subgroup of~$SL(3,\mathbb{O})$ as a one-parameter path in~$M_3(\mathbb{O})$ and frequently state that the tangent vector to this curve is the element of the Lie algebra~$sl(3,\mathbb{O})$.  We hope this abuse of language makes the following methods more intuitive.

\subsection{Constructing the algebra for~$SL(3,\mathbb{O})$}
\label{ch:E6_basic_structure.ConstructingE6algebra}

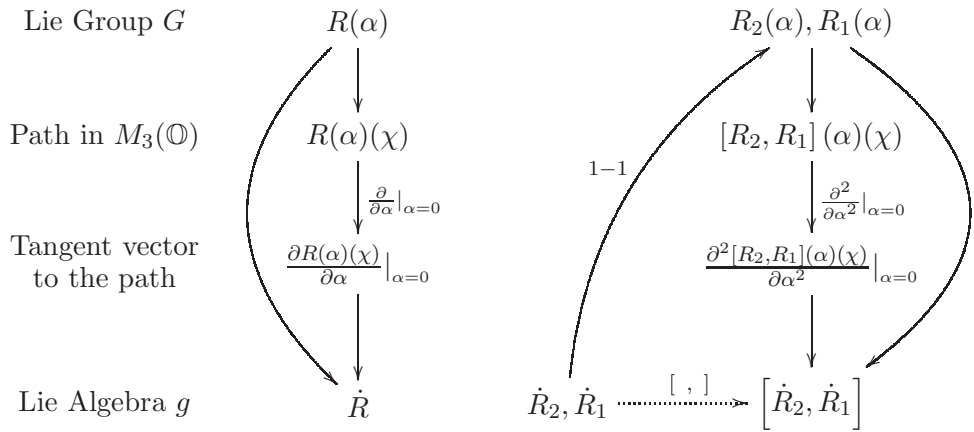
\begin{figure}[htbp]
\begin{minipage}{6in}
\begin{center}
\xymatrix@M=2pt{
*++\txt{Lie Group~$G$} &  *++\txt{$R(\alpha)$}\ar@/_3.5pc/[ddd] \ar[d] & \hspace{1.5cm} & *++\txt{$ R_2(\alpha),R_1(\alpha) $} \ar[d] \ar@/^5.0pc/[ddd]\\
*++\txt{$\textrm{Path in } M_3(\mathbb{O})$} &  *++\txt{$R(\alpha)(\chi)$} \ar[d]^{\frac{\partial }{\partial \alpha} |_{_{\alpha = 0}}} & & *++\txt{$\left[R_2,R_1\right](\alpha)(\chi)$} \ar[d]^{\frac{\partial^2 }{\partial \alpha^2}|_{_{\alpha = 0}}} \\
*++\txt{$\textrm{Tangent vector}$\\$\textrm{to the path}$} & *++\txt{$\frac{\partial R(\alpha)(\chi)}{\partial \alpha}|_{_{\alpha = 0}}$} \ar[d] & & *++\txt{$\frac{\partial^2 \left[ R_2, R_1 \right](\alpha)(\chi)}{\partial \alpha^2}|_{_{\alpha = 0}}$} \ar[d] \\
*++\txt{Lie Algebra~$g$}  &   *++\txt{$\dot R$} & *++\txt{$\dot R_2, \dot R_1 $} \ar@/^2pc/[uuur]^{1-1} \ar@{.>}[r]^{\left[ \hspace{.15cm}, \hspace{.15cm} \right]} &  *++\txt{$\left[ \dot R_2, \dot R_1 \right]$}\\
}
\caption{\noindent Calculating tangent vectors and their commutators}
\label{figure:lie_group_to_lie_algebra_commutators}
\end{center}
\end{minipage}
\end{figure}

We begin by associating each transformation in the group with a vector in the Lie algebra.  Each of the~$135$ transformations are one-parameter curves in the group.  Given a one-parameter curve~$R(\alpha)$ in a classical Lie group, the traditional method for associating it with the Lie algebra generator~$\dot R$ is to find its tangent vector~$\dot R = \frac{\partial R(\alpha)}{\partial \alpha } \mid_{_{\alpha = 0}}$ at the identity in the group.  However, our Category~$3$ transformations use nested matrices
$$  \chi \to M_4^{(a)} \left( M_3^{(a)} \left( M_2^{(a)} \left( M_1^{(a)} \chi M_1^{(a)\dagger} \right) M_2^{(a)\dagger} \right) M_3^{(a)\dagger} \right) M_4^{(a)\dagger}$$
and it is not possible to apply this technique to the one-parameter transformations~$A^{(a)}_q$, $G^{(a)}_q$, and~$S^{(a)}_q$, with~$a = 1,2,3$.  Instead, we let our one-parameter transformations~$R$ act on 
the exceptional Jordan algebra~$M_3(\mathbb{O})$, as it is a~$27$-dimensional representation of~$E_6$~\cite{corrigan}.
   When the transformation~$R$ acts on~$\chi \in M_3(\mathbb{O})$, it produces a curve~$R(\alpha)(\chi)$ in~$M_3(\mathbb{O})$.  This curve~$R(\alpha)(\chi)$ is a one-parameter curve in~$M_3(\mathbb{O})$, since~$R$ is a one-parameter transformation that is the identity,~$R(0)(\chi) = \chi$, when~$\alpha = 0$.  We use the tangent vector at the identity to this curve in~$M_3(\mathbb{O})$ to represent our tangent vector to the original one-parameter transformation~$R$.  That is, we have the association indicated in Figure~\ref{figure:lie_group_to_lie_algebra_commutators} between the group transformations and the tangent vectors.
This procedure produces~$135$ tangent vectors, one for each of the curves~$R(\alpha)(\chi)$.

We also use group orbits to construct the commutator of two tangent vectors.
  In the traditional approach to the classical matrix groups, the commutator of the tangent vectors~$\dot R_1$ and~$\dot R_2$ is defined as~\hbox{$\left[ \dot R_1, \dot R_2 \right] = \dot R_1 \dot R_2 - \dot R_2 \dot R_1$}.  However, we are not using tangent vectors to the transformations themselves but tangent vectors to curves in~$M_3(\mathbb{O})$.  To find the commutator of the tangent vectors to the curves~$R_1(\alpha)(\chi)$ and~$R_2(\alpha)(\chi)$, we create a new curve in~$M_3(\mathbb{O})$ defined by 
$$[ R_2, R_1](\alpha)(\chi) = R_2(-\frac{\alpha}{2}) \circ R_1(-\frac{\alpha}{2}) \circ R_2(\frac{\alpha}{2}) \circ R_1(\frac{\alpha}{2})(\chi)$$
where~$\circ$ denotes composition.
  This new path is not a one-parameter curve, and its first derivative is identically zero at~$\alpha = 0$, but~$\frac{\partial^2}{\partial \alpha^2}[R_2, R_1](\alpha)(\chi) |_{_{\alpha = 0}}$ is tangent to the curve~$[R_2, R_1](\alpha)(\chi)$ at~$\alpha = 0$.  Therefore, we define the commutator of~$\dot R_1$ and~$\dot R_2$ to be 
$$ \left[ \dot R_2, \dot R_1 \right] = \frac{\partial^2}{\partial \alpha^2}[R_2, R_1](\alpha)(\chi) |_{_{\alpha = 0}}$$
This formula agrees with the usual definition \cite{gilmore}.
We have a~$1 - 1$ correspondence between the group transformation~$R$, the orbit~$R(\alpha)(\chi)$ in~$M_3(\mathbb{O})$, and our tangent vector to the group orbit~$\dot R(\alpha)(\chi)$, allowing us to identify the group transformation given the tangent vector~$\dot R(\alpha)(\chi)$ to the group orbit~$R(\alpha)(\chi)$.  In particular, we can use this correspondence to identify the algebra element corresponding to~$\left[ \dot R_2, \dot R_1 \right]$ for any~$\dot R_2$ and~$\dot R_1$.  Our construction of the commutator, and this correspondence, is summarized in Figure~\ref{figure:lie_group_to_lie_algebra_commutators}.

  We note that since we are using the local action of~$SL(3,\mathbb{O})$ on~$M_3(\mathbb{O})$ to give a homomorphic image of~$sl(3,\mathbb{O})$, our construction does not lead to a readily available exponential map giving the group element corresponding to~$\left[\dot R_1, \dot R_2\right]$.  In particular, we are not always able to find the one-parameter curve whose tangent vector is~$\left[ \dot R_1, \dot R_2 \right]$. 

To illustrate our techniques, we compute the tangent vector for the group transformations~$A^{(1)}_\l$.  By our construction given in Section~\ref{ch:E6_basic_structure.Lorentz_Transformations.BetterBasis}, the finite Lorentz transformation~$A^{(1)}_\l(\alpha)$ conjugates the Jordan matrix~$\chi \in M_3(\mathbb{O})$ with the Category~$3$ matrices used in~$R^{(1)}_{i\l,i}(\alpha)$ and~$R^{(1)}_{j\l,j}(-\alpha)$.  This gives

\[ \begin{array}{ccc}
A^{(1)}_{\l}(\alpha)(\chi) & = &  R^{(1)}_{j\l,j}(-\alpha) \circ R^{(1)}_{i\l,i}(\alpha)(\chi) \\
\end{array}
\]
where
\[
\begin{array}{cl}
R^{(1)}_{i\l,i}(\alpha)(\chi) = & \left( \begin{array}{ccc} e^{\l \alpha/2 }i\l &0&0\\0&e^{\l \alpha/2 }i\l&0\\0&0&1\\ \end{array} \right) \hspace{6cm}\\
& \hspace{1cm} \left[
\left( \begin{array}{ccc} -i\l&0&0\\0&-i\l&0\\0&0&1  \end{array} \right)
\chi
\left( \begin{array}{ccc} -i\l&0&0\\0&-i\l&0\\0&0&1  \end{array} \right)^\dagger
\right] \hspace{2cm} \\
& \hspace{5cm}
\left( \begin{array}{ccc} e^{\l\alpha/2}i\l&0&0\\0&e^{\l\alpha/2}i\l&0\\0&0&1\\ \end{array} \right)^\dagger
\end{array}
\]
and
\[
\begin{array}{cl}
R^{(1)}_{j\l,j}(-\alpha)(\chi) = & \left( \begin{array}{ccc} e^{-\l \alpha/2 }j\l &0&0\\0&e^{-\l \alpha/2 }j\l&0\\0&0&1\\ \end{array} \right) \hspace{6cm} \\
& \hspace{1cm} \left[
\left( \begin{array}{ccc} -j\l&0&0\\0&-j\l&0\\0&0&1  \end{array} \right)
\chi
\left( \begin{array}{ccc} -j\l&0&0\\0&-j\l&0\\0&0&1  \end{array} \right)^\dagger
\right] \hspace{2cm}\\
& \hspace{4cm}
\left( \begin{array}{ccc} e^{-\l \alpha/2}j\l&0&0\\0&e^{-\l\alpha/2}j\l&0\\0&0&1\\ \end{array} \right)^\dagger\\
\end{array}
\]

\noindent It is useful to use a computer algebra software program to compute the result of conjugating~$\chi$ with the above matrices.  Then, the tangent vector to this group orbit is found to be

$$\dot A^{(1)}_\l = \frac{\partial A^{(1)}_\l(\alpha)(\chi)}{\partial \alpha} = 
 \left[ \begin {array}{ccc} 0& \overline{a_1}& c_1\\
\noalign{\medskip}a_1&0&\overline{b_1}\\
\noalign{\medskip}\overline{c_1}&b_1&0
\end {array} \right] 
$$
where

\begin{equation}
\label{eq:Al_tangent_coeff}
\begin{array}{c}
a_1 = -{\it a_{i\l}}\,{\it i}+{\it a_{j\l}}\,{\it j}-{\it a_j}\,{\it j\l}+{\it a_i}\,{\it i\l}\\
b_1 = -{\it b_{i\l}}\,{\it i}+{\it b_{j\l}}\,{\it j}-{\it b_j}\,{\it j\l}+{\it b_i}\,{\it i\l}\\
c_1 = -{\it c_{i\l}}\,{\it i}+{\it c_{j\l}}\,{\it j}-{\it c_j}\,{\it j\l}+{\it c_i}\,{\it i\l}
\end{array}
\end{equation}

The calculations involved in finding the commutator of two transformations are straight-forward, but time-consuming.  We used a {\tt MAPLE} program to calculate all~$3003$ commutators of pairs of the~$135$ different group generators as well as identify its tangent vector as a linear combination of the tangent vectors to our preferred~$78$ transformations.  Running this calculation on~$12$ computers running at~$1.5Ghz$, it finished in just under~6 hours.

To simplify notation, we will use the traditional notation~$\dot R$ in place of~$\frac{\partial R(\alpha)(\chi)}{\partial \alpha}|_{_{\alpha = 0}}$ for the tangent vector to~$R(\alpha)(\chi)$ at~$\alpha = 0$.  Using the isomorphism above, we can regard these tangent vectors as our Lie algebra.

\subsection{Linear Dependencies}
\label{ch:E6_basic_structure.ConstructingE6algebra.LinearDependencies}

We shall now give the dependencies among the group transformations by using linear dependencies among the Lie algebra elements.  In doing so, we will indicate which transformations can be eliminated leaving our preferred basis for the group~$SL(3,\mathbb{O})$ and the algebra~$sl(3,\mathbb{O})$.  Since we are using a homomorphic image of the Lie algebra~$sl(3,\mathbb{O})$, we check that the indicated dependencies actually do provide dependencies among the group transformations.

We begin with the Category~$3$ transformations.  Among the~$21$ transformations~$A_q$,~$G_q$, and~$S_q$ of each type, direct examination of the tangent vectors shows that
$$\dot A^{(1)}_q = \dot A^{(2)}_q = \dot A^{(3)}_q \hspace{2cm} \dot G^{(1)}_q = \dot G^{(2)}_q = \dot G^{(3)}_q$$
for each basis octonion~$q$.  That is, the~$A_q$ and~$G_q$ transformations are type independent, allowing us to drop the type designation and simply write~$\dot A_q$ and~$\dot G_q$.  These fourteen transformations generate~$G_2 = Aut(\mathbb{O})$, which is the smallest of the exceptional Lie groups.  The type independence of these transformations is an effect of {\it triality}, which will be further discussed in Section 4.3.
When added to the fourteen~$G_2$ transformations, the seven~$S^{(a)}_q$ transformations of each type produce a basis for an~$SO(7)$ of type~$a$, with~$a = 1, 2, 3$.  However, the~$S^{(a)}_q$ transformations are not independent, as their tangent vectors satisfy 
$$\dot S^{(1)}_q + \dot S^{(2)}_q + \dot S^{(3)}_q = 0$$
Hence, the union of any two of the~$SO(7)$ groups contains the third.  In particular, we may use the transformations~$S^{(a)}_q$ of type~$1$ and type~$2$ to generate the type~$3$ transformations~$S^{(3)}_q$.  These linear dependences have reduced our~$63$ Category~$3$ transformations by~$35$, trimming our original~$135$ transformations down to~$100$.

  Among the Category~$2$ transformations, we have the relation
$$ \dot R^{(1)}_{xq} + \dot R^{(2)}_{xq} + \dot R^{(3)}_{xq} = 0$$
This identity allows us to eliminate another seven transformations, and we choose to keep the type~$1$ and type~$2$ transformations.  The seven group transformations~$R^{(a)}_{xq}$ of type~$a$ are added to the group~$SO(7)$ of type~$a$ to obtain~$SO(8)$.  But the dependency among the types for the~$R^{(a)}_{xq}$ transformations indicates the group transformation~$R^{(3)}_{xq}$ can be obtained from the group transformations~$R^{(1)}_{xq}$ and~$R^{(2)}_{xq}$.  That is, the union of the groups~$SO(8)$ of type~$1$ and type~$2$ contains the group~$SO(8)$ of type~$3$.

A final pair of relations exist among the Category~$2$ and Category~$3$ transformations, however, which will eliminate an additional~$14$ transformations and condense all of the different~$SO(8)$ groups into one!  Having eliminated the~$S_q$ and~$R_{xq}$ transformations of type~$3$, these final relations allow us to write the type~$2$ transformations as linear combinations of the type~$1$ transformations.
  The tangent vectors satisfy
$$  \dot R^{(2)}_{xq} = -\frac{1}{2}\dot R^{(1)}_{xq} - \frac{1}{2}\dot S^{(1)}_q \hspace{2cm} \dot S^{(2)}_q = \frac{3}{2} \dot R^{(1)}_{xq} - \frac{1}{2}\dot S^{(1)}_q   $$
Hence, we may express any of the~$S^{(2)}_q$ and~$R^{(2)}_{xq}$ transformations in terms of~$S^{(1)}_q$ and~$R^{(1)}_{xq}$ transformations which are in the~$SO(8)$ of type~$1$.  Hence, there is only one~$SO(8)$!  This group contains three different copies of~$SO(7)$ built upon the unique group~$G_2$.  Again, this is a result of triality, which will be further discussed in Section 4.2.

Having reduced the~$135$ transformations to~$100$ and then by another~$21$ to~$79$, we expect only one additional linear dependency among the transformations.  Not surprisingly, we are left with~$52$ Category~$2$ and~$3$ transformations.  These rotations preserve the trace of~$\chi \in M_3(\mathbb{O})$ and form the Lie group~$F_4$.  Among the boost transformations in Category~$1$, we see that
$$ \dot B^{(1)}_{tz} + \dot B^{(2)}_{tz} + \dot B^{(3)}_{tz} = 0$$
From the three boosts~$B^{(1)}_{tz}$,~$B^{(2)}_{tz}$, and~$B^{(3)}_{tz}$, we choose to keep the boosts~$B^{(1)}_{tz}$ and~$B^{(2)}_{tz}$ in our preferred basis for~$SL(3,\mathbb{O})$.  We will see in the algebra that there is some evidence that we should choose the boost~$\dot B^{(2)}_{tz} - \dot B^{(3)}_{tz}$ in place of~$\dot B^{(2)}_{tz}$, but this choice makes little difference given
the identity among~$\dot B^{(1)}_{tz}, \dot B^{(2)}_{tz}$, and~$\dot B^{(3)}_{tz}$.

It can be shown that the final~$78$ tangent vectors corresponding to the remaining boosts and rotations are linearly independent.   
By looking at the tangent vectors to the group actions for the rotations, direct computation shows there are~$24$ hermitian and~$24$ anti-hermitian trace-free tangent vectors in addition to the fourteen~$G_2$ tangent vectors.  The tangent vectors for the~$26$ boosts are non-trace-free and hermitian.  It turns out that the~$64$ non-$G_2$ tangent vectors correspond to the~$64$ independent trace-free octonionic~$3 \times 3$ matrices.  Taken with the~$G_2$ transformations, this gives~$78$ independent transformations.  Hence, we do have the group~$SL(3,\mathbb{O}) = E_6$.  The commutation table for~$sl(3,\mathbb{O})$ is located online at~\cite{commutation_table_online}, and using this table, we can further identify 
$$\dot B^{(1)}_{tz}, \dot B^{(2)}_{tz}, \dot R^{1}_{x\l}, \dot A_\l, \dot G_\l, \textrm{ and }\dot S^{(1)}_{\l}$$
as the six elements corresponding to the Casimir operators in the complexified Lie algebra~$sl(3,\mathbb{O})$.  Hence, we will refer to these six elements as Casimir operators, but note that this is not a choice of orthogonal Casimir operators.

\subsection{Triplets of Subgroup Chains}
\label{ch:E6_basic_structure.Lorentz_Transformations.BetterBasisSubgroups}

\begin{figure}[htbp]
\begin{center}
\begin{minipage}{6in}
\[
\xymatrix{
               &   *++\txt{$sl(3,\mathbb{O})$\\$[E_6]$}           &             \\
               &   *++\txt{$su(3,\mathbb{O})$\\$[F_4]$}\ar[u]^{B^{(1)}_{tz}, B^{(2)}_{tz}} \\
*++\txt{$so(9) = su(2,\mathbb{O})$\\ $[B_4]$}\ar[ur] & *++\txt{$so(9) = su(2,\mathbb{O})$\\ $[B_4]$}\ar[u] & *++\txt{$so(9) = su(2,\mathbb{O})$\\ $[B_4]$}\ar[ul] \\
               &   *++\txt{$so(8)$\\$[D_4]$}\ar[u] \ar[ur] \ar[ul] \\
*++\txt{$so(7) = su(1,\mathbb{O})$\\ $[B_3]$}\ar[ur]^{R^{(1)}_{x\l}} & *++\txt{$so(7) = su(1,\mathbb{O})$\\ $[B_3]$}\ar[u]^{S^{(2)}_\l \to S^{(1)}_\l}_{R^{(1)}_{x\l}} & *++\txt{$so(7) = su(1,\mathbb{O})$\\ $[B_3]$}\ar[ul]^{S^{(3)}_\l \to S^{(1)}_\l}_{R^{(1)}_{x\l}} \\
               &   *++\txt{$Aut(\mathbb{O})$\\$[G_2]$}\ar[u]^{S^{(2)}_\l} \ar[ur]^{S^{(3)}_\l} \ar[ul]^{S^{(1)}_\l} \\
               &   *++\txt{$su(3,\mathbb{C})^C$\\$[A_2]$}\ar[u] \\
               &   *++\txt{$su(2,\mathbb{C})^C$\\$[A_1]$}\ar[u]^{G_\l} \\
               &   *++\txt{$u(1,\mathbb{C})$}\ar[u]^{A_\l}
}
\]
\caption{Chain of subgroups~$SU(3)^C \subset G_2 \subset SO(8) \subset F_4 \subset E_6$}
\label{figure:type_1_2_3_G2_so7}
\end{minipage}
\end{center}
\end{figure}
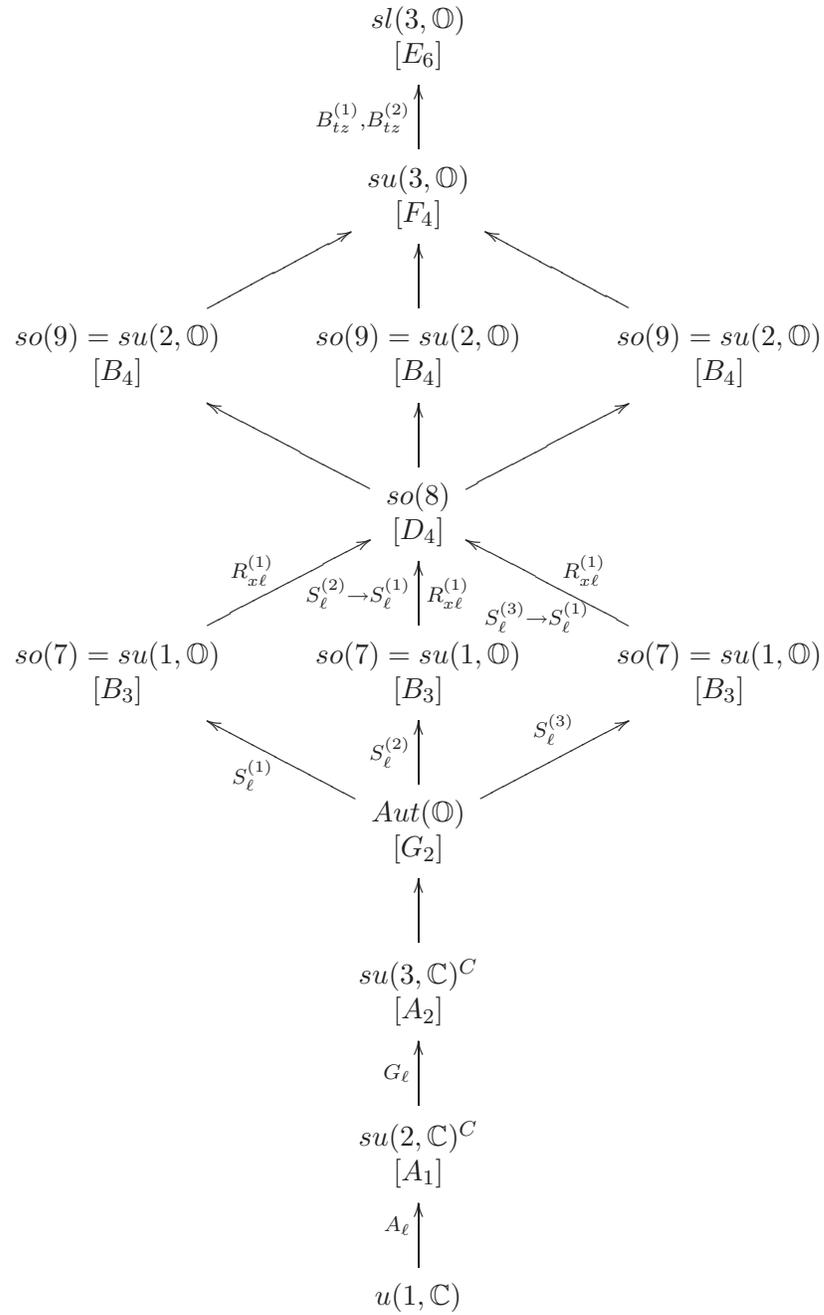

With our particular choice of basis for~$E_6$, we can identify two separate~$SO(n)$ subgroup structures within~$E_6$.
Figure~\ref{figure:type_1_2_3_so} shows the~$SO(n)$ subgroup chain of~$SO(9,1)$ of type~$1$ in~$SL(3,\mathbb{O})$, while Figure~\ref{figure:type_1_2_3_G2_so7} shows the three~$SO(9)$ subgroup chains of~$F_4$ within~$E_6$.
  In both subgroup structures, there is only one~$SO(8)$.  While~$G_2 \subset SO(7)$, it is not a subset of~$SO(6)$ in Figure~\ref{figure:type_1_2_3_so}.  Hence, we omit~$G_2$ from Figure~\ref{figure:type_1_2_3_so}, but include it in Figure~\ref{figure:type_1_2_3_G2_so7} since the~$SO(7)$ basis here may be restricted to to the basis of~$G_2$.  We also note that~$SU(3)$, consisting of~$A_i, \cdots, A_{i\l}, A_\l, G_\l$, is contained both within~$G_2$ and~$SO(6)$.  G\"unaydin denotes this~$SU(3)$ subgroup of~$G_2$ as~$SU(3)^C$ in~\cite{gunaydin_1974}, and we adopt that convention.  We use this same notation for the subgroup~$SU(2)^C \subset SU(3)^C$, which consists of~$A_k, A_{k\l}$ and~$A_\l$.  
The figures indicate which Casimir operator is added to a group when it is expanded to a larger group.  The additional Casimir operator and basis elements are also indicated in Table~\ref{table:basis_for_type_1_2_3_so} and Table~\ref{table:basis_type_1_2_3_G2_so7}.  

\begin{figure}[htbp]
\begin{center}
\begin{minipage}{6in}
\begin{center}
\[
\xymatrix{
               &   *++\txt{$sl(3,\mathbb{O})$\\$[E_6]$}           &             \\
*++\txt{$so(9,1) = sl(2,\mathbb{O})$\\ $[D_5]$}\ar[ur]^{B^{(2)}_{tz}} & *++\txt{$so(9,1) = sl(2,\mathbb{O})$\\ $[D_5]$}\ar[u]_{B^{(1)}_{tz}} & *++\txt{$so(9,1) = sl(2,\mathbb{O})$\\ $[D_5]$}\ar[ul]^{B^{(3)}_{tz} \to B^{(1)}_{tz}}_{B^{(2)}_{tz}} \\
*++\txt{$so(9) = su(2,\mathbb{O})$\\ $[B_4]$}\ar[u]^{B^{(1)}_{tz}} & *++\txt{$so(9) = su(2,\mathbb{O})$\\ $[B_4]$}\ar[u]^{B^{(2)}_{tz}} & *++\txt{$so(9) = su(2,\mathbb{O})$\\ $[B_4]$}\ar[u]^{B^{(3)}_{tz}} \\
               &   *++\txt{$so(8)$\\$[D_4]$}\ar[u] \ar[ur] \ar[ul] \\
*++\txt{$so(7) = su(1,\mathbb{O})$\\ $[B_3]$}\ar[ur]^{R^{(1)}_{x\l}} & *++\txt{$so(7) = su(1,\mathbb{O})$\\ $[B_3]$}\ar[u]^{S^{(2)}_\l \to S^{(1)}_\l}_{R^{(1)}_{x\l}} & *++\txt{$so(7) = su(1,\mathbb{O})$\\ $[B_3]$}\ar[ul]^{S^{(3)}_\l \to S^{(1)}_\l}_{R^{(1)}_{x\l}} \\
*++\txt{$so(6)$\\ $[D_3]$}\ar[u] & *++\txt{$so(6)$\\ $[D_3]$}\ar[u] & *++\txt{$so(6)$\\ $[D_3]$}\ar[u] \\
*++\txt{$so(5)$\\ $[B_2]$}\ar[u]_{G_\l-2S^{(1)}_\l} & *++\txt{$so(5)$\\ $[B_2]$}\ar[u]_{G_\l-2S^{(2)}_\l} & *++\txt{$so(5)$\\ $[B_2]$}\ar[u]_{G_\l-2S^{(3)}_\l} \\
*++\txt{$so(4) = so(3)\oplus so(3)$\\ $[D_2]$}\ar[u] & *++\txt{$so(4) = so(3)\oplus so(3)$\\ $[D_2]$}\ar[u] & *++\txt{$so(4) = so(3)\oplus so(3)$\\ $[D_2]$}\ar[u] \\
 & *++\txt{$so(3)$\\ $[B_1]$}\ar[u]_{G_\l+S^{(2)}_\l} \ar[ul]_{G_\l+2S^{(1)}_\l} \ar[ur]_{G_\l+2S^{(3)}_\l} & \\
 & *++\txt{$u(1)$}\ar[u]_{A_\l} &
}
\]
\caption{Chain of subgroups~$SO(n) \subset SO(9,1) \subset SL(3,\mathbb{O})$}
\label{figure:type_1_2_3_so}
\end{center}
\end{minipage}
\end{center}
\end{figure}
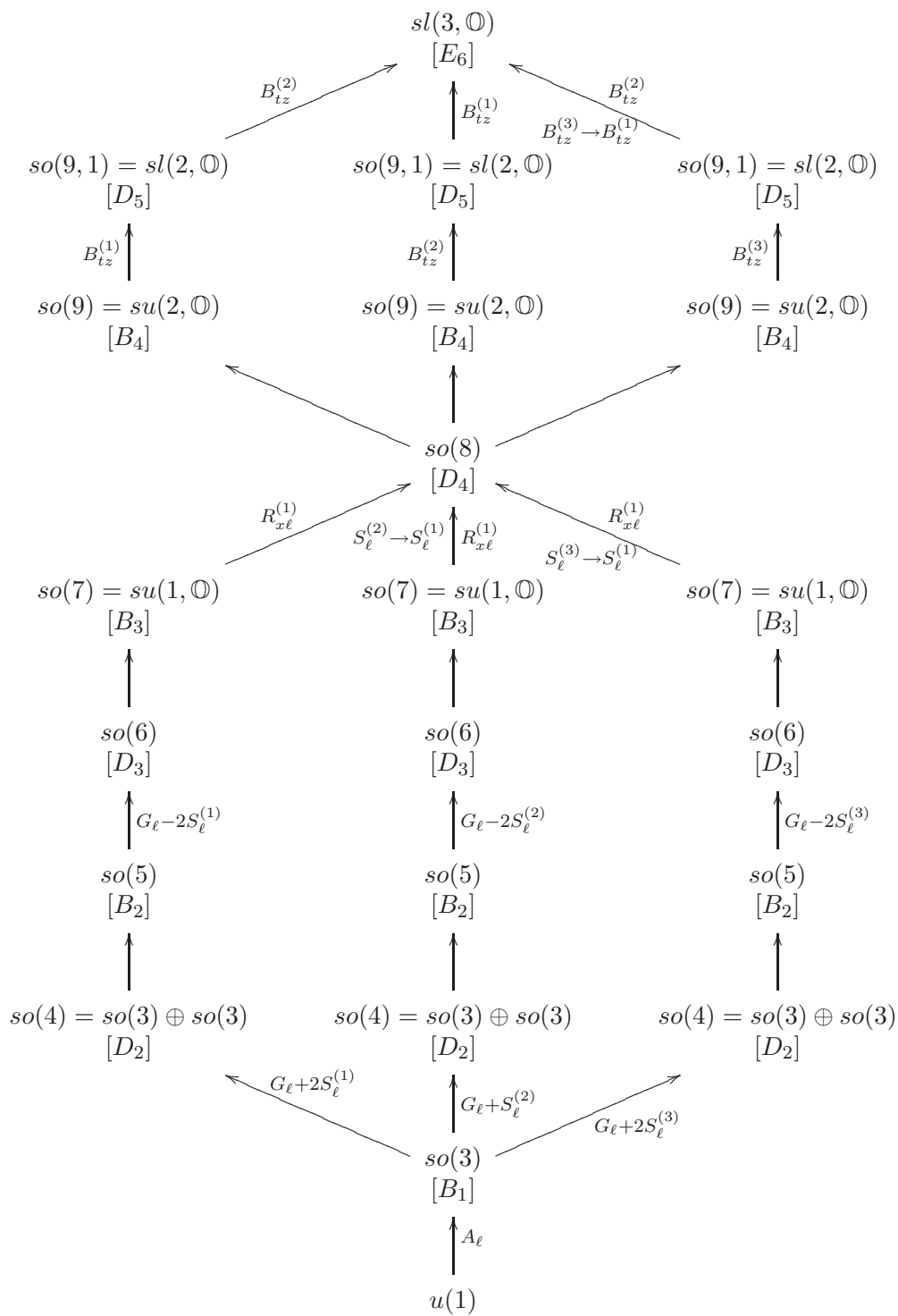

\begin{table}[htbp]
\begin{center}
\begin{tabular}{|l|l|l|l|}
\hline
Algebra & Casimir  & Basis & Set for~$q$ \\
        & Operator &       &             \\
\hline
\hline
$sl(3,\mathbb{O}) = E_6$ & $B_{tz}^{(2)}$     & $B^{(2)}_{tx}, B^{(3)}_{tx}$ & \\
                         &                      & $B^{(2)}_{tq}, B^{(3)}_{tq}$ & $i,j,k,k\l,j\l,i\l,\l$ \\
                         &                      & $R^{(2)}_{zq}, R^{(3)}_{zq}$ & $i,j,k,k\l,j\l,i\l,\l$ \\
                         &                      & $R^{(2)}_{xz}, R^{(3)}_{xz}$ & \\
\hline
$sl(2,\mathbb{O}) = so(9,1) = D_5$ & $B^{(1)}_{tz}$ & $B^{(1)}_{tx}, B^{(1)}_{tq}$ & $i,j,k,k\l,j\l,i\l,\l$ \\
                           &             & $R^{(2)}_{xz}, R^{(3)}_{xz}$ &   \\
                           &             & $R^{(2)}_{zq}, R^{(3)}_{zq}$ & $k,k\l,\l$  \\
                           &             & $B^{(1)}_{tq}$ & $i,j,i\l,j\l$ \\
\hline
$su(2,\mathbb{O}) = so(9) = B_4$ &           & $R^{(1)}_{xz}, R^{(1)}_{zq}$ & $i,j,k,k\l,j\l,i\l,\l$ \\
\hline
$so(8) = D_4$ & $R^{(1)}_{x\l}$ & $R^{(1)}_{xq}$      & $i,j,k,k\l,j\l,i\l$   \\
\hline
$su(1,\mathbb{O}) = so(7) = B_3$ &            & $G_q - S^{(1)}_q$ & $i,j,k,k\l,j\l,i\l$ \\
\hline
$so(6) = D_3$ & $G_\l-S^{(1)}_\l$ & $3A_q - G_q - 2S^{(1)}_q$ & $i,j$ \\
              &                 & $3A_q + G_q + 2S^{(1)}_q$ & $i\l,j\l$ \\
\hline
$so(5) = B_2$ &            & $3A_q + G_q + 2S^{(1)}_q$ & $i,j$ \\
              &            & $3A_q - G_q - 2S^{(1)}_q$ & $i\l,j\l$ \\
\hline
$so(4) = so(3)\oplus so(3) = D_2$ & $G_\l +2S^{(1)}_\l$ & $G_q+2S^{(1)}_q$ & $k,k\l$\\
\hline
$so(3)^C = B_1$ & &$A_q$ & $k,k\l$\\
\hline
$u(1)$ & $A_\l$ & & \\
\hline
\end{tabular}
\caption{\noindent Basis for~$SO(n) \subset SO(9,1)$ subgroup chain of Type~$1$}
\label{table:basis_for_type_1_2_3_so}
\end{center}
\end{table}

\begin{table}[htbp]
\begin{center}
\begin{tabular}{|l|l|l|l|}
\hline
Algebra & Casimir  & Basis & Set for~$q$ \\
        & Operator &       &             \\
\hline
\hline
$sl(3,\mathbb{O}) = E_6$ & $B^{(1)}_{tz}, B^{(2)}_{tz}$     & $B^{(1)}_{tx}, B^{(2)}_{tx}, B^{(3)}_{tx}$ & \\
                         &                      & $B^{(1)}_{tq}, B^{(2)}_{tq}, B^{(3)}_{tq}$ & $i,j,k,k\l,j\l,i\l,\l$ \\
\hline
$su(3,\mathbb{O}) = F_4 $ &  & $R^{(2)}_{zq}, R^{(3)}_{zq}$ & $i,j,k,k\l,j\l,i\l,\l$ \\
                          &  & $R^{(2)}_{xz}, R^{(3)}_{xz}$ & \\
\hline
$su(2,\mathbb{O}) = so(9) = B_4$ &           & $R^{(1)}_{xz}, R^{(1)}_{zq}$ & $i,j,k,k\l,j\l,i\l,\l$ \\
\hline
$so(8) = D_4$ & $R^{(1)}_{x\l}$ & $R^{(1)}_{xq}$      & $i,j,k,k\l,j\l,i\l$   \\
\hline
$su(1,\mathbb{O}) = so(7) = B_3$ &  $S^{(1)}_\l$  & $S^{(1)}_q$ & $i,j,k,k\l,j\l,i\l$ \\
\hline
$Aut(\mathbb{O}) = G_2$ &            & $G_q$ & $i,j,k,k\l,j\l,i\l$ \\
\hline
$su(3,\mathbb{C})^C$ & $G_\l$ & $A_q$ & $i,j,j\l,i\l$ \\
\hline
$su(2,\mathbb{C})^C = A_1$ & &$A_q$ & $k,k\l$\\
\hline
$u(1)$ & $A_\l$ & & \\
\hline
\end{tabular}
\caption{\noindent Basis for~$G_2 \subset SO(7) \subset F_4$ subgroup chain}
\label{table:basis_type_1_2_3_G2_so7}
\end{center}
\end{table}

\section{Triality}
\label{ch:E6_basic_structure.Triality}

In this section, we expand upon the discussion of triality given in Section~\ref{ch:Division_Algebras_Applications.Quaternions_Octonions.Triality} as it relates to~$SL(3,\mathbb{O})$ and its Lie algebra.  In particular, we note that each~$3 \times 3$ transformation in our basis of~$SL(3,\mathbb{O})$ naturally implements a~$2 \times 2$ Lorentz transformation and as well as a spinor and dual spinor version of this transformation.  However, we may also view~$\chi$ as consisting of three octonions as well as the three real components on its diagonal. The~$SO(8,\mathbb{R})$ transformations, which fix the diagonal components, only act on the octonions via left, right, or symmetric multiplication.
In this section,  we find an interesting connection between these different left, right, and symmetric transformations.
In particular, we discuss how the spinor, vector, and dual-spinor transformations embedded into our~$3 \times 3$ transformation use triality to simplify the subgroup structure of~$SO(8,\mathbb{R})$ within~$SL(3,\mathbb{O})$.  We also discuss triality in relation to our type map, and introduce a strong formulation of triality.

We quickly review the two versions of triality most pertinent to this work.  Baez~\cite{baez} mentions that each of the exterior nodes in the Dynkin diagram of~$so(8,\mathbb{R})$ may be identified with a spinor, dual-spinor, and vector representation of~$SO(8,\mathbb{R})$.  Both Baez~\cite{baez} and Conway~\cite{conway} discuss a triality between transformations of~$\mathbb{O}$ via left, right, and symmetric multiplication. Conway further states that each left multiplication may be written as a product of at most seven right or seven symmetric multiplications.  Further information regarding these claims may be found in~\cite{baez, conway} or Section~\ref{ch:Division_Algebras_Applications.Quaternions_Octonions.Triality}.

We first consider triality in relation to the transformations contained within~$G_2$.  
  Consider the one-parameter transformations~$A^{(a)}_{q}(\theta)$ and~$G^{(a)}_{q}(\theta)$, with $a = 1,2,3$, as constructed in Section~\ref{ch:E6_basic_structure.Lorentz_Transformations.BetterBasis}.  In Section~\ref{ch:E6_basic_structure.ConstructingE6algebra}, we calculated~$\dot A^{(1)}_{\l}$, which was found to be 
$$\dot A^{(1)}_\l(\theta)(\chi) = \frac{\partial A^{(1)}_\l(\theta)(\chi)}{\partial \alpha} = 
 \left[ \begin {array}{ccc} 0& a_1& \overline{c_1}\\
\noalign{\medskip}\overline{a_1}&0&b_1\\
\noalign{\medskip}c_1&\overline{b_1}&0
\end {array} \right] 
$$
$$a_1 = -{\it a_{i\l}}\,{\it i}+{\it a_{j\l}}\,{\it j}-{\it a_j}\,{\it j\l}+{\it a_i}\,{\it i\l}$$
$$b_1 = -{\it b_{i\l}}\,{\it i}+{\it b_{j\l}}\,{\it j}-{\it b_j}\,{\it j\l}+{\it b_i}\,{\it i\l}$$
$$c_1 = -{\it c_{i\l}}\,{\it i}+{\it c_{j\l}}\,{\it j}-{\it c_j}\,{\it j\l}+{\it c_i}\,{\it i\l}$$

\noindent The above expressions for~$a_1, b_1$, and~$c_1$ are very similar; Direct computation of the tangent vectors for~$A^{(2)}_\l$ or~$A^{(3)}_\l$ yields the exact same expression for~$a_1, b_1$, and~$c_1$.
  That is,
the one-parameter curves~$A^{(1)}_{\l}$,~$A^{(2)}_{\l}$, and~$A^{(3)}_{\l}$ all have this same tangent vector.  Hence, we refer to this transformation as~$A_{\l}$ without ambiguity.  Similarly, since the one-parameter curves~$G^{(1)}_\l, G^{(2)}_\l, G^{(3)}_\l$ have the same tangent vector, we also refer to this transformation as~$G_\l$ without a need to refer to type.  Similar results hold for~$A_q$ and~$G_q$ for~$q = i,j,k,k\l,j\l,i\l$.

In terms of vector, spinor, and dual spinor transformations, the~$G_2$ transformations are type independent.  This requires further explanation. Let~$q = i,j,k,k\l,j\l,i\l,\l$.   We may choose to view~$\chi$ as containing the vector~$\matX$, spinor~$\theta$, and dual spinor~$\theta^\dagger$ in any one of the three natural locations given in Table~\ref{table:three_types_of_vector_locations} of Section~\ref{ch:E6_basic_structure.Lorentz_Transformations}.   Regardless of this choice, the type~$1$ transformation~$A^{(1)}_q$ (and similarly~$G^{(1)}_q$) produces the same transformation on~$\matX, \theta, \theta^\dagger$.  On the other hand, we may fix a location for~$\matX$ in~$\chi$, and notice that the three types~$a$ of~$A^{(a)}_q$ (and similarly for~$G^{(a)}_q$) produce the same transformation on~$\matX$.  We say that the~$G_2$ transformations are type independent not only because the type~$1$ transformation may be written in terms of the type~$2$ or type~$3$ transformation, but because all three types are the same transformation!

Conway~\cite{conway} states that a general transformation of~$SO(8,\mathbb{R})$ via left-sided multiplication may be expressed in terms of at most seven right-sided multiplications or via at most seven symmetric multiplications.  Our basis~$G_2$ transformations use symmetric, left, and right-sided matrix multiplication on the vector~$\matX$, spinor~$\theta$, and dual spinor~$\theta^\dagger$ in~$\chi$ and are contained in~$SO(8,\mathbb{R})$.  By looking explicitly at the matrices involved in the type~$1$ expression for~$A_\l$, we find
\[
\begin{array}{ccl}
A_{l}(\alpha)(\chi) & = & A_{l}(\alpha)\left( \begin{array}{cc} \matX & \theta \\ \theta^\dagger & n \end{array} \right) \\
 & = & R^{(1)}_{j\l,j}(-\alpha) \circ R^{(1)}_{i\l,i}(\alpha)\left( \begin{array}{cc} \matX & \theta \\ \theta^\dagger & n \end{array} \right) \\
\end{array}
\]
where the vector~$\matX$, spinor~$\theta$, and dual spinor~$\theta^\dagger$ transform according to

$$ \matX \to \left( M_{2^\prime} \left( M_2 \left( M_{1^\prime} \left(M_1 \matX M_1^\dagger \right) M_{1^\prime}^\dagger \right) M_2^\dagger \right) M_{2^\prime}^\dagger \right)
$$

$$ \theta \to \left( M_{2^\prime} \left( M_2 \left( M_{1^\prime} \left(M_1 \theta \right) \right) \right) \right) $$

$$ \theta^\dagger \to \left( \left( \left( \left( \theta^\dagger M_1^\dagger \right) M_{1^\prime} ^\dagger \right) M_2^\dagger \right) M_{2^\prime}^\dagger \right)$$

\noindent with the~$2 \times 2$ matrices~$M_1, M_{1^\prime}, M_2$, and~$M_{2^\prime}$ given by 

$$ M_1 = -i\l I_2  \hspace{2cm} M_{1^\prime} = \left( \cos \left(\frac{\alpha}{2}\right) i\l + \sin \left(\frac{\alpha}{2}\right) i \right) I_2$$

$$ M_2 = -j\l I_2  \hspace{2cm} M_{2^\prime} = \left( \cos \left(-\frac{\alpha}{2}\right) j\l + \sin \left(-\frac{\alpha}{2}\right) j \right) I_2$$

\noindent  Indeed, we notice that the matrices involved in this~$G_2$ transformation, and all other~$G_2$ transformations, are diagonal and imaginary.  Hence, for these matrices,~$M^\dagger = -M$ and the vector transforms not only by conjugation but by symmetric multiplication!  Indeed, the spinor and dual spinor transform via four multiplications from the left or right, respectively.  We also note that our other~$G_2$ transformations of the form~$G_q$ have the vector, spinor, and dual spinor transform with expressions involving multiplication by six matrices.  Hence, there are certain transformations within~$SO(8,\mathbb{R})$ which may be expressed via symmetric, left, or right matrix multiplication using the same matrices.  We note that all fourteen of our basis transformations in~$G_2$ have this property.

Of course, we previously discussed (Section~\ref{ch:Division_Algebras_Applications.Quaternions_Octonions.Triality}) Conway's claim as it applied to transformations involving octonions.  The matrix~$\chi$ contains three octonions \hbox{$a, b, c \in \mathbb{O}$}.  
  We find by examining the transformations on~$a,b,c$ for the~$A_l$ transformation of type~$1$
$$
A_\l(\alpha)(\chi) = A_\l(\alpha)\left( \begin{array}{cc|c} t+z & a & \overline{c} \\ \overline{a} & t-z & b \\\hline c & \overline{b} & n \end{array} \right) = \left( \begin{array}{cc|c} t+z & a_1 & \overline{c_1} \\ \overline{a_1} & t-z & b_1 \\\hline c_1 & \overline{b_1} & n \end{array} \right)$$
\noindent that the octonions~$a$,~$b$, and~$c$ transform in the vector, spinor, and dual-spinor pieces via
$$ a \to a_1 = e^{-l\frac{\alpha}{2}}j\l \left( -j\l \left( e^{\l\frac{\alpha}{2}}i\l \left((-i\l) a (-i\l)^\dagger \right) (e^{\l\frac{\alpha}{2}}i\l)^\dagger \right) (-j\l)^\dagger \right) (e^{-l\frac{\alpha}{2}}j\l)^\dagger $$
$$ b \to b_1 = e^{-l\frac{\alpha}{2}}j\l \left( -j\l \left( e^{\l\frac{\alpha}{2}}i\l \left(-i\l  \left(b\right) \right) \right) \right) $$
$$ c \to c_1 = \left( \left( \left( c (-i\l)^\dagger \right) (e^{\l\frac{\alpha}{2}}i\l)^\dagger \right) (-j\l)^\dagger \right) (e^{-l\frac{\alpha}{2}}j\l)^\dagger $$
\noindent However, since all three types of the~$A_\l$ transformation are the same, we see that~$a$,~$b$, and~$c$ undergo the same transformation using conjugation, left-sided, and right-sided multiplication by octonions.  These octonions transform the same way using different forms of multiplication by the same octonions.  This is a stronger condition than Conway gives in~\cite{conway} for triality, and we call this property {\it strong triality}.  Equivalently, a transformation in~$SL(3,\mathbb{O})$ which is type independent exhibits strong triality.  Hence, every transformation in~$G_2$ gives octonionic transformations which have strong triality.

Finally, we note that~$SO(8,\mathbb{R}) \in SL(3,\mathbb{O})$ also exhibits triality.  
There is an~$SO(7)$ group of each type~$a = 1,2,3$, which consists of~$G_2$ along with the seven transformations~$S^{(a)}_{q}$.  We expand any one~$SO(7)$ group to~$SO(8)$ by including the~$R^{(a)}_{xq}$ transformations of type~$a$.  The linear dependencies in the algebra allow us to express any type~$2$ or type~$3$ transformation of the form~$\dot R^{(a)}_{xq}$ or~$\dot S^{(a)}_{q}$ in terms of type~$1$ transformations \footnote{These expressions are cyclic in~$a = 1, 2, 3$, so that we may write any~$SO(8)$ transformation of type~$a$ in terms of transformations of another type.}

$$  \dot R^{(2)}_{xq} = -\frac{1}{2}\dot R^{(1)}_{xq} - \frac{1}{2}\dot S^{(1)}_q \hspace{2cm} \dot S^{(2)}_q = \frac{3}{2} \dot R^{(1)}_{xq} - \frac{1}{2}\dot S^{(1)}_q   $$
$$  \dot R^{(3)}_{xq} = -\frac{1}{2}\dot R^{(1)}_{xq} + \frac{1}{2}\dot S^{(1)}_q \hspace{2cm} \dot S^{(3)}_q = -\frac{3}{2} \dot R^{(1)}_{xq} - \frac{1}{2}\dot S^{(1)}_q   $$

\noindent   
Since the tangent vectors for the transformations on the right-hand side of each equality commute, their corresponding group transformations commute.  Hence, we may express the corresponding type~$2$ or type~$3$ group transformation on the left-hand side in terms of type~$1$ transformations.
  For instance, the expression for~$\dot R^{(2)}_{x\l}$ given above leads to the group transformation identity 
$$ R^{(2)}_{x\l}(\alpha)(\chi) = S^{(1)}_\l(-\frac{\alpha}{2})(\chi) \circ R^{(1)}_{x\l}(-\frac{\alpha}{2})$$
  Given that it is possible in the algebra~$so(8)$ to express every transformation in terms of type~$1$ transformations or~$G_2$ transformations, we see that the group~$SO(8)$ exhibits triality.  In particular, while there are three different~$SO(7)$ groups built from~$G_2$, there is only one~$SO(8)$.  

We now show that transformations in~$SO(8,\mathbb{R}) - G_2$ exhibit the more traditional version of triality described by Baez~\cite{baez} and Conway~\cite{conway} which was described in Section~\ref{ch:Division_Algebras_Applications.Quaternions_Octonions.Triality}.  Consider the two transformations~$R^{(2)}_{x\l}(\alpha)$ and~$S^{(1)}_{\l}(-\frac{\alpha}{2}) \circ R^{(1)}_{x\l}(-\frac{\alpha}{2})$, which were just shown to be equal, and their effect on the octonions~$a,b,c \in \mathbb{O}$ which results from them being applied to~$\chi$.  
Following the construction given in Section~\ref{ch:E6_basic_structure.LieGroup}, explicit calculation shows that~$a$,~$b$, and~$c$ transform under~$R^{(2)}_{x\l}(\alpha)(\chi)$ according to

$$\begin{array}{ccccc} a \to a e^{-\l\frac{\alpha}{2}} & \hspace{1cm} & b \to e^{\l\frac{\alpha}{2}} b e^{-\l\frac{\alpha}{2}} & \hspace{1cm} &  c \to e^{\l\frac{\alpha}{2}}c \end{array}$$

\noindent while under~$S^{(1)}_{\l}(-\frac{\alpha}{2}) \circ R^{(1)}_{x\l}(-\frac{\alpha}{2})(\chi)$, we see that~$b$ and~$c$ transform via left-sided and right-sided multiplication

$$\begin{array}{ccc}
b& \to& e^{-\l\frac{\alpha}{4}}k\l \left( -k\l \left( e^{-\l\frac{\alpha}{4}}j\l \left( -j\l \left( e^{-\l\frac{\alpha}{4}}i\l \left( -i\l \left( e^{-\l\frac{\alpha}{4}} (b) \right)\right)\right)\right)\right)\right)\\
c& \to& \left( \left( \left( \left( \left( \left( (c) e^{\l\frac{\alpha}{4}} \right) i\l \right) (-i\l e^{\l\frac{\alpha}{4}}) \right) j\l \right) (-j\l e^{\l\frac{\alpha}{4}}) \right) k\l \right) (-k\l e^{\l\frac{\alpha}{4}})
\end{array}$$

\noindent and~$a$ transforms via multiplication by conjugate octonions, although we have omitted the expression for this transformation due to space constraints.  Due to the equalities of the group transformations, we see that these new expressons for~$a$,~$b$, and~$c$ are equal to the previous expressions for all values of~$\alpha$.  In particular, we have 

$$\begin{array}{ccc}
e^{\l\frac{\alpha}{2}} b e^{-\l\frac{\alpha}{2}}  &=& e^{-\l\frac{\alpha}{4}}k\l \left( -k\l \left( e^{-\l\frac{\alpha}{4}}j\l \left( -j\l \left( e^{-\l\frac{\alpha}{4}}i\l \left( -i\l \left( e^{-\l\frac{\alpha}{4}} (b) \right)\right)\right)\right)\right)\right)\\
e^{\l\frac{\alpha}{2}}c  &=&  \left( \left( \left( \left( \left( \left( (c) e^{\l\frac{\alpha}{4}} \right) i\l \right) (-i\l e^{\l\frac{\alpha}{4}}) \right) j\l \right) (-j\l e^{\l\frac{\alpha}{4}}) \right) k\l \right) (-k\l e^{\l\frac{\alpha}{4}})
\end{array}
$$

\noindent  That is, we may produce the~$b \to e^{\l \frac{\alpha}{2}} b e^{-\l \frac{\alpha}{2}}$ action by nesting seven left-sided multiplications, and we may also produce the action of multiplying the octonion~$c$ on the left by nesting seven right-sided multiplications.  The expression for~$a$ would show that multiplying~$a$ on the right by an octonion produces the same action as conjugating~$a$ with a series of octonions.  These identities hold for all values of~$\alpha \in \mathbb{R}$.  Also, these expressions are not unique.\footnote{  Note, too, that each factor~$e^{-\l\frac{\alpha}{4}}$ and~$e^{\l\frac{\alpha}{4}}$ next to~$b$ and~$c$, respectively, on the right hand side of each equality generates a rotation in the~$1\leftrightarrow \l$ plane and is perpendicular to the other three planes being by the~$S^{(1)}_{\l}$ transformation. Hence, we are free to commute the phase~$e^{-\l\frac{\alpha}{4}}$ across any pair of factors of the form~$(e^{-\l\frac{\alpha}{4}}q)(-q)$ in the expression for~$b$ and similarly in the expressions for~$c$ and~$a$.}
All of the equivalent ways of writing type~$2$ or type~$3$ transformations in~$SO(8) - G_2$ in terms of type~$1$ transformations (and cyclic permutations) will produce similar identities involving seven multiplications.

\section{Type Transformation}
\label{ch:E6_basic_structure.Type_Transformation}

Our construction of~$3 \times 3$ Lorentz transformations produces equivalent~$SO(9,1)$ groups of type~$1$, type~$2$ and type~$3$ in the following sense:  The discrete type transformation~\hbox{$\mathbf{T}: M_3(\mathbb{O}) \to M_3(\mathbb{O})$} given by 

$$\mathbf{T}(M) = \mathcal{T} M \mathcal{T}^\dagger$$
with
$$\mathcal{T} = \left( \begin{array}{ccc} 0 & 0 & 1 \\ 1 & 0 & 0 \\ 0 & 1 & 0\end{array} \right)$$
cycles the matrices used in our ~$3 \times 3$ finite Lorentz transformations from type~$1$ to type~$2$ and type~$3$ matrices
$$ \mathbf{T}(M^{(1)}) = M^{(2)} \hspace{1.5cm}  \mathbf{T}(M^{(2)}) = M^{(3)} \hspace{1.5cm} \mathbf{T}(M^{(3)}) = M^{(1)}$$
Hence, we have~$\mathbf{T}^3 = Id$ and~$\mathbf{T}^2 = \mathbf{T}^{-1}$.  In addition, the type transformation is in~$SL(3,\mathbb{O})$ since 
$$\det(\mathbf{T}(\chi)) = \det(\chi)$$
While it is {\bf not} one of our basis group transformations, we have the identities

$$ \mathbf{T} = R^{(3)}_{xz}(-\pi) \circ R^{(1)}_{xz}(-\pi) \hspace{1cm} \mathbf{T}^{-1} = R^{(2)}_{xz}(\pi) \circ R^{(1)}_{xz}(\pi) $$
We may expand~$\mathbf{T}^3 = Id$ for these expressions, showing that
$$\begin{array}{ccl}
 \mathbf{T}^{-1} &=& R^{(3)}_{xz}(-\pi) \circ R^{(1)}_{xz}(-\pi) \circ R^{(3)}_{xz}(-\pi) \circ R^{(1)}_{xz}(-\pi)\\
   \mathbf{T} &=& R^{(2)}_{xz}(\pi) \circ R^{(1)}_{xz}(\pi) \circ R^{(2)}_{xz}(\pi) \circ R^{(1)}_{xz}(\pi) 
\end{array}
$$
We also note that
$$ Id = \mathbf{T} \circ \mathbf{T}^{-1}  = R^{(3)}_{xz}(-\pi) \circ R^{(1)}_{xz}(-\pi) \circ R^{(2)}_{xz}(\pi) \circ R^{(1)}_{xz}(\pi) $$
and that
$$ Id = R^{(2)}_{xz}(\pi) \circ R^{(1)}_{xz}(\pi) \circ R^{(2)}_{xz}(\pi) \circ R^{(1)}_{xz}(\pi) \circ R^{(2)}_{xz}(\pi) \circ R^{(1)}_{xz}(\pi)$$
Combining these two expressions and solving for~$\mathbf{T}$, we obtain the expression
$$ \mathbf{T} = R^{(1)}_{xz}(\pi) \circ R^{(3)}_{xz}(\pi) \circ R^{(2)}_{xz}(\pi) \circ R^{(1)}_{xz}(\pi)$$
which involves all three types.  Similarly, many additional expressions for~$\mathbf{T}$ and~$\mathbf{T}^{-1}$ may be obtained from these expressions.

We also note that each of these expressions may be expanded into a {\it continuous type transformation} from~$\mathbb{R} \to SL(3,\mathbb{O})$ by letting the fixed angle~$\pi$ or ~$-\pi$ become an arbitrary angle~$\alpha$ which varies over~$\mathbb{R}$.  The resulting transformations are not one-parameter subgroups of~$SL(3,\mathbb{O})$, but they do connect a type~$1$ transformation with either a type~$2$ or type~$3$ transformation.

\subsection{Type Independent and Dependent Subgroups}
\label{ch:E6_basic_structure.Type_Transformation.Subgroups}
We list here some important groups which contain the type transformation.
The group~$\langle R^{(1)}_{xz}, R^{(2)}_{xz}, R^{(3)}_{xz} \rangle$ is the standard representation of~$SO(3,\mathbb{R})$, and we label this group~$SO(3,\mathbb{R})_s$ using the subscript~$s$ to designate it as the standard representation.
  This group obviously contains~$\mathbf{T}$, as does the standard representation 
$$ SL(3,\mathbb{R})_s = \langle R^{(1)}_{xz}, R^{(2)}_{xz}, R^{(3)}_{xz}, B^{(1)}_{tz}, B^{(2)}_{tz}, B^{(1)}_{tx}, B^{(2)}_{tx}, B^{(3)}_{tx} \rangle$$
Using~$\l$ as our preferred complex unit, we see that the standard representation 
$$ SU(3,\mathbb{C})_s = \langle R^{(1)}_{xz}, R^{(2)}_{xz}, R^{(3)}_{xz}, R^{(1)}_{x\l}, R^{(2)}_{x\l}, R^{(1)}_{z\l}, R^{(3)}_{z\l}, R^{(3)}_{z\l}\rangle$$
of~$SU(3,\mathbb{C})$, which consists of~$3\times 3$ special unitary matrices, contains~$SO(3,\mathbb{R})_s$.  Hence, $SU(3,\mathbb{C})_s$ also contains~$\mathbf{T}$, as does, of course, the group
$$ SL(3,\mathbb{C})_s = SU(3,\mathbb{C})_s \cup \langle B^{(1)}_{tz}, B^{(2)}_{tz}, B^{(1)}_{tx}, B^{(2)}_{tx}, B^{(3)}_{tx}, B^{(1)}_{t\l}, B^{(2)}_{t\l}, B^{(3)}_{t\l} \rangle$$
These four groups~$SO(3,\mathbb{R})_s$,~$SL(3,\mathbb{R})_s$,~$SU(3,\mathbb{C})_s$, and~$SL(3,\mathbb{C})_s$ are important because they contain the type transformation~$\mathbf{T}$.  If, for instance, the transformation~$R^{(1)}$ is in a group~$G$ which has one of these groups as a subgroup,  then the subgroup forces~$G$ to contain~$R^{(2)}$ and~$R^{(3)}$ as well.  On the other hand, if we add one of those four groups to any type-specific group~$G$, then the new group must contain~$R^{(1)}$,~$R^{(2)}$, and~$R^{(3)}$ for each transformation~$R \in G$.  The new group, in most cases, becomes much larger.

We are careful to point out that the standard representations~$SO(3,\mathbb{R})_s$ and~$SU(3,\mathbb{C})_s$ are {\bf not} our preferred representations for the groups~$SO(3,\mathbb{R})^C$ and~$SU(3,\mathbb{C})^C$, which are subgroups of~$G_2$.  For any~$R \in G_2$, we have~$R^{(1)} = R^{(2)} = R^{(3)}$.  Thus, we see that~$SU(3,\mathbb{C})_s$ and~$SU(3,\mathbb{C})^C$  are both type independent, but in~$SU(3,\mathbb{C})_s$ the transformations~$R^{(1)}, R^{(2)}$ and~$R^{(3)}$ are distinct while in in~$SU(3,\mathbb{C})^C$ the three transformations are equal.

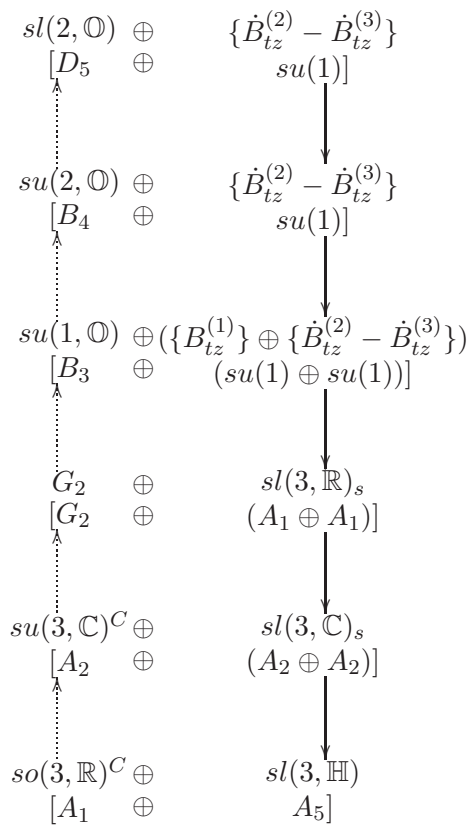
\begin{figure}[htbp]
\[
\xymatrixcolsep{1pt}
\xymatrixrowsep{30pt}
\xymatrix@M=0pt{
*++\txt{$sl(2,\mathbb{O})$  \\ $[D_5 $} & \txt{$\oplus$ \\ $\oplus$} & *++\txt{$ \lbrace \dot B^{(2)}_{tz} - \dot B^{(3)}_{tz} \rbrace $ \\ $ su(1)]$} \ar[d]<1ex> \\
*++\txt{$su(2,\mathbb{O})$ \\ $[B_4 $} \ar@{.>}[u]<1ex>  & \txt{$\oplus$ \\ $\oplus$} & *++\txt{$ \lbrace \dot B^{(2)}_{tz} - \dot B^{(3)}_{tz} \rbrace $\\ $ su(1)]$} \ar[d]<1ex> \\
*++\txt{$su(1,\mathbb{O})$ \\ $[B_3$} \ar@{.>}[u]<1ex> & \txt{$\oplus$ \\ $\oplus$} & *++\txt{$ (\lbrace B^{(1)}_{tz} \rbrace \oplus \lbrace \dot B^{(2)}_{tz} - \dot B^{(3)}_{tz} \rbrace)$ \\ $ (su(1) \oplus su(1))]$}\ar[d]<1ex> \\
*++\txt{$G_2 $ \\ $[G_2 $} \ar@{.>}[u]<1ex> & \txt{$\oplus$ \\ $\oplus$} & *++\txt{$ sl(3,\mathbb{R})_s$\\$(A_1 \oplus A_1)]$} \ar[d]<1ex> \\
*++\txt{$su(3,\mathbb{C})^C $ \\ $ [A_2$} \ar@{.>}[u]<1ex>  & \txt{$\oplus$\\ $\oplus$}  & *++\txt{$ sl(3,\mathbb{C})_s$\\$(A_2 \oplus A_2)]$} \ar[d]<1ex> \\
*++\txt{$so(3,\mathbb{R})^C $ \\ $[A_1 $} \ar@{.>}[u]<1ex> & \txt{$\oplus$ \\ $\oplus$} & *++\txt{$ sl(3,\mathbb{H})$\\ $ A_5]$} \\
}
\]
\caption{Type dependent and independent subalgebras of~$E_6$} 
\label{figure:type_subalgebras}
\end{figure}

We use the type transformation~$\mathbf{T}$ to provide insight into the structure of~$SL(3,\mathbb{O})$.  The algebras~$g$ in the left column of Figure~\ref{figure:type_subalgebras} are subalgebras of the type~$1$ copy of~$sl(2,\mathbb{O})$, while each algebra~$g^\prime$ in the right column is the largest algebra of~$sl(3,\mathbb{O})$ such that~$g \oplus g^\prime$ is still~$\mathbb{R}$-simple.  When we restrict~$g$ to a smaller subalgebra of~$sl(2,\mathbb{O})$, it is sometimes possible to exand the type independent subalgebra~$g^\prime$ to a larger subalgebra of~$sl(3,\mathbb{O})$.  Each arrow in the diagram indicates containment.

\section{Reduction of~$\mathbb{O}$ to~$\mathbb{H}$,~$\mathbb{C}$, and~$\mathbb{R}$}
\label{ch:E6_basic_structure.Reduction_of_O}
We find subalgebras of~$E_6$ by restricting our octonionic~$E_6$ generators to be quaternionic, complex, or real.

In the first case, we limit ourselves to the quaternionic case by using only the algebra generators~$\dot R_{ab}$ or~$\dot R_{aq}$ where~$a,b \in \lbrace t,x,z \rbrace$ and~$q \in \lbrace k,kl,l\rbrace$.  Of course, we must also include those generators in~$SO(7)$ which mix up these imaginary quaternionic pieces.  Given our definition of~$A_q, G_q,$ and~$S_q$, the particular combinations
$$ \dot G_k - \dot S^{(1)}_k \hspace{1.5cm} \dot G_{kl} - \dot S^{(1)}_{kl} \hspace{1.5cm} \dot G_l - \dot S^{(1)}_l$$
provide the permutations of~$\lbrace k,kl,l\rbrace$ while fixing~$\lbrace i,j,il,jl\rbrace$.  This provides~$35$ transformations, of which~$14$ are boosts.  The~$21$ compact generators form the algebra~$su(3,\mathbb{H})$, a real form of~$C_3 = sp(2\cdot 3)$, while all~$35$ together form~$sl(3,\mathbb{H})$, a real form of~\hbox{$A_5 = su(6,\mathbb{C})$}.  Restricting only to type~$1$ transformations, we obtain~$10$ rotations and~$5$ boosts.  This restriction reduces the algebra~$sl(3,\mathbb{H})$ to~$sl(2,\mathbb{H}) = so(5,1)$, a real form of~$D_3 = so(6)$.  The algebra containing only the~$10$ rotations is~$su(2,\mathbb{H}) = so(5)$, a real form of~$C_2 = sp(2\cdot 2)$.

We note that~$\dot A_k, \dot A_{kl},$ and~$\dot A_l$ fix our preferred quaternions~$\lbrace k, kl, l\rbrace$ and permute the four {\it orthogonal quaternions}~$\lbrace i, j, jl, il \rbrace$ among themselves.  Hence, we have~\hbox{$sl(3,\mathbb{H}) \oplus so(3)^C$}, where~$\lbrace \dot A_k, \dot A_kl, \dot A_l\rbrace$ form the~$so(3)^C$.  This direct sum structure~$g \oplus so(3)^C$ holds for all quaternionic subalgebras~$g$ of~$sl(3,\mathbb{H})$.

Given the octonion~$\l$, there are at least two interesting ways to break the octonions into quaternionic subalgebras.  As discussed above, we chose a preferred quaternion algebra containing~$k, k\l$ and~$\l$.  However, another choice would be to leave~$\l$ alone, instead using~$i, j,$ and~$k$.
Using this quaternionic sublagebra instead of~$\lbrace 1, k,kl,l\rbrace$, direct calculation shows that
$$\dot A_i + \frac{1}{3} \dot G_i + \frac{2}{3} \dot S^{(1)}_i \hspace{1.5cm} \dot A_j + \frac{1}{3} \dot G_j + \frac{2}{3} \dot S^{(1)}_j \hspace{1.5cm} \dot A_k + \frac{1}{3} \dot G_k + \frac{2}{3} \dot S^{(1)}_k$$
permutes our imaginary quaternionic units~$\lbrace i, j, k\rbrace$ while 
$$\dot G_i - \dot S^{(1)}_i \hspace{1.5cm}\dot G_j - \dot S^{(1)}_j \hspace{1.5cm} \dot G_k - \dot S^{(1)}_k $$
form the~$so(3)$ which fixes our quaternionic subalgebra and permutes the orthogonal quaternions~$\lbrace kl, jl, il, l \rbrace$.  However, neither of these Lie algebras contain any of the preferred Casimir operators~$\lbrace \dot B^{(1)}_{tz}, \dot B^{(2)}_{tz}, \dot R^{(1)}_{xl}, \dot A_{l}, \dot G_l, \dot S^{(1)}_{l} \rbrace $.

As previously discussed, we obtain the classical Lie algebras~$su(3,\mathbb{C})_s$ and~$sl(3,\mathbb{C})_s$ by choosing those transformations~$\dot R_{ab}$ and~$\dot R_{aq}$ of all three types, where~$a,b \in \lbrace t, x, z \rbrace$ and~$q = l$.  As there is only one octonionic unit used to form~$\mathbb{C}$, we do not need to use any of the transformations from~$SO(7)$.  While using all~$16$ transformations gives~$sl(3,\mathbb{C})_s$, a real form of~$A_2 \oplus A_2 = su(3,\mathbb{C}) \oplus su(3,\mathbb{C})$ with~$8$ boosts, we obtain~$su(3,\mathbb{C})_s$ by using only the~$8$ compact generators.  Restricting ourselves to the type~$1$ case reduces these two algebras to~$sl(2,\mathbb{C})_s = so(3,1)_s$ and~$su(3,\mathbb{C})_s$, which are real forms of~\hbox{$D_2 = su(2,\mathbb{C}) \oplus su(2,\mathbb{C})$} and~$A_1 = su(3,\mathbb{C})$.

We also note that when we restrict~$sl(3,\mathbb{C})_s$ to~$sl(2,\mathbb{C})_s$, the smaller algebra no longer contains the type transformation~$\mathbf{T}$ but it uses the octonionic direction~$\l$.  Direct calculation shows that~$sl(2,\mathbb{C})_s \oplus so(6)$, where~$so(6)$ fixes~$\l$, is an~$\mathbb{R}$-simple subalgebra of~$sl(3,\mathbb{O})$.

Finally, by restricting to real transformations, we see that in the~$3 \times 3$ case, the nine elements~$\dot R^{(T)}_{ab}, a,b \in \lbrace t,x,z \rbrace, T = 1,2,3$ along with the relation~$\dot R^{(1)}_{tz}+\dot R^{(2)}_{tz} + \dot R^{(3)}_{tz} = 0$ give a real form of~$A_2 = su(3,\mathbb{C})$ with~$5$ non-compact elements.  Note that this is {\bf not} our standard~$su(3,\mathbb{C})^C \subset G_2$, nor the standard~$su(3,\mathbb{C})_s$ as it is not complex.  It may be restricted to either~$so(3,\mathbb{R})_s$, whose group contains the type transformation, or~$so(2,1)_s$, which is a type~$1$ non-compact form of~$A_1 = so(3,\mathbb{R})$.

\begin{figure}[htbp]
\begin{center}
\begin{minipage}{6in}
\begin{center}
\xymatrixcolsep{8pt}
\xymatrix@M=1pt@H=0pt{
   & & & *++\txt{$sl(3,\mathbb{O})$\\$[E_6]$} & & & \\
   & & & *++\txt{$sl(3,\mathbb{H})$\\$[A_5]$}\ar[u] & & & \\
   & & & *++\txt{$sl(3,\mathbb{C})_s$\\$[A_2\oplus A_2]$}\ar[u] & & & \\
*++\txt{$sl(2,\mathbb{O})$\\$[D_5]$} \ar[uuurrr] & *++\txt{$sl(2,\mathbb{H})$\\$[A_3 = D_3]$} \ar[uurr] \ar[l] & *++\txt{$sl(2,\mathbb{C})_s$\\$[A_1 \oplus A_1]$} \ar[ur] \ar[l] & & *++\txt{$su(3,\mathbb{C})_s$\\$[A_2]$} \ar[ul] \ar[r]& *++\txt{$su(3,\mathbb{H})$\\$[C_3]$} \ar[uull] \ar[r]& *++\txt{$su(3,\mathbb{O})$\\$[F_4]$} \ar[uuulll]\\
   & & & *++\txt{$su(2,\mathbb{C})_s$\\$[A_1]$} \ar[d] \ar[ul] \ar[ur] & & & \\
   & & & *++\txt{$su(2,\mathbb{H})$\\$[B_2 = C_2]$}\ar[d] \ar[uurr] \ar[uull] & & & \\
   & & & *++\txt{$su(2,\mathbb{O})$\\$[B_4]$}\ar[uuurrr] \ar[uuulll] & & & \\
}
\caption{Subalgebras~$sl(n,\mathbb{F})$ and~$su(n,\mathbb{F})$ of ~$sl(3,\mathbb{O})$}
\label{fig:division_algebra_subalgebras_of_E6}
\end{center}
\end{minipage}
\end{center}
\end{figure}
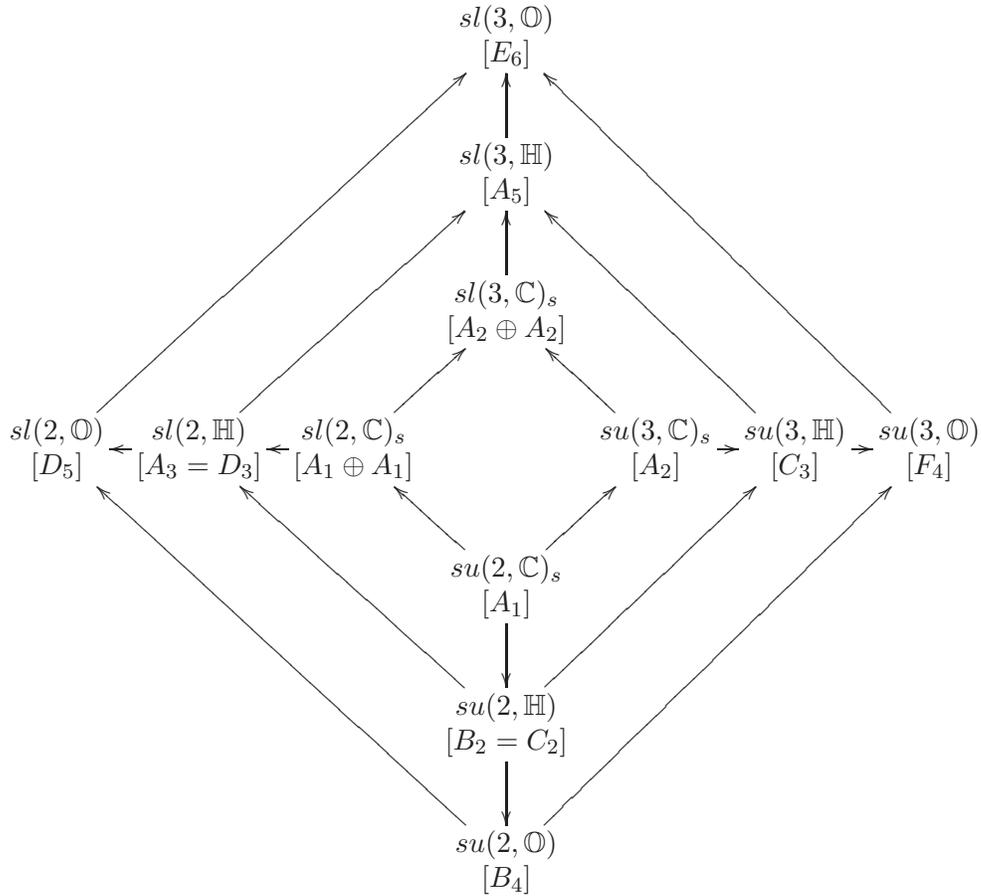

The results of restricting~$sl(3,\mathbb{O})$ to~$sl(n,\mathbb{F})$ for~$n = 1,2,3$ and~$\mathbb{F} = \mathbb{R}, \mathbb{C}, \mathbb{H}, \mathbb{O}$ is given in Figure~\ref{fig:division_algebra_subalgebras_of_E6}.  
For each algebra~$g$ in Figure~\ref{fig:division_algebra_subalgebras_of_E6}, we  list in Figure~\ref{fig:division_algebra_subalgebras_and_perp_algebras_of_E6} the maximal subalgebra~$g^\prime$ of~$E_6$ such that we have the structure~$g \oplus g^\prime$.  Here~$so(6)$ denotes the algebra which permutes~$\lbrace i,j,k,k\l,j\l,i\l\rbrace$ but fixes~$\l$.  While~$so(6) \not \subset G_2$, we do have~$su(3,\mathbb{C})^C \subset so(6)$.  The solid arrow indicates inclusion when the algebra~$g$ is expanded to a larger Lie algebra and the dotted arrow indicates the result of restricting the algebra~$g^\prime$ to a smaller algebra.

\begin{landscape}
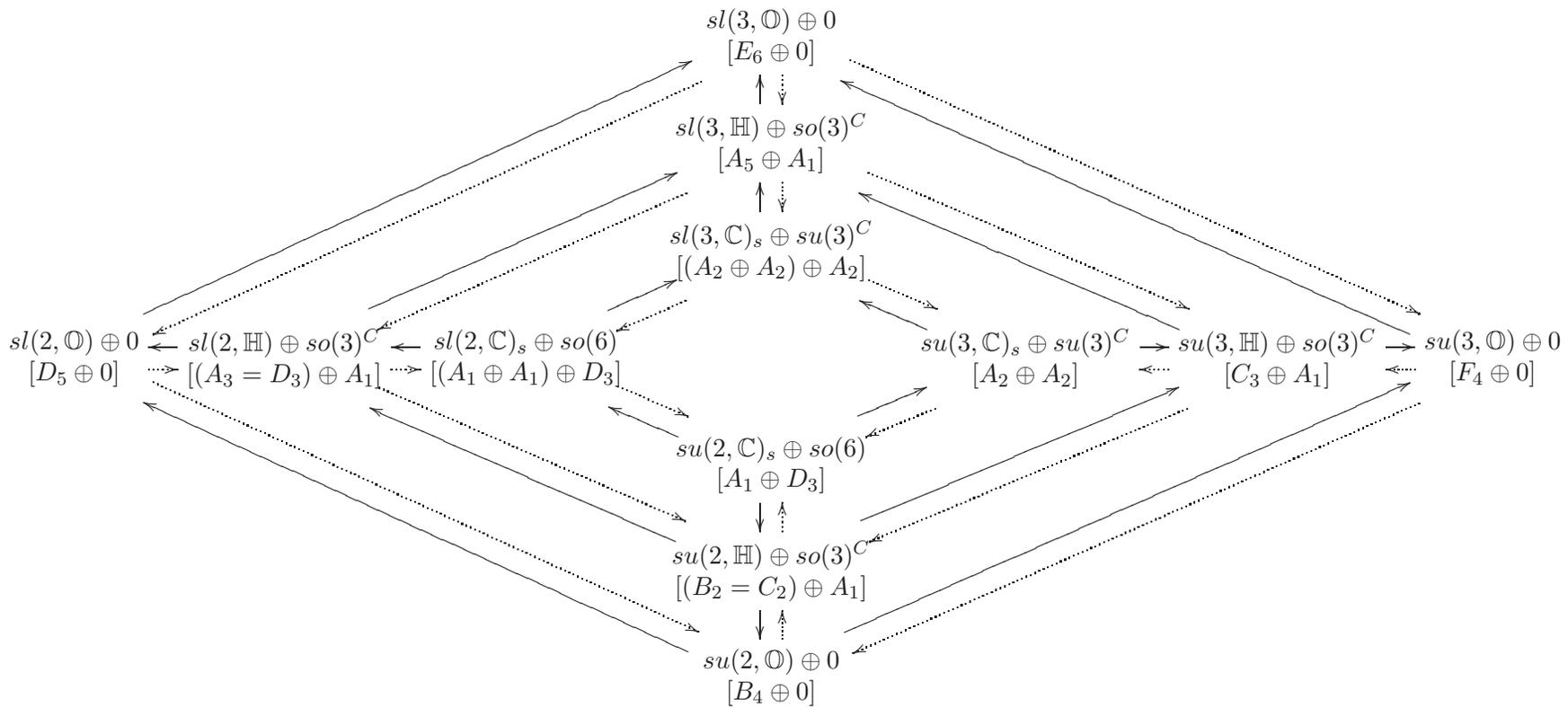
\begin{figure}[htbp]
\begin{center}
\xymatrixcolsep{12pt}
\xymatrixrowsep{12pt}
\xymatrix@M=2pt@H=2pt{
   & & & *++\txt{$sl(3,\mathbb{O})\oplus 0 $\\$[E_6 \oplus 0]$} \ar@{.>}[d]<1ex> \ar@{.>}[dddrrr]<1ex> \ar@{.>}[dddlll]<1ex>& & & \\
   & & & *++\txt{$sl(3,\mathbb{H}) \oplus so(3)^C$\\$[A_5 \oplus A_1]$}\ar[u]<1ex> \ar@{.>}[d]<1ex> \ar@{.>}[ddll]<1ex> \ar@{.>}[ddrr]<1ex>& & & \\
   & & & *++\txt{$sl(3,\mathbb{C})_s \oplus su(3)^C$\\$[(A_2\oplus A_2) \oplus A_2 ]$}\ar[u]<1ex> \ar@{.>}[dl]<1ex> \ar@{.>}[dr]<1ex>& & & \\
*++\txt{$sl(2,\mathbb{O}) \oplus 0$\\$[D_5 \oplus 0]$} \ar@{.>}[r]<-1ex> \ar[uuurrr]<1ex> \ar@{.>}[rrrddd]<1ex>& *++\txt{$sl(2,\mathbb{H}) \oplus so(3)^C$\\$[(A_3 = D_3) \oplus A_1]$} \ar@{.>}[r]<-1ex> \ar[uurr]<1ex> \ar[l]<-1ex> \ar@{.>}[ddrr]<1ex> & *++\txt{$sl(2,\mathbb{C})_s \oplus so(6)$\\$[(A_1 \oplus A_1) \oplus D_3 ]$} \ar[ur]<1ex> \ar@{.>}[dr]<1ex> \ar[l]<-1ex> & & *++\txt{$su(3,\mathbb{C})_s \oplus su(3)^C$\\$[A_2 \oplus A_2]$} \ar[ul]<1ex> \ar@{.>}[dl]<1ex> \ar[r]<1ex>& *++\txt{$su(3,\mathbb{H}) \oplus so(3)^C$\\$[C_3\oplus A_1]$} \ar[uull]<1ex> \ar@{.>}[ddll]<1ex> \ar[r]<1ex> \ar@{.>}[l]<1ex>& *++\txt{$su(3,\mathbb{O}) \oplus 0$\\$[F_4\oplus 0]$} \ar[uuulll]<1ex> \ar@{.>}[dddlll]<1ex> \ar@{.>}[l]<1ex>\\
   & & & *++\txt{$su(2,\mathbb{C})_s \oplus so(6)$\\$[A_1 \oplus D_3]$} \ar[d]<-1ex> \ar[ul]<1ex> \ar[ur]<1ex> & & & \\
   & & & *++\txt{$su(2,\mathbb{H}) \oplus so(3)^C$\\$[(B_2 = C_2) \oplus A_1]$}\ar[d]<-1ex> \ar[uurr]<1ex> \ar[uull]<1ex> \ar@{.>}[u]<-1ex> & & & \\
   & & & *++\txt{$su(2,\mathbb{O})\oplus 0$\\$[B_4 \oplus 0]$}\ar[uuurrr]<1ex> \ar[uuulll]<1ex> \ar@{.>}[u]<-1ex> & & & \\
}
\caption{Subalgebras~$sl(n,\mathbb{F}) \oplus g_1$ and~$su(n,\mathbb{F}) \oplus g_1$ of ~$sl(n,\mathbb{O})$}
\label{fig:division_algebra_subalgebras_and_perp_algebras_of_E6}
\end{center}
\end{figure}
\end{landscape}

\section{A Subalgebra fixing~$\l$}
\label{ch:E6_basic_structure.fix_l}

We examine here a subalgebra of~$sl(3,\mathbb{O})$ which fixes a preferred octonionic unit in~$\chi$.  Given 
$$\chi = \left( \begin{array}{ccc} t+z & \overline{a} & c \\ a & t-z & \overline{b} \\ \overline{c} & b & n \end{array} \right) \in M_3(\mathbb{O})$$
we choose to treat the  upper~$2 \times 2$ submatrix as a vector.  We call the unit octonion~$\l$ in this vector {\it special} and look for transformations which fix the coefficient~$a_\l$ of~$\l$ near the identity.  That is, we find~$Stab(\l)$ for this special~$\l$.  At the Lie algebra level, this implies that the coefficient of our special~$\l$ in the tangent vector must vanish, ensuring that~$\l$ is fixed by the corresponding one-parameter transformation in the group.

Examination of the~$78$ tangent vectors in our preferred basis for~$SL(3,\mathbb{O})$ results in~$52$ elements which fix our special~$\l$.  Direct computation shows these elements form a~$52$-dimensional subalgebra~$Stab(\l)$ of~$sl(3,\mathbb{O})$.
However,~$Stab(\l)$ is not a simple algebra.  It contains
$$su(3,\mathbb{C})^C = \langle A_i, A_j, A_k, A_{k\l}, A_{j\l}, A_{i\l}, A_{\l}, G_\l \rangle$$
which is a subalgebra of~$G_2$.  The algebra~$su(3,\mathbb{C})^C$ can be expanded to the algebra~$so(6,\mathbb{R})$ of type~$1$ which fixes~$\l$ via 
$$\begin{array}{ccl}
so(6,\mathbb{R}) &= & su(3,\mathbb{C})^C \cup \langle G_i + 2S^{(1)}_i, G_j + 2S^{(1)}_j, G_k + 2S^{(1)}_k, \\
      &  & \hspace{2cm} G_{k\l} + 2S^{(1)}_{k\l}, G_{j\l} + 2S^{(1)}_{j\l}, G_{i\l} + 2S^{(1)}_{i\l}, S^{(1)}_\l \rangle
\end{array}$$

We note that we may expand~$so(6,\mathbb{R})$ into~$so(8,1,\mathbb{R})$ which does {\bf not} contain our preferred~$so(8,\mathbb{R})$.  In particular, this~$so(8,1,\mathbb{R})$ is an algebra of type~$1$ and does not contain~$G_2$.  The basis of this~$so(8,1,\mathbb{R})$ is expanded from that of~$so(6,\mathbb{R})$ by including the transformations~$R^{(1)}_{xz}, B^{(1)}_{tx}, B^{(1)}_{tz}$, which form~$so(2,1,\mathbb{R})$, and 
$$ B^{(1)}_{tq}, R^{(1)}_{xq}, R^{(1)}_{zq}$$
where~$q = i,j,k,k\l,j\l,i\l$, but not~$q = \l$.
We note that~$so(2,1,\mathbb{R})$ commutes with~$so(6,\mathbb{R})$ since it contains purely real transformations.
  We note that each of the~$3 \times 3$ Lorentz transformations corresponding to the Lie algebra elements in~$so(8,1,\mathbb{R})$ are one-parameter transformations.  Hence, then the group~$SO(8,1,\mathbb{R})$ also fixes the coefficient of~$\l$ in the vector.

At the identity, the tangent vectors of our type~$2$ and type~$3$ transformations do not fix~$\l$.  However, the following real linear combinations of these tangent vectors do fix the~$\l$ of type~$1$:
$$ b_2 = \langle B^{(2)}_{tq} + R^{(2)}_{zq} | q  = i,j,k,k\l, j\l, i\l \rangle$$
$$ b_3 = \langle B^{(3)}_{tq} - R^{(3)}_{zq} | q  = i,j,k,k\l, j\l, i\l \rangle$$
$$ b_\l = \langle B^{(2)}_{tx}+R^{(2)}_{xz}, B^{(3)}_{tx}-R^{(3)}_{xz}, B^{(2)}_{t\l}+R^{(2)}_{x\l}, B^{(3)}_{t\l}-R^{(3)}_{x\l} \rangle$$
No other real linear combinations of tangent vectors fix the~$\l$ of type~$1$.
The algebra~$b = b_2 \oplus b_3 \oplus b_\l$ is abelian, and~$[ so(8,1,\mathbb{R}), b ] \subset b$.  Hence,~$b$ is an ideal, implying~$Stab(\l)$ is not simple.  We note that
$$ [su(3)^C,b_2 ] \subset b_2 \hspace{1.5cm} [su(3)^C,b_3 ] \subset b_3 \hspace{1.5cm} [su(3)^C,b_\l] = 0$$
and that~$ [so(6), b ] \subset b$ since~$[so(6), b_3] \subset b$.
The Killing form for~$E_6$ shows that the elements of~$b_2$, along with~$B^{(2)}_{tz} + R^{(2)}_{xz}$ and~$B^{(2)}_{t\l}+R^{(2)}_{x\l}$, are null transformations.

\section{Gell-Mann Matrices and~$su(3) \subset G_2$}
\label{ch:E6_basic_structure.Gell-Mann}

We list in Table~\ref{table:gell-mann_and_su3} the isomorphism between the~$3 \times 3$ Gell-Mann matrices and the subalgebra~\hbox{$su(3, \mathbb{C})^C \subset G_2$}.  This isomorphism respects~$\left[ \hspace{.15cm}, \hspace{.15cm} \right]$ in the sense that the structure constants are the same for each algebra.    We use~$\iota$ to denote a square root of~$-1$ which commutes with all octonions.  This is required, since the Gell-Mann matrices are Hermitian and use physicists' conventions, while our~$su(3,\mathbb{C})^C$ uses mathematicians' conventions.

\begin{table}[htbp]
\begin{center}
\begin{tabular}{|cl|}
\hline
Octonionic~$su(3)$ & Gell-Mann \\
Transformation & Matrix \\
\hline
$\dot A_k $ & $\lambda_1= - \iota \left(\begin{array}{ccc}0 & 1 & 0  \\1 & 0 & 0  \\0 & 0 & 0  \\\end{array} \right) $ \\
$\dot A_{k\l}$ & $\lambda_2= - \iota \left( \begin{array}{ccc}0 & -i & 0  \\i & 0 & 0  \\0 & 0 & 0  \\\end{array} \right) $\\
$\dot A_\l$ & $ \lambda_3= \iota \left( \begin{array}{ccc}1 &  0 & 0  \\0 & -1 & 0  \\0 & 0 & 0  \\\end{array} \right)$\\
$\dot A_i$ & $\lambda_4= - \iota \left( \begin{array}{ccc} 0 & 0 & 1  \\ 0 & 0 & 0  \\ 1 & 0 & 0  \\ \end{array} \right)$ \\
$\dot A_{i\l}$ & $ \lambda_5=  \iota \left( \begin{array}{ccc}0 & 0 & -i  \\0 & 0 & 0  \\i & 0 & 0  \\ \end{array} \right)$ \\
$\dot A_{j\l}$ & $ \lambda_6= - \iota \left( \begin{array}{ccc}0 & 0 & 0  \\0 & 0 & 1  \\0 & 1 & 0  \\ \end{array} \right) $ \\
$\dot A_{j}$ & $ \lambda_7= - \iota \left( \begin{array}{ccc}0 & 0 & 0  \\0 & 0 & -i  \\0 & i & 0  \\ \end{array} \right) $\\
$ \dot G_\l$ & $\lambda_8= - \iota \sqrt{3} \left( \begin{array}{ccc}\frac{1}{\sqrt{3}} & 0 & 0  \\0 & \frac{1}{\sqrt{3}} & 0  \\0 & 0 & -\frac{2}{\sqrt{3}}  \end{array} \right)$ \\
\hline
\end{tabular}
\caption{\noindent Isomorphism between~$su(3,\mathbb{C})^C$ and Gell-Mann matrices}
\label{table:gell-mann_and_su3}
\end{center}
\end{table}

\newpage{}

\part{The Further Structure of~$E_6$}
\label{ch:E6_further_structure}

In this chapter, we identify additional subalgebras and subalgebra chains of~$sl(3,\mathbb{O})$, our preferred real form of~$E_6$.  We identify real subalgebras of~$sl(3,\mathbb{O})$ by adapting techniques from the study of complex Lie algebras.  In particular, if~$g$ is a real form of a complex Lie algebra~$g^\mathbb{C}$, we use a map~$\phi$ which sends the complex Lie algebra~$g^\mathbb{C}$ to itself while sending the particular real algebra~$g$ to another real form~$g^\prime$ of~$g^\mathbb{C}$.  The algebra~$g \cap g^\prime$ is then an easily identifiable  subalgebra of~$g$.  In addition, the maximal compact subalgebra~$g^\prime_{c}$ of~$g^\prime$ can be easily identified, providing a simple way to identify the subalgebra~$\phi^{-1}\left(g^\prime_{c}\right)$ in the original algebra~$g$.

The chapter is organized as follows.  We discuss the properties of automorphisms of Lie algebras, particularly involutive automorphisms, in Section~\ref{ch:E6_further_structure.Automorphisms_of_Lie_Algebras}.  These automorphisms will later be used to find subalgebras of~$sl(3,\mathbb{O})$ in Sections~\ref{ch:E6_further_structure.three_automorphisms} and~\ref{ch:E6_further_structure.more_applications_of_three_automorphisms}.  While the automorphisms may be used to find subalgebras of~$sl(3,\mathbb{O})$, we still must identify and classify those subalgebras.  Section~\ref{ch:E6_further_structure.Automorphisms_of_Lie_Algebras.id_algebras} reviews the methods we utilize to determine whether a given subalgebra of~$sl(3,\mathbb{O})$ is semi-simple, simple, or neither, in addition to developing methods to classify it as a classical matrix algebra.  We construct three involutive automorphisms in Section~\ref{ch:E6_further_structure.three_automorphisms}, which we then use to find some subalgebras of~$sl(3,\mathbb{O})$ which we did not find in Chapter~\ref{ch:E6_basic_structure}.  These automorphisms are combined in Section~\ref{ch:E6_further_structure.more_applications_of_three_automorphisms} to give new involutive automorphisms and new subalgebras of~$sl(3,\mathbb{O})$.  These subalgebras are used to fill out the chain of subalgebras contained within~$sl(3,\mathbb{O})$ in Section~\ref{ch:E6_further_structure.chains_of_subalgebras}, where we construct a tower of subalgebras in~$sl(3,\mathbb{O})$ based upon our preferred choice of Casimir operators.

\section{Automorphisms of Lie Algebras}
\label{ch:E6_further_structure.Automorphisms_of_Lie_Algebras}
We summarize here the standard definitions and treatment of automorphisms of Lie algebras, according to~\cite{gilmore}.  These methods will be used in Sections~\ref{ch:E6_further_structure.three_automorphisms} and~\ref{ch:E6_further_structure.more_applications_of_three_automorphisms} to identify simple subalgebras of~$sl(3,\mathbb{O})$.

Let~$g$ be a real form of the complex Lie algebra~$g^\mathbb{C}$.
  An {\it automorphism of~$g^\mathbb{C}$} is a map~\hbox{$\phi : g^\mathbb{C} \to g^\mathbb{C}$} which is an isomorphism of the underlying vector space over~$\mathbb{C}$ and satisfies
$$ \phi \left(\left[X,Y\right]\right) = \left[ \phi(X), \phi(Y) \right] \hspace{1cm} X,Y \in g$$
This condition means that~$\phi$ respects the algebraic structure of~$g^\mathbb{C}$.  An {\it involutive automorphism of~$g^\mathbb{C}$} is an automorphism~$\phi$ such that 
$$\phi^2 = \textrm{Id}$$

As discussed in~\cite{gilmore}, involutive automorphisms map simple Lie algebras to simple Lie algebras.  Hence, 

{\bf Lemma: } If~$\phi : g^\mathbb{C} \to g^\mathbb{C}$ is an involutive automorphism of the associated complex Lie algebra~$g^\mathbb{C}$ of a real semi-simple Lie algebra~$g$, then the subalgebra~$g_1 \subset g$ is simple if and only if~$\phi\left( g_1 \right)$ is simple.

 Proof:  Assume that~$g_1$ is not simple.  That is, either~$g_1$ is abelian, or it possesses a proper invariant subalgebra~$h_1$.  If~$g_1$ is abelian, then for every~$x,y \in g_1$, we have 
$$\left[ \phi(x), \phi(y) \right] = \phi \left( \left[ x, y \right] \right) = 0$$
showing that~$\phi(g_1)$ is also abelian.  On the other hand, if~$h_1 \subset g_1$ is a proper invariant subalgebra of~$g_1$, then for every~$x \in g_1, y \in h_1$, we have~$[x,y] \in h_1$.
  This implies
$$\left[ \phi(x), \phi(y) \right] = \phi \left( \left[ x, y \right] \right) \in \phi(h_1)$$
Hence,~$\phi(h_1)$ is a proper invariant subalgebra of~$\phi(g_1)$.  We conclude that   if~$g_1$ is not simple, then~$\phi(g_1)$ is not simple.  Since involutive automorphisms are invertible, we can reverse the arguments above to show that~$g_1$ is not simple if~$\phi(g_1)$ is not simple, thus proving our claim. 

\myendofproof

The set of automorphisms form a group, as do the set of involutive automorphisms.  This property of involutive automorphisms will be utilized in Section~\ref{ch:E6_further_structure.more_applications_of_three_automorphisms}.  Finally, as mentioned in Section~\ref{ch:Lie_groups_Lie_algebras.Basic_Definitions.regular_rep_Killing_form}, we note that if~$\phi$ is an automorphism of a Lie algebra~$g$ with Killing form~$B$, then~$B(\phi(x), \phi(y)) = B(x,y)$ for all~$x,y \in g$.  This interaction between the Killing form of a complex Lie algebra~$g^\mathbb{C}$ and one of its involutive automorphisms allows us to naturally separate~$g^\mathbb{C}$ into two orthogonal subspaces, which we now show.  Further details may be found in~\cite{gilmore}.

Suppose~$g^\mathbb{C}$ is a semi-simple complex Lie algebra and~$\phi : g^\mathbb{C} \to g^\mathbb{C}$ is an involutive automorphism.  The following argument from~\cite{gilmore} shows that we may use~$\phi$ to naturally separate~$g^\mathbb{C}$ into two useful subspaces.  Since~$\phi^2 - \textrm{Id} = \left( \phi - \textrm{Id}\right)\left(\phi + \textrm{Id}\right) = 0$, the involutive automorphism~$\phi$ splits the vector space~$g^\mathbb{C}$ into subspaces~$\p$ and~$\n$ which are associated with the eigenvalues~$+1$ and~$-1$ of~$\phi$, respectively.  Using these eigensubspaces, we write
$$g^\mathbb{C} = \p \oplus \n$$
where we note that this is not a direct sum of Lie algebras!  However, as
$$
\begin{array}{ccl}
\phi\left(g^\mathbb{C}\right) &=& \phi\left(\p\right) \oplus \phi\left(\n\right)\\
                              &=& (+1)\p \oplus (-1)\n
\end{array}
$$
and the Killing form of the semi-simple Lie algebra~$g^\mathbb{C}$ is non-degenerate, then
$$B\left(\p, \n\right) = B\left(\phi\left(\p\right), \phi\left(\n\right)\right) = B\left(+\p,-\n\right) =-B\left(\p,\n \right) = 0$$
showing that~$\p$ and~$\n$ are orthogonal as vector spaces.

We use the Killing form and the involutive automorphism~$\phi$ to show that~$\p$ closes under commutation, that~$\left[ \n, \n \right] \subset \p$, and that~$\left[ \p, \n \right] \subset \n$.
For any~$p_1, p_2 \in \p$ and any~\hbox{$n \in \n$}, we have
$$\begin{array}{ccl}
B\left(\left[p_1,p_2\right], n \right) & = &B\left(\left[\phi(p_1),\phi(p_2)\right], \phi(n) \right) \\
                                         & = &B\left(\left[+p_1, +p_2\right], -n \right) \\
                                         & = &-B\left(\left[p_1, p_2\right], n \right) \\
                                         & = & 0 \\
\end{array}
$$
Hence,~$\left[ \p, \p \right]$ is orthogonal to~$\n$, implying~$\left[ \p, \p \right] \subset \p$ since~$B$ is non-degenerate.    Thus,~$\p$ closes under commutation and is a subalgebra of~$g$.  Note that if we had chosen~$p_1, p_2 \in \n$ instead, then the two plus signs in the commutator in the second line would become negative signs, resulting in no overall sign change.  Hence, the same calculation shows that~$\left[ \n, \n \right] \subset \p$.  Finally,
for~$n_1 \in \n$ and~$p_1,p_2 \in \p$, we have
$$\begin{array}{ccl}
B\left(\left[n_1,p_1\right], p_2 \right) &=& B\left(\left[\phi(n_1), \phi(p_1)\right], \phi(p_2)\right)  \\
                                         &=& B\left(\left[-n_1, p_1\right], p_2\right) \\
                                         &=& -B\left(\left[n_1, p_1\right], p_2\right) \\
                                         &=& 0 \\
\end{array}
$$
This shows~$\left[\n, \p\right]$ is orthogonal to~$\p$, resulting in ~$\left[ \n, \p \right] \subset \n$.
This proves~$\p$ and~$\n$ are orthogonal complementary subspaces of~$g^\mathbb{C}$.

When restricted to~$g$, a slight modification of the involutive automorphism~$\phi$ can be used to produce another real form of~$g^\mathbb{C}$.  We again use~$\phi$ to split~$g$ into subspaces~$\p \oplus \n$, and then introduce the involutive automorphism~$\phi^* : g^\mathbb{C} \to g^\mathbb{C}$ via 
$$\phi^*(p+n) = \phi(p) + \xi \phi(n)$$
where~$p \in \p$,~$n \in \n$, and~$\xi$ is a square root of~$-1$ which commutes with all imaginary units used in the representation of~$g$.
As mentioned in Section~\ref{ch:Lie_groups_Lie_algebras.Basic_Definitions.Lie_Algebras}, the structure constants of a real form~$g$ of~$g^\mathbb{C}$ are real, and since
$$ \left[ \p, \p \right] \subset \p \hspace{1cm} \left[ \p, \xi \n \right] \subset \xi \n \hspace{1cm} \left[ \xi \n, \xi \n, \right] \subset (\xi)^2 \p = (-1)\p$$
then~$\p \oplus \xi \n$ also has real structure constants.  It is a different real form of the complex Lie algebra~$g^\mathbb{C}$.  Further information regarding the interplay between involutive automorphisms, the Killing form, and real forms of a complex Lie algebra may be found in~\cite{gilmore}.

\subsection{Identification of Algebras}
\label{ch:E6_further_structure.Automorphisms_of_Lie_Algebras.id_algebras}
In Sections~\ref{ch:E6_further_structure.three_automorphisms} and~\ref{ch:E6_further_structure.more_applications_of_three_automorphisms}, we will use involutive automorphisms to find many semi-simple real subalgebras of~$sl(3,\mathbb{O})$.  Identifying these subalgebras can be fairly easy if we know they are simple.  Establishing whether a subalgebra is simple or not, however, may be time consuming.  In this section, we outline the various methods we use for identifying simple subalgebras of~$sl(3,\mathbb{O})$, and verifying that they are indeed simple.

In the following discussion,~$g$ and~$g^\prime$ are real Lie algebras.  The complexification of~$g$ is denoted~$g^\mathbb{C}$.  We are primarily interested in the real subalgebras of~$sl(3,\mathbb{O})$, and may utilize a procedure in our specific case even though it does not apply to more general Lie algebras.  We also find it useful to consider the following example to illustrate how the various methods may be used:

{\bf Example 1: }  Let~$\g$ be a subalgebra of~$sl(3,\mathbb{O})$ whose dimension is~$36$, and assume that~$\g$ contains precisely~$4$ Casimir operators.  We note that~$\g$ has rank four if it is simple, but we do not assume that~$\g$ is simple or semi-simple.

\subsubsection{Traditional Methods  }
\label{traditional_methods}

The classic approach to classifying the real Lie algebra~$g$ is to utilize its Killing form~$B$ and construct the Dynkin diagram for its associated complex Lie algebra~$g^\mathbb{C}$.  As mentioned in Section~\ref{ch:Lie_groups_Lie_algebras.Basic_Definitions.regular_rep_Killing_form}, if the Killing form is non-degenerate, then~$g$ is semi-simple.  According to~\cite{dyn2}, if~$g$ is semi-simple, then it is possible to construct the Dynkin diagram associated with~$g^\mathbb{C}$.  This diagram will be connected if and only if~$g^\mathbb{C}$ is simple~\cite{dyn2}.  In this case, the Dynkin diagram uniquely identifies the complex Lie algebra~$g^\mathbb{C}$ as one of the classical Lie algebras~$A_l, B_l, C_l, D_l$ for~$l \ge 1$ or as one of the exceptional Lie algebras~$G_2, F_4, E_6, E_7$, or~$E_8$.  Conversely, if~$g^\mathbb{C}$ is semi-simple but the Dynkin diagram is not connected, then~$g^\mathbb{C}$ is a direct sum of simple complex Lie algebras which can be identified from the resulting diagram.  Further information regarding these methods may be found in~\cite{cornwell, gilmore}.

We also note that if~$g$ is a simple Lie algebra, then the Killing form provides the signature of~$g$.  As stated in~\cite{gilmore}, except in some cases involving real forms of~$A_l$ and~$D_l$ for various values of~$l$, if~$g$ is not compact, then the signature, rank, and dimension of~$g$ identifies~$g$ as a particular real form of~$g^\mathbb{C}$.  We do not encounter any of the exceptions to this rule, due to the rank and signature of these algebras.  Tables of the various real forms of the complex Lie algebras may be found in~\cite{gilmore}.

For the algebra~$\g$ in Example 1, if we know that~$\g$ is simple, then we may classify it as either a real form of~$B_4$ or~$C_4$.  In this case, the signature of~$\g$ will usually help us identify the particular real form when we use the tables listed in~\cite{gilmore} by matching the rank, dimension, and signature.  This will fail, however, if~$\g$ is compact.  If it is not known that~$\g$ is simple, we note that it is a tedious process to calculate the Killing form of~$\g$ and to construct the Dynkin diagram of~$\g^\mathbb{C}$ for such a large algebra.  Indeed, while the traditional methods described above are effective, we prefer to find other methods which may be more helpful in identifying subalgebras of~$sl(3,\mathbb{O})$ with less effort.

\subsubsection{Counting Methods for Complex Algebras}

According to~\cite{gilmore}, the complex simple Lie algebra~$g^\mathbb{C}$ may be identified by its dimension~$n$ and rank~$l$ in all but a few cases.  The complex Lie algebra~$A_l$ has dimension~$n = l^2-1$ for~$l \ge 1$ and the dimension of~$D_l$ is~$n = l(2l-1)$ for~$l \ge 4$.  The exceptional Lie algebras~$G_2$,~$F_4$,~$E_6$,~$E_7$, and~$E_8$ have dimensions~$14, 52, 78, 133$, and~$248$.  The subscript indicates the rank of each algebra.  These algebras have different dimensions and rank.  The main exception is that both~$B_l$ and~$C_l$ have dimension~$l(2l+1)$, but that for~$l \ge 3$, these algebras are not isomorphic.  We also note that~$B_6$ and~$C_6$ have the same dimension as~$E_6$.

This counting argument can help us classify the complexified algebra of a simple real Lie algebra~$g$.  As noted in Section~\ref{traditional_methods}, we may then use the signature of~$g$ and tables found in~\cite{gilmore} to identify the particular real form of~$g$.  However, this method fails if it is not known that~$g$ is semi-simple, as is the case with the algebra~$\g$ in Example 1.

For the algebra~$\g$ in Example 1, if we know that~$\g$ is simple but do not want to determine its Killing form~$B$, then there is still a method that may identify~$\g$ as either a real form of~$B_4$ or~$C_4$.  As shown in Figure~\ref{fig:subalgebras_of_E6} in Section~\ref{ch:Joma_Paper.subalgebras_greater_than_3}, the Lie algebra~$B_4$ contains a subalgebra~$D_3$ of dimension~$28$, while the largest simple subalgebra contained within~$C_4$ is~$C_3$, which has dimension~$21$.  Hence, if we are able to find a simple subalgebra of dimension~$28$ within~$\g$, then we know that~$\g$ must be a real form of~$B_4$.  We note that the failure to find such a subalgebra does not mean~$\g$ must be a real form of~$C_4$, as not every real form of a complex algebra~$g^\mathbb{C}$ will contain real forms of the complex subalgebras of~$g^\mathbb{C}$.  

\subsubsection{Involutive Automorphisms and Compact Subalgebras}

As shown in~\cite{gilmore} and listed in its tables of real forms of complex Lie algebras, each real form~$g$ and~$g^\prime \ne g$ of~$g^\mathbb{C}$ contains a maximal compact subalgebra~$g_1 \oplus g_2 \oplus \cdots \oplus g_k \subset g$ and~$g^\prime_1 \oplus \cdots \oplus g^\prime_j \subset g^\prime$, respectively.  According to~\cite{cornwell}, there is an involutive automorphism~$\phi^* : g^\mathbb{C} \to g^\mathbb{C}$ such that~$\phi^*(g) = g^\prime$ and~$\phi^*(g) \cap g = g_1 \oplus g_2 \oplus \cdots \oplus g_k$ is the maximal compact subalgebra of~$g$.  When~$g \subset sl(3,\mathbb{O})$, this method merely identifies the maximal compact subalgebra of~$g$.  However, the pre-image of the maximal compact subalgebra of~$g^\prime = \phi^*(g)$ is also a subalgebra of~$g$ and the involutive automorphism may be used to determine the signature of this subalgebra in~$g$.  Again,~\cite{gilmore} lists the semi-simple decomposition of this non-compact subalgebra in~$g$.

Given a subalgebra~$g \subset sl(3,\mathbb{O})$, assume we know~$g$ should have the semi-simple decomposition~$g = g_1 \oplus g_2 \oplus \cdots \oplus g_k$ where~$|g_1| \le |g_2| \le \cdots \le |g_k|$ and~$|g_i|$ is the dimension of~$g_i$.  By identifying the smallest subalgebras~$g_1, \cdots, g_{k-1}$ in~$g$ which commute with each other and also which commute with~$g_k = g - g_1 - \cdots - g_{k-1}$, we automatically show that~$g_k$ is simple.  This method is most effective when~$k = 2$ or~$k=3$ and the dimension of~$g_k$ is very large compared to~$g_1, \cdots, g_{k-1}$.  Effectively, by finding some very small simple Lie algebras in~$g$, this method allows us to find a very large simple subalgebra of~$g$ with very little work.

We now show how this method may be used to show an algebra is simple by applying it the algebra~$\g$ given in Example 1.  By assumption,~$\g$ is a~$36$-dimensional subalgebra~$\g$ of~$sl(3,\mathbb{O})$.  First, we note that there is a real form~$g^\prime$ of~$E_6$ with signature~$(36,42)$, and that~\cite{gilmore} identifies its compact maximal subalgebra, denoted~$\p$, as~$C_4$.  There are involutive automorphisms~$\phi^*: E_6 \to E_6$ such that~$\phi(sl(3,\mathbb{O})) = g^\prime$.  We search for one such that the inverse image of~$\p$ is~\hbox{$\g \in sl(3,\mathbb{O})$}.  If this succeeds, we have identified~$\g$ as a real form of~$C_4$.  If it fails, however, all is not lost - we have still shown that the pre-image of~$\p$ in~$sl(3,\mathbb{O})$ is both simple and a real form of~$C_4$.

\subsubsection{Counting for Matrix Groups}

We now show how to classify the simple real classical Lie algebras, which are $su(p,q,\F)$ and~$sl(n,\F)$ for~$\F = \mathbb{R}, \mathbb{C}, \mathbb{H},$ or~$\mathbb{O}$.  While the signature of a simple algebra~$g$ was used in the past to help identify a maximal compact subalgebra, we will now see that its particular real form is identified by the algebra's signature, rank, and dimension.

For~$\F = \mathbb{R}, \mathbb{C}, \mathbb{H}$, or~$\mathbb{O}$, and for~$n \ge 1$ (with~$n \le 3$ when~$\F = \mathbb{O}$), let~$|\SO(\textrm{im}(\F))|$ denote the dimension of the special orthogonal group acting on the imaginary basis units of the division algebra~$\F$.  Hence, we have
$$
\begin{array}{ccc}
|\SO(\textrm{im}(\mathbb{R}))| = |\SO(0)| = 0 & \hspace{1cm} & |\SO(\textrm{im}(\mathbb{C}))| = |\SO(1)| = 0 \\
|\SO(\textrm{im}(\mathbb{H}))| = |\SO(3)| = 3 & &  |\SO(\textrm{im}(\mathbb{O}))| = |\SO(7)| = 21
\end{array}
$$

\noindent We note that~$Aut(\mathbb{H}) = \SO(\textrm{im}(\mathbb{H}))$, as they are both~$\SO(3)$, and that 
$$G_2 = Aut(\mathbb{O}) \subset \SO(\textrm{im}(\mathbb{O}))$$

  We now derive a formula for the dimension of~$SU(n,\F), \F \ne \mathbb{O}$ by counting the dimension of its associated Lie algebra~$su(n,\F)$.  Let~$|\F|$ denote the dimension of the division algebra~$\F$.  
  An~$n \times n$  matrix~$m \in su(n,\F)$ satisfies~$m^\dagger = -m$ and~$\tr{m} = 0$.  For the diagonal entries of~$m$, the condition~$m^\dagger = -m$ guarantees the diagonal entries are pure imaginary.  Combined with the condition~$\tr{m} = 0$, there are~$(|\F|-1)(n-1)$ real degrees of freedom on the main diagonal.  The condition~$m^\dagger = -m$ eliminates half of the~$|\F|(n(n-1))$ real degrees of freedom for the off-diagonal entries of~$m$.  Finally, we add~$|\SO(\textrm{im}(\F))|$ degrees of freedom for the transformations which mix up the imaginary units of the division algebra~$\F$.  Thus, the Lie group~$SU(n,\F)$ has dimension
$$|SU(n,\F)| = |su(n,\F)| = (|\F|-1)(n-1) + |\F|\frac{n(n-1)}{2} + |\SO(\textrm{im}(\F))|$$
which reduces to the following dimensions for the division algebras over~$\F = \mathbb{R}, \mathbb{C}$, and~$\mathbb{H}$:

$$
\begin{array}{|c|c|c|}
\hline
\F & SU(n,\F) & \textrm{Dimension} \\
\hline
\mathbb{R} & SO(n,\mathbb{R}) & \frac{n(n-1)}{2} \\
\mathbb{C} & SU(n,\mathbb{C}) & (n+1)(n-1) = n^2 - 1\\
\mathbb{H} & SU(n,\mathbb{H}) = SP(2\cdot n) & n(2n+1) \\
\hline
\end{array}
$$

Traditionally, the Lie algebra~$sp(n)$ consists of~$2n \times 2n$ complex matrices constructed in a way that resembles the representation of quaternions via~$2 \times 2$ complex matrices.  We do not use this convention, and instead use a representation which truly uses quaternions.  With our convention, the matrices in~$su(n,\mathbb{H}) = sp(n)$ are~$n \times n$ matrices.

In the case~$\F = \mathbb{O}$, the above formula produces~$|su(1,\mathbb{O})| = 21$ and~$|su(2,\mathbb{O})| = 36$, matching the dimensions of~$so(7)$ and~$so(9)$, respectively.  
For the case~$n = 3$, we must modify the formula since seven of the~$\SO(\textrm{im}(\F))$ transformations (which are not in~$G_2$) can be generated from the type~$2$ and type~$3$ off-diagonal generators.  This over-counting is a result of triality, since the three copies of~$so(7)$ (each of which contains~$G_2$) in~$so(8)$ are not independent (see Section~\ref{ch:E6_basic_structure.Triality}).  Thus, we reduce the dimension of~$su(3,\mathbb{O})$ from~$59$, as given by the formula, to~$52$, which is the dimension of~$F_4$.

$$
\begin{array}{|c|c|c|}
\hline
\F = \mathbb{O} & SU(n,\F) & \textrm{Dimension} \\
\hline
 & SU(1,\mathbb{O}) & 21 \\
 & SU(2,\mathbb{O}) & 36 \\
 & SU(3,\mathbb{O}) & 52 \\
\hline
\end{array}
$$

To produce a formula for~$SL(n,\mathbb{F})$, we again count the dimension of the associated Lie algebra~$sl(n,\mathbb{F})$.  The~$n \times n$ matrix~$m \in sl(n,\mathbb{F})$ satisfies~$m^\dagger = \pm m$ and~$\tr{m} = 0$.  When compared with the conditions for~$m \in su(n,\mathbb{F})$, there are~$n-1$ additional degrees of freedom for the diagonal entries and~$|\F|\frac{n(n-1)}{2}$ additional degrees of freedom for the off-diagonal entries of the~$sl(n,\F)$ matrices than existed for the~$su(n,\F)$ matrices.
  This produces the dimension formula

$$|sl(n,\F)| = |su(n,\F)| + (n-1) + |\F|\frac{n(n-1)}{2}$$

\noindent This formula gives the following dimensions of~$SL(n,\F)$ for~$\F = \mathbb{R}, \mathbb{C}$, and~$\mathbb{H}$:
$$
\begin{array}{|c|c|c|}
\hline
\F & SL(n,\F) & \textrm{Dimension} \\
\hline
\mathbb{R} & SL(n,\mathbb{R}) & (n+1)(n-1) \\
\mathbb{C} & SL(n,\mathbb{C}) & 2(n-1)(n+1) \\
\mathbb{H} & SL(n,\mathbb{H}) & 4n^2 - 1 \\
\hline
\end{array}
$$

In the case~$\F = \mathbb{O}$, when~$n = 1$, we find the dimension of~$sl(1,\mathbb{O})$ is~$21$.  In fact, for all~$\F$, we have~$|sl(1,\F)| = |su(1,\F)|$, 
since~$sl(n,\F)$ adds Hermitian trace-free matrices to the anti-Hermitian trace-free matrices of~$su(n,\F)$, and in the case~$n =1$, no such matrices may be added to~$su(1,\F)$.
  For~$\F = \mathbb{O}$ and~$n = 2, 3$, the formula gives the dimensions of~$sl(2,\mathbb{O})$ and~$sl(3,\mathbb{O})$ as~$45$ and~$78$, which are the dimensions of~$so(9,1)$ and~$E_6$, respectively.

$$
\begin{array}{|c|c|c|}
\hline
\F = \mathbb{O} & SU(n,\F) & \textrm{Dimension} \\
\hline
 & SL(1,\mathbb{O}) & 21 \\
 & SL(2,\mathbb{O}) & 45 \\
 & SL(3,\mathbb{O}) & 78 \\
\hline
\end{array}
$$

Finally, we note that~$su(p,q,\F)$ is another real form of the real algebra~$su(n,\F)$ when~$p+q = n$.  Thus, both~$su(p,q,\F)$ and~$su(n,\F)$ have the same dimension and rank.  The algebra~$su(n,\F)$ is compact, containing no boosts, while the algebra~$su(p,q,\F)$ contains~$|\F|pq$ boosts.  The algebra~$sl(p,q,\F)$ is~$sl(n,\F)$, since~$sl(n,\F)$ already contains the anti-hermitian and hermitian matrices which correspond to the compact and non-compact generators needed for~$su(p,q,\F)$.  Nevertheless, we will use~$sl(p,q,\F)$ in certain instances to differentiate the algebra from another (preferred) algebra~$sl(n,\F)$.  

We conclude that in all but a few cases (which are not relevant for subalgebras of~$sl(3,\mathbb{O})$), knowing the rank, dimension, and signature of a finite-dimensional simple Lie algebra~$g$ is enough information to identify~$g$ as~$su(p,q,\F)$ or as~$sl(n,\F)$, where~$\F = \mathbb{R}, \mathbb{C}, \mathbb{H}, \mathbb{O}$ and~$n = p + q$.  The exceptions to this rule are listed in~\cite{gilmore}.

This conclusion is very useful for identifying the particular real form of a simple Lie algebra contained within~$sl(3,\mathbb{O})$.  Consider the~$36$-dimensional algebra~$\g$ in Example 1, which was introduced as the start of Section~\ref{ch:E6_further_structure.Automorphisms_of_Lie_Algebras.id_algebras}.  Assuming that we have shown~$\g$ is simple, then we know~$g$ is either a real form of~$B_4$ or~$C_4$ since it has dimension~$36$ and rank~$4$.  If~$\g$ contains~$8$ boosts, then it must be~$so(8,1,\mathbb{R})$, while if it contains 12 boosts, then it must be~$su(3,1,\mathbb{H})$.  If~$\g$ is compact, we may use an involutive automorphism~$\phi^* : \g^* \to \g^*$ sending~$\g$ to another real form of its complexification~$\g^*$.  Now, depending on the number of boosts in ~$\phi^*(\g)$, we determine that~$g^*$ is either~$B_4$, in which case~$\g = so(9,\mathbb{R})$, or~$C_4$, in which case~$\g = sp(4)$.

\subsubsection{Summary}

We summarize here the methods which allow us to quickly find and identify subalgebras of~$sl(3,\mathbb{O})$.  We use an involutive automorphism~$\phi^*: g^\mathbb{C} \to g^\mathbb{C}$ not only to identify different real forms of~$g^\mathbb{C}$, but also to identify real subalgebras of our particular real form~$g$.  When applied to the compact real form~$g = \p \oplus \n$, the map~$\phi^*$ changes the signature from~$(|\p|+|\n|, 0)$ to~$(|\p|,|\n|)$.  Similar counting arguments give the signature of~$\phi^*(g)$ when~$g$ is non-compact, and we note that~$\phi^*$ changes some compact generators into non-compact generators, and vice-versa.  We then use the rank, dimension, and signature of~$\phi^*(g)$ to identify the particular real form of the algebra.  Using tables, found in~\cite{gilmore}, of real forms of~$g^\mathbb{C}$ that list the algebra's maximal compact subalgebra, we then identify~$\p = g \cap \phi^*(g)$ as a subalgebra of our original real form~$g$.  In addition, it is very easy to identify the maximal compact subalgebra~$g_{(c)}$ of~$\phi^*(g)$.  We then identify the pre-image of~$g_{(c)}$ as a non-compact subalgebra of~$g$.  The involutive automorphisms~$\phi$, and their complex extensions~$\phi^*$, provide significant insight into the structure of our real form~$g = sl(3,\mathbb{O})$ of the algebra~$E_6$.

\section{Three Important Involutive Automorphisms of~$sl(3,\mathbb{O})$}
\label{ch:E6_further_structure.three_automorphisms}

We first make some comments about our preferred basis for~$sl(3,\mathbb{O})$, which is listed in Table~\ref{table:our_basis}.  Let~$B$ and~$R$ be the vector subspaces consisting of boosts and rotations, respectively.  Our preferred choice of basis favors type~$1$ transformations, in the sense that we choose to use the linear dependencies given in Section~\ref{ch:E6_basic_structure.ConstructingE6algebra.LinearDependencies} for~$so(8)$
to express~$G_2$,~$so(7)$, and~$so(8)$ transformations in terms of  type~$1$ transformations.
  Let~$T_1$ be the subspace spanned by~$\lbrace \dotBtz \rbrace~$ and all type~$1$ transformations, and~$T_2$ and~$T_3$ the subspaces spanned by the type~$2$ and type~$3$ transformations which are not in~$T_1$.\footnote{With the condition~$\dot B^{(1)}_{tz} + \dot B^{(2)}_{tz} + \dot B^{(3)}_{tz} = 0$, we note that both~$\dot B^{(2)}_{tz}$ and~$\dot B^{(3)}_{tz}$ are viewed as elements of~$T_1$ since~$T_1$ contains~$\dot B^{(1)}_{tz}$ and~$\dotBtz$.}
  Let~$T_{(2,3)} = T_2 + T_3$.  Finally, let~$H$ be the subspace spanned by those transformations with no labels from~$\lbrace i,j,j\l,i\l \rbrace$ and~$H^\perp$ the subspace spanned by those transformations labeled with one label from~$\lbrace i, j, j\l, i\l  \rbrace$.  This subspace~$H$ contains {\it quaternionic} basis elements (using the quaternionic algebra~$\langle 1, k, k\l, \l \rangle$) such as~$\dot A_{k\l}, \dot R^{(1)}_{z\l}$, and~$\dot B^{(3)}_{tx}$, while~$H^\perp$ contains {\it orthogonal-quaternionic} basis elements such as~$\dot A_i, \dot R^{(1)}_{zi\l},$ and~$\dot B^{(3)}_{tj\l}$.

The results from Section~\ref{ch:E6_further_structure.Automorphisms_of_Lie_Algebras} show that as a vector space,~$sl(3,\mathbb{O})$ splits into the pairs of orthogonal spaces~$R \oplus B$,~$T_1 \oplus T_{(2,3)}$, and~$H \oplus H^\perp$.  Indeed, the pairs of orthogonal spaces~$R \oplus B$,~$T_1 \oplus T_{(2,3)}$, and~$H \oplus H^\perp$ satisfy the commutation relations:
$$
\begin{array}{ccc}
\left[R, R\right] \subset R  & \left[T_1, T_1 \right] \subset T_1             & \left[H, H\right] \subset H \\
\left[R, B\right] \subset B  & \hspace{1cm} \left[T_1, T_{(2,3)} \right] \subset T_{(2,3)} \hspace{1cm} &  \left[H, H^\perp\right] \subset H^\perp \\
\left[B, B\right] \subset R  & \left[T_{(2,3)}, T_{(2,3)} \right] \subset T_1 &  \left[H^\perp, H^\perp\right] \subset H \\
\end{array}
$$

We define the three maps~$\phi_{(t)}$,~$\phi_{(2,3)}$, and~$\phi_{H}$ from~$sl(3,\mathbb{O})$ to~$sl(3,\mathbb{O})$, according to the rules
$$
\begin{array}{ccc}
\phi_{(t)}(r + b) = r + \xi b & \hspace{1cm} & r \in R, b \in B \\
\phi_{(2,3)}(g_1 + g_{(2,3)}) = g_1 + \xi g_{(2,3)} & & g_1 \in T_1, g_{(2,3)} \in T_{(2,3)} \\
\phi_{(H^\perp)}(h + h^\prime) = h + \xi h^\prime & & h \in H, h^\prime \in H^\perp\\
\end{array}
$$
where~$\xi = -1$.  These maps respect the algebraic structure of~$sl(3,\mathbb{O})$, and each clearly satisfies~$\phi^2 = Id$.  Thus,~$\phi_{(t)}$,~$\phi_{(2,3)}$, and~$\phi_{(H^\perp)}$ are involutive automorphisms, inducing the structure~$sl(3,\mathbb{O}) = \p \oplus \n$ for~$\p \oplus \n = R \oplus B$,~$\p \oplus \n = T_1 \oplus T_{(2,3)}$, and~$\p \oplus \n = H \oplus H^\perp$, respectively.

{\bf Lemma:} The signatures of~$R$,~$T_1$, and~$H$ are~$(52,0)$,~$(36,10)$, and~$(24,14)$.  

{\bf Proof:}  The algebra~$R$ consists of the compact generators of~$sl(3,\mathbb{O})$, of which there are~$52$.  The algebra~$T_1$ contains~$so(9,1)$, which has signature~$(36,9)$, and the subalgebra~$\langle \dotBtz \rangle$.  The signature of~$H$ is found by counting the number of compact and non-compact generators in its basis.

\myendofproof

We extend our basis of~$sl(3,\mathbb{O})$ to a basis of the complex Lie algebra~$E_6$ by considering our basis now as a vector space over the complex numbers.  Writing~$\sqrt{-1}$ for~$\xi$, we extend our involutive automorphisms~$\phi: sl(3,\mathbb{O}) \to sl(3,\mathbb{O})$ to \hbox{$\phi^*: E_6 \to E_6$}, giving the maps
$$ \phi^*_{(t)}(R + B) = R + \sqrt{-1}B \hspace{1cm} \phi^*_{(2,3)}(T_1 + T_{(2,3)}) = T_1 + \sqrt{-1}T_{(2,3)}$$
$$\phi^*_{(H^\perp)}(H + H^\perp) = H + \sqrt{-1}H^\perp$$

These involutive automorphisms change some compact and non-compact generators of~$sl(3,\mathbb{O})$ into non-compact and compact generators, respectively, of~$\phi^*\left( sl(3,\mathbb{O})\right)$.  For instance, each boost~$b$ satisfies~$B(b, b) > 0$, where~$B( \hspace{.22cm} , \hspace{.22cm} )$ is the Killing form.  Under~$\phi^*_{(t)}$, the boosts~$b$ satisfy 
$$B\left(\phi^*_{(t)}(b), \phi^*_{(t)}(b)\right) = B(\xi b, \xi b) = \xi^2B(b,b) < 0$$
Hence,~$\phi^*_{(t)}(b)$ is a compact generator of~$\phi^*_{(t)}\left( sl(3,\mathbb{O}) \right)$.
  The involutive automorphism~$\phi^*_{(2,3)}$  changes the type~$2$ and type~$3$ rotations to non-compact generators, and the boosts in~$T_{(2,3)}$ are changed to compact generators.  Similarly, all orthogonal-quaternionic rotations and boosts in~$H^\perp$ change signature.  

The involutive automorphism~$\phi^*_{(t)}$ transforms~$sl(3,\mathbb{O})$, which has signature~$(52,26)$, into the compact real form~$\phi^*_{(t)}\left(sl(3,\mathbb{O})\right)$, which has signature~$(78,0)$.  The subalgebra
$$\p = \phi^*_{(t)}\left(sl(3,\mathbb{O})\right) \cap sl(3,\mathbb{O})$$
has dimension~$52$ and is compact.  Hence, according to the tables in~\cite{gilmore}, we see that~$\p$ is the compact real form~$su(3,\mathbb{O})$ of~$F_4$.  Of course, the compact part of~$\phi^*_{(t)}\left(sl(3,\mathbb{O})\right)$ is the entire algebra, so that in this case, its pre-image is already known.

We obtain two interesting results when we apply the involutive automorphism~$\phi^*_{(2,3)}$ to~$sl(3,\mathbb{O})$.  First, the signature of~$g^\prime = \phi^*_{(2,3)}\left(sl(3,\mathbb{O})\right)$ is~$(52,26)$, giving another real form of~$E_6$ whose maximal compact subalgebra~$g^\prime_{(c)}$ is~$F_4$.  Hence, the pre-image of~$g^\prime_{(c)}$ is a real form of~$F_4$ in~$sl(3,\mathbb{O})$.  It has signature~$(36,16)$, and the~$16$ non-compact generators identify this real form of~$F_4$ as~$su(2,1,\mathbb{O})$.  The second result comes from looking at 
$$ \p = \phi^*_{(2,3)}\left(sl(3,\mathbb{O})\right) \cap sl(3,\mathbb{O})$$
Because~$|\p| = 46$, we see from the table in~\cite{gilmore}
 that~$\p$ is a real form of~$D_5 \oplus D_1$.  Of course,~$\phi^*_{(2,3)}$ fixes the algebra~$T_1$, which has signature~$(36,9)\oplus(0,1)$.  Hence, we identify~$T_1$ as the subalgebra~$so(9,1) \oplus \langle \dot B^{(2)}_{tz} - \dot B^{(3)}_{tz} \rangle$.

We also obtain two results by applying~$\phi^*_{(H^\perp)}$ to~$sl(3,\mathbb{O})$.  We find that the signature of~$g^\prime = \phi^*_{(H^\perp)}\left(sl(3,\mathbb{O})\right)$ is~$(36,42)$.  The maximal compact subalgebra~$g^\prime_{(c)}$ of~$g^\prime$ is~$su(4,\mathbb{H})$, a real form of~$C_4$.  The pre-image of~$g^\prime_{(c)}$ has signature~$(24,12)$, and due to the~$12$ non-compact generators, must be identified as~$su(3,1,\mathbb{H})_1$.  The invariant subalgebra 
$$ \p = \phi^*_{(H^\perp)}\left(sl(3,\mathbb{O})\right) \cap sl(3,\mathbb{O})$$
has dimension~$|\p| = 38$ and signature~$(24,14)$.  According to the table in~\cite{gilmore}
~$\p$ is a real form of~$A_5 \oplus A_1$ with signature~$(21,15) \oplus (3,0)$.  The involutive automorphism~$\phi^*_{(H^\perp)}$ obviously fixes~$H$.  Of the~$24$ rotations unaffected by~$\phi^*_{(H^\perp)}$, there are~$21$ which are quaternionic and form the subalgebra~$su(3,\mathbb{H})$.  The remaining three rotations~${A_k, A_{k\l}, A_{\l}}$ are elements of~$G_2$, but are transformations involving pairs of~$i,j,i\l$, and~$j\l$ only and leave~$k,k\l$, and~$l$ fixed.  The four elements~$i,j,i\l$ and~$j\l$ can be paired into two complex pairs (of which~$i + i\l, j + j\l$ is one choice), and the~$\dot A_k, \dot A_{k\l}, \dot A_\l$ transformations act as~$su(2,\mathbb{C})$ transformations producing the other pairs of complex numbers.  Hence, the~$24$ compact elements form the algebra~$su(3,\mathbb{H})\oplus su(2,\mathbb{C})^C$, where we continue to use superscript~$C$ notation to indicate this~$su(2,\mathbb{C})$ is a subalgebra of~$G_2$ and {\bf not} the 
canonical~$su(2,\mathbb{C})_s$.  We identify~$\p = H$ as the subalgebra~$sl(3,\mathbb{H}) \oplus su(2,\mathbb{C})^C$.

The set of involutive automorphisms form a group, allowing us to consider the compositions of the involutive automorphisms~$\phi^*_{(t)}$,~$\phi^*_{(2,3)}$, and~$\phi^*_{(H^\perp)}$.  After working through one example, we provide the subalgebra structure of~$sl(3,\mathbb{O})$ that can be found using these maps in Table~\ref{table:maximal_fixed_and_compact_subalgebras}.

Consider the involutive automorphism~$\phi^*_{(2,3)} \circ \phi^*_{(H^\perp)} : E_6 \to E_6$.  We define the following subspaces~$H_1$,~$H_{(2,3)}$,~$H^\perp_1$, and~$H^\perp_{(2,3)}$ which are the intersections of pairs of spaces~$T_1, T_{(2,3)}, H$, and~$H^\perp$, as indicated in Table~\ref{table:int_of_T_H}.  For example,~$H_1 = H \cap T_1$ and~$H^\perp_{(2,3)} = H^\perp \cap T_{(2,3)}$.  The second table in Table~\ref{table:int_of_T_H} indicates the number of basis elements which are boosts and rotations (in~$sl(3,\mathbb{O})$) in each of these spaces.  These spaces satisfy the commutation rules given in Table~\ref{table_comm_of_T_H}.

\begin{table}[htbp]
\begin{center}
 \begin{tabular}{c|cccc|cc}
$\cap$  & $T_1$ & $T_{(2,3)}$ & \hspace{2cm} &$\cap$ & $T_1$ & $T_{(2,3)}$\\
\cline{1-3}\cline{5-7}
$H$ & $H_1$ & $H_{(2,3)}$  & & $H$ & $(16,6)$ & $(8,8)$ \\
$H^\perp$ & $H^\perp_1$ & $H^\perp_{(2,3)}$ & & $H^\perp$ & $(20,4)$ & $(8,8)$\\
\multicolumn{3}{c}{Subspace} & & \multicolumn{3}{c}{Signature}
\end{tabular}
\caption{Intersections of Subspaces~$T_1$,~$T_{(2,3)}$,~$H$, and~$H^\perp$}
\label{table:int_of_T_H}
\end{center}
\end{table}

\begin{table}[htbp]
\begin{center}
\begin{tabular}{|c|c|c|c|c|}
\hline
$\left[ \hspace{.15cm} , \hspace{.15cm} \right]$& $H_1$             & $H_{(2,3)}$      & $H^\perp_1$       & $H^\perp_{(2,3)}$ \\
\hline
$H_1$             & $H_1$             & $H_{(2,3)}$      & $H^\perp_1$       & $H^\perp_{(2,3)}$ \\
\hline
$H_{(2,3)}$       & $H_{(2,3)}$       & $H_1$            & $H^\perp_{(2,3)}$ & $H^\perp_1$ \\ 
\hline
$H^\perp_1 $      & $H^\perp_1$       & $H^\perp_{(2,3)}$ & $H_1$             & $H_{(2,3)}$\\
\hline
$H^\perp_{(2,3)}$ & $H^\perp_{(2,3)}$ & $H^\perp_1$      & $H_{(2,3)}$       & $H_1$ \\
\hline
\end{tabular}
\caption{Commutation structure of~$H_1, H_{(2,3)}, H^\perp_1,$ and~$H^\perp_{(2,3)}$}
\label{table_comm_of_T_H}.
\end{center}
\end{table}

We see from Table~\ref{table_comm_of_T_H} that~$H_1$,~$H_1 + H^\perp_{(2,3)}$,~$H_1 + H_{(2,3)}$, and~$H_1 + H^\perp_1$ are subalgebras of~$sl(3,\mathbb{O})$
The involutive automorphism~$\phi^*_{(2,3)} \circ \phi^*_{(H^\perp)}$ fixes 
$$
\begin{array}{ccl}
\p & = & \phi^*_{(2,3)} \circ \phi^*_{(H^\perp)}\left(sl(3,\mathbb{O})\right) \cap sl(3,\mathbb{O})\\
   & = & H_1 + H^\perp_{(2,3)} \\
\end{array}
$$
as a subalgebra, since everything in~$H_1$ is fixed by both automorphisms, while everything in~$H^\perp_{(2,3)}$ is multiplied by~$\xi^2 = -1$.    The subalgebra~$\p$ has dimension~$38$ and signature~$(24,14)$, and is thus a real form of~$A_5 \oplus A_1$.  It has a fundamentally different basis from~$sl(3,\mathbb{H})\oplus su(2,\mathbb{C})^C$, however, as it uses a mixture of the quaternionic transformations of type~$1$ with the orthogonal-quaternions transformations of type~$2$ and type~$3$, while~$sl(3,\mathbb{H})$ is comprised of only the quaternionic transformations.  We call this algebra~$sl(2,1,\mathbb{H}) \oplus su(2,\mathbb{C})_2$, where the notation~$sl(2,1,\mathbb{H})$ is used to indicate that this is {\bf not} our standard~$sl(3,\mathbb{H})$ and 
$$su(2,\mathbb{C})_2 = \langle \dot G_k + 2\dot S^1_k, \dot G_{k\l} + 2\dot S^1_{k\l}, \dot G_\l + 2\dot S^1_\l, \rangle$$
 again corresponds to permutations of ~$\lbrace i,j,j\l,i\l \rbrace$, fixes~$\lbrace k, k\l, \l \rbrace$, but is not in~$G_2$.  We also note that the compact subalgebra~$g^\prime_{(c)}$ of~$\phi^*_{(2,3)} \circ \phi^*_{(H^\perp)} \left( sl(3,\mathbb{O})\right)$ has dimension~$36$, and its preimage in~$sl(3,\mathbb{O})$ has signature~$(24,12)$.  We identify this as~$su(3,1,\mathbb{H})_2$, and include the subscript~$2$ since this real algebra has a different basis than our previously identified~$su(3,1,\mathbb{H})$, which we had called~$su(3,1,\mathbb{H})_1$.  

We summarize in Table~\ref{table:maximal_fixed_and_compact_subalgebras} the subalgebras achieved from our three involutive automorphisms as well as compositions of these automorphisms.  For each involutive automorphism~$\phi^*$, the second column lists the signature of~$\phi^*(sl(3,\mathbb{O})$.  The third column identifies the fixed subalgebra~$\p = \phi^*\left(sl(3,\mathbb{O})\right) \cap sl(3,\mathbb{O})$ as well as its signature.  In the fourth column, we identify the pre-image of the maximal compact subalgebra of~$g^\prime = \phi^*\left(sl(3,\mathbb{O})\right)$ and lists both the algebra and its signature.

\begin{table}[htbp]
\begin{center}
\begin{tabular}{|cccc|}
\hline
Involutive          &  Signature of                                    & Signature of & Signature of            \\
Automorphism        & $g^\prime = \phi^*\left(sl(3,\mathbb{O})\right)$ & $\p = g^\prime \cap sl(3,\mathbb{O})$         & $(\phi^*)^{-1}(g^\prime_{(c)})$ \\
\hline
$1$                 & $(52,26)$                                        & $(52,26)$    & $(52,0)$                \\   &   &  $sl(3,\mathbb{O})$ & $su(3,\mathbb{O})$ \\
\hline
$\phi^*_{(t)}$      & $(78,0)$                                         & $(52,0)$     & $(52,26)$               \\ & & $su(3,\mathbb{O})$ & $sl(3,\mathbb{O})$ \\
\hline
$\phi^*_{(2,3)}$    & $(52,26)$                                        & $(36,10)$    & $(36,16)$               \\ & & $sl(2,\mathbb{O}) \oplus u(1)$ & $su(2,1,\mathbb{O})$ \\
\hline
$\phi^*_{(H^\perp)}$& $(36,42)$                                        & $(24,14)$    & $(24,12)$               \\ & & $sl(3,\mathbb{H})\oplus su(2,\mathbb{C})^C$ & $su(3,1,\mathbb{H})_1$\\
\hline
$\phi^*_{(2,3)}\circ \phi^*_{(t)}$ & $(46,32)$                         & $(36,16)$             & $(36,10)$\\ & & $su(2,1,\mathbb{O})$ & $so(9,1) \oplus u(1)$\\
\hline
$\phi^*_{(H^\perp)}\circ \phi^*_{(t)}$ & $(38,40)$                         & $(24,12)$             & $(24,14)$ \\ & & $su(3,1,\mathbb{H})_1$ & $sl(3,\mathbb{H}) \oplus su(2,\mathbb{C})^C$ \\
\hline
$\phi^*_{(2,3)}\circ \phi^*_{(H^\perp)}$ & $(36,42)$                         & $(24,14)$             & $(24,12)$ \\ & & $sl(2,1,\mathbb{H})\oplus su(2,\mathbb{C})_2$ & $su(3,1,\mathbb{H})_2$ \\
\hline
$\phi^*_{(2,3)}\circ \phi^*_{(H^\perp)}\circ \phi^*_{(t)} $ & $(38,40)$                         & $(24,12)$               & $(24,14)$ \\ & & $su(3,1,\mathbb{H})_2$ & $sl(2,1,\mathbb{H})\oplus su(2,\mathbb{C})_2$ \\
\hline
\end{tabular}
\caption{Involutive automorphisms with maximal subalgebras of~$sl(3,\mathbb{O})$}
\label{table:maximal_fixed_and_compact_subalgebras}
\end{center}
\end{table}

Despite recognizing~$H_1$,~$H_1 + H^\perp_{(2,3)}$,~$H_1 + H_{(2,3)}$, and~$H_1 + H^\perp_1$ as subalgebras of~$sl(3,\mathbb{O})$ using the commutation relations in Table~\ref{table_comm_of_T_H}, we note that~$su(3,1,\mathbb{H})_2$ is not any of these subalgebras.  In the next section, we use another technique involving involutive automorphisms to give a finer refinement of subspaces of~$sl(3,\mathbb{O})$, allowing us to provide a nice basis for~$su(3,1,\mathbb{H})_2$  and other subalgebras of~$sl(3,\mathbb{O})$.

\section{Additional Subalgebras of~$sl(3,\mathbb{O})$ from Automorphisms}
\label{ch:E6_further_structure.more_applications_of_three_automorphisms}

We have already used involutive automorphisms to identify the maximal subalgebras of~$sl(3,\mathbb{O})$.  This was done by separating the algebra into two separate spaces, one of which was left invariant by the automorphism.  In this section, we use the group characteristic of the involutive automorphisms to separate~$sl(3,\mathbb{O})$ into four or more subspaces spanned by either the compact or non-compact generators, with the condition that the involutive automorphism either preserves the entire subspace or it changes the character of all the basis elements in the subspace.  We identify additional subalgebras of~$sl(3,\mathbb{O})$ by taking various combinations of these subspaces.

We continue to use the automorphisms~$\phi^*_{(t)}$,~$\phi^*_{(2,3)}$, and~$\phi^*_{(H^\perp)}$, as well as the subspaces~$R, H, T_1, T_{(2,3)}, H$, and~$H^\perp$ defined in the previous section.  

We first consider the involutive automorphism~$\phi^*_{(H^\perp)} \circ \phi^*_{(t)}: E_6 \to E_6$.  This map fixes the subspaces~$R_H = R \cap H$, consisting of quaternionic rotations, and~$B_{H^\perp} = B \cap H^\perp$, consisting of orthogonal-quaternionic boosts.  Under~$\phi^*_{(H^\perp)} \circ \phi^*_{(t)}$, the two subspaces~$B_{H}$ and~\hbox{$R_{H^\perp} = R \cap H^\perp$} change signature.  These spaces consist of orthogonal-quaternionic rotations and quaternionic boosts, respectively.  The dimensions of these four spaces are displayed in Table~\ref{table:H_T_subspaces}.

\begin{table}[htbp]
\begin{center}
\begin{tabular}{l|c|c|}
\cline{2-3}
         & \multicolumn{2}{c|}{$\phi^*_{(H^\perp)} \circ \phi^*_{(t)}$ Signature} \\
\cline{2-3}
          & $(1,0)$  &  $(0,1)$ \\
\hline
\textrm{Same Signature} &  $|R_H| = 24$    &  $|B_{H^\perp}| = 12$ \\
\hline
\textrm{Change Signature} & $|B_H| = 14$ & $|R_{H^\perp}| = 28$ \\
\hline
\end{tabular}
\caption{Splitting of $E_6$ basis under~$\phi^*_{(H^\perp)} \circ \phi^*_{(t)}$}
\label{table:H_T_subspaces}
\end{center}
\end{table}

We list the signature of these spaces under~$\phi^*_{(H^\perp)} \circ \phi^*_{(t)}$, and can thus identify subalgebras of~$\phi^*_{(H^\perp)} \circ \phi^*_{(t)}\left( sl(3,\mathbb{O}) \right)$.  However, we are primarily interested in the pre-image of these subalgebras in our preferred algebra~$sl(3,\mathbb{O})$.  The signature of these spaces in~$sl(3,\mathbb{O})$ is straight-forward, as the rotations are compact and the boosts are not compact.

We use the subspaces represented in Table~\ref{table:H_T_subspaces} to identify subalgebras of~$sl(3,\mathbb{O})$.  As previously identified, the~$24$ rotations in~$R_H$  fixed by the automorphism form the subalgebra~$su(3,\mathbb{H}) \oplus su(2,\mathbb{H}^\perp)_1$.  The entries in the first column of Table~\ref{table:H_T_subspaces} represent all quaternionic rotations and boosts, and form the subalgebra~$sl(3,\mathbb{H})_1 \oplus su(2,\mathbb{H}^\perp)_1$.  Of course, the entries on the diagonal,~$R_H$ and~$R_{H^\perp}$, form~$F_4$.  We finally consider the entries~$R_H$ and~$B_{H^\perp}$ in the top row of the table.  Two orthogonal-quaternionic boosts commute to a quaternionic rotation, and an orthogonal-quaternionic boost commuted with a quaternionic rotation is again an orthogonal-quaternionic boosts.  Hence, the subspace~$R_H + B_{H^\perp}$ closes under commutation and is a subalgebra with signature~$(24,12)$.  While both~$so(9-n, n, \mathbb{R})$ and~$su(4-n, n, \mathbb{H})$ have dimension~$36$, only~$su(3,1,\mathbb{H})$ has~$12$ boosts.  We realize~$R_H + B_{H^\perp}$ is the previously identified~$su(3,1,\mathbb{H})_1$.  

The algebra~$R_H + B_{H^\perp} = su(3,1,\mathbb{H})_1$ is a real form of the complex Lie algebra~$C_4$.  Since~$C_4$ contains~$C_3$ but not~$B_3$, we note that the~$21$-dimensional subalgebra contained within~$R_H$ is a real form of~$C_3$, not of~$B_3$.  In addition, any simple~$21$-dimensional subalgebra of~$R_H + B_{H^\perp}$ is a real from of~$C_3$.  Eliminating the boosts from~$su(3,1,\mathbb{H})_1$ leaves~$su(3,\mathbb{H}) \oplus su(2,\mathbb{H}^\perp)_1$.

We next consider the involutive automorphism~$\phi^*_{(2,3)} \circ \phi^*_{(t)}: E_6 \to E_6$.  This automorphism separates the basis for~$sl(3,\mathbb{O})$ into the subspaces
$$
\begin{array}{ccc}
R_1 = R \cap T_1 &\hspace{2cm} & B_1 = B \cap T_{(2,3)} \\
B_{(2,3)} = B \cap T_{(2,3)} & \hspace{2cm}&  R_{(2,3)} = R \cap T_{(2,3)}\\
\end{array}
$$
This automorphism leaves the signature of~$R_1$ and~$B_{(2,3)}$ alone, while it reverses the signature of the generators in~$R_{(2,3)}$ and~$B_1$.  This information is represented in Table~\ref{table:T_23_subspaces}

\begin{table}[htbp]
\begin{center}
\begin{tabular}{l|c|c|}
\cline{2-3}
         & \multicolumn{2}{c|}{$\phi^*_{(2,3)} \circ \phi^*_{(t)}$ Signature} \\
\cline{2-3}
          & $(1,0)$  &  $(0,1)$ \\
\hline
\textrm{Same Signature} &  $|R_1| = 36$    &  $|B_{(2,3)}| = 16$ \\
\hline
\textrm{Change Signature} & $|B_1| = 10$ & $|R_{(2,3)}| = 16$ \\
\hline
\end{tabular}
\caption{Splitting of $E_6$ basis under $\phi^*_{(2,3)} \circ \phi^*_{(t)}$}
\label{table:T_23_subspaces}
\end{center}
\end{table}

We again use the subspaces represented in Table~\ref{table:T_23_subspaces} to identify subalgebras of~$sl(3,\mathbb{O})$.  The subspace~$R_1$ is the subalgebra~$so(9,\mathbb{R})$, and contains all subalgebras~$so(n,\mathbb{R})$ for~$n \le 9$.  The subspace~$R_1 + B_1$ is, of course,~$so(9,1,\mathbb{R}) \oplus u(1)$, where~$u(1) = \langle \dot B^{(2)}_{tz} - \dot B^{(3)}_{tz} \rangle$.  Again, the complete set of rotations,~$R_1 + R_{(2,3)}$ is the subalgebra~$su(3,\mathbb{O})$, which is a real form of~$F_4$.  Interestingly, the subspace~$R_1 + B_{(2,3)}$ in the top row is another form of~$F_4$.  The~$16$ boosts in~$R_1 + B_{(2,3)}$ identify this version of~$F_4$ as~$su(2,1,\mathbb{O})$.

We finally consider the involutive automorphism~$\phi^*_{(2,3)} \circ \phi^*_{(H^\perp)} : E_6 \to E_6$ as the composition of the maps~$\phi^*_{(2,3)} \circ \phi^*_{(t)}$ and~$\phi^*_{(H^\perp)} \circ \phi^*_{(t)}$.  A pair of maps normally separates~$E_6$ into four different subspaces, but we can create a finer refinement of subspaces of~$sl(3,\mathbb{O})$ by using the two copies of~$\phi^*{(t)}$ to separate the boosts and rotations from each other.  These subspaces are represented in Table~\ref{table:23_H_subspaces}.  We continue with our previous conventions for designating intersections of subspaces - that is,~$R_{(2,3),H^\perp} = R \cap T_{(2,3)} \cap H^\perp$.  We use a~$4\times4$ table, rather than a~$4\times 2$ table, to help us identify relevant subalgebras.

\begin{table}[htbp]
\begin{center}
\begin{tabular}{l|cc|cc|}
\cline{2-5}
  & \multicolumn{4}{c|}{$\phi^*_{(2,3)} \circ \phi^*_{(H^\perp)}$ Signature} \\
\cline{2-5}
  & $(1,0)$  &  $(1,0)$ & $(0,1)$ & $(0,1)$ \\
\hline
$\phi^*_{(2,3)} = Id, \phi^*_{(H^\perp)} = Id$     & $|R_{1,H}|$ &         &  $|B_{1,H}|$ &   \\
                                                   & $= 16$         &         &  $=6$          &   \\
$\phi^*_{(2,3)} \ne Id, \phi^*_{(H^\perp)} \ne Id$ &          & $|R_{(2,3),H^\perp}|$ &          &  $|B_{(2,3),H^\perp)}|$ \\
                                                   &          & $= 8$                  &         &  $=8$          \\
\hline
$\phi^*_{(2,3)} \ne Id, \phi^*_{(H^\perp)} = Id$   & $|B_{(2,3),H}|$  &         & $|R_{(2,3),H}|$  &          \\
                                                   & $= 8$         &         &  $=8$          &   \\
$\phi^*_{(2,3)} = Id, \phi^*_{(H^\perp)} \ne Id$   &          & $|B_{1,H^\perp}|$ &          & $|R_{1,H^\perp}|$ \\ 
                                                   &          & $= 4$                  &         &  $=20$          \\
\hline
\end{tabular}
\caption{Splitting of~$E_6$ basis under~$\phi^*_{(2,3)} \circ \phi^*_{(H^\perp)}$}
\label{table:23_H_subspaces}
\end{center}
\end{table}

Using this division of~$sl(3,\mathbb{O})$, we find a large list of subalgebras of~$sl(3,\mathbb{O})$ simply by combining certain subspaces.  The subspace description of these algebras, as well as their identity and signature in~$sl(3,\mathbb{O})$, is listed in Table~\ref{table:finest_refinement_subalgebras}.  We note that this fine refinement of~$sl(3,\mathbb{O})$ provides a description of the basis for~$su(3,1,\mathbb{H})_1$ and~$su(3,1,\mathbb{H})_2$, as well as~$sl(3,\mathbb{H})$ and~$sl(2,1,\mathbb{H})$.  The explicit basis elements for each algebra are listed in Appendix~\ref{appendix.basis_for_various_subalgebras_of_E6}.

\begin{table}[htbp]
\begin{center}
\begin{tabular}{|l|l|c|c|}
\hline
Basis & Subalgebra of $sl(3,\mathbb{O})$        &     Our   & $\phi^*_{(2,3)}\circ \phi^*_{(H^\perp)}$ \\
      & (our signature)& Signature & Signature \\
\hline
$R_{1,H}$   &  $su(2,\mathbb{H})\oplus su(2,\mathbb{C})^C \oplus su(2)$              & $(10+3+3,0)$  &  $(16,0)$ \\
\hline
$R_{1,H}+R_{(2,3),H^\perp}$ & $su(3,\mathbb{H})_2\oplus su(2,\mathbb{C})^C$ & $(21+3,0)$ & $(21+3,0)$ \\
\hline
$R_{1,H}+B_{1,H}$ & $sl(2,\mathbb{H})\oplus su(2,\mathbb{C})^C $ & $(21+1,6)$& $(15+1, 6)$ \\
                  & \hspace{.5cm} $\oplus su(2) \oplus u(1)$ &  &  \\
\hline
$R_{1,H}+B_{(2,3),H^\perp}$ & $su(2,1,\mathbb{H})_1 \oplus su(2,\mathbb{C})^C$ & $(16,8)$  & $(16,8)$ \\
\hline
$R_{1,H}+B_{(2,3),H}$ & $su(2,1,\mathbb{H})_2\oplus su(2,\mathbb{C})^C$ & $(13+3,8)$ & $(21+3,0)$\\
\hline
$R_{1,H}+R_{(2,3),H}$ & $su(3,\mathbb{H})\oplus su(2,\mathbb{C})^C  $ & $(21+3,0)$ & $(13+3,8)$ \\
\hline
$R_{1,H}+B_{1,H^\perp}$ & $so(5,\mathbb{R})\oplus so(4,1,\mathbb{R})$ & $(10+6,4)$  & $(10+10,0)$\\
\hline
$R_{1,H}+R_{1,H^\perp}$ & $so(9) = su(2,\mathbb{O})$ & $(36,0)$      & $(16,20)$\\
\hline
$R_{1,H}+B_{1,H}$ & $sl(3,\mathbb{H}) \oplus su(2,\mathbb{C})^C $ & $(21+3,14)$ & $(21+3,14)$\\
\hspace{.5cm} $+B_{(2,3),H}+R_{(2,3),H}$  &  & & \\
\hline
$R_{1,H}+B_{1,H}$ & $sl(2,1,\mathbb{H})_1\oplus su(2,\mathbb{C})_2$ & $(21+3,14)$ & $(21+3,14)$\\
\hspace{.5cm} $+R_{(2,3),H^\perp}$ & & & \\
\hspace{.5cm} $+B_{(2,3),H^\perp}$ & & & \\
\hline
$R_{1,H}+B_{1,H}$ & $sl(2,\mathbb{O})\oplus u(1)$ & $(35+1,10)$  & $(20, 25+1)$\\
\hspace{.5cm} $ +B_{1,H^\perp}$ & & & \\
\hspace{.5cm} $ +R_{1,H^\perp}$ & & & \\
\hline
$R_{1,H}+B_{(2,3),H^\perp}$ & $su(3,1,\mathbb{H})_1$ & $(24,12)$ & $(20,16)$ \\
\hspace{.5cm}$+R_{(2,3),H}+B_{1,H^\perp}$ & & & \\
\hline
$R_{1,H}+R_{(2,3),H^\perp}$ & $su(3,1,\mathbb{H})_2$ & $(24,12)$ & $(36,0)$  \\
\hspace{.5cm} $+B_{(2,3),H}+B_{1,H^\perp}$ & & & \\
\hline
$R_{1,H}+R_{(2,3),H^\perp}$ & $su(3,\mathbb{O})$   & $(52,0)$  & $(24,28)$ \\
\hspace{.5cm} $+R_{(2,3),H}+R_{1,H^\perp}$ & & & \\
\hline
$R_{1,H}+B_{(2,3),H^\perp}$ & $su(2,1,\mathbb{O})$   & $(36,16)$ & $(24,28)$ \\
\hspace{.5cm} $+B_{(2,3),H}+R_{1,H^\perp}$ & & & \\
\hline
\end{tabular}
\caption{Subalgebras of~$sl(3,\mathbb{O})$ using~$\phi^*_{(2,3)} \circ \phi^*_{(H^\perp)}$}
\label{table:finest_refinement_subalgebras}
\end{center}
\end{table}

\section{Chains of Subalgebras of~$sl(3,\mathbb{O})$ with Compatible Bases and Casimir Operators}
\label{ch:E6_further_structure.chains_of_subalgebras}

We used the involutive automorphisms to produce large simple subalgebras of~$sl(3,\mathbb{O})$.  These subalgebras range in dimension from~$52$, for~$F_4$, to~$21$, for~$C_3$.  Each of these subalgebras have involutive automorphisms which we can also use to find even smaller subalgebras, thereby giving a catalog of subalgebras chains contained within~$sl(3,\mathbb{O})$.  However, having identified the real form of the large subalgebras of~$sl(3,\mathbb{O})$, it is not too difficult a task to find smaller algebras simply by looking for simple subalgebras of smaller dimension and/or rank, using the tables of real forms listed in~\cite{gilmore} when needed.  In addition, we choose our representatives of the smaller subalgebras so that they use a subset of the preferred Casimir operators~$\lbrace \dot B^{(1)}_{tz}, \dot B^{(2)}_{tz}, \dot R^{(1)}_{x\l}, \dot S^{(1)}_{\l}, \dot G_{\l}, \dot A_{\l} \rbrace$ which we have chosen for~$sl(3,\mathbb{O})$.  We list the basis for these subalgebras in Appendix~\ref{appendix.basis_for_various_subalgebras_of_E6}.

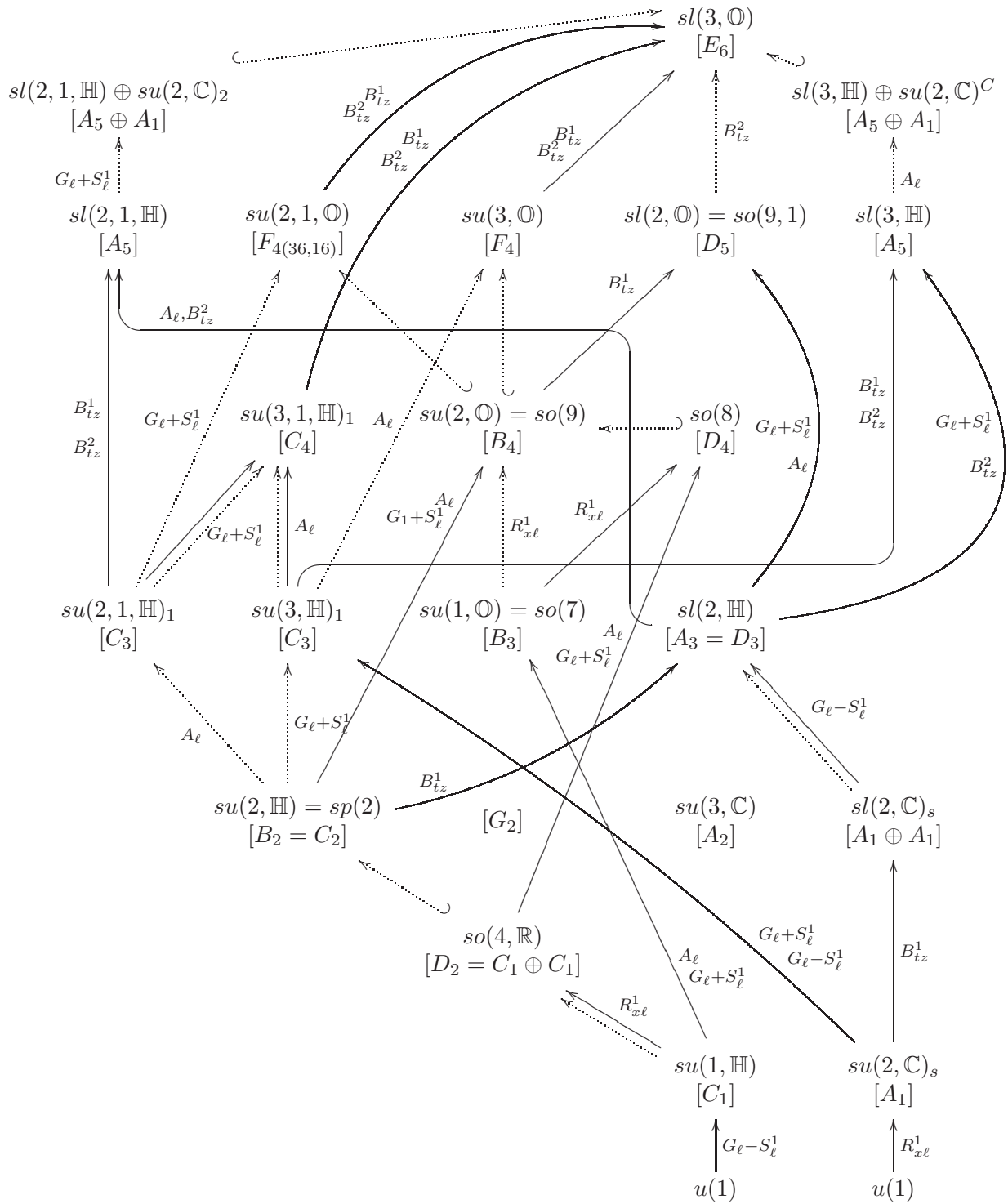
\begin{figure}[htbp]
\begin{center}
\begin{minipage}{6in}
\begin{center}
\hspace{-1in}
\xymatrixcolsep{5pt}
\xymatrix@M=3pt@H=1pt{
 & &               &   *++\txt{$sl(3,\mathbb{O})$ \\$[E_6]$}         &             \\
\save[]+<0cm,.4cm>*++\txt{$sl(2,1,\mathbb{H}) \oplus su(2,\mathbb{C})_2$\\$[A_5 \oplus A_1]$} \ar@{^{(}.>}[urrr]<2.75ex>_{} \restore & & & & \save[]+<0cm,.4cm>*++\txt{$sl(3,\mathbb{H}) \oplus su(2,\mathbb{C})^C$\\$[A_5 \oplus A_1]$}\ar@{_{(}.>}[ul]_{} \restore \\
*++\txt{$sl(2,1,\mathbb{H})$\\$[A_5]$}  \ar@{.>}[u]^{G_\l+S^1_\l} & *++\txt{$su(2,1,\mathbb{O})$\\$[F_{4(36,16)}]$} \ar@/^3pc/[uurr]^(.3){B^1_{tz}}^(.25){B^2_{tz}} & *++\txt{$su(3,\mathbb{O})$\\$[F_4]$} \ar[uur]^(.4){B^1_{tz}}^(.32){B^2_{tz}} &   *++\txt{$sl(2,\mathbb{O}) = so(9,1)$\\$[D_5]$} \ar@{.>}[uu]_{B^2_{tz}}    & *++\txt{$sl(3,\mathbb{H})$\\$[A_5]$}  \ar@{.>}[u]_{A_\l}  \\
&              &                           &                                    &              \\
& *++\txt{$su(3,1,\mathbb{H})_1$\\$[C_4]$} \ar@/^4pc/[uuuurr]^(.5){B^1_{tz}}^(.45){B^2_{tz}} & *++\txt{$su(2,\mathbb{O}) = so(9)$\\$[B_4]$} \ar@{_{(}.>}[uul]^(.4){} \ar@{_{(}.>}[uu]^(.4){} \ar[uur]^(.65){B^1_{tz}}&   *++\txt{$so(8)$ \\$[D_4]$} \ar@{_{(}.>}[l]^(.4){} &   \\
&              &                           &                                    &        \\
*++\txt{$su(2,1,\mathbb{H})_1$\\$[C_3]$} \ar[uuuu]<1ex>^(.55){B^1_{tz}}^(.45){B^2_{tz}} \ar@{.>}[uuuur]^(.5){G_\l+S^1_\l} \ar@{.>}[uur]{} \ar[uur]<1ex>_(.5){G_\l+S^1_\l} & *++\txt{$su(3,\mathbb{H})_1$\\$[C_3]$} \ar `u[r] `[rrru] [uurrruu] ^(.5){B^1_{tz}} ^(.4){B^2_{tz}} \ar@{.>}[uu]<2ex>{} \ar[uu]<1ex>_(.5){A_\l} \ar@{.>}[uuuur]^(.5){A_\l} & *++\txt{$su(1,\mathbb{O}) = so(7)$\\$[B_3]$} \ar[uur]^(.5){R^1_{x\l}} \ar@{.>}[uu]_(.5){R^1_{x\l}} &                   *++\txt{$sl(2,\mathbb{H})$\\$[A_3 = D_3]$} \ar `l[uuu] `[uuull] `[uuulll]_(.9){A_\l, B^2_{tz}} [uuuulll] \ar@/_4pc/[uuuu]^(.5){G_\l+S^1_\l}^(.4){A_\l} \ar@/_8pc/[uuuur]^(.7){G_\l+S^1_\l}^(.6){B^2_{tz}} &  \\
&               &                           &          &                \\
& *++\txt{$su(2,\mathbb{H}) = sp(2) $\\$[B_2 = C_2]$} \ar[uuuur]^(.8){A_\l}^(.75){G_1+S^1_\l} \ar@/_1.8pc/[uurr]^(.3){B^1_{tz}} \ar@{.>}[uul]^(.5){A_\l} \ar@{.>}[uu]<1ex>_(.5){G_\l+S^1_\l} & *++\txt{$[G_2]$} & *++\txt{$su(3,\mathbb{C})$\\$[A_2]$} & *++\txt{$sl(2,\mathbb{C})_s$\\ $[A_1\oplus A_1]$ } \ar@{.>}[uul]<1ex>{} \ar[uul]_(.5){G_\l-S^1_\l}\\
&               & *++\txt{$so(4,\mathbb{R})$\\$[D_2 = C_1 \oplus C_1]$} \ar@{_{(}.>}[ul]^(.4){} \ar[uuuuur]^(.6){A_\l}^(.55){G_\l+S^1_\l} &                  &  \\
&               &                           &   *++\txt{$su(1,\mathbb{H})$\\$[C_1]$} \ar[uuuul]<-1ex>_(.25){A_\l}_(.2){G_\l+S^1_\l} \ar@{.>}[ul]<1ex>{} \ar[ul]_(.45){R^1_{x\l}}      &  *++\txt{$su(2,\mathbb{C})_s$\\$[A_1]$}  \ar[uu]_(.5){B^1_{tz}}  \ar@/_1pc/[uuuulll]_(.25){G_\l+S^1_\l}_(.2){G_\l-S^1_\l} \\
& & & u(1) \ar[u]_(.4){G_\l - S^1_\l}  & u(1) \ar[u]_(.4){R^1_{x\l}}\\
}
\caption{Preferred subalgebra chains of~$E_6$ using the same basis}
\label{fig:clean_chain_of_subalgebras_of_E6}
\end{center}
\end{minipage}
\end{center}
\end{figure}

We display the chains of subalgebras of~$sl(3,\mathbb{O})$ in the following tables.  Each table is built from an~$u(1)$ algebra, which consists of merely a single Casimir operator.  We extend each algebra~$g$ to a larger algebra~$g^\prime$ by adding elements to the basis for~$g$.
  In particular, each algebra of higher rank must add compatible Casimir operators.  Figure~\ref{fig:clean_chain_of_subalgebras_of_E6} is built from the single~$u(1) = \langle G_l - S^1_l \rangle$ algebra.  This basis can be extended to 
$$u(1) \subset su(1,\mathbb{H}) \subset su(2,\mathbb{H}) \subset su(3,\mathbb{H})_1 \subset sl(3,\mathbb{H})$$
However, additional subalgebras of~$sl(3,\mathbb{O})$ can be inserted into this chain of subalgebras.  For instance, we can insert~$sl(2,\mathbb{H})$ between~$su(2,\mathbb{H})$ and~$sl(3,\mathbb{H})$, and extend~$sl(2,\mathbb{H})$ to~$sl(2,\mathbb{O})$.  We can expand~$su(1,\mathbb{H})$ to~$su(1,\mathbb{O}) = so(7)$ and insert the~$so(n,\mathbb{R})$ chain for~$n \ge 7$ into the figure.  However, we note that~$so(7)$ uses a basis in this chain that is not compatible with~$G_2 = aut(\mathbb{O})$.  We do not list all of the possible~$u(1)$ algebras, but do add the~$u(1) = \langle R^1_{x\l} \rangle$ subalgebra into the chain and extend it to~$su(2,\mathbb{C})_s$ and~$sl(2,\mathbb{C})_s$ which use the {\it standard} definition of~$su(2,\mathbb{C})$ and~$sl(2,\mathbb{C})$.  

\begin{figure}[htbp]
\begin{center}
\begin{minipage}{6in}
\begin{center}
\xymatrixcolsep{5pt}
\xymatrix@M=3pt@H=10pt{
 &  & *++\txt{$sl(3,\mathbb{O})$\\$[E_6]$}  &   &  \\
*++\txt{$sl(2,1,\mathbb{H}) \oplus su(2,\mathbb{C})_2$\\$[A_5 \oplus A_1]$} \ar@{^{(}.>}[urr]^(.5){} &  &    &   & *++\txt{$sl(3,\mathbb{H}) \oplus su(2,\mathbb{C})^C$\\$[A_5 \oplus A_1]$}  \ar@{_{(}.>}[ull]^(.5){} \\
*++\txt{$sl(2,1,\mathbb{H})$\\$[A_5]$} \ar@{.>}[u]^(.5){G_\l + S^1_\l} &  &    &   & *++\txt{$sl(3,\mathbb{H})$\\$[A_5]$} \ar@{.>}[u]_(.5){A_\l} \\
 & *++\txt{$su(2,1,\mathbb{O})$\\$[F_{4(36,16)}]$} \ar[uuur]^(.5){B^1_{tz}}^(.4){B^2_{tz}} &            & *++\txt{$su(3,\mathbb{O})$\\$[F_4]$} \ar[uuul]_(.5){B^1_{tz}}_(.4){B^2_{tz}} & \\
*++\txt{$su(3,1,\mathbb{H})_2$\\$[C_4]$} \ar@/^12.5pc/[uuuurr]^(.75){B^1_{tz}}^(.65){B^2_{tz}} & & & & *++\txt{$su(3,1,\mathbb{H})_1$\\$[C_4]$} \ar@/_12pc/[uuuull]_(.75){B^1_{tz}}_(.65){B^2_{tz}} \\
 &  &  &  &  &  \\
*++\txt{$su(2,1,\mathbb{H})_2$\\$[C_3]$} \ar@{.>}[uu]<1ex>^(.4){A_\l} \ar[uu]^(.4){} \ar@{.>}[uuur]_(.8){A_\l} \ar@/_3.5pc/[uuuurrrr]_(.55){A_\l} & *++\txt{$su(3,\mathbb{H})_2$\\$[C_3]$} \ar@{.>}[uul]<1ex>^(.4){G_\l+S^1_\l} \ar[uul]_(.4){} \ar@{.>}[uuurr]^(.75){G_\l+S^1_\l} \ar@/^1pc/[uuuul]^(.82){B^1_{tz}}^(.72){B^2_{tz}} & & *++\txt{$su(2,1,\mathbb{H})_1$\\$[C_3]$} \ar@{.>}[uur]<1ex>_(.5){} \ar[uur]^(.5){G_\l+S^1_\l}  \ar@{.>}[uuull]_(.75){G_\l+S^1_\l}  \ar@/^3pc/[uuuulll]_(.85){B^1_{tz}}_(.75){B^2_{tz}} & *++\txt{$su(3,\mathbb{H})_1$\\$[C_3]$} \ar@{.>}[uu]<1ex>_(.5){} \ar[uu]_(.4){A_\l} \ar@{.>}[uuul]_(.7){A_\l} \ar@/_4pc/[uuuu]_(.4){B^1_{tz}}_(.3){B^2_{tz}} \\
     &     &     *++\txt{$su(2,\mathbb{H})$\\$[C_2]$} \ar@{.>}[ull]<1ex>^(.5){G_\l + S^1_\l} \ar@{.>}[ul]<1ex>^(.5){A_\l} \ar@{.>}[ur]<1ex>_(.5){A_\l} \ar@{.>}[urr]<1ex>_(.5){G_\l + S^1_\l}
}
\caption{\noindent Four real forms of~$C_3$}
\label{fig:different_C3s}
\end{center}
\end{minipage}
\end{center}
\end{figure}
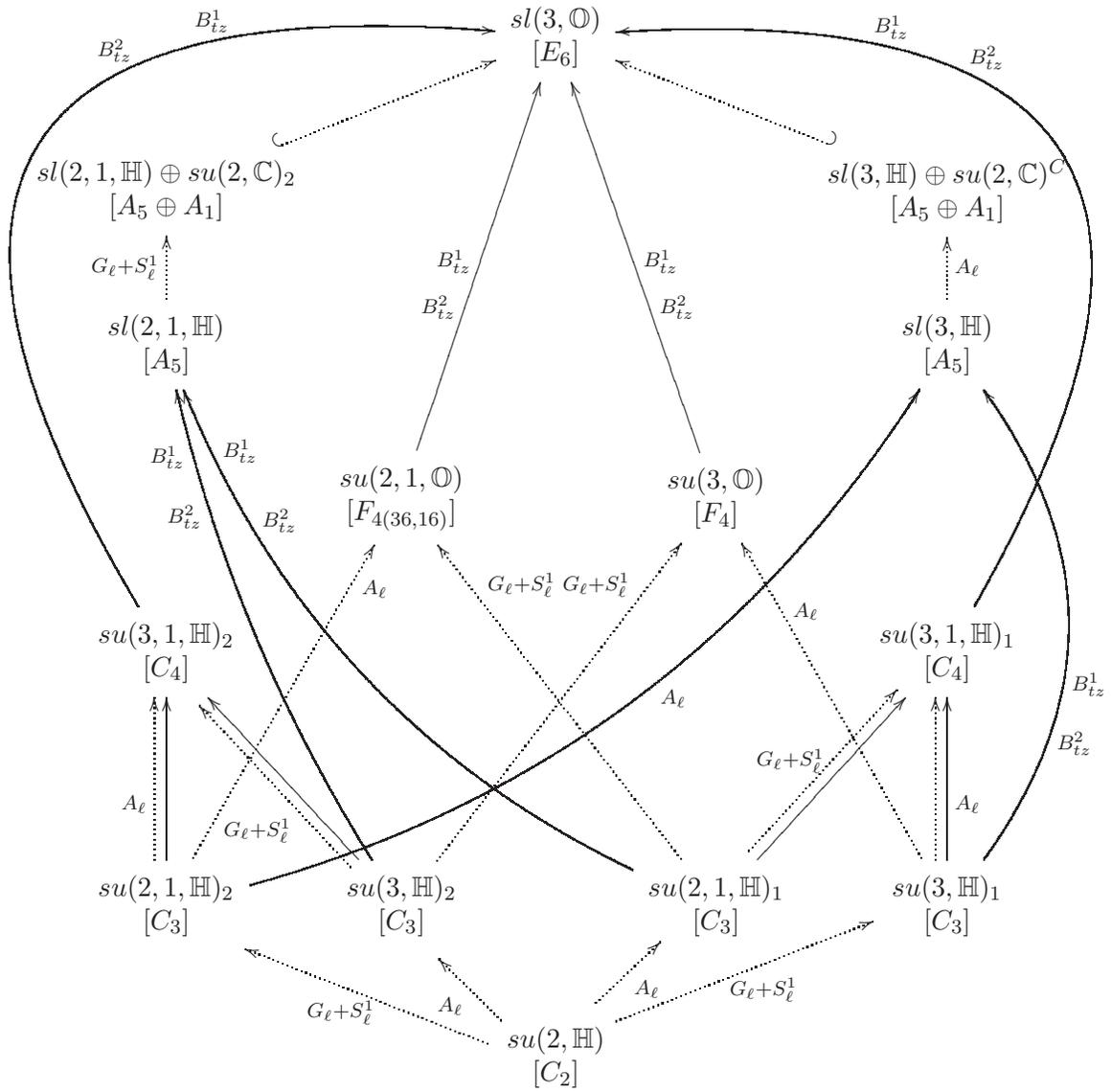

Finally, we have identified four different real forms of~$C_3$ which all contain~$su(2,\mathbb{H})$.  Space constraints limit us to listing only~$su(2,1,\mathbb{H})_1$ and~$su(3,\mathbb{H})_1$ in Figure~\ref{fig:clean_chain_of_subalgebras_of_E6}, but the algebras
$$su(2,1,\mathbb{H})_2 \hspace{1cm} su(3,\mathbb{H})_2 \hspace{1cm} su(3,1,\mathbb{H})_2$$
should also be in this table.  We list the four~$C_3$ algebras which are build from~$su(2,\mathbb{H})$ in 
Figure~\ref{fig:different_C3s}, and include all the algebras which are built from the~$C_3$ algebras.  Figure~\ref{fig:different_C3s} can be incorporated into Figure~\ref{fig:clean_chain_of_subalgebras_of_E6} without having to adjust choice of the Casimir operators.

\begin{figure}[htbp]
\begin{minipage}{6in}
\begin{center}
\xymatrixcolsep{5pt}
\hspace{-1in}
\xymatrix@M=3pt@H=10pt{
 & &               &   *++\txt{$sl(3,\mathbb{O})$ \\$[E_6]$}         &             \\
\save[]+<0cm,.4cm>*++\txt{$sl(2,1,\mathbb{H}) \oplus su(2,\mathbb{C})_2$\\$[A_5 \oplus A_1]$} \ar@{^{(}.>}[urrr]<2.75ex>_{} \restore& & & & \save[]+<0cm,.4cm>*++\txt{$sl(3,\mathbb{H}) \oplus su(2,\mathbb{C})^C$\\$[A_5 \oplus A_1]$}\ar@{_{(}.>}[ul]_{} \restore \\
*++\txt{$sl(2,1,\mathbb{H})$\\$[A_5]$}  \ar@{.>}[u]^{G_\l+S^1_\l} & *++\txt{$su(2,1,\mathbb{O})$\\$[F_{4(36,16)}]$} \ar@/^3pc/[uurr]^(.3){B^1_{tz}}^(.25){B^2_{tz}} & *++\txt{$su(3,\mathbb{O})$\\$[F_4]$} \ar[uur]^(.4){B^1_{tz}}^(.32){B^2_{tz}} &   *++\txt{$sl(2,\mathbb{O}) = so(9,1)$\\$[D_5]$} \ar@{.>}[uu]_{B^2_{tz}}    & *++\txt{$sl(3,\mathbb{H})$\\$[A_5]$}  \ar@{.>}[u]_{A_\l}  \\
&              &                           &                                    &              \\
& *++\txt{$su(3,1,\mathbb{H})_1$\\$[C_4]$} \ar@/^4pc/[uuuurr]^(.5){B^1_{tz}}^(.45){B^2_{tz}} & *++\txt{$su(2,\mathbb{O}) = so(9)$\\$[B_4]$} \ar@{_{(}.>}[uul]^(.4){} \ar@{_{(}.>}[uu]^(.4){} \ar[uur]^(.65){B^1_{tz}}&   *++\txt{$so(8)$ \\$[D_4]$} \ar@{_{(}.>}[l]^(.4){} &   \\
&              &                           &                                    &        \\
*++\txt{$su(2,1,\mathbb{H})_1$\\$[C_3]$} \ar[uuuu]<1ex>^(.55){B^1_{tz}}^(.45){B^2_{tz}} \ar@{.>}[uuuur]^(.5){G_\l+S^1_\l} \ar@{.>}[uur]{} \ar[uur]<1ex>_(.5){G_\l+S^1_\l} & *++\txt{$su(3,\mathbb{H})_1$\\$[C_3]$} \ar `u[r] `[rrru] [uurrruu] ^(.5){B^1_{tz}} ^(.4){B^2_{tz}} \ar@{.>}[uu]<2ex>{} \ar[uu]<1ex>_(.5){A_\l} \ar@{.>}[uuuur]^(.5){A_\l} & *++\txt{$su(1,\mathbb{O}) = so(7)$\\$[B_3]$} \ar[uur]^(.5){R^1_{x\l}} \ar@{.>}[uu]_(.5){R^1_{x\l}} &                   *++\txt{$sl(2,\mathbb{H})$\\$[A_3 = D_3]$} \ar `l[uuu] `[uuull] `[uuulll]_(.9){A_\l, B^2_{tz}} [uuuulll] \ar@/_4pc/[uuuu]^(.5){G_\l+S^1_\l}^(.4){A_\l} \ar@/_8pc/[uuuur]^(.7){G_\l+S^1_\l}^(.6){B^2_{tz}} &  \\
&               &                           &          &                \\
& *++\txt{$su(2,\mathbb{H}) = sp(2) $\\$[B_2 = C_2]$} \ar[uuuur]^(.8){A_\l}^(.75){G_1+S^1_\l} \ar@/_1.8pc/[uurr]^(.3){B^1_{tz}} \ar@{.>}[uul]^(.5){A_\l} \ar@{.>}[uu]<1ex>_(.5){G_\l+S^1_\l} & *++\txt{$[G_2]$} & *++\txt{$su(3,\mathbb{C})_s$\\$[A_2]$} \ar@{.>}[uull]_(.6){G_\l} & *++\txt{$sl(2,\mathbb{C})_s$\\ $[A_1\oplus A_1]$ } \ar@{.>}[uul]<1ex>{} \ar[uul]_(.5){G_\l-S^1_\l}\\
&               & *++\txt{$so(4,\mathbb{R})$\\$[D_2 = C_1 \oplus C_1]$} \ar@{_{(}.>}[ul]^(.4){} \ar[uuuuur]^(.6){A_\l}^(.55){G_\l+S^1_\l} &                  &  \\
&               &                           &   *++\txt{$su(1,\mathbb{H})$\\$[C_1]$} \ar[uuuul]<-1ex>_(.25){A_\l}_(.2){G_\l+S^1_\l} \ar@{.>}[ul]<1ex>{} \ar[ul]_(.45){R^1_{x\l}}      &  *++\txt{$su(2,\mathbb{C})_s$\\$[A_1]$}  \ar[uu]_(.5){B^1_{tz}}  \ar@{.>}[uul]<1ex>^(.5){} \ar[uul]_(.5){R^2_{x\l} \to S^1_\l } \\
& & & u(1) \ar[u]_(.4){G_\l - S^1_\l}  & u(1) \ar[u]_(.4){R^1_{x\l}}\\
}
\caption{Preferred subalgebra chains of~$E_6$ allowing a change of basis}
\label{fig:clean_change_of_basis_chain_of_subalgebras_of_E6}
\end{center}
\end{minipage}
\end{figure}
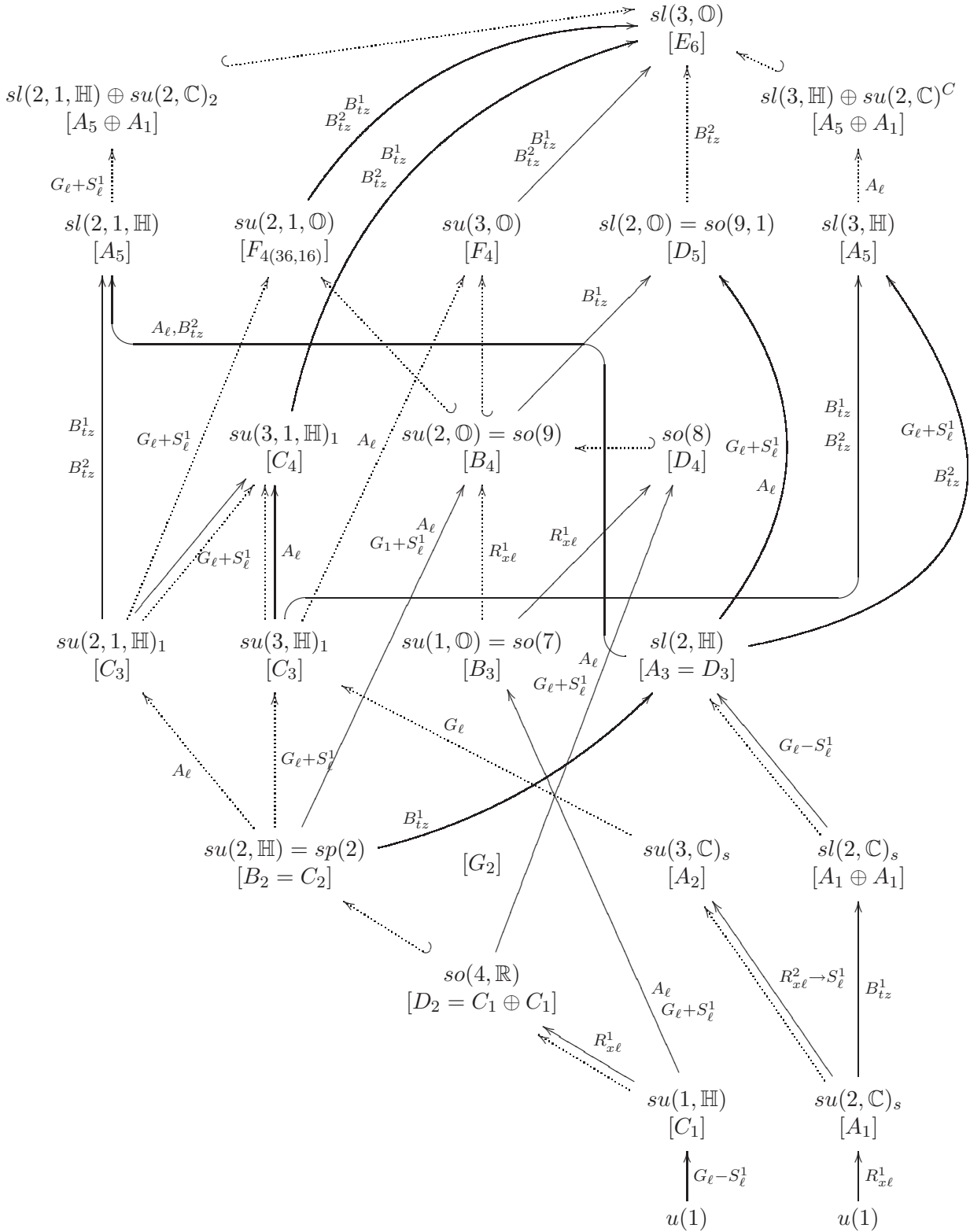

While the algebras in Figure~\ref{fig:clean_chain_of_subalgebras_of_E6} are built from its subalgebras by extending the subalgebra's basis, we do not allow a change of basis in any of the algebras.  In Figure~\ref{fig:clean_change_of_basis_chain_of_subalgebras_of_E6}, we do allow a change of basis at each subalgebra stage.  
The major difference this change allows is that the chain
$$su(2,\mathbb{C}) \subset su(3,\mathbb{C}) \subset su(3,\mathbb{H})$$
may be included with the Casimir operators~$S^1_\l$ and~$G_\l$ for each subalgebra extension.

\hspace{-1in}

\clearpage{}

\hspace{-1in}

\newpage{}

\part{Open Questions}
\label{ch:Open_Questions}

This study of $sl(3,\mathbb{O})$ has been driven by an underlying desire to find subalgebras useful for physics 
and to find a subalgebra structure which treats ``$\l$ as special.
Although we have not fully resolved these issues, this study of~$sl(3,\mathbb{O})$ has raised some additional interesting questions.

Ultimately, the interesting subalgebra structure we seek involves~$so(3,1,\mathbb{R})$ along with $su(3,\mathbb{C}) \oplus su(2,\mathbb{C}) \oplus u(1,\mathbb{C})$.  While we did not carry out an exhaustive search for this structure, we did find structures which are the direct sum of subalgebras of $sl(3,\mathbb{O})$, as indicated in Appendix~\ref{direct_sums_in_E6}.  
  This search may be more fruitful if we knew exactly how the desired algebra structure contains the direct sum $su(3,\mathbb{C}) \oplus su(2,\mathbb{C}) \oplus u(1,\mathbb{C})$ and $so(3,1,\mathbb{R})$.

There are a number of ways to interpret the phrase ``$\l$ is special'', and each leads to interesting substructures of $E_6$. In one scenario, we choose $\l$ as our preferred imaginary complex unit, leaving us free to choose the quaternionic subalgebra $\mathbb{H}$ to either contain $\l$ or to be independent of $\l$.  For instance, choosing $\mathbb{H}$ generated by $\langle 1, k, kl, l \rangle $ obviously creates the scenario $\mathbb{C} \subset \mathbb{H} \subset \mathbb{O}$, while choosing $\mathbb{H}$ to be generated by $\langle 1, i, j, k \rangle$ breaks $\mathbb{O}$ into division algebras~$\mathbb{C}, \mathbb{H}$ whose intersection is $\mathbb{R}$.  An interesting question would be to find the subalgebras of $sl(3,\mathbb{O})$ which respect these divisions of the octonions into $\mathbb{C}$ and $\mathbb{H}$.

The questions above may also be extended to the group $SL(3,\mathbb{O})$ in various ways.  Although we found in Section \ref{ch:E6_basic_structure.fix_l} the subalgebra of $sl(3,\mathbb{O})$ which fixes the $\l$ in the type $T=1$ vector of $sl(3,\mathbb{O})$, it is not clear how to give an expression for the one-parameter curves associated with the $8$ null rotations.  In addition, there may be transformations in $SL(3,\mathbb{O})$ which fix $\l$ but are not connected to the identity.  There could also be discrete transformations in $SL(3,\mathbb{O})$ which form a group and also fix $\l$.

We could use our automorphisms to study the subalgebra structures of the other real forms of $E_6$.  As a real algebra, we did not find any real form of $A_4$ in $sl(3,\mathbb{O})$, even though we can find it in $\phi^*_{H^\perp}\left( sl(3,\mathbb{O}) \right)$.  We therefore do expect that there will be slight differences in the subalgebra structures of these different real forms of $E_6$.

Finally, this work may also be extended to study the structure of the final two exceptional Lie algebras $E_7$ and $E_8$.  

\newpage{}

\part{Conclusion}
\label{ch:Conclusion}

We presented here a study of the subalgebra structure of $sl(3,\mathbb{O})$, a real form of the complex Lie algebra $E_6$.  We first examined the subalgebra structure of the complex Lie algebra~$E_6$.
  We then expanded upon Dray and Manogue's work with $SL(3,\mathbb{O})$ to provide subalgebras in the $3 \times 3$ case corresponding to the $2 \times 2$ subalgebra structure of $sl(2,\mathbb{O})$.  Finally, we used automorphisms of real algebras to provide some subalgebra structures in $sl(3,\mathbb{O})$ which are distinctly $3 \times 3$.

In Chapter \ref{ch:Lie_groups_Lie_algebras}, we presented methods which illustrated how root and weight diagrams could be used to visually identify the subalgebras of a given Lie algebra.  While the standard methods of determining subalgebras rely upon adding, removing, or folding along nodes in a Dynkin diagram, we showed here how to construct any of a Lie algebra's root or weight diagrams from its Dynkin diagram, and how to use geometric transformations to visually identify subalgebras using those weight and root diagrams.  In particular, we showed how these methods can be applied to algebras whose root and weight diagrams have dimensions four or greater.    In addition to pointing out the erroneous inclusion of ~$C_4 \subset F_4$ in ~\cite{gilmore, van_der_waerden}, we provided visual proof that ~$C_4 \subset E_6$ and listed all the complex subalgebras of ~$E_6$.  While we were primarily concerned with the subalgebras of ~$E_6$, these methods could be used to find subalgebras of any rank ~$l$ algebra.

In Chapter~\ref{ch:Division_Algebras_Applications}, we discussed the four division algebras.  We repeated the findings of Manogue and Schray who showed how multiplication in $\mathbb{H}$ and $\mathbb{O}$ can produce rotations in specific planes in $\mathbb{R}^4$ and $\mathbb{R}^8$.  Lorentz transformations in $(k+1)$ dimensions were related to determinant or trace preserving transformations of $2 \times 2$ hermitian matrices using the expressions developed in \cite{manogue_schray}.  After discussing the Albert algebra, we showed how to construct an expression for the determinant of a $3 \times 3$ octonionic hermitian matrix.  This material was needed for the construction of our real form of $E_6$.

We discussed our construction of the Lie group $SL(3,\mathbb{O})$ and its associated Lie algebra $sl(3,\mathbb{O})$ in Chapter \ref{ch:E6_basic_structure}.  The $2 \times 2$ formalism for hermitian octonionic matrices given in \cite{manogue_schray} was generalized to the $3 \times 3$ case using our notion of {\it type}.  Despite the use of nested matrices in transformations, we constructed an association between group and algebra transformations and produced the multiplication table in the Lie algebra.  We repeated many of the results shown in \cite{manogue_dray} which were needed for this work.  We identified $F_4$, $G_2$, and a preferred $su(3,\mathbb{C}) \subset G_2$.  We showed how our notion of {\it type} is related to {\it triality}, and identified a {\it strong} notion of triality among our transformations.  We also showed the existence of {\it continuous type transformations}, and identified subgroups of $SL(3,\mathbb{O})$ which contain these transformations.  We also used this type transformation to study type dependent and type independent subalgebras of $sl(3,\mathbb{O})$.  The subalgebra structures of $so(9,1)$, $sl(2,\mathbb{O})$, and $su(2,\mathbb{O})$, which were known for the $2 \times 2$ case, were generalized to the $3 \times 3$ case.  Finally, we found a non-simple subalgebra of $sl(3,\mathbb{O})$ which fixes $\l$.

In Chapter \ref{ch:E6_further_structure}, we used automorphisms of $sl(3,\mathbb{O})$ to find subalgebras of $sl(3,\mathbb{O})$.  We adapted theory related to automorphisms of real forms of complex  algebras to provide methods which helped us find subalgebras of our specific real form of $E_6$.  This search provided some surprising results. We found multiple real forms of $C_3$ and $C_4$ in $sl(3,\mathbb{O})$.  These real algebras differ both in how they reduce the octonions $\mathbb{O}$ to the quaternions $\mathbb{H}$ and how they used quaternionic and {\it orthogonal-quaternionic} spinor and dual spinor transformations.  We saw that each octonion $q$ defined both a quaternionic subalgebra which contained $q$ and a quaternionic subalgebra which was perpendicular to $q$.  Further, each of those subalgebras were used to construct real forms of $A_5$ and $C_4$ in $sl(3,\mathbb{O})$.  The automorphisms were also used to show that $sl(3,\mathbb{O})$ contains both a compact and non-compact form of $F_4$.  With the results of the previous section, we were able to construct subalgebra maps of $sl(3,\mathbb{O})$, labeled with Casimir operators, which show how its subalgebras sit in relation to each other.  When used with our basis of $sl(3,\mathbb{O})$, the automorphisms could also be used to construct these maps for the other real forms of $E_6$.

\bibliographystyle{unsrt}

\bibliography{thesis}

\newpage
\ \vspace{4.23cm}
\begin{center}
\large{APPENDICES}
\end{center}

\addcontentsline{toc}{part}{APPENDICES}
\newpage
\appendix

\section{Preferred Basis for $sl(3,\mathbb{O})$}

We list here our preferred basis for $sl(3,\mathbb{O})$.  Let $q \in \lbrace i,j,k,kl,jl,il,l \rbrace$.  In the following table, we use the letters $A_q$ and $G_q$ to denote those transverse rotations which are $G_2$ transformations, while $S_q$ denotes the transverse rotations which are in $so(7)$ but not in $G_2$.  The letter $R_{xq}$ denotes the rotation involving $x$ and $q$.  These four classes of transformations form $so(8)$, and we choose to treat each of them as a type $T = 1$ transformation.  Each of the remaining rotations $R^{(T)}_{xz}$ and $R^{(T)}_{zq}$ exist for each of the three types, as do the boosts $B^{(T)}_{tx}$ and $B^{(T)}_{tq}$.  We choose the boosts $B^{(1)}_{tz}$ and $B^{(2)}_{tz}$ to finish our basis for $sl(3,\mathbb{O})$.  Our particular choices are listed in Table \ref{table:our_basis}.

\begin{table}[htbp]
\begin{center}
\begin{tabular}{|l|ccc|}
\hline
Boosts  & \hspace{1cm}   $(2)$    \hspace{1cm}     &  \hspace{1cm}    $(3)$  \hspace{1cm} & \hspace{1cm}   $(21)$ \hspace{1cm} \\
\cline{2-4}
(Cat. 1)& $\dot B^{(1)}_{tz}$ & $\dot B^{(1)}_{tx}$ & $\dot B^{(1)}_{tq}$ \\
        & $\dot B^{(2)}_{tz}$ & $\dot B^{(2)}_{tx}$ & $\dot B^{(2)}_{tq}$ \\
        &                     & $\dot B^{(3)}_{tx}$ & $\dot B^{(3)}_{tq}$ \\
\hline
Rotations &    $(7)$            &     $(3)$      &    $(21)$ \\
\cline{2-4}
(Cat. 2)  & $\dot R^{(1)}_{xq}$ & $\dot R^{(1)}_{xz}$ & $\dot R^{(1)}_{zq}$ \\
          &                     & $\dot R^{(2)}_{xz}$ & $\dot R^{(2)}_{zq}$ \\
          &                     & $\dot R^{(3)}_{xz}$ & $\dot R^{(3)}_{zq}$ \\
\hline
Transverse &   $(7)$       &     $(7)$      &   $(7)$ \\
\cline{2-4}
Rotations  &  $\dot A_q$        &     $\dot G_q$      & $\dot S^{(1)}_q$ \\
(Cat. 3)   &   &   &  \\
\hline
\end{tabular}
\caption{Our preferred basis for $sl(3,\mathbb{O})$}
\label{table:our_basis}
\end{center}
\end{table}

\section{Preferred Basis for various subalgebras of $E_6$}
\label{appendix.basis_for_various_subalgebras_of_E6}

We list here the preferred basis for the various maximal subalgebras of $E_6$, as found in Section \ref{ch:E6_further_structure.three_automorphisms}.  We let the indices $q$ refer to $\lbrace i,j,k,kl,jl,il,l \rbrace$, and designate indices $h$ and $h^\perp$ to run over the sets $\lbrace k, kl, l \rbrace$ and $\lbrace i,j,jl,il \rbrace$, respectively.  Hence, $\dot A_q, \dot A_h$, and $\dot A_{h^\perp}$ represent $7, 3$, and $4$ basis elements, respectively.  

The tables below are arranged as follows:  We list the subalgebra of $sl(3,\mathbb{O})$, its signature, and then the basis elements.  In order to keep some consistency, an $\unde$ replaces the list of basis elements if it is not in our subalgebra.  Our preferred Casimir operators are taken as linear combinations of the elements 
$$\lbrace \dot B^{(1)}_{tz}, \dot B^{(2)}_{tz}, \dot R^{(1)}_{x\l}, \dot S^{(1)}_{\l}, \dot G_{\l}, \dot A_{\l} \rbrace$$

\begin{center}
\renewcommand{\arraystretch}{1.3}
\begin{longtable}{|c|c|ccc|}
\hline
Subalgebra & Signature & \multicolumn{3}{c|}{Basis elements} \\
\hline
$sl(3,\mathbb{O})$ & $(52,26)$ &  $\dot B^{(1)}_{tz}, \dot B^{(2)}_{tz}$ & $\dot B^{(1)}_{tx}, \dot B^{(2)}_{tx}, \dot B^{(3)}_{tx}$ & $\dot B^{(1)}_{tq}, \dot B^{(2)}_{tq}, \dot B^{(3)}_{tq}$\\
$[E_6]$            &           &  $\dot R^{(1)}_{xq}$ & $\dot R^{(1)}_{xz}, \dot R^{(2)}_{xz}, \dot R^{(3)}_{xz}$ & $\dot R^{(1)}_{zq}, \dot R^{(2)}_{zq}, \dot R^{(3)}_{zq}$\\
                   &           &  $\dot A_{q}$ & $\dot G_q$ & $\dot S^{(1)}_q$ \\
\hline
$su(3,\mathbb{O})$ & $(52,0)$ &  $\und \unde$ & $\und \und \unde$ & $\und \und \unde$\\
$[F_4]$            &           &  $\dot R^{(1)}_{xq}$ & $\dot R^{(1)}_{xz}, \dot R^{(2)}_{xz}, \dot R^{(3)}_{xz}$ & $\dot R^{(1)}_{zq}, \dot R^{(2)}_{zq}, \dot R^{(3)}_{zq}$\\
                   &           &  $\dot A_{q}$ & $\dot G_q$ & $\dot S^{(1)}_q$ \\
\hline
$su(2,1,\mathbb{O})$ & $(36,16)$ &  $\und \unde$ & $\und \dot B^{(2)}_{tx}, \dot B^{(3)}_{tx}$ & $\und \dot B^{(2)}_{tq}, \dot B^{(3)}_{tq}$\\
$[F_4]$            &           &  $\dot R^{(1)}_{xq}$ & $\dot R^{(1)}_{xz}, \und \unde$ & $\dot R^{(1)}_{zq}, \und \unde$\\
                   &           &  $\dot A_{q}$ & $\dot G_q$ & $\dot S^{(1)}_q$ \\
\hline
$sl(2,\mathbb{O})$ & $(36,9+1)$ &  $\dot B^{(1)}_{tz}, \dot B^{(2)}_{tz}$ & $\dot B^{(1)}_{tx}, \und \unde$ & $\dot B^{(1)}_{tq}, \und \unde $ \\
$\oplus \lbrace \dot B^{(2)}_{tz} -  \dot B^{(3)}_{tz} \rbrace$    &           &  $\dot R^{(1)}_{xq}$ & $\dot R^{(1)}_{xz}, \und \unde$ & $\dot R^{(1)}_{zq}, \und \unde$\\
$[D_5 \oplus u(1)]$            &           &  $\dot A_{q}$ & $\dot G_q$ & $\dot S^{(1)}_q$ \\
\cline{3-5}
   &   & \multicolumn{3}{c|}{$u(1): \dot B^{(2)}_{tz} - \dot B^{(3)}_{tz}$} \\
\hline
$sl(3,\mathbb{H})$ & $(21+3,14)$ &  $\dot B^{(1)}_{tz}, \dot B^{(2)}_{tz}$ & $\dot B^{(1)}_{tx}, \dot B^{(2)}_{tx}, \dot B^{(3)}_{tx}$ & $\dot B^{(1)}_{th}, \dot B^{(2)}_{th}, \dot B^{(3)}_{th}$\\
$\oplus su(2,\mathbb{C})^C$      &           &  $\dot R^{(1)}_{xh}$ & $\dot R^{(1)}_{xz}, \dot R^{(2)}_{xz}, \dot R^{(3)}_{xz}$ & $\dot R^{(1)}_{zh}, \dot R^{(2)}_{zh}, \dot R^{(3)}_{zh}$\\
$[A_5 \oplus C_1]$            &           & $\unde$ & $\dot G_h$ & $\dot S^{(1)}_h$ \\
\cline{3-5}
                   &           &  \multicolumn{3}{c|}{$su(2,\mathbb{C})^C: \dot A_h $} \\
\hline
$sl(2,1,\mathbb{H})$          & $(21+3,14)$ &  $\dot B^{(1)}_{tz}, \dot B^{(2)}_{tz}$ & $\dot B^{(1)}_{tx}, \und \unde$ & $\dot B^{(1)}_{th}, \dot B^{(2)}_{th^\perp}, \dot B^{(3)}_{th^\perp}$\\
$\oplus su(1,\mathbb{C})_2$    &           &  $\dot R^{(1)}_{xh}$ & $\dot R^{(1)}_{xz}, \und \unde$ & $\dot R^{(1)}_{zh}, \dot R^{(2)}_{zh^\perp}, \dot R^{(3)}_{zh^\perp}$\\
$[A_5 \oplus C_1]$            &           &  $\dot A_{h}$ & \multicolumn{2}{c|} {$\dot G_h - \dot S^{(1)}_h$} \\
\cline{3-5}
                   &           &  \multicolumn{3}{c|}{$su(2,\mathbb{C})_2 : \dot G_h + 2\dot S^{(1)}_h$} \\
\hline
$su(3,1,\mathbb{H})_1$ & $(24,12)$ &  $\und \unde$ & $\und \und \unde$ & $\dot B^{(1)}_{th^\perp}, \dot B^{(2)}_{th^\perp}, \dot B^{(3)}_{th^\perp}$\\
$[C_4]$            &           &  $\dot R^{(1)}_{xh}$ & $\dot R^{(1)}_{xz}, \dot R^{(2)}_{xz}, \dot R^{(3)}_{xz}$ & $\dot R^{(1)}_{zh}, \dot R^{(2)}_{zh}, \dot R^{(3)}_{zh}$\\
                   &           &  $\dot A_{h}$ & $\dot G_h$ & $\dot S^{(1)}_h$ \\
\hline
$su(3,1,\mathbb{H})_2$ & $(24,12)$ &  $\und \unde$ & $\und B^{(2)}_{tx}, B^{(3)}_{tx}$ & $\dot B^{(1)}_{th^\perp}, \dot B^{(2)}_{th}, \dot B^{(3)}_{th}$\\
$[C_4]$            &           &  $\dot R^{(1)}_{xh}$ & $\dot R^{(1)}_{xz}, \und \unde$ & $\dot R^{(1)}_{zh}, \dot R^{(2)}_{zh^\perp}, \dot R^{(3)}_{zh^\perp}$\\
                   &           &  $\dot A_{h}$ & $\dot G_h$ & $\dot S^{(1)}_h$ \\
\hline
$su(2,\mathbb{O}) = so(9)$ & $(36,0)$ &  $\und \unde$ & $\und \und \unde$ & $\und \und \unde$\\
$[B_4]$            &           &  $\dot R^{(1)}_{xq}$ & $\dot R^{(1)}_{xz}, \und \unde $ & $\dot R^{(1)}_{zq}, \und \unde$\\
                   &           &  $\dot A_{q}$ & $\dot G_q$ & $\dot S^{(1)}_q$ \\
\hline
$so(8)$ & $(28,0)$ &  $\und \unde$ & $\und \und \unde$ & $\und \und \unde$\\
$[D_4]$            &           &  $\dot R^{(1)}_{xq}$ & $\und \und \unde $ & $\und \und \unde$\\
                   &           &  $\dot A_{q}$ & $\dot G_q$ & $\dot S^{(1)}_q$ \\
\hline

\newpage
\hline
Subalgebra & Signature & \multicolumn{3}{c|}{Basis elements} \\
\hline
$su(3,\mathbb{H})_1$ & $(21+3,0)$ &  $\und \unde$ & $\und \und \unde$ & $ \und \und \unde$ \\
$\oplus su(2,\mathbb{C})^C$            &           &  $\dot R^{(1)}_{xh}$ & $\dot R^{(1)}_{xz}, \dot R^{(2)}_{xz}, \dot R^{(3)}_{xz} $ & $\dot R^{(1)}_{zh}, \dot R^{(2)}_{zh}, \dot R^{(3)}_{zh}$\\
$[C_3 \oplus A_1]$ &           &  $\unde$ & $\dot G_h$ & $\dot S^{(1)}_h$ \\
\cline{3-5}
                   &           &  \multicolumn{3}{c|}{$su(2,\mathbb{C})^C : \dot A_h$} \\
\hline
$su(2,1,\mathbb{H})_1$ & $(13+3,8)$ &  $\und \unde$ & $\und \und \unde$ & $ \und \dot B^{(2)}_{th^\perp} \dot B^{(3)}_{th^\perp}$ \\
$\oplus su(2,\mathbb{C})_2$            &           &  $\dot R^{(1)}_{xh}$ & $\dot R^{(1)}_{xz}, \und \unde $ & $\dot R^{(1)}_{zh}, \und \unde$\\
$[C_3 \oplus A_1]$ &           &  $\dot A_h $ & \multicolumn{2}{c|}{$\dot G_h - \dot S^{(1)}_h$} \\
\cline{3-5}
                   &           &  \multicolumn{3}{c|}{$su(2,\mathbb{C})_2: \dot G_h + 2\dot S^{(1)}_h$} \\
\hline
$su(3,\mathbb{H})_2$ & $(21+3,0)$ &  $\und \unde$ & $\und \und \unde$ & $ \und \und \unde$\\
$\oplus su(2,\mathbb{C})_2$            &           &  $\dot R^{(1)}_{xh}$ & $\dot R^{(1)}_{xz}, \und \unde $ & $\dot R^{(1)}_{zh}, \dot R^{(2)}_{zh^\perp}, \dot R^{(3)}_{zh^\perp}$ \\
$[C_3 \oplus A_1]$ &           &  $\dot A_h $ & \multicolumn{2}{c|}{$\dot G_h - \dot S^{(1)}_h$} \\
\cline{3-5}
                   &           &  \multicolumn{3}{c|}{$su(2,\mathbb{C})_2: \dot G_h + 2\dot S^{(1)}_h$} \\
\hline
$su(2,1,\mathbb{H})_2$ & $(13+3,8)$ &  $\und \unde$ & $\und \dot B^{(2)}_{tx}, \dot B^{(3)}_{tx}$ & $ \und \dot B^{(2)}_{th}, \dot B^{(3)}_{th}$ \\
$\oplus su(2,\mathbb{C})^C$            &           &  $\dot R^{(1)}_{xh}$ & $\dot R^{(1)}_{xz}, \und \unde $ & $\dot R^{(1)}_{zh}, \und \unde $\\
$[C_3 \oplus A_1]$ &           &  $\und $ & $\dot G_h$ & $\dot S^{(1)}_h$ \\
\cline{3-5}
                   &           &  \multicolumn{3}{c|}{$su(2,\mathbb{C})^C: \dot A_h $} \\
\hline

$su(1,\mathbb{O}) = so(7)$ & $(21,0)$ &  $\und \unde$ & $\und \und \unde$ & $\und \und \unde$\\
$[B_3]$            &           &  $\unde$ & $\und \und \unde $ & $\und \und \unde$\\
                   &           &  $\dot A_{q}$ & $\dot G_q$ & $\dot S^{(1)}_q$ \\
\hline
$sl(2,\mathbb{H})$ & $(15+3+3,0)$ &  $\dot B^{(1)}_{tz}, \unde$ & $\dot B^{(1)}_{tx}, \und \unde$ & $\dot B^{(1)}_{th}, \und \unde$ \\
$\oplus su(2,\mathbb{C})^C$ &           &  $\dot R^{(1)}_{xh}$ & $\dot R^{(1)}_{xz}, \und \unde $ & $\dot R^{(1)}_{zh}, \und \unde$\\
$\oplus su(2,\mathbb{C})_2$  &           &  $\unde$ & \multicolumn{2}{c|}{$\dot G_h - \dot S^{(1)}_h$} \\
\cline{3-5}
$[A_3 = D_3 \oplus A_1 \oplus A_1]$ &           &  \multicolumn{3}{c|}{$su(2,\mathbb{C})^C: \dot A_h$} \\
\cline{3-5}
                              &           &  \multicolumn{3}{c|}{$su(2,\mathbb{C})_2: \dot G_h + 2\dot S^{(1)}_h$} \\
\hline
$su(2,\mathbb{H})$ & $(10+3+3,0)$ &  $\und \unde$ & $\und \und \unde$ & $ \und \und \unde$ \\
$\oplus su(2,\mathbb{C})^C$ &           &  $\dot R^{(1)}_{xh}$ & $\dot R^{(1)}_{xz}, \und \unde $ & $\dot R^{(1)}_{zh}, \und \unde$\\
$\oplus su(2,\mathbb{C})_2$  &           &  $\unde$ & \multicolumn{2}{c|}{$\dot G_h - \dot S^{(1)}_h$} \\
\cline{3-5}
$[B_2 = C_2 \oplus A_1 \oplus A_1]$ &           &  \multicolumn{3}{c|}{$su(2,\mathbb{C})^C: \dot A_h$} \\
\cline{3-5}
                              &           &  \multicolumn{3}{c|}{$su(2,\mathbb{C})_2: \dot G_h + 2\dot S^{(1)}_h$} \\
\hline
$sl(3,\mathbb{C})_s$ & $(8,8)$ &  $\dot B^{(1)}_{tz}, \dot B^{(2)}_{tz}$ & $\dot B^{(1)}_{tx}, \dot B^{(2)}_{tx}, \dot B^{(3)}_{tx}$ & $\dot B^{(1)}_{t\l}, \dot B^{(2)}_{t\l}, \dot B^{(3)}_{t\l} $\\
$[A_2 \oplus A_2]$                      & &  $\dot R^{(1)}_{x\l}, \dot R^{(2)}_{x\l}$ & $\dot R^{(1)}_{xz}, \dot R^{(2)}_{xz}, \dot R^{(3)}_{xz} $ & $\dot R^{(1)}_{z\l}, \dot R^{(2)}_{z\l}, \dot R^{(31)}_{z\l}$ \\  
                             & &  $\unde$ & $\unde $ & $\unde$ \\
\hline
$Aut(\mathbb{O})$ & $(14,0)$ &  $\und \unde$ & $\und \und \unde$ & $\und \und \unde$\\
\nopagebreak $[G_2]$                    & &  $\unde$ & $\und \und \unde $ & $\und \und \unde$\\
\nopagebreak                           & &  $\dot A_q$ & $\dot G_q $ & $\unde$ \\
\hline
$su(3,\mathbb{C})^C$ & $(8,0)$ &  $\und \unde$ & $\und \und \unde$ & $\und \und \unde$\\
$[A_2]$                      & &  $\unde$ & $\und \und \unde $ & $\und \und \unde$\\
                             & &  $\dot A_q$ & $\dot G_\l $ & $\unde$ \\
\hline
$su(3,\mathbb{C})_s$ & $(8,0)$ &  $\und \unde$ & $\und \und \unde$ & $\und \und \unde$\\
$[A_2]$                      & &  $\dot R^{(1)}_{x\l}, \dot R^{(2)}_{x\l}$ & $\dot R^{(1)}_{xz}, \dot R^{(2)}_{xz}, \dot R^{(3)}_{xz} $ & $\dot R^{(1)}_{zh}, \dot R^{(2)}_{zh}, \dot R^{(3)}_{zh}$\\
                             & &  $\unde$ & $\unde $ & $\unde$ \\
\hline
$sl(2,\mathbb{C})_s$ & $(6,0)$ &  $\dot B^{(1)}_{tz}, \unde$ & $\dot B^{(1)}_{tx}, \und \unde$ & $\dot B^{(1)}_{t\l}, \und \unde$\\
$[A_1 \oplus A_1]$                      & &  $\dot R^{(1)}_{x\l}$ & $\dot R^{(1)}_{xz}, \und \unde $ & $\dot R^{(1)}_{z\l}, \und \unde$\\
                             & &  $\unde$ & $\unde $ & $\unde$ \\
\hline
$su(2,\mathbb{C})_s$ & $(3,0)$ &  $\und \unde$ & $\und \und \unde$ & $\und \und \unde$\\
$[A_1 ]$                      & &  $\dot R^{(1)}_{x\l}$ & $\dot R^{(1)}_{xz}, \und \unde $ & $\dot R^{(1)}_{z\l}, \und \unde$\\
                             & &  $\unde$ & $\unde $ & $\unde$ \\
\hline
$su(1,\mathbb{H})$ & $(3,0)$ &  $\und \unde$ & $\und \und \und$ & $\und \und \unde$ \\
$[A_1]$                    & &  $\und$ & $\und \und \unde$ & $\und \und \unde$ \\
                           & &  $\unde$ & \multicolumn{2}{c|}{$\dot G_h - \dot S^{(1)}_h$} \\
\hline

\caption{Preferred basis for subalgebras of $E_6$}
\label{table:subalgebra_basis}
\end{longtable}
\end{center}

\section{Basis for Small Subalgebras of $sl(3,\mathbb{O})$}

\begin{center}
\renewcommand{\arraystretch}{1.3}
\begin{longtable}{|l|l|l|l|}
\hline
Algebra & Casimir  & Basis & Set for $q$ \\
        & Operator &       &             \\
\hline
\hline
$sl(3,\mathbb{O}) = E_6$ & $B^{(1)}_{tz}, B^{(2)}_{tz}$ & $B^{(1)}_{tx}, B^{(2)}_{tx}, B^{(3)}_{tx}$ & \\
                         &                      & $B^{(1)}_{tq}, B^{(2)}_{tq}, B^{(3)}_{tq}$ & $k,k\l,\l$ \\
                         &                      & $R^{(1)}_{zq}, R^{(2)}_{zq}, R^{(3)}_{zq}$ & $i,j,i\l,j\l$ \\
                         &                      & $A_q, G_q, S^{(1)}_q$ & $i,j,i\l,j\l$ \\
\hline
$su(3,1,\mathbb{H})_1 = C_4$ & $G_\l+S^{(1)}_\l$ & $G_q + S^{(1)}_q$ & $k,k\l$ \\
                           &             & $R^{(2)}_{xz}, R^{(3)}_{xz}$ &   \\
                           &             & $R^{(2)}_{zq}, R^{(3)}_{zq}$ & $k,k\l,\l$  \\
                           &             & $B^{(1)}_{tq}$ & $i,j,i\l,j\l$ \\
\hline
$su(2,1,\mathbb{H})_1 = C_3$ &           & $B^{(2)}_{tq}, B^{(3)}_{tq}$ & $i,j,i\l,j\l$\\
\hline
$su(2,\mathbb{H})\oplus so(3)$ & $A_\l$ & $A_q$      & $k,k\l$   \\
\hline
$su(2,\mathbb{H}) = sp(2) = C_2$ &            & $R^{(1)}_{xz}, R^{(1)}_{zq}$  & $k,k\l,\l$ \\
\hline
$so(4) = so(3)\oplus so(3) = D_2$ & $R^{(1)}_{x\l}$ & $R^{(1)}_{xq}$  & $k,k\l$ \\
\hline
$su(1,\mathbb{H}) = C_1$ &            & $G_q-S^{(1)}_q$ & $k,k\l$ \\
\hline
$u(1)$ & $G_\l - S^{(1)}_\l$ & & \\
\hline
\end{longtable}
\end{center}

\newpage{}
\begin{center}
\renewcommand{\arraystretch}{1.3}
\noindent
\begin{longtable}{|l|l|l|l|l|}
\hline
Subalgebra & Algebra & Casimir  & New Basis & Choices  \\
           &         & Operator & Elements  & for $q$  \\ 
\hline
\hline
$su(2,\mathbb{H}) = C_2$ & $su(2,\mathbb{O}) = so(9)$ & $A_\l, G_\l+S^{(1)}_\l$ & $R^{(1)}_{xq}, R^{(1)}_{zq} $ & $i,j,i\l,j\l$ \\
                         & \hspace{.1in} $ = B_4$      &                       & $A_q, G_q + S^{(1)}_q$ & $k,k\l$ \\
                         &                                 &                       & $A_q, G_q, S^{(1)}_q$ & $i,j,i\l,j\l$ \\
\hline
$su(2,\mathbb{H}) = sp(2)$ & $su(3,\mathbb{H})_1 = C_3$ & $G_\l +S^{(1)}_\l$ & $R^{(2)}_{xz}, R^{(3)}_{xz}$ & \\
$ = C_2$                         &                          &              & $R^{(2)}_{zq}, R^{(3)}_{zq}$ & $k,k\l,\l$\\
                                 &                          &              & $G_q+S^{(1)}_q$ & $k,k\l$\\
\hline
$su(3,\mathbb{H})_1 = C_3$ & $sl(3,\mathbb{H}) = A_5$ & $B^{(1)}_{tz}, B^{(2)}_{tz}$ & $B^{(1)}_{tx},B^{(2)}_{tx},B^{(3)}_{tx}$ & \\
                         &                          &                      & $B^{(1)}_{tq}, B^{(2)}_{tq}, B^{(3)}_{tq}$ & $k,k\l,\l$\\
\hline
\hline
$su(2,1,\mathbb{H})_1$ & $sl(2,1,\mathbb{H})$ & $B^{(1)}_{tz}, B^{(2)}_{tz}$ & $ R^{(2)}_{zq}, R^{(3)}_{zq}$ & $i,i\l,j,j\l$\\
                     &                      &                      & $B^{(1)}_{tx}, B^{(1)}_{tq}$  & $k,k\l,\l$\\
\hline
\hline
$su(2,1,\mathbb{H})_1$ & $su(2,1,\mathbb{O})$  & $G_\l+S^{(1)}_\l$ & $G_q+S^{(1)}_q$ & $k,k\l$ \\
                     & \hspace{.1in} $= F_{4(36,16)}$&             & $B^{(2)}_{tq}, B^{(3)}_{tq}$ & $k,k\l,\l$ \\
                     &                 &             & $B^{(2)}_{tx}, B^{(3)}_{tx}$ &  \\
                     &                 &             & $R^{(1)}_{xq}, R^{(1)}_{zq}$ & $i,i\l,j,j\l$\\
                     &                 &             & $A_q, G_q, S^{(1)}_q$ & $i,i\l,j,j\l$\\
\hline
\hline
                         & $u(1)$                   & $R^{(1)}_{x\l}$ & & \\
\hline
$u(1)$                   & $su(2,\mathbb{C})_s = A_1$ &            & $R^{(1)}_{xz}, R^{(1)}_{zq}$ & $l$ \\
\hline
$su(2,\mathbb{C})_s = A_1$ & $su(3,\mathbb{C})_s = A_2$ & $R^{(2)}_{x\l}$ & $R^{(2)}_{xz}, R^{(3)}_{xz}$ & \\
                         &                          &            & $R^{(2)}_{zq}, R^{(3)}_{zq}$ & $l$ \\
\hline
$su(3,\mathbb{C})_s = A_2$ & $su(3,\mathbb{H})_1 = C_3$ & $G_\l$ & $G_q - S^{(1)}_q$ & $k,k\l$ \\
                         &                          &       & $R^{(1)}_{xq}   $ & $k,k\l$ \\
                         &                          &       & $R^{(1)}_{zq}, R^{(2)}_{zq}, R^{(3)}_{zq}$ & $k,k\l$\\
                         &                          &       & $G_q + S^{(1)}_q$ & $k,k\l$ \\
\hline
\hline
$su(3,\mathbb{H})_1 = C_3$ & $su(3,\mathbb{O}) = F_4$  & $A_\l$  & $R^{(1)}_{zq}, R^{(2)}_{zq}, R^{(3)}_{zq}$   & $i,j,i\l,j\l$ \\
                         &                          &        & $R^{(1)}_q, A_q, G_q, S^{(1)}_q$                     & $i,j,i\l,j\l$ \\
                         &                          &        & $A_q$                     & $k,k\l$ \\

\hline
\hline
$su(2,\mathbb{H}) = sp(2)$ & $sl(2,\mathbb{H})$ & $B^{(1)}_{tz}$ & $B^{(1)}_{tx}, B^{(1)}_{tq}$ & $k,k\l,\l$\\
$=C_2$                     & $=A_3 = D_3$     &              &                      & \\
\hline
$sl(2,\mathbb{H})$ & $sl(3,\mathbb{H}) = A_5$ & $B^{(2)}_{tz}, G_\l+S^{(1)}_\l$ & $R^{(2)}_{xz}, R^{(3)}_{xz}$ & \\
$=A_3=D_3$         &                          &                       & $B^{(2)}_{tx}, B^{(3)}_{tx}$ & \\
                   &                          &                       & $B^{(2)}_{tq}, B^{(3)}_{tq}$ & $k,k\l,\l$\\
                             &                          &                       & $R^{(2)}_{zq}, R^{(3)}_{zq}$ & $k,k\l,\l$\\
                             &                          &                       & $G_q+S^{(1)}_q$ & $k,k\l$\\
\hline
$sl(2,\mathbb{H})$ & $sl(2,\mathbb{O})=so(9,1)$ & $G_\l+S^{(1)}_\l, A_\l$ & $B^{(1)}_{tq}$ & $i,j,i\l,j\l$ \\
$=A_3=D_3$                     &  $= D_5$ &                  & $R^{(1)}_{xq}, R^{(1)}_{zq}$  & $i,j,i\l,j\l$ \\
                               &                                    &                  & $A_q, G_q, S^{(1)}_q$  & $i,j,i\l,j\l$ \\
                               &                                    &                  & $A_q, G_q+S^{(1)}_q$  & $k,k\l$ \\

\hline

\end{longtable}
\end{center}

\clearpage

\section{Direct sums in $sl(3,\mathbb{O})$}
\label{direct_sums_in_E6}

While the following list is probably not complete, we include it to show that we may find the following direct sums of Lie algebras in $sl(3,\mathbb{O})$:

\begin{enumerate}
    \item $u(1) \oplus u(1) \oplus su(3) \oplus so(3,1)$ where we have the following basis for the algebras:
\[
\begin{array}{|ccc|}
\hline
\textrm{Algebra} & \textrm{Casimir Operator} & \textrm{Basis} \\
\hline
u(1) & s^{(1)}_l & \\
u(1) & B^{(1)}_{tz} + 2B^{(2)}_{tz} & \\
su(3,\mathbb{C})^C & A_l, G_l & A_i, A_j A_k, A_{kl}, A_{jl}, A_{il} \\
so(3,1) & R^{(1)}_{x\l}, B^{(1)}_{tz} & R^{(1)}_{z\l}, R^{(1)}_{xz}, B^{(1)}_{tx}, B^{(1)}_{t\l} \\
\hline
\end{array}
\]
    \item $u(1) \oplus su(4) \oplus so(3,1)$ where we have the following basis for the algebras:
\[
\begin{array}{|ccc|}
\hline
\textrm{Algebra} & \textrm{Casimir Operator} & \textrm{Basis} \\
\hline
u(1) & B^{(1)}_{tz} + 2B^{(2)}_{tz} & \\
su(4) & A_l, G_l+S^{(1)}_l, G_l & A_i, \cdots, A_{il} \\
      &                     &  G_i+S^{(1)}_i, \cdots, G_{il}+S^{(1)}_{il} \\
so(3,1) & R^{(1)}_{x\l}, B^{(1)}_{tz} & R^{(1)}_{z\l}, R^{(1)}_{xz}, B^{(1)}_{tx}, B^{(1)}_{t\l} \\
\hline
\end{array}
\]
    \item $u(1) \oplus u(1) \oplus su(3) \oplus so(3,1)$ where we have the following basis for the algebras:
\[
\begin{array}{|ccc|}
\hline
\textrm{Algebra} & \textrm{Casimir Operator} & \textrm{Basis} \\
\hline
u(1) & 2G_l + s^{(1)}_l & \\
u(1) & B^{(1)}_{tz} + 2B^{(2)}_{tz} & \\
su(3) & A_l, -G_l+4S^{(1)}_l & G_i+2S^{(1)}_i, G_j+2S^{(1)}_j, A_k, A_{kl} \\
      &                  & G_{jl}+2S^{(1)}_{jl}, G_{il}+2S^{(1)}_{il} \\
so(3,1) & R^{(1)}_{x\l}, B^{(1)}_{tz} & R^{(1)}_{z\l}, R^{(1)}_{xz}, B^{(1)}_{tx}, B^{(1)}_{t\l} \\
\hline
\end{array}
\]
    \item $u(1) \oplus su(2) \oplus su(2) \oplus sl(2,\mathbb{H})$ where we have the following basis for the algebras:  
\[
\begin{array}{|ccc|}
\hline
\textrm{Algebra} & \textrm{Casimir Operator} & \textrm{Basis} \\
\hline
u(1) & B^{(1)}_{tz} + 2B^{(2)}_{tz} & \\
su(2,\mathbb{C})^C & A_l & A_k, A_{kl} \\
su(2,\mathbb{C})_2 & G_l + 2S^{(1)}_l & G_k +2S^{(1)}_k, G_{kl}+2S^{(1)}_{kl} \\
sl(2,\mathbb{H}) & B^{(1)}_{tz}, R^{(1)}_{x\l}, G_l-S^{(1)}_l & B^{(1)}_{tx}, B^{(1)}_{t\l}, R^{(1)}_{xz}, R^{(1)}_{z\l} \\
      &     & G_k-S^{(1)}_k, G_{kl}-S^{(1)}_{kl}, R^{(1)}_{xk}, R^{(1)}_{xkl} \\
      &     & R^{(1)}_{zk}, R^{(1)}_{zkl}, B^{(1)}_{tk}, B^{(1)}_{tkl}\\
\hline
\end{array}
\]
We note that $sl(2,\mathbb{H})$ has $so(3,1)\oplus u(1)$ as subalgebra.

    \item $u(1) \oplus so(5) \oplus so(4,1)$ where we have the following basis for the algebras:
\[
\begin{array}{|ccc|}
\hline
\textrm{Algebra} & \textrm{Casimir Operator} & \textrm{Basis} \\
\hline
u(1) & B^{(1)}_{tz} + 2B^{(2)}_{tz} & \\
so(5) & R^{(1)}_{x\l}, G_l-S^{(1)}_l & R^{(1)}_{xz}, R^{(1)}_{z\l}, R^{(1)}_{zk}, R^{(1)}_{zkl} \\
      &                     & R^{(1)}_{xk}, R^{(1)}_{xkl} \\
      &                     & G_k-S^{(1)}_k, G_{kl}-S^{(1)}_{kl} \\
so(4,1) & A_l, G_l+2S^{(1)}_l   & A_k, A_{kl}, G_k+2S^{(1)}_k, G_{kl}+2S^{(1)}_{kl} \\
        &                   & B^{(1)}_{ti}, B^{(1)}_{tj}, B^{(1)}_{til}, B^{(1)}_{tjl}\\
\hline
\end{array}
\]
We note that $so(4,1)$ has $so(3)\oplus so(3)$ as a subalgebra.

\end{enumerate}

\if\Book0\end{document}\fi

\addtocontents{toc}{\medskip}

\end{document}